\documentclass{article}

\usepackage[margin=1in]{geometry}
\usepackage[numbers]{natbib}

\usepackage[utf8]{inputenc} 
\usepackage[T1]{fontenc}    
\usepackage[colorlinks,linkcolor=blue,filecolor=blue,citecolor=black,urlcolor=blue]{hyperref}       
\usepackage{url}    

\usepackage{booktabs}      
\usepackage{amsfonts}       
\usepackage{nicefrac}       
\usepackage{microtype}      
\usepackage[table]{xcolor}         
\usepackage{tocvsec2}
\usepackage{titletoc}

\usepackage{xspace}

\usepackage{amsmath}
\usepackage{amsthm}
\usepackage{amssymb}

\usepackage{multicol}
\usepackage{multirow}
\usepackage{makecell}
\usepackage{enumitem}
\usepackage{wrapfig}
\usepackage{dsfont}

\usepackage{graphicx,epstopdf}
\usepackage{subfigure} 

\usepackage{algorithm}
\usepackage{algorithmic}

\usepackage{threeparttable, booktabs}
\usepackage{mathtools}
\mathtoolsset{showonlyrefs=false}
\usepackage{float}

\usepackage{hyperref}[pdfusetitle]
\hypersetup{colorlinks,urlcolor=blue}
\usepackage{xcolor,pifont}
\usepackage{MnSymbol,wasysym}
\newcommand{\cmark}{\textcolor{green}{\large\ding{52}}}
\newcommand{\xmark}{\textcolor{red}{\large\ding{55}}}

\newcommand{\rcry}{\textcolor{red}{\Large\frownie}}

\usepackage{tabularray}

\usepackage[skins,listings]{tcolorbox}

\newtcblisting{tikzexample}{
    sidebyside,
    center lower,
    bicolor,
    colbacklower=white,
    sharp corners,
    frame engine=empty
}

\usepackage{quiver}

\restylefloat{algorithm} 
\newtheorem{assumption}{Assumption}
\newtheorem{theorem}{Theorem}
\newtheorem{remark}{Remark}
\newtheorem{lemma}{Lemma}

\usepackage{diagbox}

\newcommand{\tran}{^\mathsf{T}}

\def\w{\boldsymbol{w}}

\newcommand{\ours}{\texttt{Clapping} }
\newcommand{\oursfu}{{{\texttt{Clapping}-\textbf{FU} }}}
\newcommand{\oursfc}{{{\texttt{Clapping}-\textbf{FC} }}}

\newcommand{\bv}{{{v}}}

\newcommand{\by}{{{y}}}

\def\tran{^{\mathsf{T}}}

\newcommand{\vvvert}{{\vert\kern-0.25ex\vert\kern-0.25ex\vert}}

\allowdisplaybreaks

\usepackage{array}
\usepackage{mathrsfs}
\usepackage[normalem]{ulem}

\newcommand{\oursfulast}{{{\texttt{Clapping}-\textbf{FU}}}}
\newcommand{\oursfclast}{{{\texttt{Clapping}-\textbf{FC}}}}

\title{\textbf{\texttt{Clapping}}: Removing Per-sample Storage for Pipeline Parallel Learning with Communication Compression}

\author{%
  \normalfont
  \begin{minipage}[t]{0.96\textwidth}
    \centering
    Boao Kong$^{1,*}$,\quad Xu Huang$^{1,*}$,\quad Yuqi Xu$^{1,*}$\\[2pt]
    Yixuan Liang$^{2}$,\quad Bin Wang$^{3}$,\quad Kun Yuan$^{1,\dagger}$\\[6pt]
    {\footnotesize\texttt{kongboao@stu.pku.edu.cn},\quad
     \texttt{ydove1031@gmail.com},\quad
     \texttt{xuyq@stu.pku.edu.cn}}\\[2pt]
    {\footnotesize\texttt{liangyx25@mails.tsinghua.edu.cn},\quad
     \texttt{wangbin\_bme@zju.edu.cn},\quad
     \texttt{kunyuan@pku.edu.cn}}\\[6pt]
    {\small $^{1}$ Peking University}\quad
    {\small $^{2}$ Tsinghua University}\quad
    {\small $^{3}$ MetaCarbon}\\[4pt]
    {\footnotesize
    $^{*}$ Equal contribution.\quad
    $^{\dagger}$ Corresponding author.
    }
  \end{minipage}%
}

\begin{document}

\setlength{\parindent}{0pt}
\setlength{\parskip}{0.5em}
\maketitle

\begin{abstract}
  Pipeline-parallel distributed optimization is essential for large-scale machine learning but is challenged by significant communication overhead from transmitting high-dimensional activations and gradients between workers. Existing approaches often depend on impractical unbiased gradient assumptions or incur sample-size memory overhead. This paper introduces \textbf{\texttt{Clapping}}, a \underline{\textbf{C}}ommunication compression algorithm with \underline{\textbf{LA}}zy sam\underline{\textbf{P}}ling for \underline{\textbf{P}}ipeline-parallel learn\underline{\textbf{ING}}. \texttt{Clapping} adopts a lazy sampling strategy that reuses data samples across steps, breaking sample-wise memory barrier and supporting convergence in few-epoch or online regimes. \ours comprises two variants including \texttt{Clapping}-\textbf{FC} and \texttt{Clapping}-\textbf{FU}, both of which achieve convergence without unbiased assumption for compressed gradient, effectively addressing compression error propagation in multi-worker settings. Numerical experiments validate the performance of \ours across different tasks. 
\end{abstract}

\section{Introduction}
Large-scale optimization and learning have become essential tools in numerous applications. Addressing these complex and large problems presents a substantial challenge, frequently necessitating extensive computation over many days or even months. Consequently, distributed algorithms are crucial for accelerating large-scale optimization and learning processes. In distributed optimization, multiple workers collaborate to solve a global problem with the help of communication between workers. Most existing research focuses on data-parallel distributed optimization \citep{li2014scaling,nedic2009distributed,chen2012diffusion,yuan2016convergence,konevcny2015federated,alistarh2017qsgd}. In this paradigm, each worker maintains a complete replica of the model, independently samples training data, and exchanges models or gradients at each iteration, thereby achieving significant acceleration of large-scale optimization and learning tasks.

As model parameters in contemporary optimization and learning problems have grown to hundreds of billions \citep{radford2019language, koroteev2021bert, rae2021scaling, zhang2022opt, liu2024deepseek}, these models have exceeded single-worker memory capacity, necessitating model-parallel distributed optimization. This paradigm partitions the model across multiple workers, with each worker managing a subset of parameters, enabling the training of massive models that would be intractable on a single worker. Model-parallel distributed optimization is particularly prevalent in the pre-training and fine-tuning of Large Language Models (LLMs) \citep{touvron2023LLaMA,meta2024introducing}. For instance, consider an LLM architecture comprising 24 transformer layers that exceeds the memory capacity of a single GPU cluster. The model can be efficiently segmented, with the initial 12 layers allocated to one GPU cluster and the remaining layers to another. This partitioning ensures each GPU maintains a fraction of the model's parameters, substantially reducing per-device memory requirements. 
Another prominent application of model-parallel distributed optimization is split learning \citep{thapa2022splitfed,lin2024split,wang2022fedlite}, where the entire machine learning model is divided into smaller network segments and trained independently across multiple edge computing devices. 

\begin{figure*}[t]
\centering
\vspace{-1mm}
\includegraphics[width=0.9\textwidth]{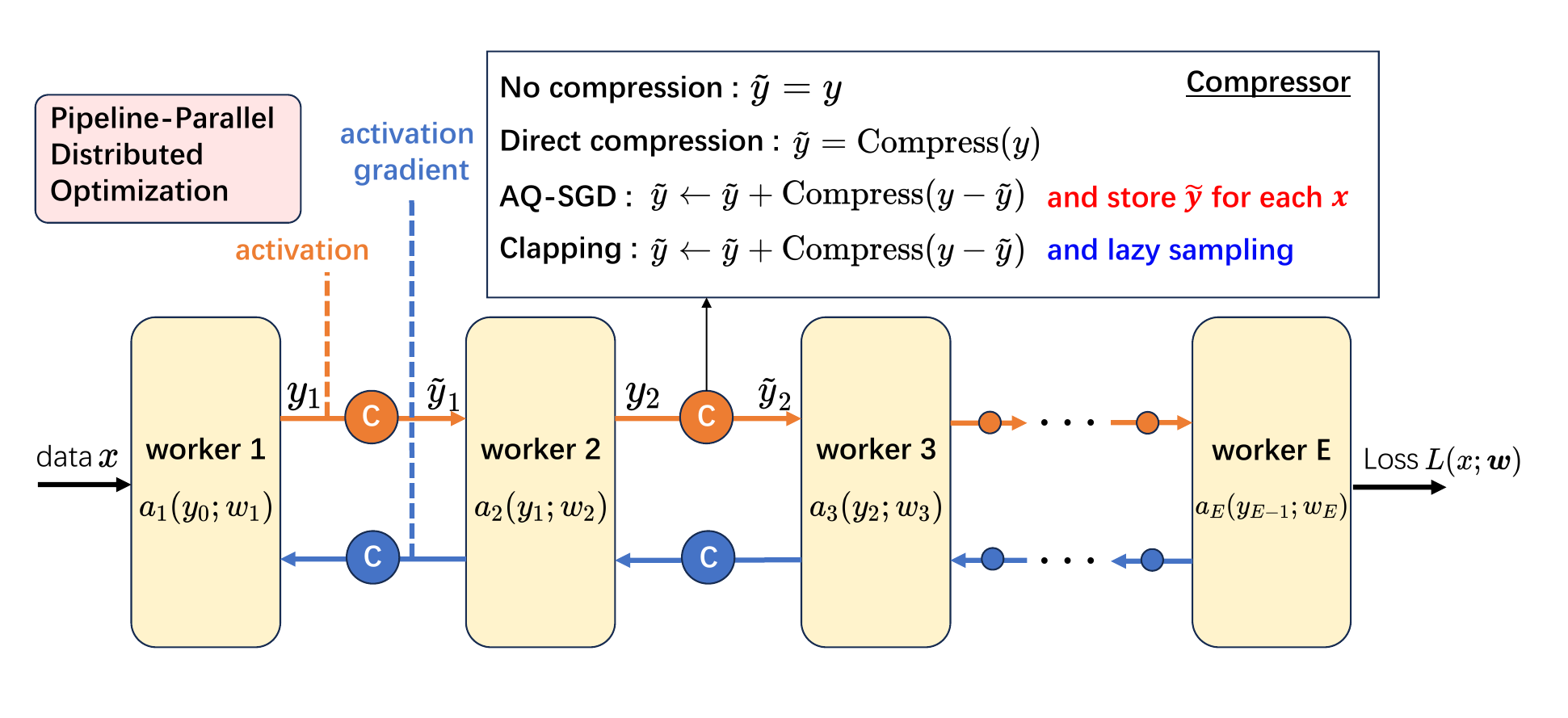}
\vspace{-1mm}
\caption{\small Illustration of the pipeline-parallel distributed optimization with communication compression.}
\label{fig:pipeline}
\vspace{-1mm}
\end{figure*}

\textbf{Problem statement.} This paper examines pipeline-parallel distributed optimization, a specific form of model-parallel  optimization that partitions model parameters in a pipeline-like manner \citep{narayanan2019pipedream, huang2019gpipe, narayanan2021efficient, ryabinin2023swarm, wan2025pipeoffload}. Consider a computing cluster consisting of $E\geq2$ workers. Model parameter $\boldsymbol{w}$ is partitioned into $E$ parts, denoted as $\w=(w_1,\cdots,w_E)\in\mathbb{R}^{d_{w_1}+\cdots+d_{w_E}}$, where each block component \( w_e \in \mathbb{R}^{d_{w_e}} \) is maintained by worker $e$. Pipeline-parallel  optimization can be  formulated as follows: 

\vspace{-2mm}
\begin{subequations}
\label{pipeline parallel}
\begin{align}
    \min_{\w \in\mathbb{R}^d}&\quad \mathbb{E}_{x\in\mathcal{D}}\Big[L(x;\w):=y_E\Big],\\
    \text{s.t.}&\quad y_e = a_e(y_{e-1}, w_e), \quad \forall e \in \{1, \cdots E\}. \label{eq:pipeline-2}
\end{align}
\end{subequations}

\vspace{-2mm}
In the above problem, the random variable $x$ denotes the data sample following distribution $\mathcal{D}$. The loss \( L(x,\boldsymbol{w} ) \) is a composite function consisting of \( E \) operators \( a_e(y_{e-1}, w_e): \mathbb{R}^{d_{e-1}} \times \mathbb{R}^{d_{w_e}} \to \mathbb{R}^{d_e} \) with $d_E = 1$, where \( y_e \), following the terminology in LLMs, refers to the activation. For initialization, we let $y_0 = x$. An illustration for problem \eqref{pipeline parallel} is shown in Fig.~\ref{fig:pipeline}, where each worker participates in the computation in a pipeline fashion. In LLMs, each operator $a_e(y_{e-1}, w_e)$ corresponds to a transformer layer. In split learning, there are two operators: \( a_1(y_0, w_1) \), which represents the client-side sub-network, and \( a_2(y_1, w_2) \), which corresponds to the server-side sub-network.

\begin{wrapfigure}{r}{0.42\textwidth}
\centering
\vspace{-6.2pt}
\includegraphics[width=0.4\textwidth]{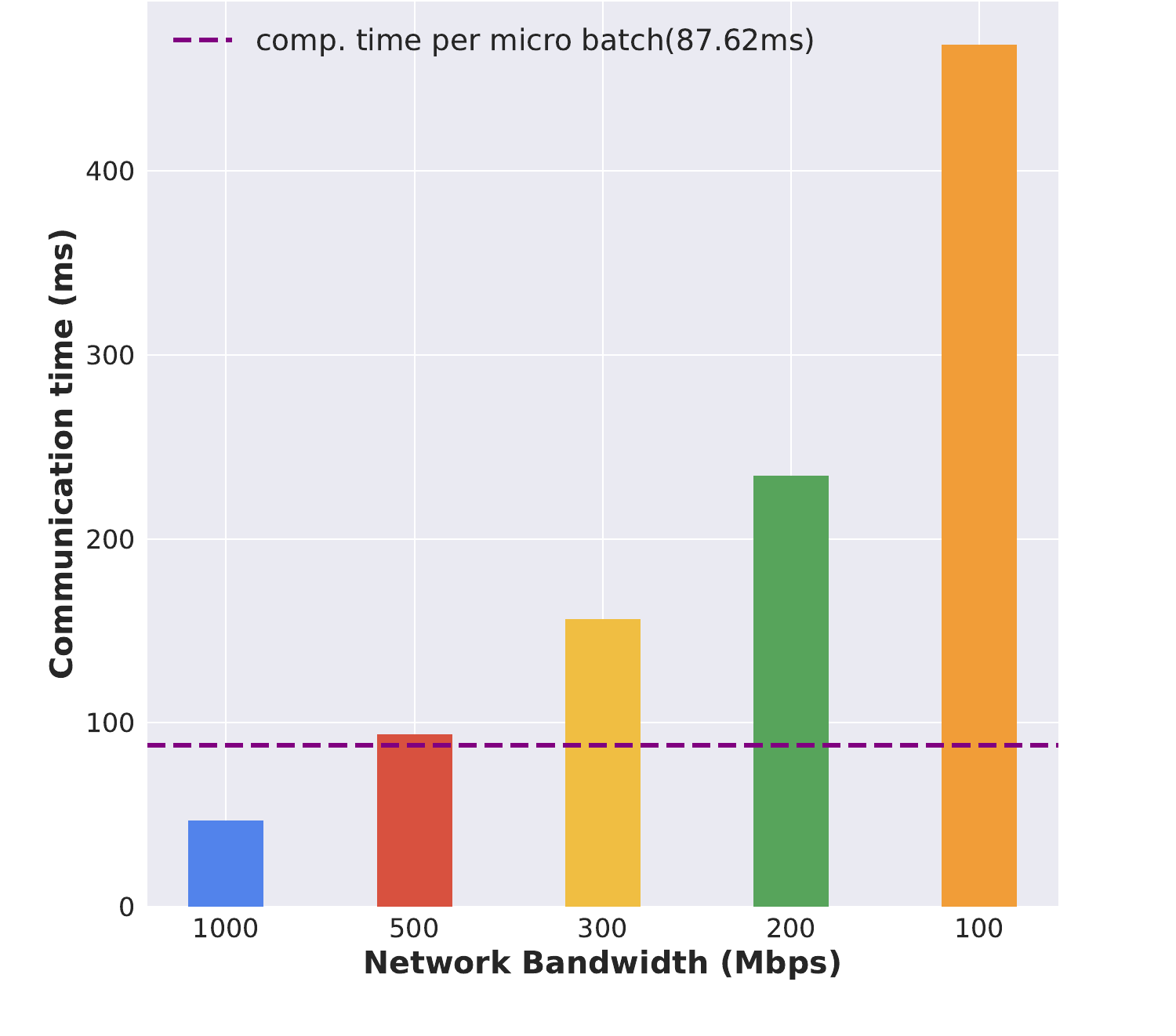}
\vspace{-10.5pt}
\caption{\small Communication time for fine-tuning GPT2-xl with different bandwidth and average computation time per micro batch on 4 Nvidia A100 GPUs.}
 \vspace{-2pt}
\label{fig: communication time}
\end{wrapfigure}

Pipeline-parallel distributed algorithms are effective for large-scale problems with multiple devices but incur significant communication overhead from transmitting high-dimensional activation and gradient vectors over low-bandwidth links \citep{diskin2021distributed,ryabinin2023swarm,wang2022fine}, as illustrated in Figure~\ref{fig: communication time}. While current LLM pre-training is typically run on tightly coupled GPU clusters with 100--400\,Gbps interconnects \citep{smith2022using}, such infrastructure is accessible only to a few industrial labs. In contrast, emerging decentralized and over-the-internet training frameworks aim to scale model parallelism across commodity GPUs connected via consumer-grade networks at only \emph{tens to hundreds of Mbps}, yet still match datacenter-level convergence \citep{ramasinghe2025protocol,wang2022fine}. This regime is fundamentally communication-bound: even a single forward–backward pass can be dominated by transmitting activations and gradients. These trends motivate studying pipeline-parallel \emph{communication compression} methods that remain stable under aggressive activation/gradient compression without assuming high-bandwidth interconnects, and thus are applicable to both datacenter and slow-network deployments.

\textbf{Fundamental challenges.} While communication compression has been extensively studied for data-parallel distributed optimization \citep{alistarh2017qsgd,stich2018sparsified,vogels2019powersgd,richtarik2021ef21,huang2022lower,fatkhullin2024momentum}, it remains largely unexplored for pipeline-parallel distributed optimization. Compressing activations and gradients, as illustrated in Fig.\ref{fig:pipeline}, presents two fundamental challenges. First, each stochastic gradient computed through forward-backward propagation requires \( 2E - 2 \) rounds of compression (\( E - 1 \) in the forward pass and \( E - 1 \) in the backward pass), introducing substantial errors in gradient estimation that may lead to non-convergence. Second, the composite structure of activations (see Eq. \eqref{eq:pipeline-2}) results in error propagation during pipeline communication, causing a non-separable entanglement of gradient information and compression errors (detailed in Section\ref{sec:direct-compression}). These challenges render most effective techniques from data-parallel optimization infeasible for pipeline-parallel optimization.

\subsection{Limitations in existing works}
Several studies have emerged to address these challenges \citep{evans2021ac,fu2020don,wang2022fine}, but they exhibit critical limitations.

\textbf{(L1) Impractical unbiased assumption for compressed gradient.}
The unbiased assumption for compressed gradient is crucial for ensuring the convergence of optimization algorithms. In line with this, existing works \citep{evans2021ac,fu2020don} assume that unbiased activation compression leads to unbiased gradient errors, but this fails due to the composite and non-linear structure of the operator \( a_e(y_{e-1}, w_e) \). For example, even if \( \mathbb{E}[\tilde{y}_{e-1}] = y_{e-1} \), it is not generally true that \( \mathbb{E}[a_e(\tilde{y}_{e-1}, w_e)] = a_e(y_{e-1}, w_e) \).

\textbf{(L2) Sample-size memory overhead.}
\cite{wang2022fine} proposes AQ-SGD, which achieves convergence without requiring the unbiased assumption for compressed gradient. Unlike approaches in \citep{evans2021ac,fu2020don} that directly compress activations, AQ-SGD compresses the change in activations for identical training samples across epochs. However, AQ-SGD requires storing activations for \textbf{all training samples}, resulting in memory overhead proportional to the sample size.

\textbf{(L3) Multiple-epoch training requirement.}
Another technique behind the convergence of AQ-SGD is error compensation \citep{seide20141,stich2018sparsified,richtarik2021ef21}, which can progressively mitigate activation compression errors. However, this approach requires extensive training epochs to eliminate errors, leading to non-convergence in few-epoch training regimes. This limitation is particularly significant for large-scale LLM training/fine-tuning, where practical workloads typically involve only a few epochs \citep{rae2021scaling,zhang2022opt,meta2024introducing,song2022hybrid,nakamoto2024cal,song2024importance,guo2024direct}.

\textbf{(L4) Limited scalability beyond two-worker setup.}
The convergence analysis in \citep{evans2021ac,fu2020don} is designed for two-worker configurations and does not extend to larger pipeline systems. While AQ-SGD offers convergence guarantees for multi-worker setups, it assumes \textbf{no error accumulation} in the multi-worker pipeline scenario, which is generally \textbf{NOT} true in practice.

\subsection{Contributions} 
This paper develops novel algorithms to address the aforementioned limitations. Our contributions are as follows: 

\textbf{(C1). }We propose {\texttt{Clapping}}, a \underline{\bf C}ommunication compression framework with \underline{\bf LA}zy sam\underline{\bf P}ling for \underline{\bf P}ipeline-parallel Learn\underline{\bf ING}. The {\texttt{Clapping}} framework is versatile and can be implemented in two variants: \oursfc and \oursfulast, which differ in their approaches to the compression strategy in forward- and backward-propagation. A core technique underpinning {\texttt{Clapping}} is its novel lazy sampling strategy. By retaining the same data sample across multiple steps, this strategy enables {\texttt{Clapping}} to overcome the sample-size memory overhead associated with AQ-SGD and achieve convergence under few-epoch or even single-epoch training regimes, effectively addressing \textbf{Limitations} \textbf{(L2)} and \textbf{(L3)}.

\textbf{(C2). }We establish convergence guarantees for two variants, \oursfc and \oursfulast, for pipeline-parallel optimization with inter-stage communication compression. Under standard smoothness assumptions and a contractive compressor model, \oursfc achieves an $O(T^{-1/3})$ stationarity rate and \oursfu achieves an $O(T^{-1/2})$ rate, without relying on the additional \emph{unbiased compressed-gradient} assumption and multi-epoch training, addressing \textbf{Limitations (L1)} and \textbf{(L4)}. Moreover, our analysis explicitly characterizes how compression errors accumulate with the pipeline depth $E$ in multi-stage settings.

\textbf{(C3). }We conduct extensive experiments across split learning, LLMs fine-tuning and pre-training. \ours is comparable to AQ-SGD in several fine-tuning with far less memory cost and has a wider applicability than the other communication compression algorithms.

\begin{table*}[t] 
\caption{Comparison between different communication compression algorithms for pipeline-parallel optimization. Notation $T$ denotes the number of iterations, $x$ denotes the training sample, $N$ is the size of training data, $B$ is the batch size, $q$ denotes the dimension of the communicated parameters. We also list the result of Momentum SGD without compression in the bottom line for reference.}
\label{Table:Comparsion}
\centering
\setlength{\tabcolsep}{5pt}
\begin{threeparttable}
{
\begin{tabular}{lccccc}
\toprule
\textbf{Algorithms} & \textbf{C. Rate}$^\Diamond$ & \multicolumn{1}{c}{\begin{tabular}[c]{@{}c@{}} \textbf{No Unbiased} \\ \textbf{Comp-Grad. Assu.}$^\dagger$\end{tabular}} & \textbf{C. Mem.}$^\ddagger$ &  \multicolumn{1}{c}{\begin{tabular}[c]{@{}c@{}} \textbf{Few-epo.} \\ \textbf{Rate}$^\triangleright$\end{tabular}} &  \multicolumn{1}{c}{\begin{tabular}[c]{@{}c@{}} \textbf{Multi.} \\ \textbf{Workers}$^\triangleleft$\end{tabular}} \\ \midrule
 \specialrule{0em}{.5pt}{.5pt}
 \multicolumn{1}{l}{\begin{tabular}[c]{@{}l@{}}AC-GC \citep{evans2021ac} \\ \end{tabular}}  & $T^{-1/2}$ & \xmark & 0 &  N.A. &  \xmark\\
 \specialrule{0em}{.5pt}{.52pt}
 \multicolumn{1}{l}{\begin{tabular}[c]{@{}l@{}}TinyScript \citep{fu2020don}\\ \end{tabular}}& $T^{-1/2}$ & \xmark & 0 &  N.A. & \xmark \\
 \specialrule{0em}{.5pt}{.5pt}
 \multicolumn{1}{l}{\begin{tabular}[c]{@{}l@{}}AQ-SGD \citep{wang2022fine}\\ \end{tabular}}   & ${N}T^{-1/2}$ & \cmark & $\mathcal{O}(Nq)$ &  \xmark &  \rcry \\
 \specialrule{0em}{.5pt}{.5pt}
\rowcolor{pink!40} \textbf{\oursfc (Ours)} & $\boldsymbol{T^{-1/3}}$ & \cmark & $\boldsymbol{\mathcal{O}(Bq)}$ & \cmark & \cmark \\
 \specialrule{0em}{.5pt}{.5pt}
\rowcolor{blue!25} \textbf{\oursfu (Ours)} & $\boldsymbol{T^{-1/2}}$ & \cmark & $\boldsymbol{\mathcal{O}(Bq)}$ & \cmark & \cmark \\ \midrule
 Momentum SGD \citep{nesterov2013introductory}   & $T^{-1/2}$ & \cmark & N.A. & \cmark & \cmark  \vspace{0.1mm}\\ \bottomrule
\end{tabular}}
\begin{tablenotes}
\small
    \item[$\Diamond$] The convergence rate with respect to total steps $T\to \infty$ (smaller is better).
    \item[$\dagger$] No additional assumptions regarding the unbiased compression-induced gradient error.
    \item[$\ddagger$] Extra memory overhead incurred during a single communication compression operation. N.A. indicates no compression. 
    \item[$\triangleright$] The algorithm can achieve convergence in the few-epoch learning tasks, i.e., $N=\mathcal{O}(T)$. 
    \item[$\triangleleft$] Convergence analysis can be extended to the multiple-worker setup.  \textcolor{red}{\frownie} indicates the analysis is with impractical assumptions.
\end{tablenotes}
\end{threeparttable}
\vspace{-2mm}
\end{table*}

All theoretical results and existing algorithms for communication compression in pipeline-parallel optimization are summarized in Table~\ref{Table:Comparsion}. Notably, \ours achieves convergence without assuming an unbiased compression-induced gradient estimator, avoids sample-wise activation storage, and applies to few-epoch workloads. Moreover, our theory characterizes how compression errors propagate and accumulate with the pipeline depth $E$, rather than implicitly assuming the depth effect away.

\textbf{Notations. }
\label{appendix: notations}
In this paper, we define $\ell(\w) := \mathbb{E}_{x\sim\mathcal{D}}[L(x;\w)]$ throughout this paper. We further denote its partial gradient with respect to the parameter \(w_e\) as \(\nabla_e \ell(\w) := {\partial \ell(\w)}/{\partial w_e} \in \mathbb{R}^{d_{w_e}}\) for \(e = 1, 2, \dots, E\). For the operator \(a_e(y_{e-1}, w_e): \mathbb{R}^{d_{e-1}} \times \mathbb{R}^{d_{w_e}} \to \mathbb{R}^{d_e}\), we define \(\nabla_1 a_e(y_{e-1}, w_e) \in \mathbb{R}^{d_e \times d_{e-1}}\) and \(\nabla_2 a_e(y_{e-1}, w_e) \in \mathbb{R}^{d_e \times d_{w_e}}\) as the Jacobian matrix with respect to \(y_{e-1}\) and \(w_e\), respectively. The notation \(\|\cdot\|\) represents the \(\ell_2\)-norm of vectors, and \(\mathbf{1}_B \in \mathbb{R}^B\) denotes a vector with all elements equal to 1. For variables \(\{y_e^{(t)}\}_{e=1,\dots,E}^{t=0,\dots,T+1}\), the subscript \(e\) (resp. superscript \(t\)) denotes the index of the worker (resp. iteration). We use \(a \lesssim b\) to indicate that there exists a constant \(C \geq 0\) such that \(a \leq Cb\), and \(a \lesssim_d b\) indicates that there exists a \(C \geq 0\) that is independent with $d$ such that \(a \leq Cb\).

    


\section{Prior Arts}

Here we display some prior works of pipeline-parallel distributed optimization. Additional related works can be found in Appendix \ref{appendix:related work}.

\subsection{Pipeline-parallel SGD}\label{sec:pp-SGD}

To solve pipeline-parallel optimization problem \eqref{pipeline parallel}, we define \( v_e := \nabla_{y_e} L(x;\w) \in \mathbb{R}^{d_e} \) as the gradient of the loss function with respect to the activation \( y_e \), referred to as the activation gradient. Similarly, we define \( u_e := \nabla_{w_e} L(x;\w) \in \mathbb{R}^{d_{w_e}} \) as the gradient with respect to the weight $w_e$, referred to as the weight gradient. Pipeline-parallel SGD performs the following steps at iteration $t$:
\begin{itemize}[leftmargin=*]
\vspace{-2mm}
    \item \textbf{Forward}: Each worker $e$ computes activation $y^{(t)}_e =a_e({y}^{(t)}_{e-1},w^{(t)}_e)$ following the forward order $e=1,\cdots, E$. We let $y^{(t)}_0$ be the random sample $x^{(t)}$ for initialization.
\vspace{-1mm}
    \item \textbf{Backward}: Each worker $e$ computes activation gradient $v^{(t)}_{e-1} = \nabla_1 a_e({y}^{(t)}_{e-1},w^{(t)}_e)\tran v^{(t)}_{e}$ following the backward order $e=E,\cdots, 2$. We initialize $v^{(t)}_{E} = 1$.
\vspace{-1mm}
    \item \textbf{Update}: Each worker $e$ computes weight gradient $u^{(t)}_{e} = \nabla_2 a_e({y}^{(t)}_{e-1},w^{(t)}_e)\tran v^{(t)}_{e}$ following the backward order $e=E,\cdots, 1$. We update parameter $w_e^{(t+1)} = w_e^{(t)} - \gamma u_e^{(t)}$.
\end{itemize}

Pipeline-parallel SGD iteratively repeats the aforementioned steps until convergence (see Fig.~\ref{fig:pipeline}). This requires transmitting all activations and their corresponding gradients between workers. To reduce communication overhead, we can apply direct compression of activations and gradients \citep{evans2020jpeg,fu2020don}. However, such a compression may cause non-convergence due to biased gradient estimation and error propagation, even with unbiased compressors (see Appendix \ref{sec:direct-compression}).

\subsection{AQ-SGD}\label{sec:aq-sgd}
To eliminate convergence bias in direct compression, \cite{wang2022fine} proposes AQ-SGD, an algorithm designed to ensure exact convergence to the stationary solution without relying on the impractical assumption of unbiased compressed gradients. The key mechanism employed by AQ-SGD is error feedback \citep{richtarik2021ef21}. For a specific random sample \( x \), we let \( y_{x,e} \) represent the activation associated with that sample on each worker \( e \), where \( y_{x,0} = x \) serves as the initialization. Rather than compressing each activation \( y_{x,e} \) directly, AQ-SGD compresses the changes in activations as follows:

\vspace{-3.5mm}
\begin{align}
\label{eq:forward-EF}
    y^{(t)}_{x,e} &= a_e(\tilde{y}^{(t)}_{x,e-1},w^{(t)}_e), \\
    \tilde{y}^{(t)}_{x,e} &= \tilde{y}^{(t-1)}_{x,e} + \mathcal{C}(y^{(t)}_{x,e} - \tilde{y}^{(t-1)}_{x,e}). 
\end{align}  
\vspace{-3.5mm}

Suppose the compressor \(\mathcal{C}(\cdot)\) is contractive, i.e., \(\|\mathcal{C}(y) - y\| \le \omega \|y\|\) for some \(\omega \in (0,1)\) (see Assumption \ref{assumption:compressor} for details), it then holds that

\vspace{-3.5mm}
\begin{align*}
\|\tilde{y}^{(t)}_{x,e} - {y}^{(t)}_{x,e}\| &= \| \mathcal{C}(y^{(t)}_{x,e} - \tilde{y}^{(t-1)}_{x,e}) - (y^{(t)}_{x,e} - \tilde{y}^{(t-1)}_{x,e})\| \le \omega \|\tilde{y}^{(t-1)}_{x,e} - {y}^{(t)}_{x,e}\|.
\end{align*}
\vspace{-3.5mm}

If we further assume that \(y^{(t)}_{x,e} \to y^\star_{x,e} := a_e(y^\star_{x,e-1}, w_e^\star)\), it follows that $\|\tilde{y}^{(t)}_{x,e} - {y}^\star_{x,e}\| \le \omega \|\tilde{y}^{(t-1)}_{x,e} - {y}^\star_{x,e}\|$, implying that \(\tilde{y}^{(t)}_{x,e}\) converges asymptotically to \(y^\star_{x,e}\). 

While error feedback progressively eliminates the compression bias in the activations, it requires repeatedly executing \eqref{eq:forward-EF} with the same sample \( x \) over multiple iterations. To enable the error feedback update, AQ-SGD passes all data samples across multiple epochs and stores each \(\tilde{y}_{x,e}\) for every data sample \( x \) and worker \( e \). When a data sample \( x \) is selected at iteration \( t \), AQ-SGD updates \(\tilde{y}^{(t)}_{x,e}\) according to \eqref{eq:forward-EF}; otherwise, it retains the previous value, ensuring \(\tilde{y}^{(t)}_{x,e} = \tilde{y}^{(t-1)}_{x,e}\). Specifically, AQ-SGD operates as follows.

\begin{itemize}[leftmargin=*]
    \item \textbf{Forward}: Suppose sample $x$ is selected at the current iteration $t$. Each worker $e$ computes activation $y^{(t)}_{x,e}$ and compresses $\tilde{y}_{x,e}^{(t)}$ according to \eqref{eq:forward-EF}, following the forward order $e=1,\cdots, E$. We let $y^{(t)}_0 = x^{(t)}$. Finally, we update $\tilde{y}_{x^\prime,e}^{(t)} = \tilde{y}_{x^\prime,e}^{(t-1)}$ for any $x^\prime \neq x$. 
\vspace{-.5mm}
    \item \textbf{Backward}: Each worker $e$ computes activation gradient $v^{(t)}_{e-1} = \nabla_1 a_e(\tilde{y}^{(t)}_{x,e-1},w^{(t)}_e)\tran \tilde{v}^{(t)}_{e}$ and compresses $\tilde{v}^{(t)}_{e-1} = \mathcal{C}(v^{(t)}_{e-1})$ following the order $e=E,\cdots, 2$. 
\vspace{-.5mm}
    \item \textbf{Update}: Each worker $e$ computes weight gradient $u^{(t)}_{e} = \nabla_2 a_e(\tilde{y}^{(t)}_{x,e-1},w^{(t)}_e)\tran \tilde{v}^{(t)}_{e}$ following the  order $e=E,\cdots, 1$. We update parameter $w_e^{(t+1)} = w_e^{(t)} - \gamma u_e^{(t)}$.
\end{itemize}

AQ-SGD operations are illustrated in Fig.~\ref{fig:pipeline}, with error feedback applied only during the forward pass (see details in \citep{wang2022fedlite}). Each worker \( e \) must store activations \( \tilde{y}_{x,e-1} \) for every sample \( x \), limiting the method to tiny datasets and making it unsuitable for online learning with infinite data streams. This requires memory proportional to sample size, resulting in significant overhead and posing a substantial challenge in large-scale optimization, especially with extensive datasets and large models.

\section{\ours Algorithm}

\begin{algorithm}[t]
  \caption{$\mbox{\sffamily LazySampling}(\mathcal{D},t,p)$ }
  \label{alg:lazy_sampling}
  \begin{algorithmic}[1]
  \IF{$t=1$}
    \STATE Sample $x^{(1)} \sim \mathcal{D}$; \ \ $f^{(1)}_{\text{FU}} \gets \texttt{True}$.
  \ELSE
    \STATE With prob. $1-p$: set $x^{(t)} \gets x^{(t-1)}$ and $f^{(t)}_{\text{FU}} \gets \texttt{False}$;
    \STATE With prob. $p$: sample $x^{(t)} \sim \mathcal{D}$ and $f^{(t)}_{\text{FU}} \gets \texttt{True}$.
  \ENDIF
  \STATE \textbf{{Return: }}$x^{(t)}$, $f^{(t)}_{\text{FU}}$.
  \end{algorithmic}
\end{algorithm}

\begin{figure}[t]
\centering
\vspace{-6.5mm}
\begin{minipage}{.49\textwidth}
{
\begin{algorithm}[H]
\caption{\footnotesize{\sffamily{Forward}$_e$$(\tilde{{y}}_e^{(t-1)},\tilde{y}_{e-1}^{(t)},w_e^{(t)},f^{(t)}_{\text{FU}})$}}
  \label{alg:forward}
  \begin{algorithmic}[1]
  \STATE $y_e^{(t)} \gets a_e(\tilde{y}_{e-1}^{(t)},w_e^{(t)})$.
  \IF{\oursfu \textbf{and} $f^{(t)}_{\text{FU}}=\texttt{True}$}
    \STATE Send $y_e^{(t)}$ from worker $e$ to $e{+}1$;
    \STATE $\tilde{{y}}_{e}^{(t)} \gets y_e^{(t)}$.
  \ELSE
    \STATE Send $\mathcal{C}({y}_e^{(t)}-\tilde{{y}}_{e}^{(t-1)})$ from worker $e$ to $e{+}1$;
    \STATE $\tilde{{y}}_{e}^{(t)} \gets \tilde{{y}}_{e}^{(t-1)}+\mathcal{C}({y}_e^{(t)}-\tilde{{y}}_{e}^{(t-1)})$.
  \ENDIF
  \STATE \textbf{{Return: }}$\tilde{{y}}_{e}^{(t)}$.
  \end{algorithmic}
\end{algorithm}}
\end{minipage}
\hfill
\begin{minipage}{.49\textwidth}
{
\begin{algorithm}[H]
  \caption{\footnotesize{\sffamily{Backward}$_{e}$$(\tilde{{v}}_{e-1}^{(t-1)},\tilde{{v}}_{e}^{(t)},\tilde{y}_{e-1}^{(t)},w_{e}^{(t)},f^{(t)}_{\text{FU}})$}}
  \label{alg:backward}
  \begin{algorithmic}[1]
  \STATE ${v}_{e-1}^{(t)} \gets \nabla_1 a_{e}(\tilde{{y}}_{e-1}^{(t)},w_{e}^{(t)})^{\!\top}\tilde{{v}}_{e}^{(t)}$.
  \IF{\oursfu \textbf{and} $f^{(t)}_{\text{FU}}=\texttt{True}$}
    \STATE Send ${v}_{e-1}^{(t)}$ from worker $e$ to $e{-}1$;
    \STATE $\tilde{{v}}_{e-1}^{(t)} \gets {v}_{e-1}^{(t)}$.
  \ELSE
    \STATE Send $\mathcal{C}({v}_{e-1}^{(t)}-\tilde{{v}}_{e-1}^{(t-1)})$ from worker $e$ to $e{-}1$;
    \STATE $\tilde{{v}}_{e-1}^{(t)} \gets \tilde{{v}}_{e-1}^{(t-1)}+\mathcal{C}({v}_{e-1}^{(t)}-\tilde{{v}}_{e-1}^{(t-1)})$.
  \ENDIF
  \STATE \textbf{{Return: }}$\tilde{{v}}_{e-1}^{(t)}$.
  \end{algorithmic}
\end{algorithm}}
\end{minipage}
\end{figure}

Here, we introduce \ours, a novel approach designed to overcome the limitations of AQ-SGD.

\textbf{Lazy sampling.} As discussed in Section~\ref{sec:aq-sgd}, error feedback is critical in communication compression for pipeline-parallel distributed optimization. The primary challenge lies in the fact that the update \eqref{eq:forward-EF} must be repeated a sufficient number of times to progressively remove the compression bias. To facilitate error feedback, AQ-SGD stores \(\tilde{y}_{x,e}\) for every data sample \( x \) on each worker \( e \). This design is the key reason for its limitations, which include a sample-size memory overhead and the inability to handle infinite datasets.

We propose a lazy sampling strategy to address this challenge. Unlike AQ-SGD, which stores \(\tilde{y}_{e,x}\) and updates it only when data \( x \) is re-sampled during multiple-epoch training, our approach employs a fixed sample \( x \) for gradient evaluation across multiple consecutive steps, combined with error compensation. Once the compression bias is sufficiently reduced, we proceed to sample the next data point. As outlined in Algorithm \ref{alg:lazy_sampling}, we retain the previous sample with probability \(1-p\); otherwise, a new sample is drawn from \(\mathcal{D}\). This strategy enables sample reuse across iterations without storing \( y_{e,x} \), thereby eliminating sample-size memory overhead. Moreover, lazy sampling facilitates error feedback updates even in settings with infinitely many data samples.

\textbf{Error feedback.} \ours employs error feedback in both forward and backward processes, as detailed in Algorithms~\ref{alg:forward} and \ref{alg:backward}. In contrast, AQ-SGD applies error feedback exclusively to the forward process. Furthermore, leveraging the lazy sampling strategy, our approach eliminates the need to maintain \(\{\tilde{y}_{e,x}\}\) and \(\{\tilde{v}_{e,x}\}\) for each sample \( x \in \mathcal{D} \). As shown in Algorithms~\ref{alg:forward} and \ref{alg:backward}, we only maintain \(\tilde{y}_{e}\) and \(\tilde{v}_{e}\) across all samples and iterations.

\textbf{Incorporation of momentum.} Momentum is a widely used technique to accelerate SGD convergence, acting as a surrogate for large-batch gradients and reducing gradient variance. Recent studies underscore its theoretical benefits: \cite{cheng2023momentum} highlights its role in mitigating data heterogeneity, while \cite{fatkhullin2024momentum} demonstrates its effectiveness in enhancing error feedback. As algorithm \ref{alg:clapping} illustrates, we can adopt the following momentum update to mitigate the gradient bias from inaccurate activation gradients:  
\begin{equation}
\begin{aligned}
\label{eq: momentum}
\tilde{u}^{(t)}_{e} &= (1-m_t)\tilde{u}^{(t-1)}_{e}    + m_t\nabla_2 a_e(\tilde{y}^{(t)}_{e-1},w^{(t)}_e)\tran \tilde{v}^{(t)}_{e}, \quad
w_e^{(t+1)} = w_e^{(t)} - \gamma  \tilde{u}^{(t)}_{e}.
\end{aligned}
\end{equation}
where the coefficient $m_t \in (0,1)$. 

\textbf{\ours algorithm.} Combining lazy sampling, error feedback, and momentum updates, the complete \ours framework is presented in Algorithm \ref{alg:clapping}. While \ours is formulated with momentum SGD and a batch size of 1, it can be extended to Adam-style optimizers~\citep{kingma2014adam,loshchilov2017decoupled}. We present \ours with Adam optimizer as well as the extended theoretical convergence analysis in Appendix \ref{section: clapping_adam}, and we evaluate \ours with Adam-based optimizers in Section \ref{section:experiment}. For large-batch scenarios, lazy sampling can be adapted batch-wise, simplifying implementation. The detailed lazy sampling strategy and algorithmic formulation for large batches are provided in Appendix \ref{Appendix: Algorithm development}.

{\small
\begin{algorithm}[t]
  \caption{\ours}
  \label{alg:clapping}
  \begin{algorithmic}[1]
  \REQUIRE Initialize $\tilde{y}^{(0)}_e = 0, \tilde{v}^{(0)}_e = 0$ for $e=1,\cdots, E-1$,
  and $\tilde{u}^{(0)}_e = 0$ for $e=1,\cdots, E$.
  Initialize dataset $\mathcal{D}$, learning rate $\gamma_t$, compressor $\mathcal{C}$, and lazy sampling rate $\{p_t\}_{t=1}^T$.
  \FOR{$t=1,\cdots,T$}        
    \STATE ${x}^{(t)},f^{(t)}_{\text{FU}} \gets {\mbox{\sffamily LazySampling}}(\mathcal{D},t,p_t)$, \ \ $\tilde{y}^{(t)}_0 \gets {x}^{(t)}$, \ \ $\tilde{v}_{E}^{(t)} \gets 1$.
    \FOR{$e=1,2,\cdots,E-1$}
      \STATE $\tilde{y}^{(t)}_e \gets {\mbox{\sffamily Forward}}_e(\tilde{{y}}_e^{(t-1)},\tilde{y}_{e-1}^{(t)},w_e^{(t)},f^{(t)}_{\text{FU}})$.
    \ENDFOR
    \FOR{$e=E,E-1,\cdots,1$}
      \STATE $\tilde{v}^{(t)}_{e-1} \gets {\mbox{\sffamily Backward}}_{e}(\tilde{{v}}_{e-1}^{(t-1)},\tilde{{v}}_{e}^{(t)},\tilde{y}_{e-1}^{(t)},w_{e}^{(t)},f^{(t)}_{\text{FU}})$ if $e\neq 1$. \ \ Update $\tilde{u}_e^{(t)}$ and $w_e^{(t+1)}$ by \eqref{eq: momentum}.
    \ENDFOR
  \ENDFOR
  \end{algorithmic}
\end{algorithm}}

\textbf{\oursfc vs. \oursfulast.}
\oursfc and \oursfu differ only on refresh iterations, i.e., iterations where a fresh sample or batch is drawn. In Algorithms~\ref{alg:forward} and~\ref{alg:backward}, \oursfc is the always-compress variant: the condition in line~2 is never activated, so both the forward and backward routines always execute the compressed error-feedback branch in lines~6--7. By contrast, \oursfu executes the raw, uncompressed branch in lines~3--4 on refresh iterations, and otherwise uses the same compressed branch as \oursfclast. Thus, at a fixed iteration budget, \oursfc sends fewer bytes than \oursfulast. For example, at \(p_t=0.4\), the fixed-budget payload reduction of \oursfc is about \(1.67\times\) that of \oursfulast. The role of \oursfu is different: raw refreshes act as resets that can improve optimization stability or time-to-target performance, rather than reducing per-iteration payload more than \oursfclast.

\textbf{Memory Overhead.} With the error feedback technique, \ours caches the current batch's activations and gradients, resulting in \(\mathcal{O}(B)\) memory overhead for batch size \( B \). In contrast, \citep{wang2022fine}'s sample-wise error compensation incurs \(\mathcal{O}(N)\) memory requirement for sample size \( N \). The memory cost reduction is significant: while \cite{wang2022fine} theoretically requires TBs for a single communication compression, \ours only needs several GBs for pre-training models with billions of parameters like LLaMA-2 7B or LLaMA-3 8B, acceptable in practice (See Appendix \ref{appendix: memory overhead} for details). This advantage is more pronounced with larger models and more extensive datasets.

\section{Theoretical Analysis}
This section presents theoretical analysis for \oursfc and \texttt{Clapping}-\textbf{FU}.

\subsection{Assumptions}
We first introduce assumptions used throughout this paper. 

\begin{assumption}
\label{assumption:smoothness}
    There exist constants $L_{\nabla \ell},C_a,L_{\nabla a},L_a$ such that:
    \begin{enumerate}[leftmargin=*]
        \item $\nabla \ell$ is $L_{\nabla \ell}$-Lipschitz continuous;
        \item For $e=1,2,\cdots,E$, the gradient of $a_e$ can be bounded by $C_a$, i.e. $||\nabla a_e(y,\w)||\leq C_a$;
        \item For $e=1,2,\cdots,E-1$, $a_e(y,\w),\nabla a_{e+1}(y,\w)$ are $L_a,L_{\nabla a}$ Lipschitz continuous with respect to $y$ and $\w$, respectively.
    \end{enumerate}
\end{assumption}

We remark that Assumption \ref{assumption:smoothness} is comparable to the smoothness assumption used in \cite{wang2022fine}.
\begin{assumption}
\label{assumption:unbiased}
    The stochastic gradient with no compression error $\nabla L(x;\w)$ is an unbiased estimate of $\nabla \ell(\w)$ with bounded variances $\sigma^2$.
\end{assumption}
\begin{remark}
    Assumption \ref{assumption:unbiased} is the standard unbiased stochastic-gradient assumption due to data sampling, which is orthogonal to the unbiased compressed-gradient assumption criticized in limitation \textbf{(L1)}. \ours removes the latter, while retaining the former as in most SGD analyses.
\end{remark}

\begin{assumption}
\label{assumption:compressor}
    For the compressor $\mathcal{C}$, there exist constants $\omega_F,\omega_B\in[0,1)$ such that:
    \begin{align*}
        \mathbb{E}\left[\left\Vert {x}-\mathcal{C}({x})\right\Vert^2\middle|{x}\right]\leq\begin{cases}\omega_F^2\left\Vert {x}\right\Vert^2,\hspace{0.5mm}\text{forward propagation},\\\omega_B^2\left\Vert {x}\right\Vert^2,\hspace{0.5mm}\text{backward propagation}.\end{cases}
    \end{align*}
    \vspace{-1.5mm}
\end{assumption}

Compressors satisfying the above assumption are referred to as contractive compressors. In general, more aggressive compression leads to greater information distortion, corresponding to a larger \(\omega\). This assumption applies to numerous compressors, including top-\(K\) and low-rank projection \citep{alistarh2017qsgd,alistarh2018convergence,vogels2019powersgd,beznosikov2023biased,stich2018sparsified}, and is widely adopted in  communication-efficient algorithms \citep{koloskova2019decentralized,richtarik2021ef21,fatkhullin2024momentum,wang2023cocktailsgd,huang2022lower}.

The assumption below is critical for lazy sampling:
\begin{assumption}
\label{assumption:sample}
    For each $x_1,x_2\sim\mathcal{D}$, there exists $\varphi>0$ such that $\mathbb{E}_{x_1,x_2\in\mathcal{D}}\left[\left\Vert x_1-x_2\right\Vert^2\right]\leq\varphi^2$.
\end{assumption}
It is important to note that Assumption \ref{assumption:sample} is not overly restrictive for most optimization and learning tasks. Indeed, Assumption \ref{assumption:sample} holds for \textbf{all finite datasets}. Moreover, it is likely to be satisfied even for infinite datasets, particularly \textbf{when normalization techniques are applied}.

\begin{table*}[t]
\caption{ The score ($\uparrow$) of all the tasks for GLUE benchmark with communication compression algorithms.}
\centering
\label{table: fine-tuning result_GLUE}
\setlength{\tabcolsep}{4pt}
\begin{threeparttable}
{
\begin{tabular}{l|cccccccc|c}
\toprule
\textbf{Algorithms}           & \textbf{MNLI} & \textbf{SST-2} & \textbf{MRPC} & \textbf{CoLA} & \textbf{QNLI} & \textbf{QQP} & \textbf{RTE} & \textbf{STS-B} & \textbf{Avg} \\ \midrule
No comp.             & \textbf{90.01} & \textbf{96.41} & 88.25 & \textbf{64.69} & \textbf{94.59} & \textbf{92.03} & 82.35 & \textbf{92.52} & 87.61 \\ \midrule
EF21                 & 89.28 & 94.79 & 87.00 & 62.97 & 93.87 & 92.01 & 81.62 & 91.30 & 86.61 \\
\oursfc                & 89.94 & 96.06 & \textbf{90.50} & 63.98 & 94.04 & 91.89 & \textbf{84.19} & 91.71 & \textbf{87.79} \\ \bottomrule
\end{tabular}}
\end{threeparttable}
\end{table*}

\subsection{\ours convergence}

\label{sec:convergence}
Firstly, we present the convergence result of \oursfc is as follows.
\begin{lemma}
\label{thm:convergence_clapping}
    Suppose $x^{(t)}$ is the sampled data at iteration $t$, if we let $p_2=1$ and $p_3=\cdots=p_T=p$ as a constant, then for \oursfclast, under Assumptions \ref{assumption:smoothness}–\ref{assumption:compressor} the following holds.
\begin{equation}
\label{equation:error analysis_clapping_fc_1}
\begin{aligned}
   \dfrac{1}{T}\sum_{t=1}^T\mathbb{E}\hspace{-0.5mm}\left[\hspace{-0.5mm}\left\Vert\nabla \ell(\w^{(t)})\right\Vert^2\right]
    \hspace{-0.5mm}\lesssim_{T,p,m}&\dfrac{1}{\gamma T}\hspace{-0.5mm}+\hspace{-0.5mm}\dfrac{1}{m T}\hspace{-0.5mm}+\hspace{-0.5mm}\dfrac{1}{T}\sum_{t=1}^T\mathbb{E}\hspace{-0.5mm}\left[\hspace{-0.5mm}\left\Vert  x^{(t+1)}\hspace{-0.5mm}-\hspace{-0.5mm} x^{(t)}\right\Vert^2\right]\hspace{-0.5mm}+\hspace{-0.5mm}\sigma^2\dfrac{(2\hspace{-0.5mm}-\hspace{-0.5mm}p)m\hspace{-0.5mm}-\hspace{-0.5mm}(1\hspace{-0.5mm}-\hspace{-0.5mm}p)m^2}{1\hspace{-0.5mm}-\hspace{-0.5mm}(1\hspace{-0.5mm}-\hspace{-0.5mm}p)(1\hspace{-0.5mm}-\hspace{-0.5mm}m)^2}\\
    &\hspace{0.5cm}+\hspace{-0.5mm}\dfrac{1}{T}\left(\dfrac{1}{m^2}\hspace{-0.5mm}-\hspace{-0.5mm}\dfrac{1}{\gamma^2}\right)\sum_{e=1}^E\sum_{t=1}^{T}\mathbb{E}\left[\left\Vert w^{(t+1)}_e\hspace{-0.5mm}-\hspace{-0.5mm}w^{(t)}_e\right\Vert^2\right]. 
\end{aligned}
\end{equation}
\end{lemma}
In inequality \eqref{equation:error analysis_clapping_fc_1}, the term \(\sum_{t=1}^T \mathbb{E}[\|x^{(t+1)} - x^{(t)}\|^2]\) arises from error feedback across iterations. If different data samples are selected at iterations \( t \) and \( t+1 \), the term \( y_e^{(t+1)} - \tilde{y}_e^{(t)} \) in Algorithm \ref{alg:forward} introduces an error due to the discrepancy between \( x^{(t+1)} \) and \( x^{(t)} \). When \(\sum_{t=1}^T \mathbb{E}[\|x^{(t+1)} - x^{(t)}\|^2] = \mathcal{O}(T)\), Algorithm \ref{alg:clapping} converges to an \(\mathcal{O}(1)\) bias, as indicated by \eqref{equation:error analysis_clapping_fc_1}. This highlights the necessity to introduce lazy sampling to mitigate this term.

As outlined in Algorithm~\ref{alg:lazy_sampling}, under a constant lazy-sampling rate $p_t\equiv p$ we have $x^{(t+1)}=x^{(t)}$ w.p.~$1-p$, and hence $\sum_{t=1}^T \mathbb{E}\!\left[\|x^{(t+1)}-x^{(t)}\|^2\right]=O(Tp)$. With this, we can present the convergence of \oursfc as follow:

\begin{theorem}
\label{thm:convergence_clapping_fc}
    Suppose $x^{(t)}$ is the sampled data at iteration $t$, there exist properly chosen constant step sizes $\gamma$, a momentum coefficient $m$, and a lazy-sampling coefficient $p$ such that, for \oursfclast, under Assumptions \ref{assumption:smoothness}–\ref{assumption:sample} the following holds (see formal version in Theorem~\ref{lemma: convergence_fc_e}).
\begin{equation}
\begin{aligned}
\label{equation:convergence_clapping_fc}
    &\dfrac{1}{T}\sum_{t=1}^T\mathbb{E}\left[\left\Vert\nabla \ell(\w^{(t)})\right\Vert^2\right]
    \lesssim\dfrac{\sigma^{\frac{4}{3}}}{T^{\frac{1}{3}}(1-\omega_B)^{\frac{4(E-1)}{3}}(1-\omega_F)^{\frac{4(E-1)}{3}}}+\dfrac{\delta}{T},  
\end{aligned}
\end{equation}
where $\delta$ is a constant only depends on $\omega_B,\omega_F,E$ as
\begin{equation}
\begin{aligned}
\label{equation:impact of compression}
\delta\lesssim\dfrac{1}{(1-\omega_F)^{E-1}(1-\omega_B)^{E-1}}
&+\dfrac{\omega_F^2+\omega_B}{(1-\omega_F)^2(1-\omega_B)^{2(E-1)}}+\dfrac{1}{(1-\omega_F)^{2(E-2)-1}}.
\end{aligned}
\end{equation}
\end{theorem}

Meanwhile, the convergence rate of \oursfu is as follows.
\begin{theorem}
\label{thm:convergence_clapping_fu}
    Suppose $x^{(t)}$ is the sampled data at iteration $t$, there exist properly chosen constant step sizes \(\gamma\), a momentum coefficient \(m\), and lazy sampling coefficient $p$ such that, for \oursfulast, under Assumptions \ref{assumption:smoothness}–\ref{assumption:compressor} the following holds (see formal version in Theorem~\ref{lemma: convergence_fu_e}).
\begin{equation}
\begin{aligned}
\label{equation:convergence_clapping_fu}
    &\dfrac{1}{T}\sum_{t=1}^T\mathbb{E}\left[\left\Vert\nabla \ell(\w^{(t)})\right\Vert^2\right]
    \lesssim\dfrac{\sigma}{\sqrt{T}}+\dfrac{1}{T(1-\omega_B)^{E-1}(1-\omega_F)^{E-1}}.
\end{aligned}
\end{equation}
\end{theorem}

As illustrated in inequality \eqref{equation:convergence_clapping_fc} and \eqref{equation:convergence_clapping_fu}, \oursfc and \oursfu can achieve an asymptotic convergence rate of $\mathcal{O}(1/\sqrt[3]{T})$ and $\mathcal{O}(1/\sqrt{T})$, respectively. This complements prior pipeline-parallel compression analyses that either require sample-wise activation storage (e.g., AQ-SGD) or impose extra unbiasedness assumptions on the compression-induced gradient estimator.

\textbf{Trade-off in the selection of $p_t$.} The lazy-sampling coefficient \(p_t\) controls the trade-off between sample refresh and sample reuse. A smaller \(p_t\) gives error feedback more consecutive updates on the same sample and, for \texttt{Clapping}-\textbf{FU}, also reduces the frequency of raw refresh communications. However, \(p_t\) cannot be arbitrarily small. For \texttt{Clapping}-\textbf{FU}, the proof uses a fixed positive refresh probability \(p=p_0\in(0,1]\) independent of \(T\); the degenerate case \(p=0\) is not covered. When \(p\to0\), the repeated-sample variance term remains of order \(\Theta(\sigma^2)\), so the proof no longer yields the same stationarity guarantee. For \texttt{Clapping}-\textbf{FC}, the \(\sigma^2\) term in \eqref{equation:error analysis_clapping_fc_1} equals \(\mathcal{O}(1)\) when \(p\to0\), leading to non-convergence. Moreover, an excessively small \(p_t\) may hurt generalization, since some samples can be over-used while others are rarely seen.

In practice, \(p_t\) should be treated as a refresh-frequency parameter rather than a delicate constant. A simple rule is to start from a moderate value such as \(p_t\in[0.4,0.5]\). If training becomes unstable or validation quality degrades, one can increase \(p_t\) to refresh data more often; if training is stable but communication remains the bottleneck, one can decrease \(p_t\) to allow more sample reuse and stronger error feedback. Existing sweeps in Appendix~\ref{app:experiment} support that moderate values are usually preferable to very small or very large refresh probabilities; see Sec.~\ref{section:experiment} for further details.

\begin{table*}[t]
\caption{The evaluation accuracy ($\uparrow$) of different communication compression algorithms when fine-tuning different models by Wikitext with Top-5\% compressor.}
\centering
\label{table: fine-tuning result_wiki}
\renewcommand{\arraystretch}{1.15}
\begin{threeparttable}
{
\begin{tabular}{c|c|cccccc}
\toprule
\multirow{2}{*}{\textbf{Model}}           & \multirow{2}{*}{\textbf{No comp.}}  & \multirow{2}{*}{\textbf{Direct comp.}} & \multirow{2}{*}{\textbf{EF21}} & \multirow{2}{*}{\textbf{AQ-SGD}} & \multicolumn{3}{c}{\textbf{\oursfc}}  \\ \cline{6-8}
          &  &  &  &  &   $p=0.3$ & $p=0.4$ & $p=0.5$ \\ \midrule
LLaMA-2 7B                 & 0.5948 & 0.5473 & 0.5678 & 0.5696 & \textbf{0.5960} & 0.5920  & 0.5877  \\
LLaMA-3 8B                & 0.5991 & 0.5671 & 0.5677 & 0.5683 & 0.5865 & \textbf{0.5969}& 0.5887  \\ \bottomrule
\end{tabular}}
\end{threeparttable}
\end{table*}

\textbf{Impact of compression error accumulation. }
Evidently, more compression entails more accumulated error and results in slower convergence. It is also noteworthy that \oursfu requires $\mathcal{O}\left(1/\varepsilon^2 + 1/\left(\varepsilon (1-\omega_B)^{E-1} (1-\omega_F)^{E-1}\right)\right)$ iterations to approach an $\varepsilon$-stationary point. Hence, the impact of compression in our proposed method can be asymptotically neglected, as this term dominates the convergence rate---a finding consistent with the result in \citep{fatkhullin2024momentum}. Moreover, the exponential accumulation of compression error with the increase in the number of workers $E$ has also been documented in prior analyses \cite{ramasinghe2025protocol,kong2026convergence}. However, the adverse effect of communication compression persists in the algorithms presented by \citep{wang2022fine}. Ablation studies on the training performance and the number of workers $E$ can be found in Appendix \ref{appendix: finetune wiki}.

\textbf{Convergence with Adam.} We provide the convergence analysis for both \oursfc and \oursfu with the Adam optimizer in Appendix \ref{appendix: convergence_clapping_adam}. It can be shown that \ours with Adam shares the same convergence rate as that with momentum SGD. Such a convergence guarantee illustrates that \ours is suitable for LLM pre-training and fine-tuning tasks, which significantly extends the applicability of our proposed algorithm.

\textbf{Convergence in large batch scenario. }We provide the convergence analysis of \ours with the batch size \( B > 1 \) in Appendix \ref{Appendix: Convergence of ours with large batch}. It can be shown that \oursfc achieves a convergence rate of \(\mathcal{O}(1/\sqrt[3]{BT})\) and \oursfu achieves an convergence rate of \(\mathcal{O}(1/\sqrt{BT})\), with the impact of compression error remaining similar to the case when \( B = 1 \). Additionally, the selection of the lazy sampling coefficient can also extended to the large-batch scenario.

\textbf{Convergence and error propagation in multi-worker scenario. }Eq. \eqref{equation:convergence_clapping_fc} and \eqref{equation:convergence_clapping_fu} illustrate the convergence result of Algorithm \ref{alg:clapping} with multiple workers with the analysis of the error accumulation in the compression. And we also present a detailed analysis of error propagation in Appendix \ref{appendix: Error propagation analysis}. Such a result remedies the shortage of \citep{wang2022fine} in multi-worker compression. See Appendix \ref{appendix:aq-sgd} for more detail. {Furthermore, we provide an empirical analysis of error propagation across multiple workers; refer to Appendix \ref{appendix: finetune wiki} for experimental details and results.}

\section{Experiments}

\label{section:experiment}
We present experiments to validate the performance of \texttt{Clapping}. Unless otherwise specified, \texttt{Clapping} refers to \texttt{Clapping}-\textbf{FC} in this section, since FC is the harsher always-compress endpoint and gives the lower fixed-budget payload. The FU variant is complementary: it uses raw refreshes to reset error-feedback states and can improve optimization stability or time-to-target performance; we therefore report FU in the communication-bound pre-training study where this distinction is most relevant. Additional experimental details, extended results, and supplementary experiments on synthetic logistic regression, ResNet-18/CIFAR-10, and LLM fine-tuning/pre-training are provided in Appendix~\ref{app:experiment}.

\textbf{Fine-tuning on GLUE benchmark. }We fine-tune pre-trained RoBERTa-large \citep{liu2019roberta} on the GLUE benchmark \citep{wang2018glue} for 10 epochs using communication compression algorithms, including \ours and direct compression with EF21 \citep{richtarik2021ef21}, and compare them with uncompressed fine-tuning on two NVIDIA A800 GPUs. A TopK compressor retains 30\% of elements at the network midpoint. As shown in Table \ref{table: fine-tuning result_GLUE}, \ours outperforms EF21 in most tasks and achieves the highest average score.

\begin{figure}[t]
    \centering
    \includegraphics[width=0.49\textwidth]{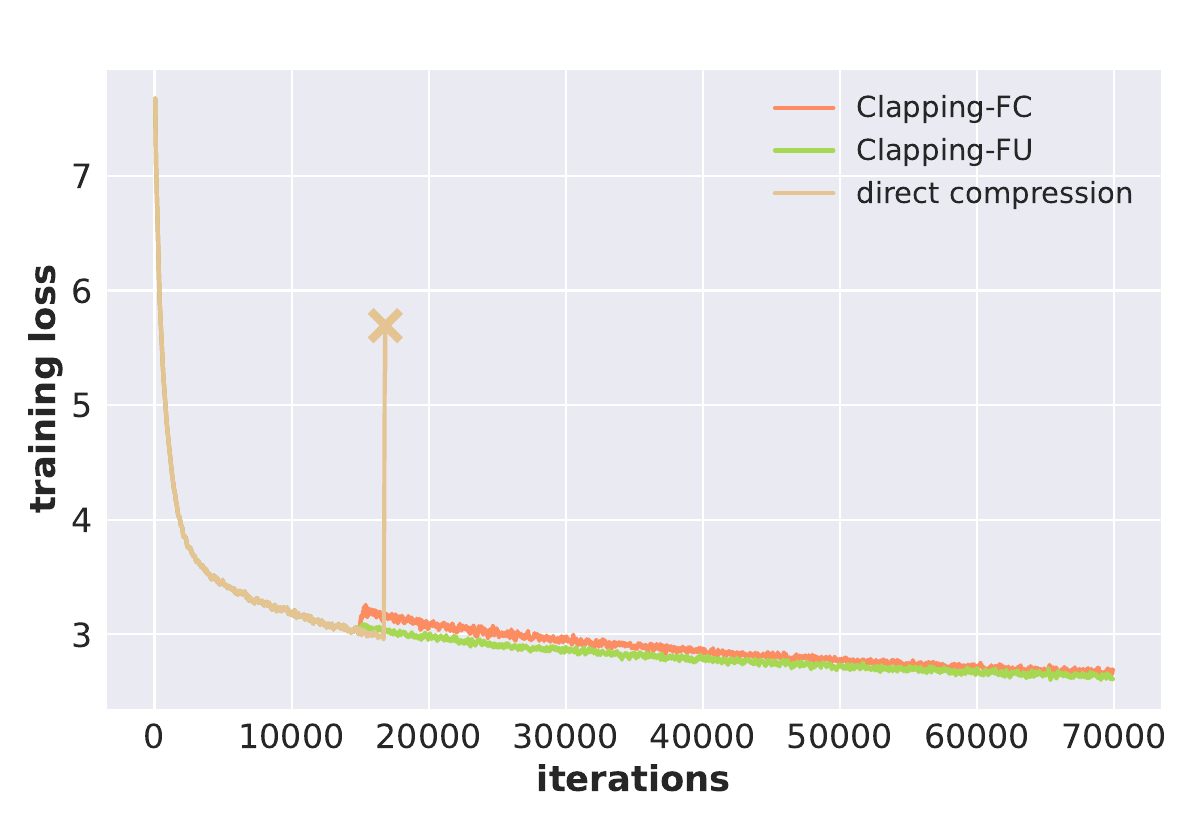}
    \hfill
    \includegraphics[width=0.49\textwidth]{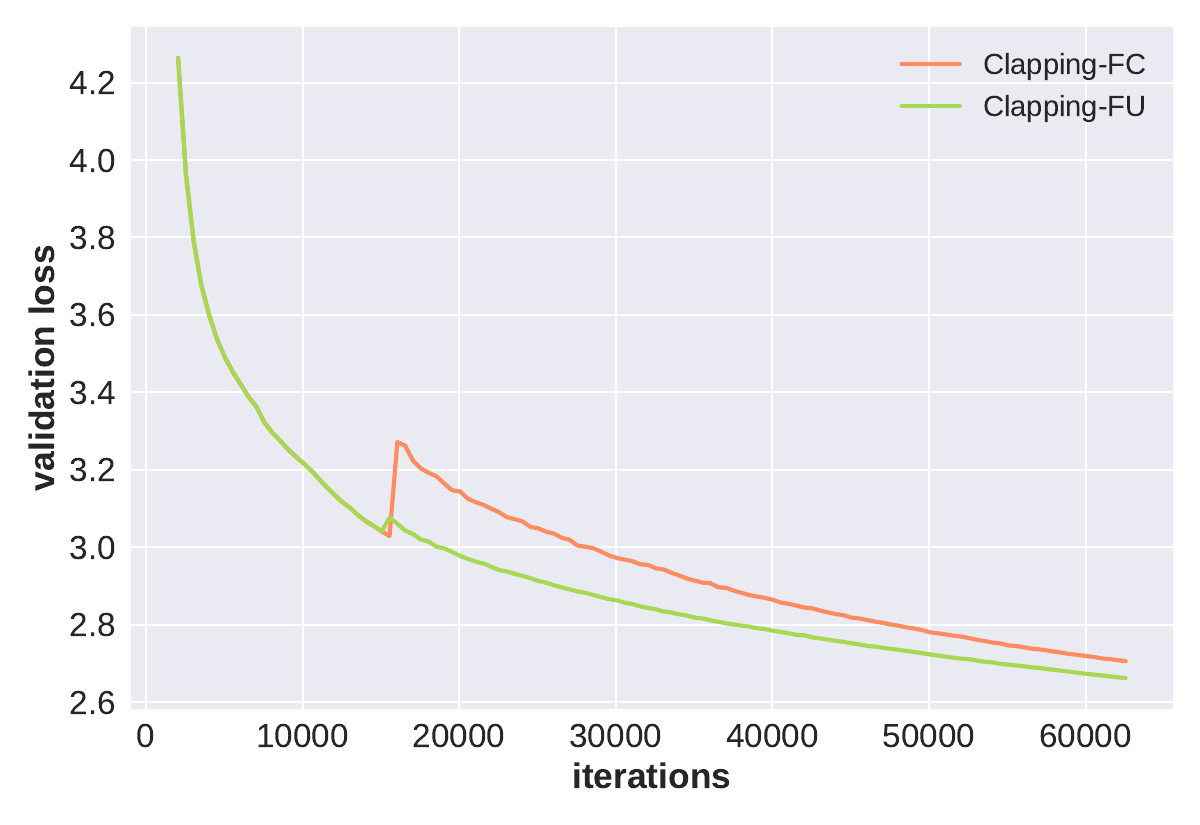}
    \caption{\small The training loss (left) and evaluation loss (right) with various communication compression algorithms for pre-training LLaMA-2 1B. Here the pre-training task under direct compression fails to converge.}
    \label{fig: LLaMA2 training-appendix1}
\end{figure}

\textbf{Fine-tuning on Wikitext. }We fine-tune pre-trained LLaMA-2-7B \citep{touvron2023LLaMA} and LLaMA-3-8B \citep{grattafiori2024llama} models on the Wikitext-2 dataset \citep{merity2016pointer}. The model is split at the midpoint, and a Top-K compressor \citep{wangni2018gradient} retains 5\% of elements. Table \ref{table: fine-tuning result_wiki} shows that \ours outperforms other algorithms, including direct compression, EF21, and AQ-SGD. By tuning $p$, \oursfc achieves $\approx95\%$ communication reduction with no more than a $0.5$ percentage-point {absolute} drop in evaluation accuracy in these fine-tuning tasks.


\textbf{Pre-training LLaMA-2 1B under a communication-bound pipeline.}
We pre-train a LLaMA-2 1B model \citep{touvron2023LLaMA} on C4-en \citep{raffel2020exploring} with pipeline parallelism under a fixed bandwidth cap (100\,MB/s) and activation/gradient compression (Top-$K$ sparsification with quantization). Following \citep{zhao2024galore,kong2025cr}'s setup (see Appendix \ref{appendix: pre-training 1B} for more details), Figure~\ref{fig: LLaMA2 training-appendix1} reports the training loss and evaluation loss \emph{versus training iterations}: direct compression fails to converge, while both \oursfu and \oursfc remain stable and track the uncompressed reference trajectory. This highlights Clapping as an \emph{optimization-stability} mechanism for pipeline-parallel communication compression in communication-bound regimes; additional implementation details are provided in Appendix~\ref{appendix: pre-training 1B}.


\section{Conclusions}
We propose \texttt{Clapping}, a communication-compression framework for pipeline-parallel distributed optimization that combines error feedback with lazy sampling. Under standard smoothness assumptions and a contractive compressor model, we establish sublinear stationarity rates for two variants: \oursfc achieves an $O(T^{-1/3})$ rate and \oursfu achieves an $O(T^{-1/2})$ rate, up to an additional term that explicitly captures compression-error accumulation with the pipeline depth $E$. Compared with prior pipeline-parallel compression methods, \ours removes the need for sample-wise activation storage and does not rely on the unbiased compressed-gradient assumption.

\section*{Acknowledgments}
This work is funded by the National Natural Science Foundation of China (No. 92370121, 12288101), and the National Key Research and Development Program of China (No. 2024YFA1012902). This research is also supported by the AI for Science Institute, Beijing, China. Part of the computational resources were supported by the Centerfor Intelligent Computing and Song-Shan Lake HPC Cen-ter (SSL-HPC) in Great Bay University, Dongguan, China.

{
\small

\bibliography{reference}
\bibliographystyle{abbrv}
}

\newpage

\appendix
\begin{center}
\Large
 \textbf{Appendix}
\end{center}
\vspace{-2mm}
\tableofcontents

\vspace{3mm}

\newpage

\section{More related works}
\subsection{Direct compression}
\label{sec:direct-compression}
To reduce communication overhead in pipeline-parallel optimization, we can directly compress the activations and their corresponding gradients, significantly reducing their size \citep{evans2020jpeg,fu2020don}. Let \(\mathcal{C}(\cdot)\) denote a compressor. The compressed pipeline-parallel SGD follows the same forward-backward procedure as described in Section~\ref{sec:pp-SGD}, with several core operations slightly modified:

\begin{itemize}[leftmargin=*]
    \item \textbf{Forward}: Each worker $e$ computes activation $y^{(t)}_e =a_e(\tilde{y}^{(t)}_{e-1},w^{(t)}_e)$ and compresses $\tilde{y}_e^{(t)} = \mathcal{C}(y^{(t)}_e)$ following the forward order $e=1,\cdots, E$. We let $y^{(t)}_0 = x^{(t)}$. 
    \item \textbf{Backward}: Each worker $e$ computes activation gradient $v^{(t)}_{e-1} = \nabla_1 a_e(\tilde{y}^{(t)}_{e-1},w^{(t)}_e)\tran \tilde{v}^{(t)}_{e}$ and compresses $\tilde{v}^{(t)}_{e-1} = \mathcal{C}(v^{(t)}_{e-1})$ following the backward order $e=E,\cdots, 2$. 
    \item \textbf{Update}: Each worker $e$ computes weight gradient $u^{(t)}_{e} = \nabla_2 a_e(\tilde{y}^{(t)}_{e-1},w^{(t)}_e)\tran \tilde{v}^{(t)}_{e}$ following the backward order $e=E,\cdots, 1$. We update parameter $w_e^{(t+1)} = w_e^{(t)} - \gamma u_e^{(t)}$.
\end{itemize}
Pipeline-parallel SGD with direct compression is also illustrated in Fig.~\ref{fig:pipeline}. Compression errors introduce unique challenges to pipeline-parallel SGD, which differ significantly from those encountered in data-parallel SGD.

\textbf{Unbiased compressor leads to biased gradient.} Direct compression in pipeline-parallel optimization naturally results in biased gradient estimates. Suppose $\mathcal{C}(\cdot)$ is a unbiased compressor such that $\mathbb{E}[\tilde{y}_e] = y_e$, it cannot be guaranteed that $\mathbb{E}[\nabla_j a_e(\tilde{y}_{e-1},w_e)] = \nabla_j a_e({y}_{e-1},w_e)\hspace{0.5mm}(j=1,2)$ due to the composite and non-linear structure of the activation operator $a_e(y_{e-1},w_e)$, leading to biased activation gradient and weight gradient estimates. However, \citep{evans2020jpeg} and \citep{fu2020don} impose the impractical assumption of an unbiased compression-induced gradient error to establish convergence guarantees. Without this assumption, direct compression cannot achieve convergence to stationary solutions. 

\textbf{Error propagation.}
Compression errors in pipeline-parallel SGD propagate in both the forward and backward processes. To illustrate this, we consider a linear neural network mapping \( a_e(y_{e-1}, w_e) = W_e y_{e-1} \), where \( W_e \) is the weight matrix reshaped from \( w_e \). Let \( \tilde{v}_e = v_e + \epsilon_e \), where \( \epsilon_e \) represents the error incurred during the compression $\tilde{v}_e = \mathcal{C}(v_e)$. From the backward step in direct compression, it holds that

\begin{align}
v_1 = W_1^\top \left( W_2^\top \left( \cdots \left( W_E^\top v_E + \epsilon_E \right)+ \cdots \right) + \epsilon_2 \right) + \epsilon_1. \nonumber 
\end{align}

It is evident that the innermost error \( \epsilon_E \) propagates through layers and can be significantly amplified by \( W_1^\top  \cdots W_{E-1}^\top \). This error propagation leads to a complex entanglement between the true gradient and the compressed one, severely impairing the performance and stability of the optimization. 

\subsection{Other related works}
\label{appendix:related work}
\textbf{Communication compression in data parallelism. }Communication compression has demonstrated significant efficacy in data-parallel distributed optimization \citep{xu2020compressed, wang2023communication}, with two core strategies underpinning its success: \textit{sparsification} and \textit{quantization}. Classical sparsification methods include Top-K \citep{wangni2018gradient, alistarh2018convergence} and Rand-K \citep{stich2018local, beznosikov2023biased}, while quantization techniques encompass Sign-SGD \citep{seide20141, bernstein2018signsgd}, TurnGrad \citep{wen2017terngrad}, and natural compression \citep{horvoth2022natural}. However, compression inevitably introduces information distortion, which can hinder convergence rates or even lead to non-convergence. To mitigate these challenges, a variety of advanced techniques have been developed, including error feedback \citep{stich2018sparsified, tang2019doublesqueeze, richtarik2021ef21, fatkhullin2024momentum}, hybrid compression \citep{wang2023cocktailsgd}, and multiple-step compression \citep{huang2022lower, he2023lower}. Furthermore, \citep{markov2023quantized} introduced weight compression as a complementary approach.
Despite these significant advancements, none of these results have been directly extended to pipeline-parallel distributed optimization.


\textbf{Activation compression. }A closely related technique is activation compression, which aims to reduce memory costs during LLM pre-training and fine-tuning \citep{jiang2018sketchml, jiang2022back, jin2021novel, evans2021ac, evans2020jpeg, fu2020don, yu2023compressing}. In communication compression for pipeline-parallel optimization, activations must also be compressed to minimize communication overhead. The key distinction lies in the fact that, in memory-efficient settings, activations are computed precisely during forward propagation, with their compressed copies stored for backward propagation \citep{jiang2022back}, thereby introducing no additional errors during forward propagation. However, in pipeline-parallel communication compression, forward propagation is inherently error-prone, as compressed activations are immediately used for subsequent computations \citep{wang2022fine}, leading to error accumulation across both forward and backward propagation. This poses a significant challenge for algorithm design and theoretical analysis. Recently, \cite{he2025tah} also presented an effective activation compression algorithm in pipeline parallelism. While favorable compression performance is achieved on the LLMs fine-tuning tasks, it still lacks theoretical convergence guarantees.
\section{Proof of the convergence rate}
\label{appendix:proof}
In this section, we prove the convergence results for Algorithm \ref{alg:clapping}. The proof separates optimization error from communication error. We first reduce convergence to the error of the evaluated gradient, then control the forward and backward compression errors, and finally combine these bounds for \oursfu and \oursfc.

\subsection{Notations and proof roadmap}
To begin with, we use different superscripts to represent variables with/without compression. Suppose the communication compression is taken in the $t$-th iteration, then we have:
\begin{itemize}[leftmargin=1em]
    \item \textbf{Variables with a hat} (like $\hat{y}_e^{(t)},\hat{v}_e^{(t)}$): The variable is obtained from standard momentum SGD algorithm \textbf{without any communication compression}. 
    \item \textbf{Variables with a tilde} (like $\tilde{y}_e^{(t)},\tilde{v}_e^{(t)}$): The variable which is used for the gradient evaluation during the $t$-th iteration. If the communication compression is taken during the $t$-th iteration, it denotes variabl obtained from Algorithm \ref{alg:lazy_sampling} - \ref{alg:clapping} and \textbf{has already been compressed}. 
    \item \textbf{Variables without any additional superscript} (like ${y}_e^{(t)},{v}_e^{(t)}$): The variable is obtained from Algorithm \ref{alg:lazy_sampling} - \ref{alg:clapping} and \textbf{has not been compressed} during the $t$-th iteration if the communication compression is taken.
\end{itemize}
If the communication compression is not taken in this iteration, same variables with different superscripts are all denote the variable without communication, such as ${y}_e^{(t)}=\hat{y}_e^{(t)}=\tilde{y}_e^{(t)}$.

Moreover, in the $t$-th iteration, we firstly denote $\mathcal{F}_e^{(t)}$ and $\mathcal{G}_e^{(t)}$ for $e=1,2,\cdots,E-1$ as follows:
\begin{itemize}[leftmargin=1em]
    \item $\mathcal{F}_e^{(t)}$: the filtration before the communication from machine $e$ to machine $e+1$ in \textbf{forward propagation},
    \item $\mathcal{G}_e^{(t)}$: the filtration before the communication from machine $e+1$ to machine $e$ in \textbf{backward propagation}.
\end{itemize}
We also define $\mathcal{F}^{(t)}$ as the $\sigma$-algebra generated by all randomness \emph{up to (but excluding)} the fresh sampling decision at iteration $t$. In particular, $x^{(t)}$ is independent of $\mathcal{F}^{(t)}$ whenever a fresh sample is drawn.

\paragraph{Proof roadmap.}
The proof has four steps. 

\textbf{Step 1. }Lemma~\ref{lem:descent} gives a descent bound and reduces the problem to controlling the evaluated-gradient error. 

\textbf{Step 2. }Lemma~\ref{lemma3} decomposes this error into stochastic noise, forward activation error, and backward adjoint error. 

\textbf{Step 3. }The technical lemmas in the next two subsections bound these forward/backward errors by parameter changes and, for \oursfclast, the sample-change term \(\|x^{(t+1)}-x^{(t)}\|^2\). 

\textbf{Step 4. }The results are combined in the master convergence bound to prove the rates for \oursfu and \oursfclast.

The key difference between the two variants is how they handle refresh iterations. \oursfu sends exact activations/adjoints on refresh iterations, so the corresponding error-feedback states are reset. \oursfc compresses every iteration, so a sample refresh can leave an additional term depending on \(x^{(t+1)}-x^{(t)}\).

\subsection{A Master Descent Bound}
We start with the optimization-level part of the proof. The next lemma is a standard smoothness descent bound: it shows that convergence follows once we can control the evaluated-gradient error produced by compression and lazy sampling.
\begin{lemma}[Descent Lemma]
\label{lem:descent}
    Suppose Assumption \ref{assumption:smoothness}-\ref{assumption:unbiased} holds and $\gamma\leq\dfrac{1}{2L_{\nabla \ell}}$, then in Algorithm \ref{alg:clapping} we have:
    \begin{equation}
    \label{descent lemma}
    \begin{aligned}
        &\dfrac{1}{T}\sum_{t=1}^T\mathbb{E}\left[\left\Vert\nabla \ell(\w^{(t)})\right\Vert^2\right]\\
        \leq&\dfrac{2}{\gamma T}\mathbb{E}\left[\ell(\w^{(1)})\hspace{-0.3mm}-\hspace{-0.3mm}\inf_{\w}\ell(W)\right]\hspace{-0.5mm}-\hspace{-0.5mm}\dfrac{1}{2\gamma^2 T}\sum_{t=1}^T\mathbb{E}\left[\left\Vert \w^{(t+1)}\hspace{-0.5mm}-\hspace{-0.3mm}\w^{(t)}\right\Vert^2\right]\hspace{-0.5mm}+\hspace{-0.5mm}\dfrac{1}{T}\sum_{t=1}^T\sum_{e=1}^E\mathbb{E}\left[\left\Vert\tilde{u}^{(t)}_e\hspace{-0.5mm}-\hspace{-0.4mm}\nabla_e\ell(\w^{(t)})\right\Vert^2\right].
    \end{aligned}
    \end{equation}
\begin{proof}
    According to the result in \citep{fatkhullin2024momentum}, we can get from Assumption \ref{assumption:unbiased} that:
    \begin{align}
        \ell(\w^{(t+1)})\leq&\ell(\w^{(t)})-\dfrac{\gamma}{2}\left\Vert\nabla \ell(\w^{(t)})\right\Vert^2\nonumber\\
        &-\left(\dfrac{1}{2\gamma}-\dfrac{L_{\nabla \ell}}{2}\right)\left\Vert \w^{(t+1)}-\w^{(t)}\right\Vert^2+\dfrac{\gamma}{2}\sum_{e=1}^E\left\Vert\tilde{u}^{(t)}_e-\nabla_e\ell(\w^{(t)})\right\Vert^2.
    \end{align}
    Then, as $\gamma\leq\dfrac{1}{2L_{\nabla \ell}}$, we have:
    \begin{equation}
    \begin{aligned}
        \ell(\w^{(t+1)})-\inf_{\w}l(\w)\leq&\ell(\w^{(t)})-\inf_{\w}\ell(\w)-\dfrac{\gamma}{2}\left\Vert\nabla \ell(\w^{(t)})\right\Vert^2\\
        &-\dfrac{1}{4\gamma}\left\Vert \w^{(t+1)}-\w^{(t)}\right\Vert^2+\dfrac{\gamma}{2}\sum_{e=1}^E\left\Vert\tilde{u}^{(t)}_e-\nabla_e\ell(\w^{(t)})\right\Vert^2.
    \end{aligned}
    \end{equation}
    Taking expectation and summation on both sides over $t=1,2,\cdots,T$, we can get:
    \begin{equation}
    \begin{aligned}
        \mathbb{E}\left[\ell(\w^{(T+1)})-\inf_{\w}\ell(\w)\right]\leq&\mathbb{E}\left[\ell(\w^{(1)})-\inf_{\w}l(\w)\right]-\dfrac{\gamma}{2}\sum_{t=1}^T\mathbb{E}\left[\left\Vert\nabla \ell(\w^{(t)})\right\Vert^2\right]\\
        &-\dfrac{1}{4\gamma}\sum_{t=1}^T\mathbb{E}\left[\left\Vert \w^{(t+1)}-\w^{(t)}\right\Vert^2\right]+\dfrac{\gamma}{2}\sum_{t=1}^T\sum_{e=1}^E\mathbb{E}\left[\left\Vert\tilde{u}^{(t)}_e-\nabla_e\ell(\w^{(t)})\right\Vert^2\right].
    \end{aligned}
    \end{equation}
    Finally, we have:
    \begin{equation*}
    \begin{aligned}
        &\dfrac{1}{T}\sum_{t=1}^T\mathbb{E}\left[\left\Vert\nabla \ell(\w^{(t)})\right\Vert^2\right]\\
        \leq&\dfrac{2}{\gamma T}\mathbb{E}\left[\ell(\w^{(1)})\hspace{-0.3mm}-\hspace{-0.3mm}\inf_{\w}\ell(W)\right]\hspace{-0.5mm}-\hspace{-0.5mm}\dfrac{1}{2\gamma^2 T}\sum_{t=1}^T\mathbb{E}\left[\left\Vert \w^{(t+1)}\hspace{-0.5mm}-\hspace{-0.3mm}\w^{(t)}\right\Vert^2\right]\hspace{-0.5mm}+\hspace{-0.5mm}\dfrac{1}{T}\sum_{t=1}^T\sum_{e=1}^E\mathbb{E}\left[\left\Vert\tilde{u}^{(t)}_e\hspace{-0.5mm}-\hspace{-0.4mm}\nabla_e\ell(\w^{(t)})\right\Vert^2\right].
    \end{aligned}
    \end{equation*}
\end{proof}
\end{lemma}

\begin{remark}
\label{remark:a_e^0}
   It is noteworthy that one can get $||\hat{v}_e^{(t)}||$ is bounded for $e=1,\cdots, E-1$ from Assumption \ref{assumption:smoothness} and the definition of $ \hat{v}_e^{(t)}$. If we let $a_e^\circ( v,\by,w_e)= \nabla_2a_e(\by,w_e)^{\tran}v$, and $L_{\nabla a}'=\mathcal{O}\left(\max\{L_{\nabla a},L_{\nabla a}||\hat{v}_e^{(t)}||\}\right)$,
then there exist $L_{\nabla a}^\circ=\mathcal{O}(C_a)$ such that
\begin{equation}
\begin{aligned}
\label{estimation of a_e^0}
    &\left\Vert a_e^\circ(\hat{v}_e^{(t)},\by,w_e)-a_e^\circ( v',\by',w_e')\right\Vert^2=\left\Vert\nabla_2a_e(\by,w_e)^{\tran}\hat{v}_e^{(t)}- \nabla_2a_e(\by',w_e')^{\tran}v'\right\Vert^2\\
    \leq&2\left\Vert\nabla_2a_e(\by,w_e)^{\tran}\hat{v}_e^{(t)}- \nabla_2a_e(\by',w_e')^{\tran}\hat{v}_e^{(t)}\right\Vert^2+2\left\Vert\nabla_2a_e(\by',w_e')^{\tran}\hat{v}_e^{(t)}- \nabla_2a_e(\by',w_e')^{\tran}v'\right\Vert^2\\
    \leq&(L_{\nabla a}')^2\left(\left\Vert \by-\by'\right\Vert^2+\left\Vert w_e-w_e'\right\Vert^2\right)+(L_{\nabla a}^\circ)^2\left\Vert  \hat{v}_e^{(t)}- v'\right\Vert^2,
\end{aligned} 
\end{equation}
where the last inequality is due to the Lipschitz continuous of $\nabla_2a_e$, the boundness of $\hat{v}_e^{(t)}$ and the boundness of $\nabla_2a_e$.
\end{remark}
Lemma~\ref{lem:descent} shows that it remains to control the evaluated-gradient
error. All subsequent lemmas are devoted to bounding this term under FC and FU.

\subsubsection{Evaluated-gradient error}
The descent bound leaves the term \(\|\tilde{u}^{(t)}_e-\nabla_e\ell(\w^{(t)})\|^2\). The next lemma decomposes this term while keeping the effects of momentum and lazy sampling explicit. After this step, the remaining work is to bound the forward activation error and the backward adjoint error.
\begin{lemma}
\label{lemma3}
Suppose Assumption \ref{assumption:smoothness} and \ref{assumption:unbiased} hold, and let $m_1=\cdots=m_T=m_{T+1}=m$ as well as $p_3=\cdots=p_T=p_{T+1}=p$. Moreover, we set $p_2=1$. Then, for all $t=2,\cdots,T+1$ we have:
    \begin{equation}
    \label{u-nablal}
        \begin{aligned}
            &\sum_{e=1}^E\sum_{t=1}^{T+1}\mathbb{E}\left[\left\Vert\tilde{u}^{(t)}_e-\nabla_e\ell(\w^{(t)})\right\Vert^2\right]\\
            \leq&32L_{\nabla \ell}^2\left(\dfrac{p+m}{m^2(1-(1-p)(1-\frac{m}{2}))}+\dfrac{1}{m^2}\right)\sum_{e=1}^E\sum_{t=1}^{T}\mathbb{E}\left[\left\Vert w^{(t+1)}_e-w^{(t)}_e\right\Vert^2\right]\\
            &+8(L_{\nabla a}^\circ)^2\sum_{e=1}^{E-1}\sum_{t=2}^{T+1}\mathbb{E}\left[\left\Vert\tilde{ v}_e^{(t)}- \hat{v}_e^{(t)}\right\Vert^2\right]+8(L_{\nabla a}')^2\sum_{e=1}^{E-1}\sum_{t=2}^{T+1}\mathbb{E}\left[\left\Vert\tilde{\by}_e^{(t)}-\hat{y}_e^{(t)}\right\Vert^2\right]\\
            &+4T\sigma^2\dfrac{(2-p)m-(1-p)m^2}{1-(1-p)(1-m)^2}+\dfrac{3}{m}\sum_{e=1}^E\mathbb{E}\left[\left\Vert\tilde{u}^{(1)}_e-\nabla_e\ell(\w^{(1)})\right\Vert^2\right].
        \end{aligned}
    \end{equation}
\begin{proof}
In the beginning, we need to claim that we set $p_1=p_2=1$ so that $\psi(t)\ge 2$ and the term $\tilde u^{(\psi(t)-1)}$ is always well-defined. This only affects $O(1)$ initial iterations and does not change the asymptotic rates.

For $t=2,3,\cdots,T+1$, we denote $\psi(t)$ as the last moment in which the sample is randomly obtained with $\mathcal{D}$ as of the $t$-th iteration. Specially,
    \begin{align*}
        \psi(t):=\max_{\tau\in\mathbb{S}_t}\tau,\quad\text{where }\mathbb{S}_t:=\{\tau=2,3,\cdots,t|\text{sampling randomly at iteration }\tau\}.
    \end{align*}
    Then, with the fact that the $p_2=1$, it holds for $\tau=2,\cdots,t$ that $\text{Pr}(\psi(t)=\tau)=\begin{cases}(1-p)^{t-2}, \text{ if }\tau=2 \\ p(1-p)^{t-\tau}, \text{ else}. \end{cases}$

    For $e=1,2,\cdots,E-1$ and $t=2,3,\cdots,T+1$, the error between the evaluated gradient and the true gradient satisfies:
    {\small
    \begin{equation}
    \label{u-nabla,new}
    \begin{aligned}
        &\tilde{u}^{(t)}_e-\nabla_e\ell(\w^{(t)})\\
        =&\sum_{\tau=\psi(t)}^tm(1-m)^{t-\tau}\nabla_2a_e(\tilde{y}_{e-1}^{(\tau)},w_e^{(\tau)})^{\tran}\tilde{v}_e^{(\tau)}+(1-m)^{t+1-\psi(t)}\tilde{u}_e^{(\psi(t)-1)}-\nabla_e\ell(\w^{(t)})\\
        =&\underbrace{\sum_{\tau=\psi(t)}^tm(1-m)^{t-\tau}\left(\nabla_2a_e(\tilde{y}_{e-1}^{(\tau)},w_e^{(\tau)})^{\tran}\tilde{v}_e^{(\tau)}-\nabla_2a_e(\hat{y}_{e-1}^{(\tau)},w_e^{(\tau)})^{\tran}\hat{v}_e^{(\tau)}\right)}_{:=\Xi_{e,1}}\\
        &+\sum_{\tau=\psi(t)}^tm(1-m)^{t-\tau}\left(\nabla_2a_e(\hat{y}_{e-1}^{(\tau)},w_e^{(\tau)})^{\tran}\hat{v}_e^{(\tau)}-\nabla_e\ell(\w^{(\tau)})\right)\\
        &+\underbrace{\sum_{\tau=\psi(t)}^tm(1-m)^{t-\tau}\left(\nabla_e\ell(\w^{(\tau)})-\nabla_e\ell(\w^{(t)})\right)}_{:=\Xi_{e,2}}+\underbrace{(1-m)^{t+1-\psi(t)}\left(\nabla_e\ell(\w^{(\psi(t)-1)})-\nabla_e\ell(\w^{(t)})\right)}_{:=\Xi_{e,3}}\\
        &+(1-m)^{t+1-\psi(t)}\left(\tilde{u}_e^{(\psi(t)-1)}-\nabla_e\ell(\w^{(\psi(t)-1)})\right),
    \end{aligned}
    \end{equation}}
    where the first equation is from the momentum update rule. Moreover, we use $\Xi_{e,1},\Xi_{e,2},\Xi_{e,3}$ to denote some complex terms, which have been shown in Eq. \eqref{u-nabla,new}.

    Additionally, we denote $\mathcal{F}^{(t)}$ as the filtration before the $t$-th iteration. Thus, $\tilde{u}_e^{(\psi(t)-1)}-\nabla_e\ell(\w^{(\psi(t)-1)})$ is measureable with respect to $\mathcal{F}^{(\psi(t))}$ for any $e=1,2,\cdots,E$. Moreover, the sampling process at the iteration $\psi(t)$ is \textbf{independent} with respect to $\mathcal{F}^{(\psi(t))}$. Thus, the non-compressed term $\nabla_2a_e(\hat{y}_{e-1}^{(\tau)},w_e^{(\tau)})\hat{v}^{(\tau)}_e$ is an unbiased estimation of the gradient $\nabla_e\ell(\w^{(\tau)})$ with bounded variance according to Assumption \ref{assumption:unbiased}. Thus, taking the $\ell_2$-norm and conditional expectation with respect to $\mathcal{F}^{(\psi(t))}$ on both sides of Eq. \eqref{u-nabla,new}, we can obtain:
    \begin{equation}
    \label{u-nabla,new2}
    \begin{aligned}
        &\mathbb{E}\left[\left\|\tilde{u}_e^{(t)}-\nabla_e\ell(\w^{(t)})\right\|^2\middle|\mathcal{F}^{(\psi(t))}\right]\\
        =&\mathbb{E}\left[\left\|\sum_{\tau=\psi(t)}^tm(1-m)^{t-\tau}\left(\nabla_2a_e(\hat{y}_{e-1}^{(\tau)},w_e^{(\tau)})^{\tran}\hat{v}_e^{(\tau)}-\nabla_e\ell(\w^{\tau})\right)\right\|^2\middle|\mathcal{F}^{(\psi(t))}\right]\\
        &+\mathbb{E}\left[\left\|(1-m)^{t+1-\psi(t)}\left(\tilde{u}_e^{(\psi(t)-1)}-\nabla_e\ell(\w^{(\psi(t)-1)})\right)+\Xi_{e,1}+\Xi_{e,2}+\Xi_{e,3}\right\|^2\middle|\mathcal{F}^{(\psi(t))}\right]\\
        &+2\mathbb{E}\Bigg[\Bigg\langle\sum_{\tau=\psi(t)}^tm(1-m)^{t-\tau}\left(\nabla_2a_e(\hat{y}_{e-1}^{(\tau)},w_e^{(\tau)})^{\tran}\hat{v}_e^{(\tau)}-\nabla_e\ell(\w^{\tau})\right),\\
        &\quad\quad\quad\quad(1-m)^{t+1-\psi(t)}\left(\tilde{u}_e^{(\psi(t)-1)}-\nabla_e\ell(\w^{(\psi(t)-1)})\right)+\Xi_{e,1}+\Xi_{e,2}+\Xi_{e,3}\Bigg\rangle\Bigg|\mathcal{F}^{(\psi(t))}\Bigg]\\
        \leq&2\mathbb{E}\left[\left\|\sum_{\tau=\psi(t)}^tm(1-m)^{t-\tau}\left(\nabla_2a_e(\hat{y}_{e-1}^{(\tau)},w_e^{(\tau)})^{\tran}\hat{v}_e^{(\tau)}-\nabla_e\ell(\w^{\tau})\right)\right\|^2\middle|\mathcal{F}^{(\psi(t))}\right]\\
        &+\mathbb{E}\left[\left\|\Xi_{e,1}+\Xi_{e,2}+\Xi_{e,3}\right\|^2\middle|\mathcal{F}^{(\psi(t))}\right]\\
        &+\mathbb{E}\left[\left\|(1-m)^{t+1-\psi(t)}\left(\tilde{u}_e^{(\psi(t)-1)}-\nabla_e\ell(\w^{(\psi(t)-1)})\right)+\Xi_{e,1}+\Xi_{e,2}+\Xi_{e,3}\right\|^2\middle|\mathcal{F}^{(\psi(t))}\right],
    \end{aligned}
    \end{equation}
    where the inequality is due to Cauchy-Schwarz inequality and Assumption \ref{assumption:unbiased}.

    For the second term of the right-hand-side of Eq. \eqref{u-nabla,new2}, it holds that:
    \begin{equation}
    \label{eq123}
    \begin{aligned}
        &\mathbb{E}\left[\left\|\Xi_{e,1}+\Xi_{e,2}+\Xi_{e,3}\right\|^2\middle|\mathcal{F}^{(\psi(t))}\right]\\
        \leq&2\sum_{\tau=\psi(t)}^tm(1-m)^{t-\tau}\mathbb{E}\left[\left\|\nabla_2a_e(\tilde{y}_{e-1}^{(\tau)},w_e^{(\tau)})^{\tran}\tilde{v}_e^{(\tau)}-\nabla_2a_e(\hat{y}_{e-1}^{(\tau)},w_e^{(\tau)})^{\tran}\hat{v}_e^{(\tau)}\right\|^2\middle|\mathcal{F}^{(\psi(t))}\right]\\
        &+2\sum_{\tau=\psi(t)}^tm(1-m)^{t-\tau}\mathbb{E}\left[\left\|\nabla_e\ell(\w^{(\tau)})-\nabla_e\ell(\w^{(t)})\right\|^2\middle|\mathcal{F}^{(\psi(t))}\right]\\
        &+(1-m)^{t+1-\psi(t)}\mathbb{E}\left[\left\|\nabla_e\ell(\w^{(\psi(t)-1)})-\nabla_e\ell(\w^{(t)})\right\|^2\middle|\mathcal{F}^{(\psi(t))}\right],
    \end{aligned}
    \end{equation}
    where the inequality holds is due to the convexity of the $\ell_2$-norm.

    Moreover, for the last term, it also holds that:
    {\small
    \begin{equation}
    \label{psi_i+123}
    \begin{aligned}
        &\mathbb{E}\left[\left\|(1-m)^{t+1-\psi(t)}\left(\tilde{u}_e^{(\psi(t)-1)}-\nabla_e\ell(\w^{(\psi(t)-1)})\right)+\Xi_{e,1}+\Xi_{e,2}+\Xi_{e,3}\right\|^2\middle|\mathcal{F}^{(\psi(t))}\right]\\
        \leq&2\sum_{\tau=\psi(t)}^tm(1-m)^{t-\tau}\mathbb{E}\left[\left\|\nabla_2a_e(\tilde{y}_{e-1}^{(\tau)},w_e^{(\tau)})^{\tran}\tilde{v}_e^{(\tau)}-\nabla_2a_e(\hat{y}_{e-1}^{(\tau)},w_e^{(\tau)})^{\tran}\hat{v}_e^{(\tau)}\right\|^2\middle|\mathcal{F}^{(\psi(t))}\right]\\
        &+2\sum_{\tau=\psi(t)}^tm(1-m)^{t-\tau}\mathbb{E}\left[\left\|\nabla_e\ell(\w^{(\tau)})-\nabla_e\ell(\w^{(t)})\right\|^2\middle|\mathcal{F}^{(\psi(t))}\right]\\
        &+(1-m)^{t+1-\psi(t)}\mathbb{E}\left[\left\|\left(\tilde{u}_e^{(\psi(t)-1)}-\nabla_e\ell(\w^{(\psi(t)-1)})\right)+\left(\nabla_e\ell(\w^{(\psi(t)-1)})-\nabla_e\ell(\w^{(t)})\right)\right\|^2\middle|\mathcal{F}^{(\psi(t))}\right].
    \end{aligned}
    \end{equation}}
    Form Young's inequality, for any $u>0$ the last term of \eqref{psi_i+123} holds that:
    \begin{equation}
    \label{psi_i+123_2}
    \begin{aligned}
        &\mathbb{E}\left[\left\|\left(\tilde{u}_e^{(\psi(t)-1)}-\nabla_e\ell(\w^{(\psi(t)-1)})\right)+\left(\nabla_e\ell(\w^{(\psi(t)-1)})-\nabla_e\ell(\w^{(t)})\right)\right\|^2\middle|\mathcal{F}^{(\psi(t))}\right]\\
        \leq&(1+u)\mathbb{E}\left[\left\|\tilde{u}_e^{(\psi(t)-1)}-\nabla_e\ell(\w^{(\psi(t)-1)})\right\|^2\middle|\mathcal{F}^{(\psi(t))}\right]\\
        &+\left(1+\dfrac{1}{u}\right)\mathbb{E}\left[\left\|\nabla_e\ell(\w^{(\psi(t)-1)})-\nabla_e\ell(\w^{(t)})\right\|^2\middle|\mathcal{F}^{(\psi(t))}\right].
    \end{aligned}
    \end{equation}
    If we take $u=\dfrac{\left(1-\frac{m}{2}\right)^{t+1-\psi(t)}-\left(1-m\right)^{t+1-\psi(t)}}{\left(1-m\right)^{t+1-\psi(t)}}$, then $1+u=\dfrac{\left(1-\frac{m}{2}\right)^{t+1-\psi(t)}}{\left(1-m\right)^{t+1-\psi(t)}}$ and 
    \begin{align*}
        1+\dfrac{1}{u}=&1+\dfrac{\left(1-m\right)^{t+1-\psi(t)}}{\left(1-\frac{m}{2}\right)^{t+1-\psi(t)}-\left(1-m\right)^{t+1-\psi(t)}}=1+\dfrac{\left(1-m\right)^{t+1-\psi(t)}}{\dfrac{m}{2}\left(\sum_{j=0}^{t-\psi(t)}\left(1-\frac{m}{2}\right)^{j}\left(1-m\right)^{t-\psi(t)-j}\right)}\\
        \leq&1+\dfrac{2}{m(t+1-\psi(t))}.
    \end{align*}

    Substituting \eqref{psi_i+123_2} into Eq. \eqref{psi_i+123} and taking the value of $u$, we can obtain that:
    {\small
    \begin{equation}
    \label{psi_i+123_3}
    \begin{aligned}
        &\mathbb{E}\left[\left\|(1-m)^{t+1-\psi(t)}\left(\tilde{u}_e^{(\psi(t)-1)}-\nabla_e\ell(\w^{(\psi(t)-1)})\right)+\Xi_{e,1}+\Xi_{e,2}+\Xi_{e,3}\right\|^2\middle|\mathcal{F}^{(\psi(t))}\right]\\
        \leq&2\sum_{\tau=\psi(t)}^tm(1-m)^{t-\tau}\mathbb{E}\left[\left\|\nabla_2a_e(\tilde{y}_{e-1}^{(\tau)},w_e^{(\tau)})^{\tran}\tilde{v}_e^{(\tau)}-\nabla_2a_e(\hat{y}_{e-1}^{(\tau)},w_e^{(\tau)})^{\tran}\hat{v}_e^{(\tau)}\right\|^2\middle|\mathcal{F}^{(\psi(t))}\right]\\
        &+2\sum_{\tau=\psi(t)}^tm(1-m)^{t-\tau}\mathbb{E}\left[\left\|\nabla_e\ell(\w^{(\tau)})-\nabla_e\ell(\w^{(t)})\right\|^2\middle|\mathcal{F}^{(\psi(t))}\right]\\
        &+\left(1-\frac{m}{2}\right)^{t+1-\psi(t)}\mathbb{E}\left[\left\|\tilde{u}_e^{(\psi(t)-1)}-\nabla_e\ell(\w^{(\psi(t)-1)})\right\|^2\middle|\mathcal{F}^{(\psi(t))}\right]\\
        &+(1-m)^{t+1-\psi(t)}\left(1+\dfrac{2}{m(t+1-\psi(t))}\right)\mathbb{E}\left[\left\|\nabla_e\ell(\w^{(\psi(t)-1)})-\nabla_e\ell(\w^{(t)})\right\|^2\middle|\mathcal{F}^{(\psi(t))}\right].
    \end{aligned}
    \end{equation}}

    Substituting Eq. \eqref{eq123} and Eq. \eqref{psi_i+123_3} into Eq.\eqref{u-nabla,new2}, we can obtain that:
    \begin{equation}
    \label{u-nabla,new3}
    \begin{aligned}
        &\mathbb{E}\left[\left\|\tilde{u}_e^{(t)}-\nabla_e\ell(\w^{(t)})\right\|^2\middle|\mathcal{F}^{(\psi(t))}\right]\\
        \leq&2\mathbb{E}\left[\left\|\sum_{\tau=\psi(t)}^tm(1-m)^{t-\tau}\left(\nabla_2a_e(\hat{y}_{e-1}^{(\tau)},w_e^{(\tau)})^{\tran}\hat{v}_e^{(\tau)}-\nabla_e\ell(\w^{\tau})\right)\right\|^2\middle|\mathcal{F}^{(\psi(t))}\right]\\
        &+4\sum_{\tau=\psi(t)}^tm(1-m)^{t-\tau}\mathbb{E}\left[\left\|\nabla_2a_e(\tilde{y}_{e-1}^{(\tau)},w_e^{(\tau)})^{\tran}\tilde{v}_e^{(\tau)}-\nabla_2a_e(\hat{y}_{e-1}^{(\tau)},w_e^{(\tau)})^{\tran}\hat{v}_e^{(\tau)}\right\|^2\middle|\mathcal{F}^{(\psi(t))}\right]\\
        &+4\sum_{\tau=\psi(t)}^tm(1-m)^{t-\tau}\mathbb{E}\left[\left\|\nabla_e\ell(\w^{(\tau)})-\nabla_e\ell(\w^{(t)})\right\|^2\middle|\mathcal{F}^{(\psi(t))}\right]\\
        &+\left(1-\frac{m}{2}\right)^{t+1-\psi(t)}\mathbb{E}\left[\left\|\tilde{u}_e^{(\psi(t)-1)}-\nabla_e\ell(\w^{(\psi(t)-1)})\right\|^2\middle|\mathcal{F}^{(\psi(t))}\right]\\
        &+2(1-m)^{t+1-\psi(t)}\left(1+\dfrac{1}{m(t+1-\psi(t))}\right)\mathbb{E}\left[\left\|\nabla_e\ell(\w^{(\psi(t)-1)})-\nabla_e\ell(\w^{(t)})\right\|^2\middle|\mathcal{F}^{(\psi(t))}\right].
    \end{aligned}
    \end{equation}

    With Assumption \ref{assumption:unbiased}, we can get:
    {\small
    \begin{equation}
    \begin{aligned}
        &\mathbb{E}\left[\sum_{e=1}^E\left\|\sum_{\tau=\psi(t)}^tm(1-m)^{t-\tau}\left(\nabla_2a_e(\hat{y}_{e-1}^{(\tau)},w_e^{(\tau)})^{\tran}\hat{v}_e^{(\tau)}-\nabla_e\ell(\w^{\tau})\right)\right\|^2\middle|\mathcal{F}^{(\psi(t))}\right]
        \leq\left(\sum_{\tau=\psi(t)}^{t}m(1-m)^{t-\tau}\right)^2\sigma^2.
    \end{aligned}
    \end{equation}}

    We denote $\Lambda^{(\tau)}$ as the vector resulting from connecting $\nabla_2a_e(\tilde{y}_{e-1}^{(\tau)},w_e^{(\tau)})^{\tran}\tilde{v}_e^{(\tau)}-\nabla_2a_e(\hat{y}_{e-1}^{(\tau)},w_e^{(\tau)})^{\tran}\hat{v}_e^{(\tau)}$ for all $e=1,2,\cdots,E$ end-to-end. Then, taking the summation on both sides of Eq. \eqref{u-nabla,new3}, it holds that:
    \begin{equation}
    \label{u-nabla,new4}
    \begin{aligned}
        &\mathbb{E}\left[\sum_{e=1}^E\left\|\tilde{u}_e^{(t)}-\nabla_e\ell(\w^{(t)})\right\|^2\middle|\mathcal{F}^{(\psi(t))}\right]\\
        \leq&2\left(\sum_{\tau=\psi(t)}^{t}m(1-m)^{t-\tau}\right)^2\sigma^2+4\sum_{\tau=\psi(t)}^tm(1-m)^{t-\tau}\mathbb{E}\left[\left\|\Lambda^{(\tau)}\right\|^2\middle|\mathcal{F}^{(\psi(t))}\right]\\
        &+4\sum_{\tau=\psi(t)}^tm(1-m)^{t-\tau}\mathbb{E}\left[\left\|\nabla\ell(\w^{(\tau)})-\nabla\ell(\w^{(t)})\right\|^2\middle|\mathcal{F}^{(\psi(t))}\right]\\
        &+\left(1-\frac{m}{2}\right)^{t+1-\psi(t)}\mathbb{E}\left[\sum_{e=1}^E\left\|\tilde{u}_e^{(\psi(t)-1)}-\nabla_e\ell(\w^{(\psi(t)-1)})\right\|^2\middle|\mathcal{F}^{(\psi(t))}\right]\\
        &+2(1-m)^{t+1-\psi(t)}\left(1+\dfrac{1}{m(t+1-\psi(t))}\right)\mathbb{E}\left[\left\|\nabla\ell(\w^{(\psi(t)-1)})-\nabla\ell(\w^{(t)})\right\|^2\middle|\mathcal{F}^{(\psi(t))}\right].
    \end{aligned}
    \end{equation}

    Then, taking the conditional expectation with respect to $\psi(t)$ on both sides of \eqref{u-nabla,new4}, it holds that:
    \begin{equation}
    \label{u-nabla,new5}
    \begin{aligned}
        &\mathbb{E}\left[\sum_{e=1}^E\left\|\tilde{u}_e^{(t)}-\nabla_e\ell(\w^{(t)})\right\|^2\middle|\psi(t)\right]\\
        \leq&2\mathbb{E}\left[\left(\sum_{\tau=\psi(t)}^{t}m(1-m)^{t-\tau}\right)^2\sigma^2\middle|\psi(t)\right]+4\sum_{\tau=\psi(t)}^tm(1-m)^{t-\tau}\mathbb{E}\left[\left\|\Lambda^{(\tau)}\right\|^2\middle|\psi(t)\right]\\
        &+4\sum_{\tau=\psi(t)}^tm(1-m)^{t-\tau}\mathbb{E}\left[\left\|\nabla\ell(\w^{(\tau)})-\nabla\ell(\w^{(t)})\right\|^2\middle|\psi(t)\right]\\
        &+\left(1-\frac{m}{2}\right)^{t+1-\psi(t)}\mathbb{E}\left[\sum_{e=1}^E\left\|\tilde{u}_e^{(\psi(t)-1)}-\nabla_e\ell(\w^{(\psi(t)-1)})\right\|^2\middle|\psi(t)\right]\\
        &+2(1-m)^{t+1-\psi(t)}\left(1+\dfrac{1}{m(t+1-\psi(t))}\right)\mathbb{E}\left[\left\|\nabla\ell(\w^{(\psi(t)-1)})-\nabla\ell(\w^{(t)})\right\|^2\middle|\psi(t)\right].
    \end{aligned}
    \end{equation}

    Furthermore, taking the expectation over $\psi(t)$, it holds that:
    \begin{equation}
    \label{u-nabla,new61}
    \begin{aligned}
        &\mathbb{E}\left[\sum_{e=1}^E\left\|\tilde{u}_e^{(t)}-\nabla_e\ell(\w^{(t)})\right\|^2\right]\\
        \leq&\sum_{\kappa=2}^t\text{Pr}(\psi(t)=\kappa)\left(1-\frac{m}{2}\right)^{t+1-\kappa}\mathbb{E}\left[\sum_{e=1}^E\left\|\tilde{u}_e^{(\kappa-1)}-\nabla_e\ell(\w^{(\kappa-1)})\right\|^2\right]\\
        &+4\sum_{\kappa=2}^t\text{Pr}(\psi(t)=\kappa)\sum_{\tau=\kappa}^tm(1-m)^{t-\tau}\mathbb{E}\left[\left\|\nabla\ell(\w^{(\tau)})-\nabla\ell(\w^{(t)})\right\|^2\right]\\
        &+2\sum_{\kappa=2}^t\text{Pr}(\psi(t)=\kappa)\left(1-(1-m)^{t-\kappa+1}\right)^2\sigma^2\hspace{-0.7mm}+\hspace{-0.7mm}4\sum_{\kappa=2}^t\text{Pr}(\psi(t)=\kappa)\sum_{\tau=\kappa}^tm(1-m)^{t-\tau}\mathbb{E}\left[\left\|\Lambda^{(\tau)}\right\|^2\right]\\
        &+2\sum_{\kappa=2}^t\text{Pr}(\psi(t)=\kappa)(1-m)^{t+1-\kappa}\left(1+\dfrac{1}{m(t+1-\kappa)}\right)\mathbb{E}\left[\left\|\nabla\ell(\w^{(\kappa-1)})-\nabla\ell(\w^{(t)})\right\|^2\right].
    \end{aligned}
    \end{equation}

    Thus,
    \begin{equation}
    \label{u-nabla,new6}
    \begin{aligned}
        &\mathbb{E}\left[\sum_{e=1}^E\left\|\tilde{u}_e^{(t)}-\nabla_e\ell(\w^{(t)})\right\|^2\right]\\
        \leq&\left(1-\frac{m}{2}\right)\sum_{\kappa=2}^t\text{Pr}(\psi(t)=\kappa)\left(1-\frac{m}{2}\right)^{t-\kappa}\mathbb{E}\left[\sum_{e=1}^E\left\|\tilde{u}_e^{(\kappa-1)}-\nabla_e\ell(\w^{(\kappa-1)})\right\|^2\right]\\
        &+4\sum_{\kappa=2}^t\left(\text{Pr}(\psi(t)=\kappa)(1-m)^{t+1-\kappa}\left(1+\dfrac{1}{m(t+1-\kappa)}\right)\right.\\
        &\quad\quad\left.+\sum_{\tau=2}^{\kappa-1}\text{Pr}(\psi(t)=\tau)m(1-m)^{t+1-\kappa}\right)\mathbb{E}\left[\left\|\nabla\ell(\w^{(\kappa-1)})-\nabla\ell(\w^{(t)})\right\|^2\right]\\
        &+2\sum_{\kappa=2}^t\hspace{-0.5mm}\text{Pr}(\psi(t)=\kappa)\hspace{-0.3mm}\left(1-(1-m)^{t-\kappa+1}\right)^2\hspace{-0.5mm}\sigma^2\hspace{-0.5mm}+\hspace{-0.5mm}4\sum_{\kappa=2}^t\hspace{-0.5mm}\left(\sum_{\tau=2}^{\kappa}\text{Pr}(\psi(t)=\tau)\right)\hspace{-0.5mm}m(1-m)^{t-\kappa}\mathbb{E}\hspace{-0.8mm}\left[\left\|\Lambda^{(\kappa)}\right\|^2\right].
    \end{aligned}
    \end{equation}

    Taking summatation over $t=2,\cdots,T+1$, we can obtain that:
    \begin{equation}
    \label{u-nabla,new7}
    \begin{aligned}
        &\sum_{t=2}^{T+1}\mathbb{E}\left[\sum_{e=1}^E\left\|\tilde{u}_e^{(t)}-\nabla_e\ell(\w^{(t)})\right\|^2\right]\\
        \leq&\left(1-\frac{m}{2}\right)\sum_{t=2}^{T+1}\sum_{\kappa=2}^t\text{Pr}(\psi(t)=\kappa)\left(1-\frac{m}{2}\right)^{t-\kappa}\mathbb{E}\left[\sum_{e=1}^E\left\|\tilde{u}_e^{(\kappa-1)}-\nabla_e\ell(\w^{(\kappa-1)})\right\|^2\right]\\
        &+4\sum_{t=2}^{T+1}\sum_{\kappa=2}^t\left(\text{Pr}(\psi(t)=\kappa)(1-m)^{t+1-\kappa}\left(1+\dfrac{1}{m(t+1-\kappa)}\right)\right.\\
        &\quad\quad\quad\quad\quad\quad\left.+\sum_{\tau=2}^{\kappa-1}\text{Pr}(\psi(t)=\tau)m(1-m)^{t+1-\kappa}\right)\mathbb{E}\left[\left\|\nabla\ell(\w^{(\kappa-1)})-\nabla\ell(\w^{(t)})\right\|^2\right]\\
        &+2\sum_{t=2}^{T+1}\sum_{\kappa=2}^t\text{Pr}(\psi(t)=\kappa)\left(1-(1-m)^{t-\kappa+1}\right)^2\sigma^2\\
        &+4\sum_{t=2}^{T+1}\sum_{\kappa=2}^t\left(\sum_{\tau=2}^{\kappa}\text{Pr}(\psi(t)=\tau)\right)m(1-m)^{t-\kappa}\mathbb{E}\left[\left\|\Lambda^{(\kappa)}\right\|^2\right].
    \end{aligned}
    \end{equation}

    Here we consider each term of Eq. \eqref{u-nabla,new7}. Firstly, we can obtain that:
    {\small
    \begin{equation}
    \label{huanxu_27}
    \begin{aligned}
        &\sum_{\kappa=2}^t(1-m)^{t+1-\kappa}\hspace{-1mm}\left(\text{Pr}(\psi(t)=\kappa)\hspace{-0.8mm}\left(1+\dfrac{1}{m(t+1-\kappa)}\right)\hspace{-1mm}+\hspace{-0.5mm}m\hspace{-0.3mm}\sum_{\tau=2}^{\kappa-1}\text{Pr}(\psi(t)=\tau)\hspace{-0.8mm}\right)\hspace{-0.3mm}\mathbb{E}\hspace{-0.8mm}\left[\hspace{-0.3mm}\left\|\nabla\ell(\w^{(\kappa-1)})\hspace{-0.5mm}-\hspace{-0.7mm}\nabla\ell(\w^{(t)})\hspace{-0.3mm}\right\|^2\right]\\
        \leq&\sum_{\kappa=2}^t(t+1-\kappa)(1-m)^{t+1-\kappa}\left(\text{Pr}(\psi(t)=\kappa)\left(1+\dfrac{1}{m(t+1-\kappa)}\right)+m\sum_{\tau=2}^{\kappa-1}\text{Pr}(\psi(t)=\tau)\right)\\
        &\quad\quad\quad\quad\quad\quad\quad\quad\quad\quad\quad\quad\quad\quad\quad\quad\quad\quad\quad\quad\quad\quad\quad\quad\quad\cdot\sum_{\tau=\kappa}^t\mathbb{E}\left[\left\|\nabla\ell(\w^{(\tau-1)})-\nabla\ell(\w^{(\tau)})\right\|^2\right]\\
        =&\sum_{\kappa=1}^{t-1}(1-m)^{t-\kappa}\left(\text{Pr}(\psi(t)=\kappa+1)(t-\kappa+\frac{1}{m})+m(t-\kappa)\sum_{\tau=2}^{\kappa}\text{Pr}(\psi(t)=\tau)\right)\\
        &\quad\quad\quad\quad\quad\quad\quad\quad\quad\quad\quad\quad\quad\quad\quad\quad\quad\quad\quad\quad\quad\quad\quad\quad\quad\cdot\sum_{\tau=\kappa}^{t-1}\mathbb{E}\left[\left\|\nabla\ell(\w^{(\tau)})-\nabla\ell(\w^{(\tau+1)})\right\|^2\right]\\
        =&\sum_{\kappa=1}^{t-1}\left[\sum_{\tau=1}^{\kappa}(1-m)^{t-\tau}\left(\text{Pr}(\psi(t)=\tau+1)(t-\tau+\frac{1}{m})+m(t-\tau)\sum_{\iota=2}^{\tau}\text{Pr}(\psi(t)=\iota)\right)\right]\\
        &\quad\quad\quad\quad\quad\quad\quad\quad\quad\quad\quad\quad\quad\quad\quad\quad\quad\quad\quad\quad\quad\quad\quad\quad\quad\cdot\mathbb{E}\left[\left\|\nabla\ell(\w^{(\kappa)})-\nabla\ell(\w^{(\kappa+1)})\right\|^2\right].
    \end{aligned}
    \end{equation}}

    let $s:=t-\tau$, then the coefficient of the term $\mathbb{E}\left[\left\|\nabla\ell(\w^{(\kappa)})-\nabla\ell(\w^{(\kappa+1)})\right\|^2\right]$ in Eq. \eqref{huanxu_27} holds that:
    {\small
    \begin{equation}
    \label{huanxu_28}
    \begin{aligned}
        &\sum_{\tau=1}^{\kappa}(1-m)^{t-\tau}\hspace{-0.6mm}\left(\hspace{-0.6mm}\text{Pr}(\psi(t)\hspace{-0.4mm}=\hspace{-0.4mm}\tau\hspace{-0.4mm}+\hspace{-0.4mm}1)(t\hspace{-0.4mm}-\hspace{-0.4mm}\tau\hspace{-0.4mm}+\hspace{-0.4mm}\frac{1}{m})\hspace{-0.4mm}+\hspace{-0.4mm}m(t\hspace{-0.4mm}-\hspace{-0.4mm}\tau)\sum_{\iota=2}^{\tau}\text{Pr}(\psi(t)\hspace{-0.4mm}=\hspace{-0.4mm}\iota)\hspace{-0.6mm}\right)\\
        =&\sum_{\tau=2}^{\kappa}(1-m)^{t-\tau}\hspace{-0.6mm}\left(\hspace{-0.6mm}\text{Pr}(\psi(t)\hspace{-0.4mm}=\hspace{-0.4mm}\tau\hspace{-0.4mm}+\hspace{-0.4mm}1)(t\hspace{-0.4mm}-\hspace{-0.4mm}\tau\hspace{-0.4mm}+\hspace{-0.4mm}\frac{1}{m})\hspace{-0.4mm}+\hspace{-0.4mm}m(t\hspace{-0.4mm}-\hspace{-0.4mm}\tau)\sum_{\iota=2}^{\tau}\text{Pr}(\psi(t)\hspace{-0.4mm}=\hspace{-0.4mm}\iota)\hspace{-0.6mm}\right)\hspace{-0.6mm}+\hspace{-0.4mm}(1\hspace{-0.4mm}-\hspace{-0.4mm}m)^{t-1}\hspace{-0.2mm}(1\hspace{-0.4mm}-\hspace{-0.4mm}p)^{t-2}(t\hspace{-0.4mm}-\hspace{-0.4mm}1\hspace{-0.4mm}+\hspace{-0.4mm}\frac{1}{m})\\
        =&\sum_{\tau=2}^{\kappa}(1-m)^{t-\tau}\hspace{-0.6mm}\left(\hspace{-0.6mm}p(1-p)^{t-\tau-1}(t\hspace{-0.4mm}-\hspace{-0.4mm}\tau\hspace{-0.4mm}+\hspace{-0.4mm}\frac{1}{m})\hspace{-0.4mm}+\hspace{-0.4mm}m(t\hspace{-0.4mm}-\hspace{-0.4mm}\tau)\hspace{-0.4mm}(1-p)^{t-\tau}\hspace{-0.6mm}\right)\hspace{-0.6mm}+\hspace{-0.4mm}(1\hspace{-0.4mm}-\hspace{-0.4mm}m)^{t-1}\hspace{-0.2mm}(1\hspace{-0.4mm}-\hspace{-0.4mm}p)^{t-2}(t\hspace{-0.4mm}-\hspace{-0.4mm}1\hspace{-0.4mm}+\hspace{-0.4mm}\frac{1}{m})\\
        =&\sum_{s=t-\kappa}^{t-2}(1-m)^{s}\left(\hspace{-0.6mm}p(1-p)^{s-1}(s\hspace{-0.4mm}+\hspace{-0.4mm}\frac{1}{m})\hspace{-0.4mm}+\hspace{-0.4mm}ms\hspace{-0.4mm}(1-p)^{s}\hspace{-0.6mm}\right)\hspace{-0.6mm}+\hspace{-0.4mm}(1\hspace{-0.4mm}-\hspace{-0.4mm}m)^{t-1}\hspace{-0.2mm}(1\hspace{-0.4mm}-\hspace{-0.4mm}p)^{t-2}(t\hspace{-0.4mm}-\hspace{-0.4mm}1\hspace{-0.4mm}+\hspace{-0.4mm}\frac{1}{m})\\
        =&\sum_{s=t-\kappa}^{t-2}\left[(1-m)^s(1-p)^{s-1}s(p+m(1-p))+\dfrac{p}{m}(1-m)^s(1-p)^{s-1}\right]\hspace{-0.6mm}+\hspace{-0.4mm}(1\hspace{-0.4mm}-\hspace{-0.4mm}m)^{t-1}\hspace{-0.2mm}(1\hspace{-0.4mm}-\hspace{-0.4mm}p)^{t-2}(t\hspace{-0.4mm}-\hspace{-0.4mm}1\hspace{-0.4mm}+\hspace{-0.4mm}\frac{1}{m}),
    \end{aligned}
    \end{equation}}
    where the first equation is due to the fact that $\sum_{\iota=2}^{\tau}\text{Pr}(\psi(t)=\iota)=(1-p)^{t-\tau}$ for $\tau=2,\cdots,T$.
    
    Substituting Eq. \eqref{huanxu_28} into Eq. \eqref{huanxu_27}, it holds that:
    Thus, we can obtain that:
    {\small
    \begin{equation}
    \begin{aligned}
        &\sum_{\kappa=2}^t(1-m)^{t+1-\kappa}\hspace{-1mm}\left(\text{Pr}(\psi(t)=\kappa)\hspace{-0.8mm}\left(1+\dfrac{1}{m(t+1-\kappa)}\right)\hspace{-1mm}+\hspace{-0.5mm}m\hspace{-0.3mm}\sum_{\tau=2}^{\kappa-1}\text{Pr}(\psi(t)=\tau)\hspace{-0.8mm}\right)\hspace{-0.3mm}\mathbb{E}\hspace{-0.8mm}\left[\hspace{-0.3mm}\left\|\nabla\ell(\w^{(\kappa-1)})\hspace{-0.5mm}-\hspace{-0.7mm}\nabla\ell(\w^{(t)})\hspace{-0.3mm}\right\|^2\right]\\
        \leq&\sum_{\kappa=1}^{t-1}\left[\sum_{s=t-\kappa}^{t-2}\hspace{-0.3mm}\left[\hspace{-0.3mm}(1-m)^s\hspace{-0.3mm}(1-p)^{s-1}\hspace{-0.3mm}s(p\hspace{-0.3mm}+\hspace{-0.3mm}m(1\hspace{-0.3mm}-\hspace{-0.3mm}p))\hspace{-0.5mm}+\hspace{-0.5mm}\dfrac{p}{m}(1-m)^s(1-p)^{s-1}\hspace{-0.4mm}\right]\hspace{-0.5mm}+\hspace{-0.5mm}(1-m)^{t-1}\hspace{-0.3mm}(1-p)^{t-2}\hspace{-0.3mm}(t\hspace{-0.4mm}-\hspace{-0.4mm}1\hspace{-0.4mm}+\hspace{-0.4mm}\frac{1}{m})\hspace{-0.5mm}\right]\\
        &\quad\quad\quad\quad\quad\quad\quad\quad\quad\quad\quad\quad\quad\quad\quad\quad\quad\quad\quad\quad\quad\quad\quad\quad\quad\cdot\mathbb{E}\left[\left\|\nabla\ell(\w^{(\kappa)})-\nabla\ell(\w^{(\kappa+1)})\right\|^2\right]\\
        \leq&(p+m(1-p))\sum_{\kappa=1}^{t-1}\left[\left(\sum_{\tau=t-\kappa}^{t-2}(1-m)^{\tau}(1-p)^{\tau-1}\tau\right)\mathbb{E}\left[\left\|\nabla\ell(\w^{(\kappa)})-\nabla\ell(\w^{(\kappa+1)})\right\|^2\right]\right]\\
        &+\dfrac{p}{m}\sum_{\kappa=1}^{t-1}\left[\left(\sum_{\tau=t-\kappa}^{t-2}(1-m)^{\tau}(1-p)^{\tau-1}\right)\mathbb{E}\left[\left\|\nabla\ell(\w^{(\kappa)})-\nabla\ell(\w^{(\kappa+1)})\right\|^2\right]\right]\\
        &+(1-m)^{t-1}(1-p)^{t-2}(t-1+\frac{1}{m})\sum_{\kappa=1}^{t-1}\mathbb{E}\left[\left\|\nabla\ell(\w^{(\kappa)})-\nabla\ell(\w^{(\kappa+1)})\right\|^2\right].
    \end{aligned}
    \end{equation}}
    Taking summation on both sides over $t=1,2,\cdots,T$ and use the fact that $p\leq m\leq1$, it holds that:
    {\small
    \begin{equation}
    \label{huanxu30}
    \begin{aligned}
        &\sum_{t=1}^T\Bigg[\sum_{\kappa=2}^t(1-m)^{t+1-\kappa}\left(\text{Pr}(\psi(t)=\kappa)\left(1+\dfrac{1}{m(t+1-\kappa)}\right)+m\sum_{\tau=2}^{\kappa-1}\text{Pr}(\psi(t)=\tau)\right)\\
        &\quad\quad\quad\quad\quad\quad\quad\quad\quad\quad\quad\quad\quad\quad\quad\quad\quad\quad\quad\quad\quad\quad\quad\quad\quad\cdot\mathbb{E}\left[\left\|\nabla\ell(\w^{(\kappa-1)})-\nabla\ell(\w^{(t)})\right\|^2\right]\Bigg]\\
        \leq&\sum_{t=1}^T\hspace{-0.5mm}\left(\hspace{-0.5mm}(p\hspace{-0.5mm}+\hspace{-0.5mm}m(1\hspace{-0.5mm}-\hspace{-0.5mm}p))\sum_{\tau=1}^{+\infty}(1\hspace{-0.5mm}-\hspace{-0.5mm}m)^{\tau}(1\hspace{-0.5mm}-\hspace{-0.5mm}p)^{\tau-1}\tau^2\hspace{-0.5mm}+\hspace{-0.5mm}\dfrac{p\hspace{-0.5mm}+\hspace{-0.5mm}m}{m}\hspace{-0.5mm}\sum_{\tau=1}^{+\infty}(1\hspace{-0.5mm}-\hspace{-0.5mm}m)^{\tau}(1\hspace{-0.5mm}-\hspace{-0.5mm}p)^{\tau-1}\tau\hspace{-0.5mm}+\hspace{-0.5mm}\dfrac{1}{m}\sum_{\tau=1}^{+\infty}(1\hspace{-0.5mm}-\hspace{-0.5mm}m)^{\tau}(1\hspace{-0.5mm}-\hspace{-0.5mm}p)^{\tau-1}\hspace{-0.5mm}\right)\\
        &\quad\quad\quad\quad\quad\quad\quad\quad\quad\quad\quad\quad\quad\quad\quad\quad\quad\quad\quad\quad\quad\quad\quad\quad\quad\cdot\mathbb{E}\left[\left\|\nabla\ell(\w^{(t)})-\nabla\ell(\w^{(t+1)})\right\|^2\right]\\
        \leq&\sum_{t=1}^T\left(\dfrac{(1-m)(1+(1-p)(1-m)}{(1-(1-p)(1-m))^2}+\dfrac{p+m}{m}\cdot\dfrac{1-m}{(1-(1-p)(1-m))^2}+\dfrac{1}{m}\cdot\dfrac{1-m}{1-(1-p)(1-m)}\right)\\
        &\quad\quad\quad\quad\quad\quad\quad\quad\quad\quad\quad\quad\quad\quad\quad\quad\quad\quad\quad\quad\quad\quad\quad\quad\quad\cdot\mathbb{E}\left[\left\|\nabla\ell(\w^{(t)})-\nabla\ell(\w^{(t+1)})\right\|^2\right]\\
        \leq&\left(\dfrac{4(p+m)}{m(1-(1-p)(1-m))^2}+\dfrac{4}{m(1-(1-p)(1-m))}\right)\sum_{t=1}^T\mathbb{E}\left[\left\|\nabla\ell(\w^{(t)})-\nabla\ell(\w^{(t+1)})\right\|^2\right].
    \end{aligned}
    \end{equation}}

    Moreover, it also holds that:
    \begin{equation}
    \label{huanxu31}
    \begin{aligned}
        &\sum_{\kappa=2}^t\text{Pr}(\psi(t)=\kappa)\left(1-(1-m)^{t-\kappa+1}\right)^2\\
        =&(1-p)^{t-2}\left(1-(1-m)^{t-1}\right)^2+\sum_{\kappa=3}^tp(1-p)^{t-\kappa}\left(1-(1-m)^{t-\kappa+1}\right)^2\\
        \leq&\sum_{s=t-2}^{+\infty}p(1-p)^s\left(1-(1-m)^{s+1}\right)^2+\sum_{\kappa=3}^tp(1-p)^{t-\kappa}\left(1-(1-m)^{t-\kappa+1}\right)^2\\
        =&\sum_{s=0}^{+\infty}p(1-p)^s\left(1-(1-m)^{s+1}\right)^2=\dfrac{(2-p)m^2-(1-p)m^3}{(1-(1-p)(1-m))(1-(1-p)(1-m)^2)}.
    \end{aligned}
    \end{equation}

    Finally, we can obtain that:
    \begin{equation}
    \label{huanxu32}
    \begin{aligned}
        &(1-\frac{m}{2})\sum_{t=2}^{T+1}\left(\sum_{\kappa=2}^t\text{Pr}(\psi(t)=\kappa)(1-\frac{m}{2})^{t-\kappa}\mathbb{E}\left[\sum_{e=1}^E\left\Vert\tilde{u}_e^{(\kappa-1)}-\nabla_e\ell(\w^{(\kappa-1)})\right\Vert^2\right]\right)\\
        =&(1-\frac{m}{2})\sum_{t=2}^{T+1}\left(\sum_{\kappa=t}^{T+1}\text{Pr}(\psi(\kappa)=t)(1-\frac{m}{2})^{\kappa-t}\mathbb{E}\left[\sum_{e=1}^E\left\Vert\tilde{u}_e^{(t-1)}-\nabla_e\ell(\w^{(t-1)})\right\Vert^2\right]\right).
    \end{aligned}
    \end{equation}
    We note that
    \begin{align*}
        \sum_{\kappa=t}^{T+1}(1-m)^{\kappa-t}(1-p)^{\kappa-t}\leq\dfrac{1}{1-(1-m)(1-p)}.
    \end{align*}

    For $t\geq3$, it holds that
    \begin{align*}
        \sum_{\kappa=t}^{T+1}\text{Pr}(\psi(\kappa)=t)(1-\frac{m}{2})^{\kappa-t}=p\sum_{\kappa=t}^{T+1}(1-p)^{\kappa-t}(1-\frac{m}{2})^{\kappa-t}\leq\dfrac{p}{1-(1-p)(1-\frac{m}{2})}.
    \end{align*}
    For $t=2$, it holds that
    \begin{align*}
        \sum_{\kappa=2}^{T+1}\text{Pr}(\psi(\kappa)=2)(1-\frac{m}{2})^{\kappa-2}=\sum_{\kappa=2}^{T+1}(1-p)^{\kappa-2}(1-\frac{m}{2})^{\kappa-2}\leq\dfrac{1}{1-(1-p)(1-\frac{m}{2})}.
    \end{align*}

    Combining Eq. \eqref{u-nabla,new7}, \eqref{huanxu30}, \eqref{huanxu31}, \eqref{huanxu32} together, we can obtain that:
    {\small
    \begin{equation}
    \begin{aligned}
        &\sum_{t=2}^{T+1}\mathbb{E}\left[\sum_{e=1}^E\left\Vert\tilde{u}_{e}^{(t)}\hspace{-0.5mm}-\hspace{-0.5mm}\nabla_e\ell(\w^{(t)})\right\Vert^2\right]\\
        \leq&2T\sigma^2\dfrac{(2-p)m^2-(1-p)m^3}{(1-(1-p)(1-m))(1-(1-p)(1-m)^2)}+\dfrac{4m}{1-(1-p)(1-\frac{m}{2})}\sum_{t=2}^{T+1}\mathbb{E}\left[\left\Vert\Lambda^{(t)}\right\Vert^2\right]\\
        &+\dfrac{(1\hspace{-0.5mm}-\hspace{-0.5mm}\frac{m}{2})p}{1\hspace{-0.5mm}-\hspace{-0.5mm}(1\hspace{-0.5mm}-\hspace{-0.5mm}p)(1\hspace{-0.5mm}-\hspace{-0.5mm}\frac{m}{2})}\sum_{t=2}^T\mathbb{E}\left[\sum_{e=1}^E\left\Vert\tilde{u}_{e}^{(t)}\hspace{-0.5mm}-\hspace{-0.5mm}\nabla_e\ell(\w^{(t)})\right\Vert^2\right]+\dfrac{1}{1\hspace{-0.5mm}-\hspace{-0.5mm}(1\hspace{-0.5mm}-\hspace{-0.5mm}p)(1\hspace{-0.5mm}-\hspace{-0.5mm}\frac{m}{2})}\mathbb{E}\left[\sum_{e=1}^E\left\Vert\tilde{u}_{e}^{(1)}\hspace{-0.5mm}-\hspace{-0.5mm}\nabla_e\ell(\w^{(1)})\right\Vert^2\right]\\
        &+16\left(\dfrac{1}{m(1-(1-p)(1-\frac{m}{2}))}+\dfrac{p+m}{m(1-(1-p)(1-\frac{m}{2}))^2}\right)\sum_{t=1}^T\mathbb{E}\left[\left\Vert\nabla\ell(\w^{(t)})-\nabla\ell(\w^{(t+1)})\right\Vert^2\right].
    \end{aligned}
    \end{equation}}

    Thus, we can obtain that:
    {\small
    \begin{equation}
    \begin{aligned}
        &\dfrac{\frac{m}{2}}{1-(1-p)(1-\frac{m}{2})}\sum_{t=1}^{T+1}\mathbb{E}\left[\sum_{e=1}^E\left\Vert\tilde{u}_{e}^{(t)}\hspace{-0.5mm}-\hspace{-0.5mm}\nabla_e\ell(\w^{(t)})\right\Vert^2\right]\\
        \leq&2T\sigma^2\dfrac{(2-p)m^2-(1-p)m^3}{(1-(1-p)(1-\frac{m}{2}))(1-(1-p)(1-m)^2)}+\dfrac{4m}{1-(1-p)(1-\frac{m}{2})}\sum_{t=2}^{T+1}\mathbb{E}\left[\left\Vert\Lambda^{(t)}\right\Vert^2\right]\\
        &+\dfrac{\frac{m}{2}+1}{1-(1-p)(1-\frac{m}{2})}\mathbb{E}\left[\sum_{e=1}^E\left\Vert\tilde{u}_{e}^{(1)}-\nabla_e\ell(\w^{(1)})\right\Vert^2\right]\\
        &+16\left(\dfrac{1}{m(1-(1-p)(1-\frac{m}{2}))}+\dfrac{p+m}{m(1-(1-p)(1-\frac{m}{2}))^2}\right)\sum_{t=1}^T\mathbb{E}\left[\left\Vert\nabla\ell(\w^{(t)})-\nabla\ell(\w^{(t+1)})\right\Vert^2\right].
    \end{aligned}
    \end{equation}}

    Then, we can get that:
    {\small
    \begin{equation}
    \begin{aligned}
        &\sum_{t=1}^{T+1}\mathbb{E}\left[\sum_{e=1}^E\left\Vert\tilde{u}_{e}^{(t)}\hspace{-0.5mm}-\hspace{-0.5mm}\nabla_e\ell(\w^{(t)})\right\Vert^2\right]\\
        \leq&4T\sigma^2\dfrac{(2-p)m-(1-p)m^2}{1-(1-p)(1-m)^2}+\dfrac{3}{m}\mathbb{E}\left[\sum_{e=1}^E\left\Vert\tilde{u}_{e}^{(1)}-\nabla_e\ell(\w^{(1)})\right\Vert^2\right]+8\sum_{t=2}^{T+1}\mathbb{E}\left[\left\Vert\Lambda^{(t)}\right\Vert^2\right]\\
        &+32\left(\dfrac{p+m}{m^2(1-(1-p)(1-\frac{m}{2}))}+\dfrac{1}{m^2}\right)\sum_{t=1}^T\mathbb{E}\left[\left\Vert\nabla\ell(\w^{(t)})-\nabla\ell(\w^{(t+1)})\right\Vert^2\right].
    \end{aligned}
    \end{equation}}

    Finally, with Eq. \eqref{estimation of a_e^0}, it holds that
    \begin{equation}
    \begin{aligned}
        \mathbb{E}\left[\left\Vert\Lambda^{(t)}\right\Vert^2\right]=&\mathbb{E}\left[\sum_{e=1}^E\left\Vert\nabla_2a_e(\tilde{y}_{e-1}^{(t)},w_e^{(t)})^{\tran}\tilde{v}_e^{(t)}-\nabla_2a_e(\hat{y}_{e-1}^{(t)},w_e^{(t)})^{\tran}\hat{v}_e^{(t)}\right\Vert^2\right]\\
        \leq&(L_{\nabla a}^{\circ})^2\sum_{e=1}^{E}\mathbb{E}\left[\left\Vert\tilde{v}_e^{(t)}-\hat{v}_e^{(t)}\right\Vert^2\right]+(L_{\nabla a}')^2\sum_{e=1}^{E}\mathbb{E}\left[\left\Vert\tilde{y}_{e-1}^{(t)}-\hat{y}_{e-1}^{(t)}\right\Vert^2\right]\\
        =&(L_{\nabla a}^{\circ})^2\sum_{e=1}^{E-1}\mathbb{E}\left[\left\Vert\tilde{v}_e^{(t)}-\hat{v}_e^{(t)}\right\Vert^2\right]+(L_{\nabla a}')^2\sum_{e=1}^{E-1}\mathbb{E}\left[\left\Vert\tilde{y}_{e}^{(t)}-\hat{y}_{e}^{(t)}\right\Vert^2\right],
    \end{aligned}
    \end{equation}
    where the last equation is from the fact that $\hat{y}_0=\tilde{y}_0=x_0$ and $\hat{v}_E=\tilde{v}_E=1$.

    Thus, with Assumption \ref{assumption:smoothness}, it holds that:
    \begin{equation*}
        \begin{aligned}
            &\sum_{e=1}^E\sum_{t=1}^{T+1}\mathbb{E}\left[\left\Vert\tilde{u}^{(t)}_e-\nabla_e\ell(\w^{(t)})\right\Vert^2\right]\\
            \leq&32L_{\nabla \ell}^2\left(\dfrac{p+m}{m^2(1-(1-p)(1-\frac{m}{2}))}+\dfrac{1}{m^2}\right)\sum_{e=1}^E\sum_{t=1}^{T}\mathbb{E}\left[\left\Vert w^{(t+1)}_e-w^{(t)}_e\right\Vert^2\right]\\
            &+8(L_{\nabla a}^\circ)^2\sum_{e=1}^{E-1}\sum_{t=2}^{T+1}\mathbb{E}\left[\left\Vert\tilde{ v}_e^{(t)}- \hat{v}_e^{(t)}\right\Vert^2\right]+8(L_{\nabla a}')^2\sum_{e=1}^{E-1}\sum_{t=2}^{T+1}\mathbb{E}\left[\left\Vert\tilde{\by}_e^{(t)}-\hat{y}_e^{(t)}\right\Vert^2\right]\\
            &+4T\sigma^2\dfrac{(2-p)m-(1-p)m^2}{1-(1-p)(1-m)^2}+\dfrac{3}{m}\sum_{e=1}^E\mathbb{E}\left[\left\Vert\tilde{u}^{(1)}_e-\nabla_e\ell(\w^{(1)})\right\Vert^2\right].
        \end{aligned}
    \end{equation*}
    Thus, we finish the proof of this lemma.
\end{proof}
\end{lemma}
\begin{remark}
    Throughout our analysis, we set $p_2=1$, i.e., the second iteration always draws a fresh sample from $\mathcal{D}$. This acts as a benign initialization that keeps the subsequent recursion clean: in particular, it ensures that the ``last refresh time'' $\psi(t)$ used in the proofs is well-defined for all $t\ge 2$ and avoids carrying arbitrary initial error-feedback states into the main induction. Since this affects only a constant number of iterations, it does not change the asymptotic stationarity rates and incurs negligible additional communication compared with $T$.
\end{remark}

\subsection{Forward-side technical lemmas}
We next bound the errors introduced by error-feedback compression. This subsection focuses on the forward pass and on temporal changes of compressed activations. These estimates are also used later when bounding the backward recursion.

\subsubsection{One-step compression error}
We first isolate one compression operation at one pipeline boundary. The following lemma uses the contractive compressor assumption to relate the current compression error to the change of the uncompressed variable between two consecutive iterations.
\begin{lemma}[Compress error of forward and backward propagation]
\label{lemma:compress error}
Suppose Assumption \ref{assumption:compressor} holds, then for $e=1,2,\cdots,E-1$ and $T_1>T_0\geq1$ we have:
\begin{subequations}
    \begin{align}
        \label{expection_v_chi_2}
        \sum_{t=T_0+1}^{T_1+1}\mathbb{E}\left[\left\Vert\tilde{ v}^{(t)}_e-v^{(t)}_e\right\Vert^2\right]
        &\leq\dfrac{\omega_B^2}{(1-\omega_B)^2}\sum_{t=T_0}^{T_1}\mathbb{E}\left[\left\Vert v^{(t+1)}_e-v^{(t)}_e\right\Vert^2\right]+\dfrac{\omega_B}{1-\omega_B}\mathbb{E}\left[\left\Vert\tilde{ v}^{(T_0)}_{e}-v^{(T_0)}_e\right\Vert^2\right],\\
        \label{expection_y_theta_2}
        \sum_{t=T_0+1}^{T_1+1}\mathbb{E}\left[\left\Vert\tilde{\by}^{(t)}_e-\by^{(t)}_e\right\Vert^2\right]
        &\leq\dfrac{\omega_F^2}{(1-\omega_F)^2}\sum_{t=T_0}^{T_1}\mathbb{E}\left[\left\Vert\by^{(t+1)}_e-\by^{(t)}_e\right\Vert^2\right]+\dfrac{\omega_F}{1-\omega_F}\mathbb{E}\left[\left\Vert\tilde{\by}^{(T_0)}_{e}-\by^{(T_0)}_e\right\Vert^2\right].
    \end{align}
\end{subequations}
\begin{proof}
    Firstly we consider the term $\left\Vert\tilde{ v}^{(t)}_e-\bv^{(t)}_e\right\Vert^2$ for $e=1,2,\cdots,E-1$. According to Algorithm \ref{alg:backward}, we have:
    \begin{equation}
    \label{expection_v_chi}
    \begin{aligned}
        \mathbb{E}\left[\left\Vert\tilde{ v}^{(t)}_e-\bv^{(t)}_e\right\Vert^2\middle|\mathcal{G}_e^{(t)}\right]=&\mathbb{E}\left[\left\Vert\tilde{ v}^{(t-1)}_{e}+\mathcal{C}(\bv^{(t)}_e-\tilde{ v}^{(t-1)}_{e})-\bv^{(t)}_e\right\Vert^2\middle|\mathcal{G}_e^{(t)}\right]\leq\omega_B^2\left\Vert\tilde{ v}^{(t-1)}_{e}-\bv^{(t)}_e\right\Vert^2\\
        \leq&\omega_B\left\Vert\tilde{ v}^{(t-1)}_{e}-\bv^{(t-1)}_e\right\Vert^2+\dfrac{\omega_B^2}{1-\omega_B}\left\Vert\bv^{(t)}_e-\bv^{(t-1)}_e\right\Vert^2,
    \end{aligned}
    \end{equation}
    where the first inequality is due to Assumption \ref{assumption:compressor} and the second inequality uses Young's inequality. Then, taking expectation on both sides and then taking summation over $t=T_0+1,\cdots,T_1+1$, we have:
    \begin{align}
        \sum_{t=T_0+1}^{T_1+1}\mathbb{E}\left[\left\Vert\tilde{ v}^{(t)}_e-\bv^{(t)}_e\right\Vert^2\right]
        \leq\omega_B\sum_{t=T_0}^{T_1}\mathbb{E}\left[\left\Vert\tilde{ v}^{(t)}_{e}-\bv^{(t)}_e\right\Vert^2\right]+\dfrac{\omega_B^2}{1-\omega_B}\sum_{t=T_0}^{T_1}\mathbb{E}\left[\left\Vert\bv^{(t+1)}_e-\bv^{(t)}_e\right\Vert^2\right].
    \end{align}
    Then, we can get:
    \begin{align*}
         \sum_{t=T_0+1}^{T_1+1}\mathbb{E}\left[\left\Vert\tilde{ v}^{(t)}_e-\bv^{(t)}_e\right\Vert^2\right]
        \leq\dfrac{\omega_B^2}{(1-\omega_B)^2} \sum_{t=T_0}^{T_1}\mathbb{E}\left[\left\Vert\bv^{(t+1)}_e-\bv^{(t)}_e\right\Vert^2\right]+\dfrac{\omega_B}{1-\omega_B}\mathbb{E}\left[\left\Vert\tilde{ v}^{(T_0)}_{e}-\bv^{(T_0)}_e\right\Vert^2\right].
    \end{align*}
    Thus, Eq. \eqref{expection_v_chi_2} holds.

    Next, we consider the term $\left\Vert\tilde{\by}^{(t)}_e-\by^{(t)}_e\right\Vert^2$ for $e=1,2,\cdots,E-1$. According to Algorithm \ref{alg:forward}, we have:
    \begin{equation}
    \label{expection_y_chi}
    \begin{aligned}
        \mathbb{E}\left[\left\Vert\tilde{\by}^{(t)}_e-\by^{(t)}_e\right\Vert^2\middle|\mathcal{F}_{e}^{(t)}\right]=&\mathbb{E}\left[\left\Vert\tilde{\by}^{(t-1)}_{e}+\mathcal{C}(\by^{(t)}_e-\tilde{\by}^{(t-1)}_{e})-\by^{(t)}_e\right\Vert^2\middle|\mathcal{F}_{e}^{(t)}\right]\leq\omega_F^2\left\Vert\tilde{\by}^{(t-1)}_{e}-\by^{(t)}_e\right\Vert^2\\
        \leq&\omega_F\left\Vert\tilde{\by}^{(t-1)}_{e}-\by^{(t-1)}_e\right\Vert^2+\dfrac{\omega_F^2}{1-\omega_F}\left\Vert\by^{(t)}_e-\by^{(t-1)}_e\right\Vert^2,
    \end{aligned}
    \end{equation}
    where the first inequality is due to Assumption \ref{assumption:compressor} and the second inequality uses Young's inequality. Then, taking expectation on both sides and then taking summation over $t=T_0+1,\cdots,T_1+1$, we have:
    \begin{align}
        \sum_{t=T_0+1}^{T_1+1}\mathbb{E}\left[\left\Vert\tilde{\by}^{(t)}_e-\by^{(t)}_e\right\Vert^2\right]
        \leq\omega_F\sum_{t=T_0}^{T_1}\mathbb{E}\left[\left\Vert\tilde{\by}^{(t)}_{e}-\by^{(t)}_e\right\Vert^2\right]+\dfrac{\omega_F^2}{1-\omega_F}\sum_{t=T_0}^{T_1}\mathbb{E}\left[\left\Vert\by^{(t+1)}_e-\by^{(t)}_e\right\Vert^2\right].
    \end{align}
    Then, we can get:
    \begin{align*}
        \sum_{t=T_0+1}^{T_1+1}\mathbb{E}\left[\left\Vert\tilde{\by}^{(t)}_e-\by^{(t)}_e\right\Vert^2\right]
        \leq\dfrac{\omega_F^2}{(1-\omega_F)^2}\sum_{t=T_0}^{T_1}\mathbb{E}\left[\left\Vert\by^{(t+1)}_e-\by^{(t)}_e\right\Vert^2\right]+\dfrac{\omega_F}{1-\omega_F}\mathbb{E}\left[\left\Vert\tilde{\by}^{(T_0)}_{e}-\by^{(T_0)}_e\right\Vert^2\right].
    \end{align*}
    Thus, Eq. \eqref{expection_y_theta_2} holds.
\end{proof}
\end{lemma}

\subsubsection{Forward activation error}
We now propagate the one-step forward compression error across pipeline stages. The goal is to compare the activation used by \ours with the activation that would have been obtained by the uncompressed pipeline.
\begin{lemma}
    Suppose Assumption \ref{assumption:smoothness} and \ref{assumption:compressor} holds, then for any $T_1>T_0\geq1$ we have:
    \begin{equation}
    \label{y_tilde-y}
    \begin{aligned}
        \sum_{e=1}^{E-1}\sum_{t=T_0+1}^{T_1+1}\mathbb{E}\left[\left\Vert\tilde{\by}_e^{(t)}-\hat{y}_e^{(t)}\right\Vert^2\right]\leq&\sum_{e=1}^{E-1}\sum_{\iota=e}^{E-1}2(2L_a^2)^{\iota-e}\dfrac{\omega_F^2}{(1-\omega_F)^2}\sum_{t=T_0}^{T_1}\mathbb{E}\left[\left\Vert\by^{(t+1)}_e-\by^{(t)}_e\right\Vert^2\right]\\
        &+\sum_{e=1}^{E-1}\sum_{\iota=e}^{E-1}2(2L_a^2)^{\iota-e}\dfrac{\omega_F}{1-\omega_F}\mathbb{E}\left[\left\Vert\tilde{\by}^{(T_0)}_{e}-\by^{(T_0)}_e\right\Vert^2\right].
    \end{aligned}
    \end{equation}
\begin{proof}
    For $1\leq e\leq E-1$, we have:
    \begin{equation}
    \label{zhankai_y}
    \begin{aligned}
        \left\Vert\tilde{\by}_e^{(t)}-\hat{y}_e^{(t)}\right\Vert^2\leq&2\left\Vert\tilde{\by}_e^{(t)}-\by_e^{(t)}\right\Vert^2+2\left\Vert\by_e^{(t)}-\hat{y}_e^{(t)}\right\Vert^2\\
        =&2\left\Vert\tilde{\by}_e^{(t)}-\by_e^{(t)}\right\Vert^2+2\left\Vert a_e(\tilde{\by}_{e-1}^{(t)},w_e^{(t)})-a_e(\hat{y}_{e-1}^{(t)},w_e^{(t)})\right\Vert^2\\
        \leq&2\left\Vert\tilde{\by}_e^{(t)}-\by_e^{(t)}\right\Vert^2+2L_{a}^2\left\Vert\tilde{\by}_{e-1}^{(t)}-\hat{y}_{e-1}^{(t)}\right\Vert^2\leq\cdots\leq\sum_{\iota=1}^{e}2(2L_a^2)^{e-\iota}\left\Vert\tilde{\by}_\iota^{(t)}-\by_\iota^{(t)}\right\Vert^2,
    \end{aligned}
    \end{equation}
    where the second inequality is due to Assumption \ref{assumption:smoothness}.

    Taking expectation and then taking summation on both sides of \eqref{zhankai_y} over $t=T_0+1,\cdots,T_1+1$, then we have:
    \begin{equation}
    \label{zhankai_y_21}
    \begin{aligned}
        &\sum_{t=T_0+1}^{T_1+1}\mathbb{E}\left[\left\Vert\tilde{\by}_e^{(t)}-\hat{y}_e^{(t)}\right\Vert^2\right]\leq\sum_{\iota=1}^{e}2(2L_a^2)^{e-\iota}\sum_{t=T_0+1}^{T_1+1}\mathbb{E}\left[\left\Vert\tilde{\by}_\iota^{(t)}-\by_\iota^{(t)}\right\Vert^2\right]\\
        \leq&\sum_{\iota=1}^{e}2(2L_a^2)^{e-\iota}\dfrac{\omega_F^2}{(1-\omega_F)^2}\sum_{t=T_0}^{T_1}\mathbb{E}\left[\left\Vert\by^{(t+1)}_{\iota}-\by^{(t)}_{\iota}\right\Vert^2\right]+\sum_{\iota=1}^{e}2(2L_a^2)^{e-\iota}\dfrac{\omega_F}{1-\omega_F}\mathbb{E}\left[\left\Vert\tilde{\by}^{(T_0)}_{\iota}-\by^{(T_0)}_{\iota}\right\Vert^2\right],
    \end{aligned}
    \end{equation}
    where the last inequality is due to the result of \eqref{expection_y_theta_2}.

    Taking summation on both sides of \eqref{zhankai_y_21} over $e$, we can get:
    \begin{equation*}
    \label{zhankai_y_2}
    \begin{aligned}
        \sum_{e=1}^{E-1}\sum_{t=T_0+1}^{T_1+1}\mathbb{E}\left[\left\Vert\tilde{\by}_e^{(t)}-\hat{y}_e^{(t)}\right\Vert^2\right]\leq&\sum_{e=1}^{E-1}\sum_{\iota=e}^{E-1}2(2L_a^2)^{\iota-e}\dfrac{\omega_F^2}{(1-\omega_F)^2}\sum_{t=T_0}^{T_1}\mathbb{E}\left[\left\Vert\by^{(t+1)}_e-\by^{(t)}_e\right\Vert^2\right]\\
        &+\sum_{e=1}^{E-1}\sum_{\iota=e}^{E-1}2(2L_a^2)^{\iota-e}\dfrac{\omega_F}{1-\omega_F}\mathbb{E}\left[\left\Vert\tilde{\by}^{(T_0)}_{e}-\by^{(T_0)}_e\right\Vert^2\right].
    \end{aligned}
    \end{equation*}
\end{proof}
\end{lemma}

\subsubsection{Forward temporal variation}
The backward analysis needs to compare compressed activations from two consecutive iterations. The next lemma converts this temporal variation into the corresponding change of the uncompressed activation, plus the initial compression error.
\begin{lemma}[Temporal variation of compressed activations]
\label{lemma6}
    Suppose Assumption \ref{assumption:compressor} and Eq. \eqref{expection_y_theta_2} holds, then for any $T_1>T_0\geq1$ we have:
    \begin{equation}
    \label{yt+1-yt_main}
        \begin{aligned}
            \sum_{t=T_0}^{T_1}\mathbb{E}\left[\left\Vert\tilde{\by}_e^{(t+1)}-\tilde{\by}_e^{(t)}\right\Vert^2\right]
            \leq\dfrac{8}{(1-\omega_F)^2}\sum_{t=T_0}^{T_1}\mathbb{E}\left[\left\Vert\by_e^{(t+1)}-\by_e^{(t)}\right\Vert^2\right]+\dfrac{8}{1-\omega_F}\mathbb{E}\left[\left\Vert\tilde{\by}_e^{(T_0)}-\by_e^{(T_0)}\right\Vert^2\right].
        \end{aligned}
    \end{equation}
\begin{proof}
    For $t=T_0+1,\cdots,T_1+1$, we have:
    \begin{equation}
    \label{yt+1-yt}
        \begin{aligned}
            &\mathbb{E}\left[\left\Vert\tilde{\by}_e^{(t)}-\tilde{\by}_e^{(t-1)}\right\Vert^2\middle|\mathcal{F}^{(t)}\right]=\mathbb{E}\left[\left\Vert\mathcal{C}\left(\by_e^{(t)}-\tilde{\by}_e^{(t-1)}\right)\right\Vert^2\middle|\mathcal{F}^{(t)}\right]\\
            \leq&\left(1+\dfrac{1}{\omega_F}\right)\mathbb{E}\left[\left\Vert\mathcal{C}\left(\by_e^{(t)}-\tilde{\by}_e^{(t-1)}\right)-\left(\by_e^{(t)}-\tilde{\by}_e^{(t-1)}\right)\right\Vert^2\middle|\mathcal{F}^{(t)}\right]+(1+\omega_F)\left\Vert\by_e^{(t)}-\tilde{\by}_e^{(t-1)}\right\Vert^2\\
            \leq&(1+\omega_F)^2\left\Vert\by_e^{(t)}-\tilde{\by}_e^{(t-1)}\right\Vert^2\leq8\left\Vert\by_e^{(t)}-\by_e^{(t-1)}\right\Vert^2+8\left\Vert\by_e^{(t-1)}-\tilde{\by}_e^{(t-1)}\right\Vert^2,
        \end{aligned}
    \end{equation}
    where the first inequality is due to Young's inequality, the second inequality is due to Assumption \ref{assumption:compressor}, and the last inequality is due to $\omega_F\in[0,1)$ in Assumption \ref{assumption:compressor}.

    Then, taking expectation and taking summation on both sides of \eqref{yt+1-yt} over $t=T_0+1,\cdots,T_1+1$, we have:
    \begin{equation*}
        \begin{aligned}
            \sum_{t=T_0}^{T_1}\mathbb{E}\left[\left\Vert\tilde{\by}_e^{(t+1)}-\tilde{\by}_e^{(t)}\right\Vert^2\right]
            \leq&8\sum_{t=T_0}^{T_1}\mathbb{E}\left[\left\Vert\by_e^{(t+1)}-\by_e^{(t)}\right\Vert^2\right]+8\sum_{t=T_0}^{T_1}\mathbb{E}\left[\left\Vert\by_e^{(t)}-\tilde{\by}_e^{(t)}\right\Vert^2\right]\\
            \leq&\dfrac{8}{(1-\omega_F)^2}\sum_{t=T_0}^{T_1}\mathbb{E}\left[\left\Vert\by_e^{(t+1)}-\by_e^{(t)}\right\Vert^2\right]+\dfrac{8}{1-\omega_F}\mathbb{E}\left[\left\Vert\tilde{\by}_e^{(T_0)}-\by_e^{(T_0)}\right\Vert^2\right],
        \end{aligned}
    \end{equation*}
    where the last inequality is due to Eq. \eqref{expection_y_theta_2} and the fact that $1+\dfrac{\omega_F^2}{(1-\omega_F)^2}\leq\dfrac{1}{(1-\omega_F)^2}$ as $\omega_F\in[0,1)$.
\end{proof}
\end{lemma}

For intuition, Appendix~\ref{app:two-worker-warmup} provides a simplified proof for \(E=2\). The general proof below does not rely on the $E=2$ case.

\subsection{Variation Bounds}
\label{app:variation-bounds}
We now turn to the backward pass and to the variation terms created by the pipeline recursion. The following bounds translate activation and adjoint changes into quantities that can later be absorbed by the descent inequality.

\subsubsection{Forward variation}
We first relate the backward adjoint error \(\|\tilde v_e^{(t)}-\hat v_e^{(t)}\|^2\) to two simpler sources: local backward compression error and the forward activation error bounded above.

\begin{lemma}
\label{lemma7}
Suppose Assumption \ref{assumption:smoothness} and \ref{assumption:compressor} holds, then for any $T_1>T_0\geq1$ and $e=1,2,\cdots,E-1$ we have:
    {\small
    \begin{equation}
    \label{v_tilde-v}
    \begin{aligned}
        &\sum_{t=T_0+1}^{T_1+1}\mathbb{E}\left[\left\Vert\tilde{ v}^{(t)}_e- \hat{v}^{(t)}_e\right\Vert^2\right]\leq2\dfrac{\omega_B^2}{(1-\omega_B)^2}\sum_{\iota=e}^{E-1}(2(L_{\nabla a}^\circ)^2)^{\iota-e}\sum_{t=T_0}^{T_1}\mathbb{E}\left[\left\Vert\bv^{(t+1)}_{\iota}-\bv^{(t)}_{\iota}\right\Vert^2\right]\\
        +&(2(L_{\nabla a}')^2)\sum_{\iota=e}^{E-1}(2(L_{\nabla a}^\circ)^2)^{\iota-e}\sum_{t=T_0+1}^{T_1+1}\mathbb{E}\left[\left\Vert \tilde{\by}_{\iota}^{(t)}-\hat{y}_{\iota}^{(t)}\right\Vert^2\right]+2\dfrac{\omega_B}{1-\omega_B}\sum_{\iota=e}^{E-1}(2(L_{\nabla a}^\circ)^2)^{\iota-e}\mathbb{E}\left[\left\Vert\tilde{ v}^{(T_0)}_{\iota}-\bv^{(T_0)}_{\iota}\right\Vert^2\right].
    \end{aligned}
    \end{equation}}
    \begin{proof}
        For $1\leq e\leq E-1$, we have:
        \begin{equation}
            \begin{aligned}
                \left\Vert\tilde{ v}^{(t)}_e- \hat{v}^{(t)}_e\right\Vert^2\leq&2\left\Vert\tilde{ v}^{(t)}_e-\bv^{(t)}_e\right\Vert^2+2\left\Vert\bv^{(t)}_e- \hat{v}^{(t)}_e\right\Vert^2\\
                =&2\left\Vert\tilde{ v}^{(t)}_e-\bv^{(t)}_e\right\Vert^2+2
                \left\Vert\nabla_1a_{e+1}(\tilde{\by}_{e}^{(t)},w_{e+1}^{(t)})^{\tran}\tilde{ v}_{e+1}^{(t)}- \nabla_1a_{e+1}(\hat{y}_{e}^{(t)},w_{e+1}^{(t)})^{\tran}\hat{v}_{e+1}^{(t)}\right\Vert^2\\
                \leq&2\left\Vert\tilde{ v}^{(t)}_e-\bv^{(t)}_e\right\Vert^2+2(L_{\nabla a}')^2\left\Vert \tilde{\by}_{e}^{(t)}-\hat{y}_{e}^{(t)}\right\Vert^2+2(L_{\nabla a}^\circ)^2\left\Vert \tilde{ v}_{e+1}^{(t)}- \hat{v}_{e+1}^{(t)}\right\Vert^2,
            \end{aligned}
        \end{equation}
    where the second inequality is due to \eqref{estimation of a_e^0}. Then, we can get:
        \begin{equation}
        \label{zhankai_v}
            \begin{aligned}
                &\left\Vert\tilde{ v}^{(t)}_e- \hat{v}^{(t)}_e\right\Vert^2\leq\sum_{\iota=e}^{E-1}2(2(L_{\nabla a}^\circ)^2)^{\iota-e}\left[\left\Vert\tilde{ v}^{(t)}_\iota-\bv^{(t)}_\iota\right\Vert^2+(L_{\nabla a}')^2\left\Vert \tilde{\by}_{\iota}^{(t)}-\hat{y}_{\iota}^{(t)}\right\Vert^2\right]\\
                \leq&2\sum_{\iota=e}^{E-1}(2(L_{\nabla a}^\circ)^2)^{\iota-e}\left\Vert\tilde{ v}^{(t)}_\iota-\bv^{(t)}_\iota\right\Vert^2+(2(L_{\nabla a}')^2)\sum_{\iota=e}^{E-1}(2(L_{\nabla a}^\circ)^2)^{\iota-e}\left\Vert \tilde{\by}_{\iota}^{(t)}-\hat{y}_{\iota}^{(t)}\right\Vert^2.
            \end{aligned}
        \end{equation}
    
    Taking expectation and then taking summation on both sides of \eqref{zhankai_v} over $t=T_0+1,\cdots,T_1+1$, then we gan get
        {\small
        \begin{equation}
        \label{v_tilde-v,chubu_0}
            \begin{aligned}
                &\sum_{t=T_0+1}^{T_1+1}\mathbb{E}\left[\left\Vert\tilde{ v}^{(t)}_e- \hat{v}^{(t)}_e\right\Vert^2\right]\\
                \leq&2\sum_{\iota=e}^{E-1}(2(L_{\nabla a}^\circ)^2)^{\iota-e}\sum_{t=T_0+1}^{T_1+1}\mathbb{E}\left[\left\Vert\tilde{ v}^{(t)}_\iota-\bv^{(t)}_\iota\right\Vert^2\right]+(2(L_{\nabla a}')^2)\sum_{\iota=e}^{E-1}(2(L_{\nabla a}^\circ)^2)^{\iota-e}\sum_{t=T_0+1}^{T_1+1}\mathbb{E}\left[\left\Vert \tilde{\by}_{\iota}^{(t)}-\hat{y}_{\iota}^{(t)}\right\Vert^2\right]\\
                \leq&2\dfrac{\omega_B^2}{(1\hspace{-0.5mm}-\hspace{-0.5mm}\omega_B)^2}\hspace{-0.5mm}\sum_{\iota=e}^{E-1}(2(L_{\nabla a}^\circ)^2)^{\iota-e}\sum_{t=T_0}^{T_1}\mathbb{E}\left[\left\Vert\bv^{(t+1)}_{\iota}\hspace{-0.5mm}-\hspace{-0.5mm}\bv^{(t)}_{\iota}\right\Vert^2\right]\hspace{-0.5mm}+\hspace{-0.5mm}2\dfrac{\omega_B}{1-\omega_B}\sum_{\iota=e}^{E-1}(2(L_{\nabla a}^\circ)^2)^{\iota-e}\mathbb{E}\left[\left\Vert\tilde{ v}^{(T_0)}_{\iota}\hspace{-0.5mm}-\hspace{-0.5mm}\bv^{(T_0)}_{\iota}\right\Vert^2\right]\\
                &+(2(L_{\nabla a}')^2)\sum_{\iota=e}^{E-1}(2(L_{\nabla a}^\circ)^2)^{\iota-e}\sum_{t=T_0+1}^{T_1+1}\mathbb{E}\left[\left\Vert \tilde{\by}_{\iota}^{(t)}-\hat{y}_{\iota}^{(t)}\right\Vert^2\right],\\
            \end{aligned}
        \end{equation}}
    where the last inequality is due to Eq. \eqref{expection_v_chi_2}.
    \end{proof}
\end{lemma}

\subsubsection{Exact-adjoint variation}
The previous bound contains the temporal difference \(\|\bv_e^{(t+1)}-\bv_e^{(t)}\|^2\). We next control this term by changes in the next-stage adjoint, the compressed activation, and the local model block.
\begin{lemma}
\label{lemma8}
Suppose Assumption \ref{assumption:smoothness} and \eqref{yt+1-yt_main} hold, then for any $T_1>T_0\geq1$ and $e=1,2,\cdots,E-1$ we have:
{\small
\begin{equation}
\label{e_e_chi}
\begin{aligned}
    &\sum_{t=T_0}^{T_1}\mathbb{E}\left[\left\Vert\bv^{(t+1)}_{e}-\bv^{(t)}_{e}\right\Vert^2\right]\leq5(L_{\nabla a}^\circ)^2\sum_{t=T_0+1}^{T_1+1}\mathbb{E}\left[\left\Vert\tilde{ v}_{e+1}^{(t)}- \hat{v}_{e+1}^{(t)}\right\Vert^2\right]+5(L_{\nabla a}^\circ)^2\sum_{t=T_0}^{T_1}\mathbb{E}\left[\left\Vert \hat{v}_{e+1}^{(t+1)}-\hat{v}_{e+1}^{(t)}\right\Vert^2\right]\\
    +&5(L_{\nabla a}')^2\sum_{t=T_0}^{T_1}\mathbb{E}\left[\left\Vert w_{e+1}^{(t+1)}-w_{e+1}^{(t)}\right\Vert^2\right]+\dfrac{40(L_{\nabla a}')^2}{(1-\omega_F)^2}\sum_{t=T_0}^{T_1}\mathbb{E}\left[\left\Vert\by_e^{(t+1)}-\by_e^{(t)}\right\Vert^2\right]\\
    +&\dfrac{40(L_{\nabla a}')^2}{1-\omega_F}\mathbb{E}\left[\left\Vert\tilde{\by}_e^{(T_0)}-\by_e^{(T_0)}\right\Vert^2\right]+\dfrac{5}{2}(L_{\nabla a}^\circ)^2\mathbb{E}\left[\left\Vert\tilde{ v}_{e+1}^{(T_0)}-\hat{v}_{e+1}^{(T_0)}\right\Vert^2\right].
\end{aligned}
\end{equation}}
\begin{proof}
Firstly, we consider the case of $e=1,\cdots,E-1$, we have:
\begin{equation}
\begin{aligned}
    &\sum_{t=T_0}^{T_1}\mathbb{E}\left[\left\Vert\bv^{(t+1)}_{e}-\bv^{(t)}_{e}\right\Vert^2\right]\\
    =&\sum_{t=T_0}^{T_1}\mathbb{E}\left[\left\Vert\nabla_1a_{e+1}(\tilde{\by}_e^{(t+1)},w_{e+1}^{(t+1)})^{\tran}\tilde{ v}_{e+1}^{(t+1)}-\nabla_1a_{e+1}(\tilde{\by}_e^{(t)},w_{e+1}^{(t)})^{\tran}\tilde{ v}_{e+1}^{(t)}\right\Vert^2\right]\\
    \leq&\dfrac{5}{2}\sum_{t=T_0}^{T_1}\mathbb{E}\left[\left\Vert\nabla_1a_{e+1}(\tilde{\by}_e^{(t+1)},w_{e+1}^{(t+1)})^{\tran}\tilde{ v}_{e+1}^{(t+1)}-\nabla_1a_{e+1}(\tilde{\by}_e^{(t+1)},w_{e+1}^{(t+1)})^{\tran}\hat{ v}_{e+1}^{(t+1)}\right\Vert^2\right]\\
    &+\dfrac{5}{2}\sum_{t=T_0}^{T_1}\mathbb{E}\left[\left\Vert\nabla_1a_{e+1}(\tilde{\by}_e^{(t)},w_{e+1}^{(t)})^{\tran}\tilde{ v}_{e+1}^{(t)}-\nabla_1a_{e+1}(\tilde{\by}_e^{(t)},w_{e+1}^{(t)})^{\tran}\hat{ v}_{e+1}^{(t)}\right\Vert^2\right]\\
    &+5\sum_{t=T_0}^{T_1}\mathbb{E}\left[\left\Vert \nabla_1a_{e+1}(\tilde{\by}_e^{(t+1)},w_{e+1}^{(t+1)})^{\tran}\hat{ v}_{e+1}^{(t+1)}-\nabla_1a_{e+1}(\tilde{\by}_e^{(t)},w_{e+1}^{(t)})^{\tran}\hat{ v}_{e+1}^{(t)}\right\Vert^2\right].
\end{aligned}
\end{equation}

Then, we have:
\begin{equation}
\begin{aligned}
    &\sum_{t=T_0}^{T_1}\mathbb{E}\left[\left\Vert\bv^{(t+1)}_{e}\hspace{-0.5mm}-\hspace{-0.5mm}\bv^{(t)}_{e}\right\Vert^2\right]\\
    \leq&\dfrac{5}{2}(L_{\nabla a}^\circ)^2\sum_{t=T_0}^{T_1}\mathbb{E}\left[\left\Vert\tilde{ v}_{e+1}^{(t+1)}\hspace{-0.5mm}-\hspace{-0.5mm}\hat{v}_{e+1}^{(t+1)}\right\Vert^2\right]\hspace{-0.5mm}+\hspace{-0.5mm}\dfrac{5}{2}(L_{\nabla a}^\circ)^2\sum_{t=T_0}^{T_1}\mathbb{E}\left[\left\Vert\tilde{ v}_{e+1}^{(t)}\hspace{-0.5mm}-\hspace{-0.5mm}\hat{v}_{e+1}^{(t)}\right\Vert^2\right]\\
    &+\hspace{-0.5mm}5(L_{\nabla a}^\circ)^2\sum_{t=T_0}^{T_1}\mathbb{E}\left[\left\Vert \hat{v}_{e+1}^{(t+1)}\hspace{-0.5mm}-\hspace{-0.5mm}\hat{v}_{e+1}^{(t)}\right\Vert^2\right]\hspace{-0.5mm}+\hspace{-0.5mm}5(L_{\nabla a}')^2\sum_{t=T_0}^{T_1}\left(\mathbb{E}\left[\left\Vert\tilde{\by}_{e}^{(t+1)}\hspace{-0.5mm}-\hspace{-0.5mm}\tilde{\by}_{e}^{(t)}\right\Vert^2\right]\hspace{-0.5mm}+\hspace{-0.5mm}\mathbb{E}\left[\left\Vert w_{e+1}^{(t+1)}\hspace{-0.5mm}-\hspace{-0.5mm}w_{e+1}^{(t)}\right\Vert^2\right]\right)\\
    \leq&5(L_{\nabla a}^\circ)^2\sum_{t=T_0+1}^{T_1+1}\mathbb{E}\left[\left\Vert\tilde{ v}_{e+1}^{(t)}\hspace{-0.5mm}-\hspace{-0.5mm}\hat{v}_{e+1}^{(t)}\right\Vert^2\right]\hspace{-0.5mm}+\hspace{-0.5mm}5(L_{\nabla a}^\circ)^2\sum_{t=T_0}^{T_1}\mathbb{E}\left[\left\Vert \hat{v}_{e+1}^{(t+1)}\hspace{-0.5mm}-\hspace{-0.5mm}\hat{v}_{e+1}^{(t)}\right\Vert^2\right]\\
    &+\hspace{-0.5mm}5(L_{\nabla a}')^2\sum_{t=T_0}^{T_1}\mathbb{E}\left[\left\Vert w_{e+1}^{(t+1)}\hspace{-0.5mm}-\hspace{-0.5mm}w_{e+1}^{(t)}\right\Vert^2\right]\hspace{-0.5mm}+\hspace{-0.5mm}\dfrac{40(L_{\nabla a}')^2}{(1-\omega_F)^2}\sum_{t=T_0}^{T_1}\mathbb{E}\left[\left\Vert\by_e^{(t+1)}\hspace{-0.5mm}-\hspace{-0.5mm}\by_e^{(t)}\right\Vert^2\right]\\
    &+\hspace{-0.5mm}\dfrac{40(L_{\nabla a}')^2}{1-\omega_F}\mathbb{E}\left[\left\Vert\tilde{\by}_e^{(T_0)}\hspace{-0.5mm}-\hspace{-0.5mm}\by_e^{(T_0)}\right\Vert^2\right]\hspace{-0.5mm}+\hspace{-0.5mm}\dfrac{5}{2}(L_{\nabla a}^\circ)^2\mathbb{E}\left[\left\Vert\tilde{ v}_{e+1}^{(T_0)}\hspace{-0.5mm}-\hspace{-0.5mm}\hat{v}_{e+1}^{(T_0)}\right\Vert^2\right],
\end{aligned}
\end{equation}
where the first inequality is due to Eq. \eqref{estimation of a_e^0} and the last inequality is due to Eq. \eqref{yt+1-yt_main}.
\end{proof}
\end{lemma}

\subsubsection{Aggregated backward error bound}
Combining the two previous recursions gives a pipeline-level bound for all backward adjoint errors. This is where the compression constants begin to accumulate with the number of pipeline stages.
\begin{lemma}
Suppose Lemma \ref{lemma7} and Lemma \ref{lemma8} are all satisfied, then there exist coefficients $C_{y,e}^\circ,C_{W,e}^\circ,C_{v,e}^\circ,C_{v,e}^1,C_{v,e}^2,C_{\theta,e}^\circ,C_{\theta,e}^1\geq0$ (which has been defined by Eq. \eqref{quzhi_1}) for each $e=1,2,\cdots,E-1$ and $T_1>T_0\geq1$ such that:
{\small\begin{equation}
\label{v_tilde-v_3}
\begin{aligned}
    &\sum_{e=1}^{E-1}\sum_{t=T_0+1}^{T_1+1}\mathbb{E}\left[\left\Vert\tilde{ v}^{(t)}_e\hspace{-0.5mm}-\hspace{-0.5mm}\hat{v}^{(t)}_e\right\Vert^2\right]\\
    \leq&\sum_{e=1}^{E-1}C_{y,e}^{\circ}\sum_{t=T_0+1}^{T_1+1}\mathbb{E}\left[\left\Vert \tilde{\by}_{e}^{(t)}\hspace{-0.5mm}-\hspace{-0.5mm}\hat{y}_{e}^{(t)}\right\Vert^2\right]+\sum_{e=1}^{E}C_{W,e}^{\circ}\sum_{t=T_0}^{T_1}\mathbb{E}\left[\left\Vert w_{e}^{(t+1)}\hspace{-0.5mm}-\hspace{-0.5mm}w_{e}^{(t)}\right\Vert^2\right]\\
    &+\sum_{e=1}^{E-1}C_{v,e}^{\circ}\sum_{t=T_0}^{T_1}\mathbb{E}\left[\left\Vert\hat{v}_{e}^{(t+1)}\hspace{-0.5mm}-\hspace{-0.5mm}\hat{v}_{e}^{(t)}\right\Vert^2\right]+\sum_{e=1}^{E-1}C_{\theta,e}^{\circ}
    \sum_{t=T_0}^{T_1}\mathbb{E}\left[\left\Vert\by_e^{(t+1)}\hspace{-0.5mm}-\hspace{-0.5mm}\by_e^{(t)}\right\Vert^2\right]+\sum_{e=1}^{E-1}C_{\theta,e}^{1}\mathbb{E}\left[\left\Vert\tilde{\by}_e^{(T_0)}\hspace{-0.5mm}-\hspace{-0.5mm}\by_e^{(T_0)}\right\Vert^2\right]\\
    &+\sum_{e=1}^{E-1}C_{v,e}^{1}\mathbb{E}\left[\left\Vert\tilde{ v}_{e}^{(T_0)}\hspace{-0.5mm}-\hspace{-0.5mm}\hat{v}_{e}^{(T_0)}\right\Vert^2\right]+\sum_{e=1}^{E-1}C_{v,e}^{2}\mathbb{E}\left[\left\Vert\tilde{ v}_{e}^{(T_0)}\hspace{-0.5mm}-\hspace{-0.5mm}\bv_{e}^{(T_0)}\right\Vert^2\right].
\end{aligned}
\end{equation}}
\begin{proof}
Combing \eqref{v_tilde-v} and \eqref{e_e_chi} together, we have:
{\small\begin{equation}
\label{v_tilde-v1}
\begin{aligned}
    &\sum_{t=T_0+1}^{T_1+1}\mathbb{E}\left[\left\Vert\tilde{ v}^{(t)}_e-\hat{v}^{(t)}_e\right\Vert^2\right]\\
    \leq&\dfrac{5(L_{\nabla a}^*)^2\omega_B^2}{(1-\omega_B)^2}\sum_{\iota=e}^{E-1}\sum_{t=T_0+1}^{T_1+1}\mathbb{E}\left[\left\Vert\tilde{ v}_{\iota+1}^{(t)}-\hat{v}_{\iota+1}^{(t)}\right\Vert^2\right]+\dfrac{5(L_{\nabla a}^*)^2\omega_B^2}{(1-\omega_B)^2}\sum_{\iota=e}^{E-1}\sum_{t=T_0}^{T_1}\mathbb{E}\left[\left\Vert\hat{v}_{\iota+1}^{(t+1)}-\hat{v}_{\iota+1}^{(t)}\right\Vert^2\right]\\
    &+(L_{\nabla a}^*)^2(2(L_{\nabla a}')^2)\sum_{\iota=e}^{E-1}\sum_{t=T_0+1}^{T_1+1}\mathbb{E}\left[\left\Vert \tilde{\by}_{\iota}^{(t)}-\hat{y}_{\iota}^{(t)}\right\Vert^2\right]+\dfrac{20(L_{\nabla a}^*)^2(L_{\nabla a}')^2\omega_B^2}{(1-\omega_B)^2}\sum_{\iota=e}^{E-1}\sum_{t=T_0}^{T_1}\mathbb{E}\left[\left\Vert w_{\iota+1}^{(t+1)}-w_{\iota+1}^{(t)}\right\Vert^2\right]\\
    &+\dfrac{80(L_{\nabla a}^*)^2(L_{\nabla a}')^2\omega_B^2}{(1-\omega_B)^2(1-\omega_F)^2}\sum_{\iota=e}^{E-1}\sum_{t=T_0}^{T_1}\mathbb{E}\left[\left\Vert\by_\iota^{(t+1)}-\by_\iota^{(t)}\right\Vert^2\right]+\dfrac{80(L_{\nabla a}^*)^2(L_{\nabla a}')^2\omega_B^2}{(1-\omega_B)^2(1-\omega_F)}\sum_{\iota=e}^{E-1}\mathbb{E}\left[\left\Vert\tilde{\by}_\iota^{(T_0)}-\by_\iota^{(T_0)}\right\Vert^2\right]\\
    &+\dfrac{5(L_{\nabla a}^*)^2\omega_B^2}{(1-\omega_B)^2}\sum_{\iota=e}^{E-2}\mathbb{E}\left[\left\Vert\tilde{ v}_{\iota+1}^{(T_0)}-\hat{v}_{\iota+1}^{(T_0)}\right\Vert^2\right]+2(L_{\nabla a}^*)^2\dfrac{\omega_B}{1-\omega_B}\sum_{\iota=e}^{E-1}\mathbb{E}\left[\left\Vert\tilde{ v}^{(T_0)}_{\iota}-\bv^{(T_0)}_{\iota}\right\Vert^2\right]
\end{aligned}
\end{equation}}
holds for $e=1,2,\cdots,E-1$, where $(L_{\nabla a}^*)^2:=\max\{1,(2(L_{\nabla a}^{\circ})^2)^{E}\}\geq0$.

For $e=1,2,\cdots,E-1$, define $\{C_{v,e}\}$ as:
\begin{align}
\label{def_cve}
C_{v,1}=1,\quad C_{v,e}=1+\sum_{\iota=1}^{e-1}C_{v,\iota}\dfrac{5(L_{\nabla a}^*)^2{\omega_B^2}}{(1-\omega_B)^2}\quad(e=2,\cdots,E-1).
\end{align}
Then taking summation on both sides of \eqref{v_tilde-v1} over $e=1,2,\cdots,E-1$, we have:
{\small
\begin{equation}
\label{v_tilde-v2}
\begin{aligned}
    &\sum_{e=1}^{E-1}C_{v,e}\sum_{t=T_0+1}^{T_1+1}\mathbb{E}\left[\left\Vert\tilde{ v}^{(t)}_e-\hat{v}^{(t)}_e\right\Vert^2\right]\leq\sum_{e=1}^{E-1}\left(\sum_{\iota=1}^{e-1}C_{v,\iota}\dfrac{5(L_{\nabla a}^*)^2\omega_B^2}{(1-\omega_B)^2}\right)\sum_{t=T_0+1}^{T_1+1}\mathbb{E}\left[\left\Vert\tilde{ v}_{e}^{(t)}-\hat{v}_{e}^{(t)}\right\Vert^2\right]\\
    &\hspace{2.5mm}+\sum_{e=1}^{E-1}\left(\sum_{\iota=1}^{e-1}C_{v,\iota}\dfrac{5(L_{\nabla a}^*)^2\omega_B^2}{(1-\omega_B)^2}\right)\sum_{t=T_0}^{T_1}\mathbb{E}\left[\left\Vert\hat{v}_{e}^{(t+1)}-\hat{v}_{e}^{(t)}\right\Vert^2\right]\\
    &\hspace{2.5mm}+\sum_{e=1}^{E-1}\left(\sum_{\iota=1}^{e}C_{v,\iota}(L_{\nabla a}^*)^2(2(L_{\nabla a}')^2)\right)\sum_{t=T_0+1}^{T_1+1}\mathbb{E}\left[\left\Vert \tilde{\by}_{e}^{(t)}-\hat{y}_{e}^{(t)}\right\Vert^2\right]\\
    &\hspace{2.5mm}+\sum_{e=1}^{E}\left(\sum_{\iota=1}^{e-1}C_{v,\iota}\dfrac{20(L_{\nabla a}^*)^2(L_{\nabla a}')^2\omega_B^2}{(1-\omega_B)^2}\right)\sum_{t=T_0}^{T_1}\mathbb{E}\left[\left\Vert w_{e}^{(t+1)}-w_{e}^{(t)}\right\Vert^2\right]\\
    &\hspace{2.5mm}+\sum_{e=1}^{E-1}\left(\sum_{\iota=1}^{e}C_{v,\iota}\dfrac{80(L_{\nabla a}^*)^2(L_{\nabla a}')^2\omega_B^2}{(1-\omega_B)^2(1-\omega_F)^2}\right)\sum_{t=T_0}^{T_1}\mathbb{E}\left[\left\Vert\by_e^{(t+1)}-\by_e^{(t)}\right\Vert^2\right]\\
    &\hspace{2.5mm}+\hspace{-0.5mm}\sum_{e=1}^{E-1}\left(\sum_{\iota=1}^{e-1}C_{v,\iota}\dfrac{5(L_{\nabla a}^*)^2\omega_B^2}{(1-\omega_B)^2}\right)\mathbb{E}\left[\left\Vert\tilde{ v}_{e}^{(T_0)}\hspace{-0.5mm}-\hspace{-0.5mm}\hat{v}_{e}^{(T_0)}\right\Vert^2\right]\hspace{-0.5mm}+\hspace{-0.5mm}\sum_{e=1}^{E-1}\left(\sum_{\iota=1}^{e}C_{v,\iota}2(L_{\nabla a}^*)^2\dfrac{\omega_B}{1-\omega_B}\right)\mathbb{E}\left[\left\Vert\tilde{ v}^{(T_0)}_{e}\hspace{-0.5mm}-\hspace{-0.5mm}\bv^{(T_0)}_{e}\right\Vert^2\right]\\
    &\hspace{2.5mm}+\sum_{e=1}^{E-1}\left(\sum_{\iota=1}^{e}C_{v,\iota}\dfrac{80(L_{\nabla a}^*)^2(L_{\nabla a}')^2\omega_B^2}{(1-\omega_B)^2(1-\omega_F)}\right)\mathbb{E}\left[\left\Vert\tilde{\by}_e^{(T_0)}-\by_e^{(T_0)}\right\Vert^2\right].
\end{aligned}
\end{equation}}
Then, with the definition of $C_{v,e}$ in Eq. \eqref{def_cve}, we have:
{\small
\begin{equation*}
\begin{aligned}
    &\sum_{e=1}^{E-1}\sum_{t=T_0+1}^{T_1+1}\mathbb{E}\left[\left\Vert\tilde{ v}^{(t)}_e\hspace{-0.5mm}-\hspace{-0.5mm}\hat{v}^{(t)}_e\right\Vert^2\right]\\
    \leq&\sum_{e=1}^{E-1}C_{y,e}^{\circ}\sum_{t=T_0+1}^{T_1+1}\mathbb{E}\left[\left\Vert \tilde{\by}_{e}^{(t)}\hspace{-0.5mm}-\hspace{-0.5mm}\hat{y}_{e}^{(t)}\right\Vert^2\right]+\sum_{e=1}^{E}C_{W,e}^{\circ}\sum_{t=T_0}^{T_1}\mathbb{E}\left[\left\Vert w_{e}^{(t+1)}\hspace{-0.5mm}-\hspace{-0.5mm}w_{e}^{(t)}\right\Vert^2\right]\\
    &+\sum_{e=1}^{E-1}C_{v,e}^{\circ}\sum_{t=T_0}^{T_1}\mathbb{E}\left[\left\Vert\hat{v}_{e}^{(t+1)}\hspace{-0.5mm}-\hspace{-0.5mm}\hat{v}_{e}^{(t)}\right\Vert^2\right]+\sum_{e=1}^{E-1}C_{\theta,e}^{\circ}
    \sum_{t=T_0}^{T_1}\mathbb{E}\left[\left\Vert\by_e^{(t+1)}\hspace{-0.5mm}-\hspace{-0.5mm}\by_e^{(t)}\right\Vert^2\right]+\sum_{e=1}^{E-1}C_{\theta,e}^{1}\mathbb{E}\left[\left\Vert\tilde{\by}_e^{(T_0)}\hspace{-0.5mm}-\hspace{-0.5mm}\by_e^{(T_0)}\right\Vert^2\right]\\
    &+\sum_{e=1}^{E-1}C_{v,e}^{1}\mathbb{E}\left[\left\Vert\tilde{ v}_{e}^{(T_0)}\hspace{-0.5mm}-\hspace{-0.5mm}\hat{v}_{e}^{(T_0)}\right\Vert^2\right]+\sum_{e=1}^{E-1}C_{v,e}^{2}\mathbb{E}\left[\left\Vert\tilde{ v}_{e}^{(T_0)}\hspace{-0.5mm}-\hspace{-0.5mm}\bv_{e}^{(T_0)}\right\Vert^2\right],
\end{aligned}
\end{equation*}}
where the coefficients are defined as follows:
\begin{subequations}
\label{quzhi_1}
\begin{align}
    C_{v,e}^{\circ}&=\sum_{\iota=1}^{e-1}C_{v,\iota}\dfrac{5(L_{\nabla a}^*)^2{\omega_B^2}}{(1-\omega_B)^2}\leq\mathcal{O}\left(\dfrac{{\omega_B^{2}}}{(1-\omega_B)^{2(e-1)}}\right),\\
    C_{v,e}^{1}&=\sum_{\iota=1}^{e-1}C_{v,\iota}\dfrac{5(L_{\nabla a}^*)^2{\omega_B^2}}{(1-\omega_B)^2}\leq\mathcal{O}\left(\dfrac{{\omega_B^{2}}}{(1-\omega_B)^{2(e-1)}}\right),\\
    C_{v,e}^{2}&=\sum_{\iota=1}^{e}C_{v,\iota}\dfrac{2(L_{\nabla a}^*)^2{\omega_B}}{1-\omega_B}\leq\mathcal{O}\left(\dfrac{{\omega_B}}{(1-\omega_B)^{2e-1}}\right),\\
    C_{W,e}^{\circ}&=\sum_{\iota=1}^{e-1}C_{v,\iota}\dfrac{20(L_{\nabla a}^*)^2(L_{\nabla a}')^2\omega_B^2}{(1-\omega_B)^2}\leq\mathcal{O}\left(\dfrac{{\omega_B^{2}}}{(1-\omega_B)^{2(e-1)}}\right),\\
    C_{y,e}^{\circ}&=\sum_{\iota=1}^{e}C_{v,\iota}(L_{\nabla a}^*)^2(2(L_{\nabla a}')^2)\leq\mathcal{O}\left(\dfrac{{1}}{(1-\omega_B)^{2(e-1)}}\right),\quad\hspace{0.3mm}\\
    C_{\theta,e}^{\circ}&=\sum_{\iota=1}^{e}C_{v,\iota}\dfrac{80(L_{\nabla a}^*)^2(L_{\nabla a}')^2\omega_B^2}{(1-\omega_B)^2(1-\omega_F)^2}\leq\mathcal{O}\left(\dfrac{{\omega_B^2}}{(1-\omega_B)^{2e}(1-\omega_F)^2}\right),\\
    C_{\theta,e}^{1}&=\sum_{\iota=1}^{e}C_{v,\iota}\dfrac{80(L_{\nabla a}^*)^2(L_{\nabla a}')^2\omega_B^2}{(1-\omega_B)^2(1-\omega_F)}\leq\mathcal{O}\left(\dfrac{{\omega_B^{2}}}{(1-\omega_B)^{2e}(1-\omega_F)}\right).
\end{align}
\end{subequations}
\end{proof}
\end{lemma}

Now, we try to combine the result of Eq. \eqref{v_tilde-v_3} with Eq. \eqref{u-nablal}. For \oursfulast, it holds for any $r=1,2,\cdots,r_0$ that:

\begin{equation}
\label{v_tilde-v1_1}
\begin{aligned}
    &\sum_{e=1}^{E-1}\sum_{t=Q_r+1}^{Q_{r+1}}\mathbb{E}\left[\left\Vert\tilde{ v}^{(t)}_e\hspace{-0.5mm}-\hspace{-0.5mm}\hat{v}^{(t)}_e\right\Vert^2\right]=\sum_{e=1}^{E-1}\sum_{t=Q_r+1}^{Q_{r+1}-1}\mathbb{E}\left[\left\Vert\tilde{ v}^{(t)}_e\hspace{-0.5mm}-\hspace{-0.5mm}\hat{v}^{(t)}_e\right\Vert^2\right]\\
    \leq&\sum_{e=1}^{E-1}C_{y,e}^{\circ}\sum_{t=Q_r+1}^{Q_{r+1}-1}\mathbb{E}\left[\left\Vert \tilde{\by}_{e}^{(t)}\hspace{-0.5mm}-\hspace{-0.5mm}\hat{y}_{e}^{(t)}\right\Vert^2\right]+\sum_{e=1}^{E}C_{W,e}^{\circ}\sum_{t=Q_r}^{Q_{r+1}-2}\mathbb{E}\left[\left\Vert w_{e}^{(t+1)}\hspace{-0.5mm}-\hspace{-0.5mm}w_{e}^{(t)}\right\Vert^2\right]\\
    &+\sum_{e=1}^{E-1}C_{v,e}^{\circ}\sum_{t=Q_r}^{Q_{r+1}-2}\mathbb{E}\left[\left\Vert\hat{v}_{e}^{(t+1)}\hspace{-0.5mm}-\hspace{-0.5mm}\hat{v}_{e}^{(t)}\right\Vert^2\right]+\sum_{e=1}^{E-1}C_{\theta,e}^{\circ}
    \sum_{t=Q_r}^{Q_{r+1}-2}\mathbb{E}\left[\left\Vert\by_e^{(t+1)}\hspace{-0.5mm}-\hspace{-0.5mm}\by_e^{(t)}\right\Vert^2\right],
\end{aligned}
\end{equation}
where it holds due to the fact that in \oursfu at $Q_r$-iteration it satisfies that
\begin{align*}
    \tilde{v}_e^{(Q_r)}={v}_e^{(Q_r)}=\hat{v}_e^{(Q_r)},\quad \tilde{y}_e^{(Q_r)}={y}_e^{(Q_r)}=\hat{y}_e^{(Q_r)}.
\end{align*}
Taking summation over $r$, it holds that:
\begin{equation}
\label{v_tilde-v1_1_fu}
\begin{aligned}
    &\sum_{e=1}^{E-1}\sum_{t=2}^{T+1}\mathbb{E}\left[\left\Vert\tilde{ v}^{(t)}_e\hspace{-0.5mm}-\hspace{-0.5mm}\hat{v}^{(t)}_e\right\Vert^2\right]\\
    \leq&\sum_{e=1}^{E-1}C_{y,e}^{\circ}\sum_{t=2}^{T+1}\mathbb{E}\left[\left\Vert \tilde{\by}_{e}^{(t)}\hspace{-0.5mm}-\hspace{-0.5mm}\hat{y}_{e}^{(t)}\right\Vert^2\right]+\sum_{e=1}^{E}C_{W,e}^{\circ}\sum_{t=1}^{T}\mathbb{E}\left[\left\Vert w_{e}^{(t+1)}\hspace{-0.5mm}-\hspace{-0.5mm}w_{e}^{(t)}\right\Vert^2\right]\\
    &+\sum_{e=1}^{E-1}C_{v,e}^{\circ}\sum_{r=1}^{r_0}\sum_{t=Q_r}^{Q_{r+1}-2}\mathbb{E}\left[\left\Vert\hat{v}_{e}^{(t+1)}\hspace{-0.5mm}-\hspace{-0.5mm}\hat{v}_{e}^{(t)}\right\Vert^2\right]+\sum_{e=1}^{E-1}C_{\theta,e}^{\circ}
    \sum_{r=1}^{r_0}\sum_{t=Q_r}^{Q_{r+1}-2}\mathbb{E}\left[\left\Vert\by_e^{(t+1)}\hspace{-0.5mm}-\hspace{-0.5mm}\by_e^{(t)}\right\Vert^2\right].
\end{aligned}
\end{equation}

Plugging Eq. \eqref{v_tilde-v1_1_fu} into \eqref{u-nablal}, we can find that:
{\small
\begin{equation}
\label{u-nablal_1_fu}
\begin{aligned}
    &\sum_{e=1}^E\sum_{t=1}^{T+1}\mathbb{E}\left[\left\Vert\tilde{u}^{(t)}_e-\nabla_e\ell(\w^{(t)})\right\Vert^2\right]\\
    \leq&\sum_{e=1}^E\left(32L_{\nabla \ell}^2\left(\dfrac{p+m}{m^2(1-(1-p)(1-\frac{m}{2}))}+\dfrac{1}{m^2}\right)+8(L_{\nabla a}^\circ)^2C_{W,e}^{\circ}\right)\sum_{t=1}^{T}\mathbb{E}\left[\left\Vert w^{(t+1)}_e-w^{(t)}_e\right\Vert^2\right]\\
    &+\sum_{e=1}^{E-1}\left(8(L_{\nabla a}')^2+8(L_{\nabla a}^\circ)^2C_{y,e}^{\circ}\right)\sum_{t=2}^{T+1}\mathbb{E}\left[\left\Vert\tilde{\by}_e^{(t)}-\hat{y}_e^{(t)}\right\Vert^2\right]\\
    &+\sum_{e=1}^{E-1}8(L_{\nabla a}^\circ)^2C_{v,e}^{\circ}\sum_{r=1}^{r_0}\sum_{t=Q_r}^{Q_{r+1}-2}\mathbb{E}\left[\left\Vert\hat{v}_{e}^{(t+1)}\hspace{-0.5mm}-\hspace{-0.5mm}\hat{v}_{e}^{(t)}\right\Vert^2\right]+\sum_{e=1}^{E-1}8(L_{\nabla a}^\circ)^2C_{\theta,e}^{\circ}
    \sum_{r=1}^{r_0}\sum_{t=Q_r}^{Q_{r+1}-2}\mathbb{E}\left[\left\Vert\by_e^{(t+1)}\hspace{-0.5mm}-\hspace{-0.5mm}\by_e^{(t)}\right\Vert^2\right]\\
    &+4T\sigma^2\dfrac{(2-p)m-(1-p)m^2}{1-(1-p)(1-m)^2}+\dfrac{3}{m}\sum_{e=1}^E\mathbb{E}\left[\left\Vert\tilde{u}^{(1)}_e-\nabla_e\ell(\w^{(1)})\right\Vert^2\right].
\end{aligned}
\end{equation}}

Moreover for \oursfclast, it holds from Eq. \eqref{v_tilde-v_3} that: 
{\small\begin{equation}
\label{v_tilde-v1_1_fc}
\begin{aligned}
    &\sum_{e=1}^{E-1}\sum_{t=2}^{T+1}\mathbb{E}\left[\left\Vert\tilde{ v}^{(t)}_e\hspace{-0.5mm}-\hspace{-0.5mm}\hat{v}^{(t)}_e\right\Vert^2\right]\\
    \leq&\sum_{e=1}^{E-1}C_{y,e}^{\circ}\sum_{t=2}^{T+1}\mathbb{E}\left[\left\Vert \tilde{\by}_{e}^{(t)}\hspace{-0.5mm}-\hspace{-0.5mm}\hat{y}_{e}^{(t)}\right\Vert^2\right]+\sum_{e=1}^{E}C_{W,e}^{\circ}\sum_{t=1}^{T}\mathbb{E}\left[\left\Vert w_{e}^{(t+1)}\hspace{-0.5mm}-\hspace{-0.5mm}w_{e}^{(t)}\right\Vert^2\right]\\
    &+\sum_{e=1}^{E-1}C_{v,e}^{\circ}\sum_{t=1}^{T}\mathbb{E}\left[\left\Vert\hat{v}_{e}^{(t+1)}\hspace{-0.5mm}-\hspace{-0.5mm}\hat{v}_{e}^{(t)}\right\Vert^2\right]+\sum_{e=1}^{E-1}C_{\theta,e}^{\circ}
    \sum_{t=1}^{T}\mathbb{E}\left[\left\Vert\by_e^{(t+1)}\hspace{-0.5mm}-\hspace{-0.5mm}\by_e^{(t)}\right\Vert^2\right]+\sum_{e=1}^{E-1}C_{\theta,e}^{1}\mathbb{E}\left[\left\Vert\tilde{\by}_e^{(1)}\hspace{-0.5mm}-\hspace{-0.5mm}\by_e^{(1)}\right\Vert^2\right]\\
    &+\sum_{e=1}^{E-1}C_{v,e}^{1}\mathbb{E}\left[\left\Vert\tilde{ v}_{e}^{(1)}\hspace{-0.5mm}-\hspace{-0.5mm}\hat{v}_{e}^{(1)}\right\Vert^2\right]+\sum_{e=1}^{E-1}C_{v,e}^{2}\mathbb{E}\left[\left\Vert\tilde{ v}_{e}^{(1)}\hspace{-0.5mm}-\hspace{-0.5mm}\bv_{e}^{(1)}\right\Vert^2\right].
\end{aligned}
\end{equation}}

Plugging Eq. \eqref{v_tilde-v1_1_fc} into \eqref{u-nablal}, we can find that:
{\small
\begin{equation}
\label{u-nablal_1_fc}
\begin{aligned}
    &\sum_{e=1}^E\sum_{t=2}^{T+1}\mathbb{E}\left[\left\Vert\tilde{u}^{(t)}_e-\nabla_e\ell(\w^{(t)})\right\Vert^2\right]\\
    \leq&\sum_{e=1}^E\left(32L_{\nabla \ell}^2\left(\dfrac{p+m}{m^2(1-(1-p)(1-\frac{m}{2}))}+\dfrac{1}{m^2}\right)+8(L_{\nabla a}^\circ)^2C_{W,e}^{\circ}\right)\sum_{t=1}^{T}\mathbb{E}\left[\left\Vert w^{(t+1)}_e-w^{(t)}_e\right\Vert^2\right]\\
    &+\sum_{e=1}^{E-1}\left(8(L_{\nabla a}')^2+8(L_{\nabla a}^\circ)^2C_{y,e}^{\circ}\right)\sum_{t=2}^{T+1}\mathbb{E}\left[\left\Vert\tilde{\by}_e^{(t)}-\hat{y}_e^{(t)}\right\Vert^2\right]\\
    &+\sum_{e=1}^{E-1}8(L_{\nabla a}^\circ)^2C_{v,e}^{\circ}\sum_{t=1}^{T}\mathbb{E}\left[\left\Vert\hat{v}_{e}^{(t+1)}\hspace{-0.5mm}-\hspace{-0.5mm}\hat{v}_{e}^{(t)}\right\Vert^2\right]+\sum_{e=1}^{E-1}8(L_{\nabla a}^\circ)^2C_{\theta,e}^{\circ}
    \sum_{t=1}^{T}\mathbb{E}\left[\left\Vert\by_e^{(t+1)}\hspace{-0.5mm}-\hspace{-0.5mm}\by_e^{(t)}\right\Vert^2\right]\\
    &+4T\sigma^2\dfrac{(2\hspace{-0.5mm}-\hspace{-0.5mm}p)m\hspace{-0.5mm}-\hspace{-0.5mm}(1\hspace{-0.5mm}-\hspace{-0.5mm}p)m^2}{1\hspace{-0.5mm}-\hspace{-0.5mm}(1\hspace{-0.5mm}-\hspace{-0.5mm}p)(1\hspace{-0.5mm}-\hspace{-0.5mm}m)^2}\hspace{-0.5mm}+\hspace{-0.5mm}\dfrac{3}{m}\sum_{e=1}^E\mathbb{E}\left[\left\Vert\tilde{u}^{(1)}_e\hspace{-0.5mm}-\hspace{-0.5mm}\nabla_e\ell(\w^{(1)})\right\Vert^2\right]+\sum_{e=1}^{E-1}8(L_{\nabla a}^\circ)^2C_{\theta,e}^{1}\mathbb{E}\left[\left\Vert\tilde{\by}_e^{(1)}\hspace{-0.5mm}-\hspace{-0.5mm}\by_e^{(1)}\right\Vert^2\right]\\
    &+\sum_{e=1}^{E-1}8(L_{\nabla a}^\circ)^2C_{v,e}^{1}\mathbb{E}\left[\left\Vert\tilde{ v}_{e}^{(1)}\hspace{-0.5mm}-\hspace{-0.5mm}\hat{v}_{e}^{(1)}\right\Vert^2\right]+\sum_{e=1}^{E-1}8(L_{\nabla a}^\circ)^2C_{v,e}^{2}\mathbb{E}\left[\left\Vert\tilde{ v}_{e}^{(1)}\hspace{-0.5mm}-\hspace{-0.5mm}\bv_{e}^{(1)}\right\Vert^2\right].
\end{aligned}
\end{equation}}

\subsubsection{Forward variation to parameter/sample variation}
The bounds above still contain activation-variation terms. To close the proof, we convert these terms into parameter increments and, when samples change, the lazy-sampling difference \(\|x^{(t+1)}-x^{(t)}\|^2\). The precise combination that needs to be bounded is
{\small
\begin{align*}
&\sum_{e=1}^{E-1}\left(8(L_{\nabla a}')^2+8(L_{\nabla a}^\circ)^2C_{y,e}^{\circ}\right)\sum_{t=T_0+1}^{T_1+1}\mathbb{E}\left[\left\Vert \tilde{\by}_{e}^{(t)}-\hat{y}_{e}^{(t)}\right\Vert^2\right]+\sum_{e=1}^{E-1}8(L_{\nabla a}^\circ)^2C_{\theta,e}^{\circ}\sum_{t=T_0}^{T_1}\mathbb{E}\left[\left\Vert\by_e^{(t+1)}-\by_e^{(t)}\right\Vert^2\right]\\
+&\sum_{e=1}^{E-1}8(L_{\nabla a}^\circ)^2C_{v,e}^{\circ}\sum_{t=1}^{T}\mathbb{E}\left[\left\Vert\hat{v}_{e}^{(t+1)}\hspace{-0.5mm}-\hspace{-0.5mm}\hat{v}_{e}^{(t)}\right\Vert^2\right].
\end{align*}}
The next lemma provides this conversion for any \(T_1>T_0\geq1\).
\begin{lemma}
Suppose Assumption \ref{assumption:smoothness} and \ref{assumption:compressor} holds, then there exist coefficients $C_{y,\theta,e}^{\circ}$ and $C_{y,\theta,e}^{1}$ defined by \eqref{quzhi2} for each $e=1,2,\cdots,E-1$ and $T_1>T_0\geq1$ such that:
{\small
\begin{equation}
\label{yidui_y_theta_1}
\begin{aligned}
&\sum_{e=1}^{E-1}\left(8(L_{\nabla a}')^2+8(L_{\nabla a}^\circ)^2C_{y,e}^{\circ}\right)\sum_{t=T_0+1}^{T_1+1}\mathbb{E}\left[\left\Vert \tilde{\by}_{e}^{(t)}-\hat{y}_{e}^{(t)}\right\Vert^2\right]+\sum_{e=1}^{E-1}8(L_{\nabla a}^\circ)^2C_{\theta,e}^{\circ}\sum_{t=T_0}^{T_1}\mathbb{E}\left[\left\Vert\by_e^{(t+1)}-\by_e^{(t)}\right\Vert^2\right]\\
\leq&\sum_{e=1}^{E-1}L_a^2\left(\sum_{\iota=e}^{E-1}C_{y,\theta,\iota}^{\circ}\left(\dfrac{8L_a^2}{(1-\omega_F)^2}\right)^{\iota-e}\right)\sum_{t=T_0}^{T_1}\mathbb{E}\left[\left\Vert w_{e}^{(t+1)}-w_{e}^{(t)}\right\Vert^2\right]\\
&+\left(\sum_{\iota=1}^{E-1}C_{y,\theta,\iota}^{\circ}L_a^2\left(\dfrac{8L_a^2}{(1-\omega_F)^2}\right)^{\iota-1}\right) \sum_{t=T_0}^{T_1}\mathbb{E}\left[\left\Vert  x^{(t+1)}- x^{(t)}\right\Vert^2\right]\\
&+\sum_{e=1}^{E-1}\left(C_{y,\theta,e}^{1}+\sum_{\iota=e+1}^{E-1}\dfrac{8L_a^2}{1-\omega_F}\left(\dfrac{8L_a^2}{(1-\omega_F)^2}\right)^{\iota-e-1}\right)\mathbb{E}\left[\left\Vert\tilde{\by}^{(T_0)}_e-\by^{(T_0)}_e\right\Vert^2\right].
\end{aligned}
\end{equation}}
\begin{proof}
For $e=2,\cdots,E-1$ and $t=T_0,\cdots,T_1$, we have:
\begin{equation}
\begin{aligned}
&\left\Vert\by_e^{(t+1)}-\by_e^{(t)}\right\Vert^2= \left\Vert a_e(\tilde{\by}_{e-1}^{(t+1)},w_e^{(t+1)})-a_e(\tilde{\by}_{e-1}^{(t)},w_e^{(t)})\right\Vert^2\\
\leq&L_a^2\left\Vert \tilde{\by}_{e-1}^{(t+1)}-\tilde{\by}_{e-1}^{(t)}\right\Vert^2+L_a^2\left\Vert w_e^{(t+1)}-w_e^{(t)}\right\Vert^2.
\end{aligned}
\end{equation}
Taking expectation and then taking summation on boths over $t$, we can get:
{\small\begin{equation}
\begin{aligned}
&\sum_{t=T_0}^{T_1}\mathbb{E}\left[\left\Vert\by_e^{(t+1)}-\by_e^{(t)}\right\Vert^2\right] 
\leq L_a^2\sum_{t=T_0}^{T_1}\mathbb{E}\left[\left\Vert \tilde{\by}_{e-1}^{(t+1)}-\tilde{\by}_{e-1}^{(t)}\right\Vert^2\right]+L_a^2\sum_{t=T_0}^{T_1}\mathbb{E}\left[\left\Vert w_e^{(t+1)}-w_e^{(t)}\right\Vert^2\right]\\
\leq&\dfrac{8L_a^2}{(1-\omega_F)^2}\sum_{t=T_0}^{T_1}\mathbb{E}\left[\left\Vert\by_{e-1}^{(t+1)}-\by_{e-1}^{(t)}\right\Vert^2\right]+\dfrac{8L_a^2}{1-\omega_F}\mathbb{E}\left[\left\Vert\tilde{\by}_{e-1}^{(T_0)}-\by_{e-1}^{(T_0)}\right\Vert^2\right]\\
&+L_a^2\sum_{t=T_0}^{T_1}\mathbb{E}\left[\left\Vert w_e^{(t+1)}-w_e^{(t)}\right\Vert^2\right],
\end{aligned}
\end{equation}}
where the first inequality is due to Assumption \ref{assumption:smoothness} and the second inequality is due to \eqref{yt+1-yt_main}.

Then, we have:
{\small\begin{equation}
\begin{aligned}
\label{theta_t+1-theta_t_1}
&\sum_{t=T_0}^{T_1}\mathbb{E}\left[\left\Vert\by_e^{(t+1)}-\by_e^{(t)}\right\Vert^2\right]\\
\leq&\left(\dfrac{8L_a^2}{(1-\omega_F)^2}\right)^{e-1}\sum_{t=T_0}^{T_1}\mathbb{E}\left[\left\Vert\by_{1}^{(t+1)}-\by_{1}^{(t)}\right\Vert^2\right]+\sum_{\iota=1}^{e-1}\dfrac{8L_a^2}{1-\omega_F}\left(\dfrac{8L_a^2}{(1-\omega_F)^2}\right)^{e-1-\iota}\mathbb{E}\left[\left\Vert\tilde{\by}_{\iota}^{(T_0)}-\by_{\iota}^{(T_0)}\right\Vert^2\right]\\
&+\sum_{\iota=2}^{e}L_a^2\left(\dfrac{8L_a^2}{(1-\omega_F)^2}\right)^{e-\iota}\sum_{t=T_0}^{T_1}\mathbb{E}\left[\left\Vert w_{\iota}^{(t+1)}-w_{\iota}^{(t)}\right\Vert^2\right].
\end{aligned}
\end{equation}}

Moreover, from Assumption \ref{assumption:smoothness} we get:
{\small\begin{equation}
\begin{aligned}
\label{theta_1-theta_1}
&\left\Vert\by_1^{(t+1)}-\by_1^{(t)}\right\Vert^2= \left\Vert a_1( x^{(t+1)},w_1^{(t+1)})-a_1( x^{(t)},w_1^{(t)})\right\Vert^2\\
\leq&L_a^2\left\Vert  x^{(t+1)}- x^{(t)}\right\Vert^2+L_a^2\left\Vert w_1^{(t+1)}-w_1^{(t)}\right\Vert^2.
\end{aligned}
\end{equation}}

Plugging \eqref{theta_1-theta_1} into \eqref{theta_t+1-theta_t_1}, then we have:
{\small\begin{equation}
\begin{aligned}
\label{theta_t+1-theta_t_2}
&\sum_{t=T_0}^{T_1}\mathbb{E}\left[\left\Vert\by_e^{(t+1)}-\by_e^{(t)}\right\Vert^2\right]\\
\leq&\sum_{\iota=1}^{e}L_a^2\left(\dfrac{8L_a^2}{(1\hspace{-0.5mm}-\hspace{-0.5mm}\omega_F)^2}\right)^{e-\iota}\sum_{t=T_0}^{T_1}\mathbb{E}\left[\left\Vert w_{\iota}^{(t+1)}\hspace{-0.5mm}-\hspace{-0.5mm}w_{\iota}^{(t)}\right\Vert^2\right]\hspace{-0.5mm}+\hspace{-0.5mm}L_a^2\left(\dfrac{8L_a^2}{(1\hspace{-0.5mm}-\hspace{-0.5mm}\omega_F)^2}\right)^{e-1}\sum_{t=T_0}^{T_1}\mathbb{E}\left[\left\Vert  x^{(t+1)}\hspace{-0.5mm}-\hspace{-0.5mm}x^{(t)}\right\Vert^2\right]\\
&+\sum_{\iota=1}^{e-1}\dfrac{8L_a^2}{1\hspace{-0.5mm}-\hspace{-0.5mm}\omega_F}\left(\dfrac{8L_a^2}{(1\hspace{-0.5mm}-\hspace{-0.5mm}\omega_F)^2}\right)^{e-1-\iota}\mathbb{E}\left[\left\Vert\tilde{\by}_{\iota}^{(T_0)}\hspace{-0.5mm}-\hspace{-0.5mm}\by_{\iota}^{(T_0)}\right\Vert^2\right].
\end{aligned}
\end{equation}}

Then, consider Eq. \eqref{zhankai_y_21}, we have:
{\small
\begin{equation}
\label{yidui_y_theta}
\begin{aligned}
&\sum_{e=1}^{E-1}\left(8(L_{\nabla a}')^2+8(L_{\nabla a}^\circ)^2C_{y,e}^{\circ}\right)\sum_{t=T_0+1}^{T_1+1}\mathbb{E}\left[\left\Vert \tilde{\by}_{e}^{(t)}-\hat{y}_{e}^{(t)}\right\Vert^2\right]+\sum_{e=1}^{E-1}8(L_{\nabla a}^\circ)^2C_{\theta,e}^{\circ}\sum_{t=T_0}^{T_1}\mathbb{E}\left[\left\Vert\by_e^{(t+1)}-\by_e^{(t)}\right\Vert^2\right]\\
\leq&\sum_{e=1}^{E-1}\left(8(L_{\nabla a}')^2+8(L_{\nabla a}^\circ)^2C_{y,e}^{\circ}\right)\sum_{\iota=1}^{e}2(2L_a^2)^{e-\iota}\dfrac{\omega_F^2}{(1-\omega_F)^2}\sum_{t=T_0}^{T_1}\mathbb{E}\left[\left\Vert\by^{(t+1)}_\iota-\by^{(t)}_\iota\right\Vert^2\right]\\
&+\sum_{e=1}^{E-1}\left(8(L_{\nabla a}')^2+8(L_{\nabla a}^\circ)^2C_{y,e}^{\circ}\right)\sum_{\iota=1}^{e}2(2L_a^2)^{e-\iota}\dfrac{\omega_F}{1-\omega_F}\mathbb{E}\left[\left\Vert\tilde{\by}^{(T_0)}_{\iota}-\by^{(T_0)}_\iota\right\Vert^2\right]\\
&+\sum_{e=1}^{E-1}8(L_{\nabla a}^\circ)^2C_{\theta,e}^{\circ}\sum_{t=T_0}^{T_1}\mathbb{E}\left[\left\Vert\by_e^{(t+1)}-\by_e^{(t)}\right\Vert^2\right]\\
\leq&\sum_{e=1}^{E-1}C_{y,\theta,e}^{\circ}\sum_{t=T_0}^{T_1}\mathbb{E}\left[\left\Vert\by_e^{(t+1)}-\by_e^{(t)}\right\Vert^2\right]+\sum_{e=1}^{E-1}C_{y,\theta,e}^{1}\mathbb{E}\left[\left\Vert\tilde{\by}^{(T_0)}_e-\by^{(T_0)}_e\right\Vert^2\right],
\end{aligned}
\end{equation}}
where
\begin{subequations}
\label{quzhi2}
\begin{align}
    C_{y,\theta,e}^{\circ}&=8(L_{\nabla a}^\circ)^2C_{\theta,e}^{\circ}+\dfrac{2\omega_F^2\sum_{\iota=e}^{E-1}\left(8(L_{\nabla a}')^2+8(L_{\nabla a}^\circ)^2C_{y,\iota}^{\circ}\right)(2L_a^2)^{\iota-e}}{(1-\omega_F)^2},\\
    C_{y,\theta,e}^{1}&=\dfrac{2\omega_F\sum_{\iota=e}^{E-1}\left(8(L_{\nabla a}')^2+8(L_{\nabla a}^\circ)^2C_{y,\iota}^{\circ}\right)(2L_a^2)^{\iota-e}}{1-\omega_F}.
\end{align}
\end{subequations}

From \eqref{quzhi_1} we know that:
\begin{align*}
C_{y,e}^{\circ}&\leq\mathcal{O}\left(\dfrac{{1}}{(1-\omega_B)^{2(e-1)}}\right),\quad C_{\theta,e}^{\circ}\leq\mathcal{O}\left(\dfrac{{\omega_B^{2}}}{(1-\omega_B)^{2e}(1-\omega_F)^2}\right).\\
C_{y,\theta,e}^{\circ}&\leq\mathcal{O}\left(\dfrac{{\omega_B^{2}}}{(1-\omega_B)^{2}(1-\omega_F)^2}+\dfrac{{\omega_F^2}}{(1-\omega_B)^{2(E-2)}(1-\omega_F)^2}\right),\\
C_{y,\theta,e}^{1}&\leq\mathcal{O}\left(\dfrac{{\omega_F^2}}{(1-\omega_B)^{2(E-2)}(1-\omega_F)^2}\right).
\end{align*}
Then, plugging \eqref{theta_t+1-theta_t_2} into \eqref{yidui_y_theta}, we can get:
{\small
\begin{equation*}
\begin{aligned}
&\sum_{e=1}^{E-1}\left(8(L_{\nabla a}')^2+8(L_{\nabla a}^\circ)^2C_{y,e}^{\circ}\right)\sum_{t=T_0+1}^{T_1+1}\mathbb{E}\left[\left\Vert \tilde{\by}_{e}^{(t)}-\hat{y}_{e}^{(t)}\right\Vert^2\right]+\sum_{e=1}^{E-1}8(L_{\nabla a}^\circ)^2C_{\theta,e}^{\circ}\sum_{t=T_0}^{T_1}\mathbb{E}\left[\left\Vert\by_e^{(t+1)}-\by_e^{(t)}\right\Vert^2\right]\\
\leq&\sum_{e=1}^{E-1}L_a^2\left(\sum_{\iota=e}^{E-1}C_{y,\theta,\iota}^{\circ}\left(\dfrac{8L_a^2}{(1-\omega_F)^2}\right)^{\iota-e}\right)\sum_{t=T_0}^{T_1}\mathbb{E}\left[\left\Vert w_{e}^{(t+1)}-w_{e}^{(t)}\right\Vert^2\right]\\
&+\left(\sum_{\iota=1}^{E-1}C_{y,\theta,\iota}^{\circ}L_a^2\left(\dfrac{8L_a^2}{(1-\omega_F)^2}\right)^{\iota-1}\right) \sum_{t=T_0}^{T_1}\mathbb{E}\left[\left\Vert  x^{(t+1)}- x^{(t)}\right\Vert^2\right]\\
&+\sum_{e=1}^{E-1}\left(C_{y,\theta,e}^{1}+\sum_{\iota=e+1}^{E-1}\dfrac{8L_a^2}{1-\omega_F}\left(\dfrac{8L_a^2}{(1-\omega_F)^2}\right)^{\iota-e-1}\right)\mathbb{E}\left[\left\Vert\tilde{\by}^{(T_0)}_e-\by^{(T_0)}_e\right\Vert^2\right].
\end{aligned}
\end{equation*}}

\end{proof}
\end{lemma}

\subsubsection{Exact-adjoint variation to parameter/sample variation}
The last remaining variation term is the change of the no-compression adjoint. Since it follows the exact pipeline recursion, it can be bounded directly by sample changes and parameter increments.
\begin{lemma}
\label{lemma12}
Suppose Assumption \ref{assumption:smoothness} and \ref{assumption:compressor} holds, then there exist coefficients $C_{v,x}^{\circ}$ and $C_{v,w,e}^{\circ}$ defined by \eqref{quzhi3} for each $e=1,2,\cdots,E$ and $T_1>T_0\geq1$ such that:
\begin{equation}
\label{y_t+1-y_t_notilde_final}
\begin{aligned}
&\sum_{e=1}^{E-1}8(L_{\nabla a}^\circ)^2C_{v,e}^{\circ}\sum_{t=T_0}^{T_1}\mathbb{E}\left[\left\Vert\hat{v}_{e}^{(t+1)}-\hat{v}_{e}^{(t)}\right\Vert^2\right]\\
\leq&C_{v,x}^{\circ}\sum_{t=T_0}^{T_1}\mathbb{E}\left[\left\Vert x^{(t+1)}- x^{(t)}\right\Vert^2\right]+\sum_{e=1}^EC_{v,w,e}^{\circ}\sum_{t=T_0}^{T_1}\mathbb{E}\left[\left\Vert w_e^{(t+1)}-w_e^{(t)}\right\Vert^2\right].
\end{aligned}
\end{equation}
\begin{proof}
For every $e=1,2,\cdots,E-2$ and $t=1,2,\cdots,T$ and $T=T_0,\cdots,T_1$, we have:
\begin{equation}
\label{v_t+1-v_t_notilde}
\begin{aligned}
    \left\Vert\hat{v}_{e}^{(t+1)}-\hat{v}_{e}^{(t)}\right\Vert^2=&\left\Vert\nabla_1a_{e+1}(\hat{y}_e^{(t+1)},W_{e+1}^{(t+1)})^{\tran}\hat{v}_{e+1}^{(t+1)}-\nabla_1a_{e+1}(\hat{y}_e^{(t)},W_{e+1}^{(t)})^{\tran}\hat{v}_{e+1}^{(t)}\right\Vert^2\\
    \leq&(L_{\nabla a}')^2\left(\left\Vert\hat{y}_e^{(t+1)}-\hat{y}_e^{(t)}\right\Vert^2+\left\Vert w_{e+1}^{(t+1)}-w_{e+1}^{(t)}\right\Vert^2\right)+(L_{\nabla a}^{\circ})^2\left\Vert\hat{v}_{e+1}^{(t+1)}-\hat{v}_{e+1}^{(t)}\right\Vert^2,
\end{aligned}
\end{equation}
where the first inequality is due to \eqref{estimation of a_e^0}.

Then we can get:
\begin{equation}
\label{v_t+1-v_t_notilde}
\begin{aligned}
    \left\Vert\hat{v}_{e}^{(t+1)}-\hat{v}_{e}^{(t)}\right\Vert^2
    \leq&\sum_{\iota=e}^{E-1}(L_{\nabla a}')^2((L_{\nabla a}^{\circ})^2)^{\iota-e}\left\Vert\hat{y}_\iota^{(t+1)}-\hat{y}_\iota^{(t)}\right\Vert^2\\
    &+\sum_{\iota=e+1}^{E}(L_{\nabla a}')^2((L_{\nabla a}^{\circ})^2)^{\iota-e-1}\left\Vert w_\iota^{(t+1)}-w_\iota^{(t)}\right\Vert^2,
\end{aligned}
\end{equation}
where the second inequality is due to the definition of $ \hat{v}_{E-1}^{(t)}$ and Assumption \ref{assumption:smoothness}. And we can know that \eqref{v_t+1-v_t_notilde} also holds in the case of $e=E-1$.

Then we consider the term $\left\Vert\hat{y}_{e}^{(t+1)}-\hat{y}_{e}^{(t)}\right\Vert^2$. For $e=1,2,\cdots,E-1$, we obtain from Assumption \ref{assumption:smoothness} that:
\begin{equation}
\label{y_t+1-y_t_notilde}
\begin{aligned}
    \left\Vert\hat{y}_{e}^{(t+1)}-\hat{y}_{e}^{(t)}\right\Vert^2=&\left\Vert a_e(\hat{y}_{e-1}^{(t+1)},w_e^{(t+1)})-a_e(\hat{y}_{e-1}^{(t)},w_e^{(t)})\right\Vert^2\\
    \leq& L_a^2\left\Vert \hat{y}_{e-1}^{(t+1)}-\hat{y}_{e-1}^{(t)}\right\Vert^2+L_a^2\left\Vert{w}_{e}^{(t+1)}-{w}_{e}^{(t)}\right\Vert^2.
\end{aligned}
\end{equation}
Then we have:
\begin{equation}
\label{y_t+1-y_t_notilde_1}
\begin{aligned}
    \left\Vert\hat{y}_{e}^{(t+1)}-\hat{y}_{e}^{(t)}\right\Vert^2\leq \sum_{\iota=1}^e(L_a^2)^{e-\iota+1}\left\Vert {w}_{\iota}^{(t+1)}-{w}_{\iota}^{(t)}\right\Vert^2+(L_a^2)^e\left\Vert{ x}^{(t+1)}- x^{(t)}\right\Vert^2.
\end{aligned}
\end{equation}

Combining \eqref{v_t+1-v_t_notilde} and \eqref{y_t+1-y_t_notilde_1} together, then we can get:
\begin{equation}
\label{y_t+1-y_t_notilde_2}
\begin{aligned}
&\sum_{e=1}^{E-1}8(L_{\nabla a}^\circ)^2C_{v,e}^{\circ}\sum_{t=1}^{T}\mathbb{E}\left[\left\Vert\hat{v}_{e}^{(t+1)}-\hat{v}_{e}^{(t)}\right\Vert^2\right]\\
\leq&\sum_{e=1}^{E-1}\left(\sum_{\iota=1}^e8C_{v,\iota}^{\circ}(L_{\nabla a}')^2((L_{\nabla a}^{\circ})^2)^{e-\iota+1}\right)\sum_{t=1}^{T}\mathbb{E}\left[\left\Vert\hat{y}_{e}^{(t+1)}-\hat{y}_{e}^{(t)}\right\Vert^2\right]\\
&+\sum_{e=2}^{E}\left(\sum_{\iota=1}^{e-1}8C_{v,\iota}^{\circ}(L_{\nabla a}')^2((L_{\nabla a}^{\circ})^2)^{e-\iota}\right)\sum_{t=1}^{T}\mathbb{E}\left[\left\Vert w_e^{(t+1)}-w_e^{(t)}\right\Vert^2\right]\\
\leq&\sum_{e=1}^{E-1}\left(\sum_{\iota=1}^e8C_{v,\iota}^{\circ}(L_{\nabla a}')^2((L_{\nabla a}^{\circ})^2)^{e-\iota+1}\right)\\
&\hspace{3.5cm}\cdot\sum_{t=1}^{T}\mathbb{E}\left[\sum_{\eta=1}^e(L_a^2)^{e-\eta+1}\left\Vert {w}_{\eta}^{(t+1)}-{w}_{\eta}^{(t)}\right\Vert^2+(L_a^2)^e\left\Vert{ x}^{(t+1)}- x^{(t)}\right\Vert^2\right]\\
&+\sum_{e=2}^{E}\left(\sum_{\iota=1}^{e-1}8C_{v,\iota}^{\circ}(L_{\nabla a}')^2((L_{\nabla a}^{\circ})^2)^{e-\iota}\right)\sum_{t=1}^{T}\mathbb{E}\left[\left\Vert w_e^{(t+1)}-w_e^{(t)}\right\Vert^2\right]\\
\leq&C_{v,x}^{\circ}\sum_{t=1}^T\mathbb{E}\left[\left\Vert x^{(t+1)}- x^{(t)}\right\Vert^2\right]+\sum_{e=1}^EC_{v,w,e}^{\circ}\sum_{t=1}^T\mathbb{E}\left[\left\Vert w_e^{(t+1)}-w_e^{(t)}\right\Vert^2\right],
\end{aligned}
\end{equation}
where
\begin{subequations}
\label{quzhi3}
\begin{align}
C_{v,x}^{\circ} &= \sum_{e=1}^{E-1}\left(\sum_{\iota=1}^e8C_{v,\iota}^{\circ}(L_{\nabla a}')^2((L_{\nabla a}^{\circ})^2)^{e-\iota+1}\right)(L_a^2)^e,\\
C_{v,w,e}^{\circ} &= 
\sum_{\iota=1}^{e-1}8C_{v,\iota}^{\circ}(L_{\nabla a}')^2((L_{\nabla a}^{\circ})^2)^{e-\iota}+\sum_{\eta=e}^{E-1}\left(\sum_{\iota=1}^\eta8C_{v,\iota}^{\circ}(L_{\nabla a}')^2((L_{\nabla a}^{\circ})^2)^{\eta-\iota+1}\right)(L_a^2)^{\eta-e+1}.
\end{align}
\end{subequations}
And from \eqref{quzhi_1} we know that:
\begin{align*}
C_{v,x}^{\circ}\leq\mathcal{O}\left(\dfrac{{\omega_B^{2}}}{(1-\omega_B)^{2(E-2)}}\right),\quad C_{v,w,e}^{\circ}\leq\mathcal{O}\left(\dfrac{{\omega_B^{2}}}{(1-\omega_B)^{2(E-2)}}\right).   
\end{align*}
\end{proof}
\end{lemma}

\subsection{Master Convergence Bound}
\label{app:master-convergence-bound}
We now combine the descent bound, the evaluated-gradient decomposition, and the compression-error bounds. The two variants use the same ingredients, but refresh iterations affect their error-feedback states differently.
\subsubsection{Convergence rate of \oursfu}
\label{app:proof-fu-rate}
For \oursfu, a refresh iteration sends the exact activation/adjoint and resets the corresponding error-feedback state. Thus the initial compression errors vanish at the beginning of each refresh interval, and the sample-change term does not appear inside such intervals.
\begin{theorem}[Formal convergence rate of \textbf{\texttt{Clapping}-FU}]
\label{lemma: convergence_fu_e}
Suppose Assumption \ref{assumption:smoothness}-\ref{assumption:sample} hold. Then for Algorithm \ref{alg:clapping} there exist $\gamma,m,p>0$ such that: 
\begin{equation}
\label{descent lemma__GENERAL_e3}
\begin{aligned}
    \dfrac{1}{T}\sum_{t=1}^T\mathbb{E}\left[\left\Vert\nabla \ell(\w^{(t)})\right\Vert^2\right]\lesssim&\dfrac{\sigma}{\sqrt{T}}+\dfrac{1}{T(1-\omega_B)^{E-1}(1-\omega_F)^{E-1}}.
\end{aligned}
\end{equation}
\begin{proof}
Combining \eqref{yidui_y_theta_1}, and \eqref{y_t+1-y_t_notilde_final} together, it holds that:
{\small
\begin{equation}
\label{u-nablal_3_more_1}
\begin{aligned}
&\sum_{e=1}^{E-1}8(L_{\nabla a}^\circ)^2C_{v,e}^{\circ}\sum_{r=1}^{r_0}\sum_{t=Q_r}^{Q_{r+1}-2}\mathbb{E}\left[\left\Vert\hat{v}_{e}^{(t+1)}\hspace{-0.5mm}-\hspace{-0.5mm}\hat{v}_{e}^{(t)}\right\Vert^2\right]+\sum_{e=1}^{E-1}8(L_{\nabla a}^\circ)^2C_{\theta,e}^{\circ}\sum_{r=1}^{r_0}\sum_{t=Q_r}^{Q_{r+1}-2}\mathbb{E}\left[\left\Vert\by_e^{(t+1)}\hspace{-0.5mm}-\hspace{-0.5mm}\by_e^{(t)}\right\Vert^2\right]\\
&+\sum_{e=1}^{E-1}\left(8(L_{\nabla a}')^2+8(L_{\nabla a}^\circ)^2C_{y,e}^{\circ}\right)\sum_{t=2}^{T+1}\mathbb{E}\left[\left\Vert\tilde{\by}_e^{(t)}-\hat{y}_e^{(t)}\right\Vert^2\right]\\
\leq&\sum_{e=1}^{E}\left(L_a^2\left(\sum_{\iota=e}^{E-1}C_{y,\theta,\iota}^{\circ}\left(\dfrac{8L_a^2}{(1-\omega_F)^2}\right)^{\iota-e}\right)+C_{v,w,e}^{\circ}\right)\sum_{t=1}^{T}\mathbb{E}\left[\left\Vert w_{e}^{(t+1)}-w_{e}^{(t)}\right\Vert^2\right].
\end{aligned}
\end{equation}}
Eq. \eqref{u-nablal_3_more_1} holds is also due to the fact that $x^{(t+1)}= x^{(t)}$ for any $t=Q_r,\cdots,Q_{r+1}-2$, ${y}_e^{(Q_r)}=\hat{y}_e^{(Q_r)}=\tilde{y}_e^{(Q_r)}$, and ${v}_e^{(Q_r)}=\hat{v}_e^{(Q_r)}=\tilde{v}_e^{(Q_r)}$ for any $r=1,2,\cdots,r_0$.

Combining \eqref{u-nablal_1_fu}, \eqref{yidui_y_theta_1}, and \eqref{y_t+1-y_t_notilde_final} together, we can obtain:
\begin{equation}
\label{u-nablal_3_more}
    \begin{aligned}
    &\sum_{e=1}^E\sum_{t=1}^{T+1}\mathbb{E}\left[\left\Vert\tilde{u}^{(t)}_e-\nabla_e\ell(\w^{(t)})\right\Vert^2\right]\\
    \leq&\sum_{e=1}^E\left(32L_{\nabla \ell}^2\left(\dfrac{p+m}{m^2(1-(1-p)(1-\frac{m}{2}))}+\dfrac{1}{m^2}\right)+C_{w,e}\right)\sum_{t=1}^{T}\mathbb{E}\left[\left\Vert w^{(t+1)}_e-w^{(t)}_e\right\Vert^2\right]\\
    &+4T\sigma^2\dfrac{(2-p)m-(1-p)m^2}{1-(1-p)(1-m)^2}+\dfrac{3}{m}\sum_{e=1}^E\mathbb{E}\left[\left\Vert\tilde{u}^{(1)}_e-\nabla_e\ell(\w^{(1)})\right\Vert^2\right],
    \end{aligned}
\end{equation}
where
\begin{equation}
\begin{aligned}
C_{w,e}=&8(L_{\nabla a}^\circ)^2C_{W,e}^{\circ}+L_a^2\left(\sum_{\iota=e}^{E-1}C_{y,\theta,\iota}^{\circ}\left(\dfrac{8L_a^2}{(1-\omega_F)^2}\right)^{\iota-e}\right)+C_{v,w,e}^{\circ}\\
\lesssim&\dfrac{1}{(1-\omega_B)^{2(E-1)}(1-\omega_F)^{2(E-1)}}.
\end{aligned}
\end{equation}

Let $p=p_0=\mathcal{O}(1)$ as a constant with respect to $\sigma,T,\omega_F,\omega_B$. Then let 
\begin{equation*}
\begin{aligned}
    m \sim&\left(\dfrac{1}{(1-\omega_B)^{E-1}(1-\omega_F)^{E-1}} +\sigma\sqrt{T}\right)^{-1},\ m\leq 1,\\
    \gamma\sim&\left(\dfrac{1}{(1-\omega_B)^{E-1}(1-\omega_F)^{E-1}} +\sigma\sqrt{T}\right)^{-1},\ \gamma\leq 1.
\end{aligned}
\end{equation*}
At this time, if $\gamma/m$ is sufficiently small, we can obtain that \[32L_{\nabla \ell}^2\left(\dfrac{p+m}{m^2(1-(1-p)(1-\frac{m}{2}))}+\dfrac{1}{m^2}\right)+C_{w,e}-\dfrac{1}{2\gamma^2}\leq0.\]  Then we have from Eq. \eqref{descent lemma}:
\begin{equation*}
\begin{aligned}
    \dfrac{1}{T}\sum_{t=1}^T\mathbb{E}\left[\left\Vert\nabla \ell(\w^{(t)})\right\Vert^2\right]\lesssim&\dfrac{\sigma}{\sqrt{T}}+\dfrac{1}{T(1-\omega_B)^{E-1}(1-\omega_F)^{E-1}}.
\end{aligned}
\end{equation*}
\end{proof}
\end{theorem}

\subsubsection{Convergence rate of \oursfc}
\label{app:proof-fc-rate}
For \oursfc, compression is applied at every iteration, including refresh iterations. Therefore the proof keeps an additional sample-change term \(\sum_t\mathbb{E}\|x^{(t+1)}-x^{(t)}\|^2\), which is then controlled by lazy sampling and Assumption~\ref{assumption:sample}.
\begin{theorem}[Formal convergence rate of \oursfclast]
\label{lemma: convergence_fc_e}
Suppose Assumption \ref{assumption:smoothness}-\ref{assumption:sample} hold. Then for Algorithm \ref{alg:clapping} there exist $\gamma,m>0$ such that: 
\begin{equation}
\label{descent lemma__GENERAL_e3}
\begin{aligned}
    &\dfrac{1}{T}\sum_{t=1}^T\mathbb{E}\left[\left\Vert\nabla \ell(\w^{(t)})\right\Vert^2\right]\lesssim\dfrac{\sigma^{\frac{4}{3}}}{T^{\frac{1}{3}}(1-\omega_B)^{\frac{4(E-1)}{3}}(1-\omega_F)^{\frac{4(E-1)}{3}}}\\
    &\quad\quad+\dfrac{1}{T}\left(\dfrac{1}{(1-\omega_B)^{E-1}(1-\omega_F)^{E-1}}+\dfrac{\omega_F^2+\omega_B}{(1-\omega_B)^{2(E-1)}(1-\omega_F)^2}+\dfrac{1}{(1-\omega_F)^{2(E-2)-1}}\right).
\end{aligned}
\end{equation}
\begin{proof}
Combining \eqref{u-nablal_1_fc}, \eqref{yidui_y_theta_1}, and \eqref{y_t+1-y_t_notilde_final} together, we can obtain:
\begin{equation}
\label{u-nablal_3_more1}
    \begin{aligned}
    &\sum_{e=1}^E\sum_{t=2}^{T+1}\mathbb{E}\left[\left\Vert\tilde{u}^{(t)}_e-\nabla_e\ell(\w^{(t)})\right\Vert^2\right]\\
    \leq&\sum_{e=1}^E\left(32L_{\nabla \ell}^2\left(\dfrac{p+m}{m^2(1-(1-p)(1-\frac{m}{2}))}+\dfrac{1}{m^2}\right)+C_{w,e}\right)\sum_{t=1}^{T}\mathbb{E}\left[\left\Vert w^{(t+1)}_e-w^{(t)}_e\right\Vert^2\right]\\
    &+C_{x} \sum_{t=1}^T\mathbb{E}\left[\left\Vert  x^{(t+1)}-x^{(t)}\right\Vert^2\right]+\sum_{e=1}^{E-1}C_{y,\theta,e}\mathbb{E}\left[\left\Vert\tilde{\by}^{(1)}_e-\by^{(1)}_e\right\Vert^2\right]+\sum_{e=1}^{E-1}C_{v,e}\mathbb{E}\left[\left\Vert\tilde{ v}_{e}^{(1)}-\hat{v}_{e}^{(1)}\right\Vert^2\right]\\
    &+\sum_{e=1}^{E-1}C_{v,\chi,e}\mathbb{E}\left[\left\Vert\tilde{ v}_{e}^{(1)}-\bv_{e}^{(1)}\right\Vert^2\right]+4T\sigma^2\dfrac{(2-p)m-(1-p)m^2}{1-(1-p)(1-m)^2}\\
    &+\dfrac{3}{m}\sum_{e=1}^E\mathbb{E}\left[\left\Vert\tilde{u}^{(1)}_e-\nabla_e\ell(\w^{(1)})\right\Vert^2\right],
    \end{aligned}
\end{equation}
where
\begin{subequations}
\begin{align}
C_{w,e}=&8(L_{\nabla a}^\circ)^2C_{W,e}^{\circ}+L_a^2\left(\sum_{\iota=e}^{E-1}C_{y,\theta,\iota}^{\circ}\left(\dfrac{8L_a^2}{(1-\omega_F)^2}\right)^{\iota-e}\right)+C_{v,w,e}^{\circ}\nonumber\\
\lesssim&\dfrac{1}{(1-\omega_B)^{2(E-1)}(1-\omega_F)^{2(E-1)}},\\
C_{x}=&C_{v,x}^{\circ}+\sum_{\iota=1}^{E-1}C_{y,\theta,\iota}^{\circ}L_a^2\left(\dfrac{8L_a^2}{(1-\omega_F)^2}\right)^{\iota-1}\nonumber\\
\lesssim&\dfrac{\omega_F^2+\omega_B^2}{(1-\omega_B)^{2(E-1)}(1-\omega_F)^{2(E-1)}},\\
C_{y,\theta,e}=&8(L_{\nabla a}^\circ)^2C_{\theta,e}^{1}+C_{y,\theta,e}^{1}+\sum_{\iota=e+1}^{E-1}\dfrac{8L_a^2}{1-\omega_F}\left(\dfrac{8L_a^2}{(1-\omega_F)^2}\right)^{\iota-e-1} \nonumber\\
\lesssim&\dfrac{\omega_B^2}{(1-\omega_B)^{2e}(1-\omega_F)^2}+\dfrac{{\omega_F^2}}{(1-\omega_B)^{2(E-2)}(1-\omega_F)^2}+\dfrac{1}{(1-\omega_F)^{2(E-e-2)+1}},\\
C_{v,e}=&8(L_{\nabla a}^\circ)^2C_{v,e}^{1},\\
C_{v,\chi,e}=&8(L_{\nabla a}^\circ)^2C_{v,e}^{2}.
\end{align}
\end{subequations}

Thus, we know that:
\begin{equation}
\begin{aligned}
&\sum_{e=1}^{E-1}C_{y,\theta,e}\mathbb{E}\left[\left\Vert\tilde{\by}^{(1)}_e-\by^{(1)}_e\right\Vert^2\right]+\sum_{e=1}^{E-1}C_{v,e}\mathbb{E}\left[\left\Vert\tilde{ v}_{e}^{(1)}- \hat{v}_{e}^{(1)}\right\Vert^2\right]+\sum_{e=1}^{E-1}C_{v,\chi,e}\mathbb{E}\left[\left\Vert\tilde{ v}_{e}^{(1)}-\bv_{e}^{(1)}\right\Vert^2\right]\\
\lesssim&\dfrac{\omega_B}{(1-\omega_B)^{2(E-1)}(1-\omega_F)^2}+\dfrac{{\omega_F^2}}{(1-\omega_B)^{2(E-2)}(1-\omega_F)^2}+\dfrac{1}{(1-\omega_F)^{2(E-2)-1}}.
\end{aligned}
\end{equation}
Plugging \eqref{u-nablal_3_more} into \eqref{descent lemma}, and use the fact that $\dfrac{p+m}{1-(1-p)(1-\frac{m}{2})}\leq2$ when $p,m\leq1$, we can get:
\begin{equation}
\label{descent lemma_more}
\begin{aligned}
    &\dfrac{1}{T}\sum_{t=1}^T\mathbb{E}\left[\left\Vert\nabla \ell(\w^{(t)})\right\Vert^2\right]\\
    \leq&\dfrac{2}{\gamma T}\mathbb{E}\left[\ell(\w^{(1)})-\inf_{\w}l(\w)\right]+\dfrac{C_{x}}{T}\sum_{t=1}^T\mathbb{E}\left[\left\Vert  x^{(t+1)}- x^{(t)}\right\Vert^2\right]\\
    &+\dfrac{1}{T}\sum_{e=1}^E\left(96L_{\nabla \ell}^2\dfrac{1}{m^2}+C_{w,e}-\dfrac{1}{2\gamma^2}\right)\sum_{t=1}^{T}\mathbb{E}\left[\left\Vert w^{(t+1)}_e-w^{(t)}_e\right\Vert^2\right]\\
    &+\sum_{e=1}^{E-1}\dfrac{C_{y,\theta,e}}{T}\mathbb{E}\left[\left\Vert\tilde{\by}^{(1)}_e-\by^{(1)}_e\right\Vert^2\right]+\sum_{e=1}^{E-1}\dfrac{C_{v,e}}{T}\mathbb{E}\left[\left\Vert\tilde{ v}_{e}^{(1)}- \hat{v}_{e}^{(1)}\right\Vert^2\right]+\sum_{e=1}^{E-1}\dfrac{C_{v,\chi,e}}{T}\mathbb{E}\left[\left\Vert\tilde{ v}_{e}^{(1)}-\bv_{e}^{(1)}\right\Vert^2\right]\\
    &+4\sigma^2\dfrac{(2-p)m-(1-p)m^2}{1-(1-p)(1-m)^2}+\dfrac{3}{mT}\sum_{e=1}^E\mathbb{E}\left[\left\Vert\tilde{u}^{(1)}_e-\nabla_e\ell(\w^{(1)})\right\Vert^2\right]+\dfrac{1}{T}\mathbb{E}\left[\left\Vert\nabla \ell(\w^{(1)})\right\Vert^2\right].
\end{aligned}
\end{equation}
Taking $p=\sqrt{m}$ where $h$ is an undetermined function with respectively to $\omega_F,\omega_B$, then the term of stochastic noise satisfies:
\begin{align}
    4\sigma^2\dfrac{(2-p)m-(1-p)m^2}{1-(1-p)(1-m)^2}\leq4\sigma^2\dfrac{(2-\sqrt{m})m-(1-\sqrt{m})m^2}{1-(1-\sqrt{m})(1-m)^2}\leq4\sigma^2\cdot 2\sqrt{m}.
\end{align}

Let
\begin{equation*}
\begin{aligned}
    m \sim&\left(\dfrac{1}{(1-\omega_B)^{E-1}(1-\omega_F)^{E-1}} +\dfrac{\sigma^{\frac{4}{3}}T^{\frac{2}{3}}}{(1-\omega_B)^{\frac{4(E-1)}{3}}(1-\omega_F)^{\frac{4(E-1)}{3}}}\right)^{-1},\ m\leq 1,\\
    \gamma\sim&\left(\dfrac{1}{(1-\omega_B)^{E-1}(1-\omega_F)^{E-1}} +\dfrac{\sigma^{\frac{4}{3}} T^{\frac{2}{3}}}{(1-\omega_B)^{\frac{4(E-1)}{3}}(1-\omega_F)^{\frac{4(E-1)}{3}}}\right)^{-1},\ \gamma\leq 1.
\end{aligned}
\end{equation*}
At this time, if $\gamma/m$ is sufficiently small, we can obtain that \[96L_{\nabla \ell}^2\dfrac{1}{m^2}+C_{w,e}-\dfrac{1}{2\gamma^2}\leq0.\] 

Then with Assumption \ref{assumption:sample}, we have:
{\small
\begin{equation*}
\begin{aligned}
    &\dfrac{1}{T}\sum_{t=1}^T\mathbb{E}\left[\left\Vert\nabla l(W^{(t)})\right\Vert^2\right]\\
    \leq&\dfrac{2}{\gamma T}\mathbb{E}\left[\ell(\w^{(1)})-\inf_{\w}l(\w)\right]+C_{x}\sqrt{m}\varphi^2+8\sigma^2\sqrt{m}+\dfrac{3}{mT}\sum_{e=1}^E\mathbb{E}\left[\left\Vert\tilde{u}^{(1)}_e-\nabla_e\ell(\w^{(1)})\right\Vert^2\right]\\
    &+\sum_{e=1}^{E-1}\dfrac{C_{y,\theta,e}}{T}\mathbb{E}\left[\left\Vert\tilde{\by}^{(1)}_e-\by^{(1)}_e\right\Vert^2\right]+\sum_{e=1}^{E-1}\dfrac{C_{v,e}}{T}\mathbb{E}\left[\left\Vert\tilde{ v}_{e}^{(1)}- \hat{v}_{e}^{(1)}\right\Vert^2\right]+\sum_{e=1}^{E-1}\dfrac{C_{v,\chi,e}}{T}\mathbb{E}\left[\left\Vert\tilde{ v}_{e}^{(1)}-\bv_{e}^{(1)}\right\Vert^2\right]\\
    &+\dfrac{1}{T}\mathbb{E}\left[\left\Vert\nabla \ell(\w^{(1)})\right\Vert^2\right]\\
    \lesssim&\dfrac{\sigma^{\frac{4}{3}}}{T^{\frac{1}{3}}(1-\omega_B)^{\frac{4(E-1)}{3}}(1-\omega_F)^{\frac{4(E-1)}{3}}}\\
    &+\dfrac{1}{T}\left(\dfrac{1}{(1-\omega_B)^{E-1}(1-\omega_F)^{E-1}}+\dfrac{\omega_F^2+\omega_B}{(1-\omega_B)^{2(E-1)}(1-\omega_F)^2}+\dfrac{1}{(1-\omega_F)^{2(E-2)-1}}\right).
\end{aligned}
\end{equation*}}
\end{proof}
\end{theorem}

\begin{remark}
    It can be observed that the conclusion in Lemma \ref{thm:convergence_clapping} can be directly obtained from \eqref{descent lemma_more}.
\end{remark}

\subsection{A simplified case: convergence rate of the two-worker case}
\label{app:two-worker-warmup}
The proof above covers the general \(E\)-worker case. For intuition, we also include the two-worker case, where there is only one forward communication and one backward communication. This simplified proof makes the role of error feedback and lazy sampling easier to see.

\begin{lemma}[Convergence rate of {\oursfu} in the case $E=2$]
\label{lemma: FU_e_2}
Suppose Assumption \ref{assumption:smoothness}, \ref{assumption:unbiased}, and \ref{assumption:compressor} hold. Then for {\oursfu}, there exist $\gamma,m>0$ such that: 
\begin{equation}
\label{descent lemma_e3}
\begin{aligned}
    \dfrac{1}{T}\sum_{t=1}^T\mathbb{E}\left[\left\Vert\nabla \ell(\w^{(t)})\right\Vert^2\right]\lesssim\dfrac{\sigma}{\sqrt{T}}+\dfrac{1}{T(1-\omega_B)(1-\omega_F)}.
\end{aligned}
\end{equation}
\begin{proof}

Substituting $E=2$ into \eqref{u-nablal}, we have:
\begin{equation}
\label{u-nablal_2_1}
    \begin{aligned}
        &\sum_{t=1}^{T+1}\sum_{e=1}^2\mathbb{E}\left[\left\Vert\tilde{u}^{(t)}_e-\nabla_e\ell(\w^{(t)})\right\Vert^2\right]\\
        \leq&32L_{\nabla \ell}^2\left(\dfrac{p+m}{m^2(1-(1-p)(1-\frac{m}{2}))}+\dfrac{1}{m^2}\right)\sum_{t=1}^{T}\sum_{e=1}^2\mathbb{E}\left[\left\Vert w^{(t+1)}_e-w^{(t)}_e\right\Vert^2\right]\\
        &+8(L_{\nabla a}^\circ)^2\sum_{t=2}^{T+1}\mathbb{E}\left[\left\Vert\tilde{ v}_1^{(t)}- \hat{v}_1^{(t)}\right\Vert^2\right]+8(L_{\nabla a}')^2\sum_{t=2}^{T+1}\mathbb{E}\left[\left\Vert\tilde{\by}_1^{(t)}-\hat{y}_1^{(t)}\right\Vert^2\right]\\
        &+4T\sigma^2\dfrac{(2-p)m-(1-p)m^2}{1-(1-p)(1-m)^2}+\dfrac{3}{m}\sum_{e=1}^2\mathbb{E}\left[\left\Vert\tilde{u}^{(1)}_e-\nabla_e\ell(\w^{(1)})\right\Vert^2\right].
    \end{aligned}
\end{equation}

For the first $T-1$ iterations, suppose $1=Q_1<Q_2<\cdots<Q_{r_0}\leq T$ are the all moments at which the sample is randomly obtained from $\mathcal{D}$. We also denote $Q_{r_0+1}=T+1$.

For any $1\leq r\leq r_0$, we note that $\tilde{\by}^{(Q_{r})}_{1}=\by^{(Q_{r})}_1=\hat{\by}^{(Q_{r})}_1$ as we do not compress the activation and gradients in \oursfulast. Thus, substituting $E=2$ into \eqref{y_tilde-y}, we can get:
\begin{equation}
\label{y_tilde-y2_1}
\begin{aligned}
    \sum_{t=Q_r+1}^{Q_{r+1}}\mathbb{E}\left[\left\Vert\tilde{\by}_1^{(t)}-\hat{y}_1^{(t)}\right\Vert^2\right]&=\sum_{t=Q_r+1}^{Q_{r+1}-1}\mathbb{E}\left[\left\Vert\tilde{\by}_1^{(t)}-\hat{y}_1^{(t)}\right\Vert^2\right] 
    \leq C_y^2\dfrac{\omega_F^2}{(1-\omega_F)^2}\sum_{t=Q_r}^{Q_{r+1}-2}\mathbb{E}\left[\left\Vert\by^{(t+1)}_1-\by^{(t)}_1\right\Vert^2\right].
\end{aligned}
\end{equation}

Taking summation over $r=1,2,\cdots,r_0$, it holds that:
\begin{equation}
\label{y_tilde-y2}
\begin{aligned}
    \sum_{t=2}^{T+1}\mathbb{E}\left[\left\Vert\tilde{\by}_1^{(t)}-\hat{y}_1^{(t)}\right\Vert^2\right]&=\sum_{r=1}^{r_0}\sum_{t=Q_r+1}^{Q_{r+1}}\mathbb{E}\left[\left\Vert\tilde{\by}_1^{(t)}-\hat{y}_1^{(t)}\right\Vert^2\right]\\
    &\leq C_y^2\dfrac{\omega_F^2}{(1-\omega_F)^2}\sum_{r=1}^{r_0}\sum_{t=Q_r}^{Q_{r+1}-2}\mathbb{E}\left[\left\Vert\by^{(t+1)}_1-\by^{(t)}_1\right\Vert^2\right].
\end{aligned}
\end{equation}



Then, use the fact that $\nabla a_2$ is Lipschitz continous, we have:
    \begin{align}
    \label{lipschitz_a_2}
    \left\Vert\bv^{(t)}_{1}- \hat{v}^{(t)}_{1}\right\Vert^2=\left\Vert\nabla a_2(\tilde{\by}^{(t)}_{1},w_2^{(t)})-\nabla a_2(\hat{y}^{(t)}_{1},w_2^{(t)})\right\Vert^2\leq L_{\nabla a}^2\left\Vert\tilde{\by}^{(t)}_{1}-\hat{y}^{(t)}_{1}\right\Vert^2.
    \end{align}

From $\tilde{ v}^{(Q_r)}_{1}=\bv^{(Q_r)}_1=\hat{ v}^{(Q_r)}_{1}$, we can obtain that:
\begin{equation}
\label{v_tilde-v_21}
\begin{aligned}
    &\sum_{t=Q_r+1}^{Q_{r+1}}\mathbb{E}\left[\left\Vert \tilde{ v}_1^{(t)}-\hat{v}_1^{(t)}\right\Vert^2\right]\leq2\sum_{t=Q_r+1}^{Q_{r+1}}\mathbb{E}\left[\left\Vert \tilde{ v}_1^{(t)}-{\bv}_1^{(t)}\right\Vert^2\right]+2\sum_{t=Q_r+1}^{Q_{r+1}}\mathbb{E}\left[\left\Vert \bv_1^{(t)}-\hat{v}_1^{(t)}\right\Vert^2\right]\\
    =&2\sum_{t=Q_r+1}^{Q_{r+1}-1}\mathbb{E}\left[\left\Vert \tilde{ v}_1^{(t)}-{\bv}_1^{(t)}\right\Vert^2\right]+2\sum_{t=Q_r+1}^{Q_{r+1}-1}\mathbb{E}\left[\left\Vert \bv_1^{(t)}-\hat{v}_1^{(t)}\right\Vert^2\right]\\
    \leq&2\dfrac{\omega_B^2}{(1-\omega_B)^2}\sum_{t=Q_r}^{Q_{r+1}-2}\mathbb{E}\left[\left\Vert\bv^{(t+1)}_1-\bv^{(t)}_1\right\Vert^2\right]+2L_{\nabla a}^2\sum_{t=Q_r+1}^{Q_{r+1}-1}\mathbb{E}\left[\left\Vert\tilde{\by}^{(t)}_{1}-\hat{y}^{(t)}_{1}\right\Vert^2\right],
\end{aligned}
\end{equation}
where the second inequality is due to \eqref{expection_v_chi_2} and \eqref{lipschitz_a_2}.
Taking summation over $r$, it holds that:
\begin{equation}
\label{v_tilde-v_2}
\begin{aligned}
    &\sum_{t=2}^{T+1}\mathbb{E}\left[\left\Vert \tilde{ v}_1^{(t)}-\hat{v}_1^{(t)}\right\Vert^2\right]\\
    \leq&2\dfrac{\omega_B^2}{(1-\omega_B)^2}\sum_{r=1}^{r_0}\sum_{t=Q_r}^{Q_{r+1}-2}\mathbb{E}\left[\left\Vert\bv^{(t+1)}_1-\bv^{(t)}_1\right\Vert^2\right]+2L_{\nabla a}^2\sum_{t=2}^{T+1}\mathbb{E}\left[\left\Vert\tilde{\by}^{(t)}_{1}-\hat{y}^{(t)}_{1}\right\Vert^2\right],
\end{aligned}
\end{equation}

Plugging Eq. \eqref{y_tilde-y2} and Eq. \eqref{v_tilde-v_2} into Eq. \eqref{u-nablal_2_1}, it holds that:
\begin{equation}
\label{u-nablal_2_2}
    \begin{aligned}
        &\sum_{t=1}^{T+1}\sum_{e=1}^2\mathbb{E}\left[\left\Vert\tilde{u}^{(t)}_e-\nabla_e\ell(\w^{(t)})\right\Vert^2\right]\\
        \leq&32L_{\nabla \ell}^2\left(\dfrac{p+m}{m^2(1-(1-p)(1-\frac{m}{2}))}+\dfrac{1}{m^2}\right)\sum_{t=1}^{T}\sum_{e=1}^2\mathbb{E}\left[\left\Vert w^{(t+1)}_e-w^{(t)}_e\right\Vert^2\right]\\
        &+16(L_{\nabla a}^\circ)^2\dfrac{\omega_B^2}{(1-\omega_B)^2}\sum_{r=1}^{r_0}\sum_{t=Q_r}^{Q_{r+1}-2}\mathbb{E}\left[\left\Vert\bv^{(t+1)}_1-\bv^{(t)}_1\right\Vert^2\right]\\
        &+\left(8(L_{\nabla a}')^2+16(L_{\nabla a}^\circ)^2L_{\nabla a}^2\right)C_y^2\dfrac{\omega_F^2}{(1-\omega_F)^2}\sum_{r=1}^{r_0}\sum_{t=Q_r}^{Q_{r+1}-2}\mathbb{E}\left[\left\Vert\by^{(t+1)}_1-\by^{(t)}_1\right\Vert^2\right]\\
        &+4T\sigma^2\dfrac{(2-p)m-(1-p)m^2}{1-(1-p)(1-m)^2}+\dfrac{3}{m}\sum_{e=1}^2\mathbb{E}\left[\left\Vert\tilde{u}^{(1)}_e-\nabla_e\ell(\w^{(1)})\right\Vert^2\right].
    \end{aligned}
\end{equation}

Note that $E=2$, thus for any $1\leq r\leq r_0$ we have:
\begin{equation}
\label{e_2_chi1}
\begin{aligned}
    &\sum_{t=Q_r}^{Q_{r+1}-2}\mathbb{E}\left[\left\Vert\bv^{(t+1)}_1-\bv^{(t)}_1\right\Vert^2\right]=\sum_{t=Q_r}^{Q_{r+1}-2}\mathbb{E}\left[\left\Vert\nabla_1a_2(\tilde{\by}^{(t+1)}_1,w_2^{(t+1)})-\nabla_1a_2(\tilde{\by}^{(t)}_1,w_2^{(t)})\right\Vert^2\right]\\
    \leq&L_{\nabla a}^2\sum_{t=Q_r}^{Q_{r+1}-2}\mathbb{E}\left[\left\Vert\tilde{\by}^{(t+1)}_1-\tilde{\by}^{(t)}_1\right\Vert^2+\left\Vert w_2^{(t+1)}-w_2^{(t)}\right\Vert^2\right]\\
    \leq&L_{\nabla a}^2\dfrac{8}{(1-\omega_F)^2}\sum_{t=Q_r}^{Q_{r+1}-2}\mathbb{E}\left[\left\Vert\by_1^{(t+1)}-\by_1^{(t)}\right\Vert^2\right]+L_{\nabla a}^2\sum_{t=Q_r}^{Q_{r+1}-2}\mathbb{E}\left[\left\Vert w_2^{(t+1)}-w_2^{(t)}\right\Vert^2\right],
\end{aligned}
\end{equation}
where the first inequality is due to Assumption \ref{assumption:smoothness} and the second inequality is due to \eqref{yt+1-yt_main}. 

Taking summation over $r$, it holds that:
\begin{equation}
\label{e_2_chi}
\begin{aligned}
    &\sum_{r=1}^{r_0}\sum_{t=Q_r}^{Q_{r+1}-2}\mathbb{E}\left[\left\Vert\bv^{(t+1)}_1-\bv^{(t)}_1\right\Vert^2\right]\\
    \leq&L_{\nabla a}^2\dfrac{8}{(1-\omega_F)^2}\sum_{r=1}^{r_0}\sum_{t=Q_r}^{Q_{r+1}-2}\mathbb{E}\left[\left\Vert\by_1^{(t+1)}-\by_1^{(t)}\right\Vert^2\right]+L_{\nabla a}^2\sum_{t=1}^{T}\mathbb{E}\left[\left\Vert w_2^{(t+1)}-w_2^{(t)}\right\Vert^2\right],
\end{aligned}
\end{equation}

Then as $x^{(t+1)}= x^{(t)}$ for any $t=Q_r,\cdots,Q_{r+1}-2$, we also have
\begin{equation}
\label{e_2_theta1}
\begin{aligned}
    &\sum_{t=Q_r}^{Q_{r+1}-2}\mathbb{E}\left[\left\Vert\by^{(t+1)}_1-\by^{(t)}_1\right\Vert^2\right]=\sum_{t=Q_r}^{Q_{r+1}-2}\mathbb{E}\left[\left\Vert\hat{y}^{(t+1)}_1-\hat{y}^{(t)}_1\right\Vert^2\right]\\
    =&\sum_{t=Q_r}^{Q_{r+1}-2}\mathbb{E}\left[\left\Vert a_1( x^{(t+1)},w_1^{(t+1)})-a_1( x^{(t)},w_1^{(t)})\right\Vert^2\right]\leq L_a^2\sum_{t=Q_r}^{Q_{r+1}-2}\mathbb{E}\left[\left\Vert w_1^{(t+1)}- w_1^{(t)}\right\Vert^2\right],
\end{aligned}
\end{equation}
where the inequality is due to Assumption \ref{assumption:smoothness}.

Taking summation over $r$, it holds that:
\begin{equation}
\label{e_2_theta}
\begin{aligned}
    \sum_{r=1}^{r_0}\sum_{t=Q_r}^{Q_{r+1}-2}\mathbb{E}\left[\left\Vert\by^{(t+1)}_1-\by^{(t)}_1\right\Vert^2\right]\leq L_a^2\sum_{t=1}^{T}\mathbb{E}\left[\left\Vert w^{(t+1)}- w^{(t)}\right\Vert^2\right],
\end{aligned}
\end{equation}

Plugging \eqref{e_2_chi} and \eqref{e_2_theta} into \eqref{u-nablal_2_2}, we can get:
\begin{equation}
\label{u-nablal_3}
    \begin{aligned}
        &\sum_{t=1}^{T+1}\sum_{e=1}^2\mathbb{E}\left[\left\Vert\tilde{u}^{(t)}_e-\nabla_e\ell(\w^{(t)})\right\Vert^2\right]\\
        \leq&32L_{\nabla \ell}^2\left(\dfrac{p+m}{m^2(1-(1-p)(1-\frac{m}{2}))}+\dfrac{1}{m^2}\right)\sum_{t=1}^{T}\sum_{e=1}^2\mathbb{E}\left[\left\Vert w^{(t+1)}_e-w^{(t)}_e\right\Vert^2\right]\\
        &+16L_{\nabla_a}^2(L_{\nabla a}^\circ)^2\dfrac{\omega_B^2}{(1-\omega_B)^2}\sum_{t=1}^{T}\mathbb{E}\left[\left\Vert w^{(t+1)}_2-w^{(t)}_2\right\Vert^2\right]\\
        &+\left(\left(8(L_{\nabla a}')^2+16(L_{\nabla a}^\circ)^2L_{\nabla a}^2\right)C_y^2\dfrac{\omega_F^2}{(1-\omega_F)^2}+128L_{\nabla_a}^2(L_{\nabla a}^\circ)^2\dfrac{\omega_B^2}{(1-\omega_B)^2(1-\omega_F)^2}\right)\\
        &\hspace{8cm}\cdot L_a^2 \sum_{t=1}^{T}\mathbb{E}\left[\left\Vert w^{(t+1)}-w^{(t)}\right\Vert^2\right]\\
        &+4T\sigma^2\dfrac{(2-p)m-(1-p)m^2}{1-(1-p)(1-m)^2}+\dfrac{3}{m}\sum_{e=1}^2\mathbb{E}\left[\left\Vert\tilde{u}^{(1)}_e-\nabla_e\ell(\w^{(1)})\right\Vert^2\right]\\
        \leq&C_w\sum_{t=1}^{T}\mathbb{E}\left[\left\Vert \w^{(t+1)}-\w^{(t)}\right\Vert^2\right]+4T\sigma^2\dfrac{(2-p)m-(1-p)m^2}{1-(1-p)(1-m)^2}+\dfrac{3}{m}\sum_{e=1}^2\mathbb{E}\left[\left\Vert\tilde{u}^{(1)}_e-\nabla_e\ell(\w^{(1)})\right\Vert^2\right].
    \end{aligned}
\end{equation}
where:
\begin{align*}
C_w=&32L_{\nabla \ell}^2\left(\dfrac{p+m}{m^2(1-(1-p)(1-\frac{m}{2}))}+\dfrac{1}{m^2}\right)+16L_{\nabla_a}^2(L_{\nabla a}^\circ)^2\dfrac{\omega_B^2}{(1-\omega_B)^2}\\
&+\left(8(L_{\nabla a}')^2+16(L_{\nabla a}^\circ)^2L_{\nabla a}^2\right)L_a^2 C_y^2\dfrac{\omega_F^2}{(1-\omega_F)^2}+128L_{\nabla_a}^2(L_{\nabla a}^\circ)^2L_a^2 \dfrac{\omega_B^2}{(1-\omega_B)^2(1-\omega_F)^2}.  
\end{align*}

Plugging $E=2$ into \eqref{descent lemma}, combining it with \eqref{u-nablal_3}, and then taken $p=m$, we can get
\begin{equation}
\label{descent lemma_e2}
\begin{aligned}
    &\dfrac{1}{T}\sum_{t=1}^T\mathbb{E}\left[\left\Vert\nabla \ell(\w^{(t)})\right\Vert^2\right]\\
    \leq&\dfrac{2}{\gamma T}\mathbb{E}\left[\ell(\w^{(1)})-\inf_{\w}\ell(\w)\right]+\dfrac{1}{T}\left(C_w-\dfrac{1}{2\gamma^2}\right)\sum_{t=1}^{T}\mathbb{E}\left[\left\Vert \w^{(t+1)}-\w^{(t)}\right\Vert^2\right]\\
    &+4\sigma^2\dfrac{(2-p)m-(1-p)m^2}{1-(1-p)(1-m)^2}+\dfrac{3}{mT}\sum_{e=1}^2\mathbb{E}\left[\left\Vert\tilde{u}^{(1)}_e-\nabla_e\ell(\w^{(1)})\right\Vert^2\right].
\end{aligned}
\end{equation}

Let $p=p_0$ as a constant with the order of $\mathcal{O}(1)$, and let:
\begin{align*}
    m&\sim\left(\dfrac{1}{(1-\omega_B)(1-\omega_F)}+\sigma\sqrt{T}\right)^{-1},\quad
    \gamma\sim\left(\dfrac{1}{(1-\omega_B)(1-\omega_F)}+\sigma\sqrt{T}\right)^{-1}\ \text{and}\ m,\gamma\leq 1.
\end{align*}
Then, as $\gamma$ has the same order as $m$ with respect to $\omega_B,\omega_F,\sigma,T$, it holds that $C_w-\dfrac{1}{2\gamma^2}\leq0$ if $\gamma/m$ is sufficiently small.

Thus, we have:
\begin{equation*}
\begin{aligned}
    \dfrac{1}{T}\sum_{t=1}^T\mathbb{E}\left[\left\Vert\nabla  \ell(\w^{(t)})\right\Vert^2\right]\lesssim&\dfrac{\sigma}{\sqrt{T}}+\dfrac{1}{T(1-\omega_B)(1-\omega_F)}.
\end{aligned}
\end{equation*}
\end{proof}
\end{lemma}

\section{Error propagation analysis for a single forward-backward step}
\label{appendix: Error propagation analysis}
The following lemma states how the compressed error propagates for one forward and backward step.


\begin{lemma}
\label{lemma: error accumulation and propagation_proof}
    Suppose Assumptions \ref{assumption:smoothness} and \ref{assumption:compressor} hold. Then, for \( e = 1, \dots, E-1 \), the error of the forward activation in Algorithm \ref{alg:clapping} can be bounded above as follows:
    \begin{equation}
    \label{forward_compression_error_proof}
    \begin{aligned}
        \left\Vert\tilde{\by}_e^{(t)}-\hat{y}_e^{(t)}\right\Vert^2\leq\sum_{\iota=1}^{e}2(2L_a^2)^{e-\iota}\left\Vert\tilde{\by}_\iota^{(t)}-\by_\iota^{(t)}\right\Vert^2.
    \end{aligned}
    \end{equation}
    For the error of backward gradients in Algorithm \ref{alg:clapping}, there exist constants $L_{\nabla a}^{\circ},L_{\nabla a}'>0$ such that:
    \begin{equation}
    \label{backward_compression_error_proof}
    \begin{aligned}
        \left\Vert\tilde{ v}^{(t)}_e-\hat{v}^{(t)}_e\right\Vert^2\leq&2\sum_{\iota=e}^{E-1}(2(L_{\nabla a}^\circ)^2)^{\iota-e}\left\Vert\tilde{ v}^{(t)}_\iota-\bv^{(t)}_\iota\right\Vert^2\\
        &+4(L_{\nabla a}')^2\sum_{\iota=1}^{E-1}\sum_{s=\max\{e,\iota\}}^{E-1}(2(L_{\nabla a}^\circ)^2)^{s-e}(2L_{a}^2)^{s-\iota}\left\Vert \tilde{\by}_{\iota}^{(t)}-y_{\iota}^{(t)}\right\Vert^2.
    \end{aligned}
    \end{equation}
\begin{proof}
    Eq. \eqref{forward_compression_error_proof} is actually the same as Eq. \eqref{zhankai_y}. Then from \eqref{zhankai_v} we can get:
        \begin{equation}
        \label{zhankai_v_proof}
            \begin{aligned}
                &\left\Vert\tilde{ v}^{(t)}_e- \hat{v}^{(t)}_e\right\Vert^2\\
                \leq&2\sum_{\iota=e}^{E-1}(2(L_{\nabla a}^\circ)^2)^{\iota-e}\left\Vert\tilde{ v}^{(t)}_\iota-\bv^{(t)}_\iota\right\Vert^2+(2(L_{\nabla a}')^2)\sum_{\iota=e}^{E-1}(2(L_{\nabla a}^\circ)^2)^{\iota-e}\left\Vert \tilde{\by}_{\iota}^{(t)}-\hat{y}_{\iota}^{(t)}\right\Vert^2,
            \end{aligned}
        \end{equation}
where the last inequality uses the fact that $L_{\nabla a}'\geq L_{\nabla a}$ in Remark \ref{remark:a_e^0}.

Then, plugging \eqref{forward_compression_error_proof} into \eqref{zhankai_v_proof}, we can obtain:
\begin{equation}
    \begin{aligned}
        &\left\Vert\tilde{ v}^{(t)}_e-\hat{v}^{(t)}_e\right\Vert^2\\
        \leq&2\sum_{\iota=e}^{E-1}(2(L_{\nabla a}^\circ)^2)^{\iota-e}\left\Vert\tilde{ v}^{(t)}_\iota-\bv^{(t)}_\iota\right\Vert^2+(2(L_{\nabla a}')^2)\sum_{\iota=e}^{E-1}(2(L_{\nabla a}^\circ)^2)^{\iota-e}\sum_{s=1}^{\iota}2(2L_a^2)^{\iota-s}\left\Vert\tilde{\by}_s^{(t)}-\by_s^{(t)}\right\Vert^2\\
        \leq&2\sum_{\iota=e}^{E-1}(2(L_{\nabla a}^\circ)^2)^{\iota-e}\left\Vert\tilde{ v}^{(t)}_\iota-\bv^{(t)}_\iota\right\Vert^2+4(L_{\nabla a}')^2\sum_{\iota=1}^{E-1}\sum_{s=\max\{e,\iota\}}^{E-1}(2(L_{\nabla a}^\circ)^2)^{s-e}(2L_{a}^2)^{s-\iota}\left\Vert \tilde{\by}_{\iota}^{(t)}-y_{\iota}^{(t)}\right\Vert^2.
    \end{aligned}
\end{equation}
Then we know that Eq. \eqref{backward_compression_error_proof} holds for $e=1,\cdots,E-1$.
\end{proof}
\end{lemma}
As indicated by Eq. \eqref{forward_compression_error_proof}, the error in forward activation is directly accumulated due to the compression operations of the preceding machines. Meanwhile, according to Eq. \eqref{backward_compression_error_proof}, the error in backward gradient evaluation, denoted as $\tilde{v}_e^{(t)}=\tilde{v}_{e+1}^{(t)}\nabla_1a_{e+1}(\tilde{y}_e^{(t)},W_e^{(t)})$, arises from two aspects. One aspect pertains to the accumulated errors during the backward compression in preceding machines, while the other relates to the error of the forward activation $\tilde{y}_e^{(t)}$. Moreover, the mathematical formulations presented in Eq. \eqref{forward_compression_error_proof} and Eq. \eqref{backward_compression_error_proof} jointly demonstrate that compression errors exhibit exponential amplification across distributed computing nodes. Furthermore, this analysis reveals a positive correlation between error magnitude and system complexity: as model dimensionality and parallelism scale increase, communication-induced compression errors emerge as a critical bottleneck in distributed training architectures.

\section{\ours with Adam optimizer}
\label{section: clapping_adam}
In this section, we present \ours equipped with Adam \cite{kingma2014adam} optimizer, which can especially fit the need of pre-training and fine-tuning tasks for LLMs.

\subsection{Algorithm design}
Similar to the standard back-propagation algorithms with Adam optimizer, we introuduce two coefficients $\beta_1,\beta_2\in(0,1)$ and use $\nu_e^{(t)}$ and $\upsilon_e^{(t)}$ to record the first- and second-order optimizer states for the parameter $w_e^{(t)}$ for $e=1,2,\cdots,E$, respectively. Then we use the following update rules to optimize $w_e^{(t)}$:
\begin{subequations}
\label{eq: adam_update}
\begin{align}
\tilde{u}_e^{(t)}& = (1-\beta_1)\tilde{u}_e^{(t-1)}+ \beta_1\nabla_2 a_e(\tilde{y}^{(t)}_{e-1},w^{(t)}_e)\tran \tilde{v}^{(t)}_{e},\\
\upsilon_e^{(t)}& = (1-\beta_2)\upsilon_e^{(t-1)}+ \beta_2\left(\nabla_2 a_e(\tilde{y}^{(t)}_{e-1},w^{(t)}_e)\tran \tilde{v}^{(t)}_{e}\right)^{\odot 2},\\
\nu^{(t)}_{e}& = \dfrac{\tilde{u}_e^{(t)}}{\sqrt{\upsilon_e^{(t)}+\varepsilon}},\label{eq: adam_update_chu}\\
w_e^{(t+1)} &= w_e^{(t)} - \gamma  \nu^{(t)}_{e},
\end{align}
\end{subequations}
where $\odot 2$ denotes the second moment and $\varepsilon>0$ is a fixed constant. The division and addition operate in \eqref{eq: adam_update_chu} are all sample-wised. Then we can present \ours with Adam optimizer as Algorithm \ref{alg:clapping_adam}.

\begin{algorithm}[t]
  \caption{\ours with Adam optimizer}
  \label{alg:clapping_adam}
  \begin{algorithmic}[1]
  \REQUIRE Initialize $\tilde{y}^{(0)}_e = 0, \tilde{v}^{(0)}_e = 0$ for $e=1,\cdots, E-1$,
  and $\upsilon^{(0)}_e = 0, \tilde{u}^{(0)}_e = 0$ for $e=1,\cdots, E$.
  Initialize dataset $\mathcal{D}$, learning rate $\gamma_t$, compressor $\mathcal{C}$, and lazy sampling rate $\{p_t\}_{t=1}^T$.
  \FOR{$t=1,\cdots,T$}        
    \STATE ${x}^{(t)},f^{(t)}_{\text{FU}} \gets {\mbox{\sffamily LazySampling}}(\mathcal{D},t,p_t)$;
    \STATE $\tilde{y}^{(t)}_0 \gets {x}^{(t)}$, \ \ $\tilde{v}_{E}^{(t)} \gets 1$.
    \FOR{$e=1,2,\cdots,E-1$}
      \STATE $\tilde{y}^{(t)}_e \gets {\mbox{\sffamily Forward}}_e(\tilde{{y}}_e^{(t-1)},\tilde{y}_{e-1}^{(t)},w_e^{(t)},f^{(t)}_{\text{FU}})$.
    \ENDFOR
    \FOR{$e=E,E-1,\cdots,1$}
      \IF{$e\neq 1$}
        \STATE $\tilde{v}^{(t)}_{e-1} \gets {\mbox{\sffamily Backward}}_{e}(\tilde{{v}}_{e-1}^{(t-1)},\tilde{{v}}_{e}^{(t)},\tilde{y}_{e-1}^{(t)},w_{e}^{(t)},f^{(t)}_{\text{FU}})$.
      \ENDIF
      \STATE Update $\upsilon^{(t)}_e$, $\tilde{u}^{(t)}_e$ and $w_e^{(t+1)}$ by \eqref{eq: adam_update}.
    \ENDFOR
  \ENDFOR
  \end{algorithmic}
\end{algorithm}

\subsection{Convergence of \ours with Adam optimizer}
\label{appendix: convergence_clapping_adam}
In this subsection we present the convergence analysis of \ours with Adam optimizer. Firstly, we need to present an additional assumption for the bounded gradient estimation:
\begin{assumption}
\label{assumption: bounded gradient estimation}
There exist $M_1\geq0$ such that for each $e=1,2,\cdots,E$, the gradient estimation can be bonuded as:
\begin{align*}
\left\Vert\nabla_2 a_e(\tilde{y}^{(t)}_{e-1},w^{(t)}_e)\tran \tilde{v}^{(t)}_{e}\right\Vert\leq M_1.
\end{align*}
\end{assumption}
\begin{remark}
We assume the gradient of $a_e(y,\w)$ is bounded according to Assumption \ref{assumption:smoothness}. Thus, if there exists a constant \(\omega_0\) such that the compressor \(\mathcal{C}\) satisfies \(\left\Vert x-\mathcal{C}(x)\right\Vert^2\leq\omega_0\left\Vert x\right\Vert^2\) for any input vector $x$, then Assumption \ref{assumption: bounded gradient estimation} can be satisfied. Such a bounded property is common for traditional compressors like TopK and quantization.
\end{remark}

With Assumption \ref{assumption: bounded gradient estimation}, one can obtain that the second-order optimizer state $\upsilon_e^{(t)}$ satisfy that $\left\Vert\upsilon_e^{(t)}\right\Vert\leq M_1^2$. Then we can directly obtain the following lemma
\begin{lemma}
    There exists $M>0$ such that for $t=1,2,\cdots,T$ and $e=1,2,\cdots,E$ it holds that: $$\varepsilon\leq\left\Vert\sqrt{\upsilon_e^{(t)}+\varepsilon}\right\Vert_1\leq M,$$
    where $\Vert\cdot\Vert_1$ denote the $\ell_1$-norm.
\end{lemma}

\begin{lemma}
Suppose Assumption \ref{assumption:unbiased} and Assumption \ref{assumption: bounded gradient estimation} hold. And the step-size $\gamma$ satisfies that $\gamma\leq\min\left\{\dfrac{1}{4L_{\nabla \ell}},\dfrac{1}{4C_{\nu}}\right\}$, where $C_{\nu}^2\leq\max\left\{(M-1)^2,(\varepsilon-1)^2\right\}$ is a constant. Then in Algorithm \ref{alg:clapping_adam} we have:
    \begin{equation}
    \begin{aligned}
    \label{descent lemma_adam}
        &\dfrac{1}{T}\sum_{t=1}^T\mathbb{E}\left[\left\Vert\nabla \ell(\w^{(t)})\right\Vert^2\right]\\
        \leq&\dfrac{2}{\gamma T}\mathbb{E}\left[\ell(\w^{(1)})\hspace{-0.3mm}-\hspace{-0.3mm}\inf_{\w}\ell(W)\right]\hspace{-0.5mm}-\hspace{-0.5mm}\dfrac{1}{2\gamma^2 T}\sum_{t=1}^T\mathbb{E}\left[\left\Vert \w^{(t+1)}\hspace{-0.5mm}-\hspace{-0.3mm}\w^{(t)}\right\Vert^2\right]\hspace{-0.5mm}+\hspace{-0.5mm}\dfrac{2}{T}\sum_{t=1}^T\sum_{e=1}^E\mathbb{E}\left[\left\Vert\tilde{u}^{(t)}_e\hspace{-0.5mm}-\hspace{-0.4mm}\nabla_e\ell(\w^{(t)})\right\Vert^2\right].
    \end{aligned}
    \end{equation}
\begin{proof}
As $\ell$ is $L_{\nabla\ell}$-smooth, then it holds that:
\begin{equation}
\begin{aligned}
    &\ell(\w^{(t+1)})\\
    \leq&\ell(\w^{(t)})+\langle\nabla\ell(\w^{(t)}),\w^{(t+1)}-\w^{(t)}\rangle+\dfrac{L_{\nabla\ell}}{2}\left\Vert\w^{(t+1)}-\w^{(t)}\right\Vert^2\\
    \leq&\ell(\w^{(t)})-\dfrac{\gamma}{2}\left\Vert\nabla\ell(\w^{(t)})\right\Vert^2-\left(\dfrac{\gamma}{2}-\dfrac{L_{\nabla\ell}}{2}\right)\left\Vert\w^{(t+1)}-\w^{(t)}\right\Vert^2+\dfrac{\gamma}{2}\sum_{e=1}^E\left\Vert\nu_e^{(t)}-\nabla\ell(\w^{(t)})\right\Vert^2\\
    \leq&\ell(\w^{(t)})-\dfrac{\gamma}{2}\left\Vert\nabla\ell(\w^{(t)})\right\Vert^2-\left(\dfrac{\gamma}{2}-\dfrac{L_{\nabla\ell}}{2}\right)\left\Vert\w^{(t+1)}-\w^{(t)}\right\Vert^2+\gamma\sum_{e=1}^E\left\Vert{\tilde{u}_e^{(t)}}-\nabla\ell(\w^{(t)})\right\Vert^2\\
    &+\gamma\sum_{e=1}^E\left\Vert\left(\sqrt{\upsilon_e^{(t)}+\varepsilon}-1\right)\odot\nu_e^{(t)}\right\Vert^2\\
    \leq&\ell(\w^{(t)})-\dfrac{\gamma}{2}\left\Vert\nabla\ell(\w^{(t)})\right\Vert^2-\left(\dfrac{\gamma}{2}-\dfrac{L_{\nabla\ell}}{2}-\dfrac{C_{\nu}^2}{\gamma}\right)\left\Vert\w^{(t+1)}-\w^{(t)}\right\Vert^2+\gamma\sum_{e=1}^E\left\Vert{\tilde{u}_e^{(t)}}-\nabla\ell(\w^{(t)})\right\Vert^2,
\end{aligned}
\end{equation}
where $C_{\nu}^2\leq\max\left\{(M-1)^2,(\varepsilon-1)^2\right\}$.

Then as $\gamma\leq\min\left\{\dfrac{1}{4L_{\nabla \ell}},\dfrac{1}{4C_{\nu}}\right\}$, we obtain that:
\begin{equation}
\begin{aligned}
\label{eq: descent_adam_1}
    &\ell(\w^{(t+1)})\\
    \leq&\ell(\w^{(t)})-\dfrac{\gamma}{2}\left\Vert\nabla\ell(\w^{(t)})\right\Vert^2-\dfrac{\gamma}{4}\left\Vert\w^{(t+1)}-\w^{(t)}\right\Vert^2+\gamma\sum_{e=1}^E\left\Vert{\tilde{u}_e^{(t)}}-\nabla\ell(\w^{(t)})\right\Vert^2.
\end{aligned}
\end{equation}
Similar to the proof of Lemma \ref{lem:descent}, we can obtain from Eq. \eqref{eq: descent_adam_1} that:
    \begin{equation*}
    \begin{aligned}
        &\dfrac{1}{T}\sum_{t=1}^T\mathbb{E}\left[\left\Vert\nabla \ell(\w^{(t)})\right\Vert^2\right]\\
        \leq&\dfrac{2}{\gamma T}\mathbb{E}\left[\ell(\w^{(1)})\hspace{-0.3mm}-\hspace{-0.3mm}\inf_{\w}\ell(W)\right]\hspace{-0.5mm}-\hspace{-0.5mm}\dfrac{1}{2\gamma^2 T}\sum_{t=1}^T\mathbb{E}\left[\left\Vert \w^{(t+1)}\hspace{-0.5mm}-\hspace{-0.3mm}\w^{(t)}\right\Vert^2\right]\hspace{-0.5mm}+\hspace{-0.5mm}\dfrac{2}{T}\sum_{t=1}^T\sum_{e=1}^E\mathbb{E}\left[\left\Vert\tilde{u}^{(t)}_e\hspace{-0.5mm}-\hspace{-0.4mm}\nabla_e\ell(\w^{(t)})\right\Vert^2\right].
    \end{aligned}
    \end{equation*}
Then we finish the proof of Eq. \eqref{descent lemma_adam}.
\end{proof}
\end{lemma}

As the analysis of compress error and error accumulation process are independent with the optimizer, we can directly use the other lemmas in Appendix \ref{appendix:proof}. Thus, \eqref{u-nablal_3_more} and \eqref{u-nablal_3_more1} also hold for \oursfu and \oursfclast, respectively. Finally, we can present the convergence rate of \ours with Adam optimizer as follows:

\begin{lemma}[Convergence rate of \oursfu with Adam optimizer]
Suppose Assumption \ref{assumption:smoothness}-\ref{assumption: bounded gradient estimation} hold. Then for Algorithm \ref{alg:clapping_adam} there exist $\gamma,\beta_1,p>0$ such that: 
\begin{equation}
\label{descent lemma__GENERAL_e3_adam}
\begin{aligned}
    \dfrac{1}{T}\sum_{t=1}^T\mathbb{E}\left[\left\Vert\nabla \ell(\w^{(t)})\right\Vert^2\right]\lesssim&\dfrac{\sigma}{\sqrt{T}}+\dfrac{1}{T(1-\omega_B)^{E-1}(1-\omega_F)^{E-1}}.
\end{aligned}
\end{equation}
\end{lemma}
\begin{lemma}[Convergence rate of \oursfc with Adam optimizer]
Suppose Assumption \ref{assumption:smoothness}-\ref{assumption: bounded gradient estimation} hold. Then for Algorithm \ref{alg:clapping_adam} there exist $\gamma,\beta_1,p>0$ such that: 
\begin{equation}
\label{descent lemma__GENERAL_e3_adam}
\begin{aligned}
    &\dfrac{1}{T}\sum_{t=1}^T\mathbb{E}\left[\left\Vert\nabla \ell(\w^{(t)})\right\Vert^2\right]\lesssim\dfrac{\sigma^{\frac{4}{3}}}{T^{\frac{1}{3}}(1-\omega_B)^{\frac{4(E-1)}{3}}(1-\omega_F)^{\frac{4(E-1)}{3}}}\\
    &\quad\quad+\dfrac{1}{T}\left(\dfrac{1}{(1-\omega_B)^{E-1}(1-\omega_F)^{E-1}}+\dfrac{\omega_F^2+\omega_B}{(1-\omega_B)^{2(E-1)}(1-\omega_F)^2}+\dfrac{1}{(1-\omega_F)^{2(E-2)-1}}\right).
\end{aligned}
\end{equation}
\end{lemma}

\section{\ours with large batch}
\subsection{Algorithm development}
\label{Appendix: Algorithm development}
Here we establish the convergence of \ours with large-batch gradients. Before presenting the theoretical analysis, we firstly present the detailed algorithm of \ours with large-batch gradients.

\textbf{Notations. }Suppose the batch size is $B\geq1$, and we use bold type like $\boldsymbol{y}, \boldsymbol{v}$ to represent a high-dimensional matrix formed by variables computed by the current batch. And we use standard type with subscribe $b=1,2,\cdots,B$ denote the activation and gradient of activation in the corresponding batch. For example, we can denote $\boldsymbol{y}_e^{(t)}:=(y_{e,1}^{(t)},y_{e,2}^{(t)},\cdots,y_{e,B}^{(t)})^{\tran}$ as the activation of the $e$-th layer. 

\textbf{Lazy sampling strategies of \ours with large batch. }The lazy sampling strategies of \ours with large batch are shown as Algorithm \ref{alg:lazy_sampling_large}. The strategies in the large-batch are including {\sffamily{\colorbox{green!30}{Sample-wise rule}}} and {\sffamily{\colorbox{orange!30}{Batch-wise rule}}}. The sample-wise rule means taking lazy sampling process \textbf{sample by sample}. Meanwhile, the batch-wise rule means one can keep the whole batch with probability $1-p$ and use a new batch with probability $p$. We will show that both strategies can achieve the convergence later in Appendix \ref{Appendix: Convergence of ours with large batch}.

\textbf{Algorithm formulation of \ours with large batch. }In large batch scenario, the forward and backward function for large-bach scenario are easy to obtain as they are all sample-wise. The detail is summarized as Algorithm \ref{alg:forwardef_large} and \ref{alg:backwardef_large}. With the lazy sampling strategy as well as the forward/backward process, the \ours algorithm in with large batch can be summarized as Algorithm \ref{alg:clapping_large}.

\paragraph{Practical implementation of sample-wise FU.}
The sample-wise rule in Algorithm~\ref{alg:lazy_sampling_large} does not
require literal per-sample branching inside a communication kernel. In practice,
one maintains a boolean refresh mask \(f_{\mathrm{FU},b}\) over batch slots. Before
each boundary transmission, the batch is stably partitioned into two contiguous
sub-batches: refreshed samples with \(f_{\mathrm{FU},b}=\texttt{True}\), and reused
samples with \(f_{\mathrm{FU},b}=\texttt{False}\). The refreshed group takes the
raw path, while the reused group takes the compressed error-feedback path. The
outputs are then scattered back to the original batch order. Thus the
sample-wise rule can be implemented using standard batched gather/compute/send
and scatter operations. The batch-wise rule is the simplest default, while the
sample-wise rule is the more flexible fine-grained version.

\begin{algorithm}[t!]
  \caption{{\sffamily{LazySampling\_LargeBatch}}($\mathcal{D},t,p$) with 
  {\sffamily{\colorbox{green!30}{Sample-wise rule}}} / {\sffamily{\colorbox{orange!30}{Batch-wise rule}}}}
  \label{alg:lazy_sampling_large}
  \begin{algorithmic}[1]
  \IF{$t=1$}
    \STATE Get \colorbox{green!30}{samples $x^{(1)}_{1},\cdots,x^{(1)}_{B}$} /
           \colorbox{orange!30}{the first batch $\boldsymbol{x}^{(1)}$} from $\mathcal{D}$.
    \STATE Set $f^{(1)}_{\text{FU},b}\gets \texttt{True}$ for all $b=1,\cdots,B$.
  \ELSE
    \STATE \colorbox{green!30}{\textbf{for} $b=1,\cdots,B$ independently \textbf{do}}
    \STATE \colorbox{green!30}{\quad With prob. $1-p$: $x_{b}^{(t)}\gets x_{b}^{(t-1)}$, $f^{(t)}_{\text{FU},b}\gets\texttt{False}$.}
    \STATE \colorbox{green!30}{\quad With prob. $p$: \ \ $x_{b}^{(t)}\sim\mathcal{D}$, $f^{(t)}_{\text{FU},b}\gets\texttt{True}$.}
    \STATE \colorbox{green!30}{\textbf{end for}}
    \STATE \colorbox{orange!30}{With prob. $1-p$: $\boldsymbol{x}^{(t)}\gets \boldsymbol{x}^{(t-1)}$ and set $f^{(t)}_{\text{FU},b}\gets\texttt{False}$ for all $b$.}
    \STATE \colorbox{orange!30}{With prob. $p$: $\boldsymbol{x}^{(t)}\sim\mathcal{D}$ and set $f^{(t)}_{\text{FU},b}\gets\texttt{True}$ for all $b$.}
  \ENDIF
  \STATE \textbf{{Return: }}$\boldsymbol{x}^{(t)}=(x^{(t)}_{1},\cdots,x^{(t)}_{B})^\top$,
  $f^{(t)}_{\text{FU}}=\left(f^{(t)}_{\text{FU},1},\cdots,f^{(t)}_{\text{FU},B}\right)$.
  \end{algorithmic}
\end{algorithm}

\begin{algorithm}[t!]
  \caption{{\sffamily ForwardEF\_LargeBatch}$_e(\tilde{y}_{e,b}^{(t-1)},\tilde{y}_{e-1,b}^{(t)},w_e^{(t)},f^{(t)}_{\mathrm{FU},b})$}
  \label{alg:forwardef_large}
  \begin{algorithmic}
    \STATE In worker $e$: $y_{e,b}^{(t)} = a_e(\tilde{y}_{e-1,b}^{(t)}, w_e^{(t)})$.
    \IF{\oursfu \textbf{and} $f^{(t)}_{\mathrm{FU},b}=\texttt{True}$}
      \STATE Send $y_{e,b}^{(t)}$ from worker $e$ to worker $e+1$.
      \STATE $\tilde{y}_{e,b}^{(t)} \leftarrow y_{e,b}^{(t)}$.
    \ELSE
      \STATE Send $\mathcal{C}\!\left(y_{e,b}^{(t)}-\tilde{y}_{e,b}^{(t-1)}\right)$ from worker $e$ to worker $e+1$.
      \STATE $\tilde{y}_{e,b}^{(t)} \leftarrow \tilde{y}_{e,b}^{(t-1)} + \mathcal{C}\!\left(y_{e,b}^{(t)}-\tilde{y}_{e,b}^{(t-1)}\right)$.
    \ENDIF
    \STATE \textbf{Return:} $\tilde{y}_{e,b}^{(t)}$.
  \end{algorithmic}
\end{algorithm}

\begin{algorithm}[t!]
  \caption{{\sffamily BackwardEF\_LargeBatch}$_e(\tilde{v}_{e-1,b}^{(t-1)},\tilde{v}_{e,b}^{(t)},\tilde{y}_{e-1,b}^{(t)},w_e^{(t)},f^{(t)}_{\mathrm{FU},b})$}
  \label{alg:backwardef_large}
  \begin{algorithmic}
    \STATE In worker $e$: $v_{e-1,b}^{(t)} = \nabla_1 a_e(\tilde{y}_{e-1,b}^{(t)}, w_e^{(t)})^{\top}\tilde{v}_{e,b}^{(t)}$.
    \IF{\oursfu \textbf{and} $f^{(t)}_{\mathrm{FU},b}=\texttt{True}$}
      \STATE Send $v_{e-1,b}^{(t)}$ from worker $e$ to worker $e-1$.
      \STATE $\tilde{v}_{e-1,b}^{(t)} \leftarrow v_{e-1,b}^{(t)}$.
    \ELSE
      \STATE Send $\mathcal{C}\!\left(v_{e-1,b}^{(t)}-\tilde{v}_{e-1,b}^{(t-1)}\right)$ from worker $e$ to worker $e-1$.
      \STATE $\tilde{v}_{e-1,b}^{(t)} \leftarrow \tilde{v}_{e-1,b}^{(t-1)} + \mathcal{C}\!\left(v_{e-1,b}^{(t)}-\tilde{v}_{e-1,b}^{(t-1)}\right)$.
    \ENDIF
    \STATE \textbf{Return:} $\tilde{v}_{e-1,b}^{(t)}$.
  \end{algorithmic}
\end{algorithm}

\begin{algorithm}[t!]
  \caption{\sffamily \ours with large-batch gradients.}
  \label{alg:clapping_large}
  \begin{algorithmic}
    \REQUIRE Initialize $\{\tilde{y}_{e,b}^{(0)}=0\}_{e=1}^{E-1}{}_{b=1}^{B}$,
    $\{\tilde{v}_{e,b}^{(0)}=0\}_{e=1}^{E-1}{}_{b=1}^{B}$,
    and $\{\tilde{u}_{e}^{(0)}=0\}_{e=1}^{E}$.
    Dataset $\mathcal{D}$, stepsize $\{\gamma_t\}$, compressor $\mathcal{C}$, lazy-sampling $\{p_t\}$.
    \FOR{$t=1,\cdots,T$}
      \STATE $\boldsymbol{x}^{(t)}, \boldsymbol{f}^{(t)}_{\mathrm{FU}}
      \leftarrow \text{\sffamily LazySampling\_LargeBatch}(\mathcal{D},t,p_t)$.
      \STATE Set $\tilde{y}_{0,b}^{(t)} \leftarrow x_{b}^{(t)}$ for all $b=1,\cdots,B$, and set $\tilde{v}_{E,b}^{(t)} \leftarrow 1$.
      \FOR{$e=1,2,\cdots,E-1$}
        \FOR{$b=1,2,\cdots,B$}
          \STATE $\tilde{y}_{e,b}^{(t)} \leftarrow \text{\sffamily ForwardEF\_LargeBatch}_e(\tilde{y}_{e,b}^{(t-1)},\tilde{y}_{e-1,b}^{(t)},w_e^{(t)},f^{(t)}_{\mathrm{FU},b})$.
        \ENDFOR
      \ENDFOR
      \FOR{$e=E,E-1,\cdots,1$}
        \STATE $g_e^{(t)} \leftarrow \frac{1}{B}\sum_{b=1}^B \nabla_2 a_e(\tilde{y}_{e-1,b}^{(t)},w_e^{(t)})^\top \tilde{v}_{e,b}^{(t)}$.
        \STATE $\tilde{u}_e^{(t)} \leftarrow (1-m_t)\tilde{u}_e^{(t-1)} + m_t g_e^{(t)}$.
        \STATE $w_e^{(t+1)} \leftarrow w_e^{(t)} - \gamma_t \tilde{u}_e^{(t)}$.
        \IF{$e\neq 1$}
          \FOR{$b=1,2,\cdots,B$}
            \STATE $\tilde{v}_{e-1,b}^{(t)} \leftarrow \text{\sffamily BackwardEF\_LargeBatch}_e(\tilde{v}_{e-1,b}^{(t-1)},\tilde{v}_{e,b}^{(t)},\tilde{y}_{e-1,b}^{(t)},w_e^{(t)},f^{(t)}_{\mathrm{FU},b})$.
          \ENDFOR
        \ENDIF
      \ENDFOR
    \ENDFOR
  \end{algorithmic}
\end{algorithm}

\subsection{Convergence of \ours with large batch}
\label{Appendix: Convergence of ours with large batch}
In this subsection, we present the convergence analysis of \ours with large batch. 
Firstly, it is worth noting that the Descent Lemma (Lemma \ref{lem:descent}) remains applicable. Moreover, since the forward and backward propagations are, in fact, sample-wise, the propagations of activations and the gradients of activations during the forward and backward passes of Algorithm \ref{alg:clapping_large} are identical to those of Algorithm \ref{alg:clapping}. The sole distinction between Algorithm \ref{alg:clapping_large} and Algorithm \ref{alg:clapping} lies in the acquisition of $\tilde{u}$. Consequently, the lemmas regarding error accumulation, specifically Lemmas \ref{lemma:compress error} - \ref{lemma6} and Lemmas \ref{lemma7} - \ref{lemma12}, can be directly employed in the analysis of Algorithm \ref{alg:clapping_large}.

The following lemma is the large-batch version of Lemma \ref{lemma3}:
\begin{lemma}[Large-batch version of Lemma \ref{lemma3}]
\label{lemma3_large}
Suppose Assumption \ref{assumption:smoothness} and \ref{assumption:unbiased} holds, and let $m_1=\cdots=m_T=m_{T+1}=m$ as well as $p_3=\cdots=p_T=p_{T+1}=p$. Moreover, we set $p_2=1$. Then, for all $t=2,\cdots,T+1$ we have:
    \begin{equation}
    \label{u-nablal_large}
        \begin{aligned}
            &\sum_{e=1}^E\sum_{t=1}^{T+1}\mathbb{E}\left[\left\Vert\tilde{u}^{(t)}_e-\nabla_e\ell(\w^{(t)})\right\Vert^2\right]\\
            \leq&32L_{\nabla \ell}^2\left(\dfrac{p+m}{m^2(1-(1-p)(1-\frac{m}{2}))}+\dfrac{1}{m^2}\right)\sum_{e=1}^E\sum_{t=1}^{T}\mathbb{E}\left[\left\Vert w^{(t+1)}_e-w^{(t)}_e\right\Vert^2\right]\\
            &+\dfrac{8(L_{\nabla a}^\circ)^2}{B}\sum_{e=1}^{E-1}\sum_{b=1}^B\sum_{t=2}^{T+1}\mathbb{E}\left[\left\Vert\tilde{ v}_{e,b}^{(t)}- \hat{v}_{e,b}^{(t)}\right\Vert^2\right]+
            \dfrac{8(L_{\nabla a}')^2}{B}\sum_{e=1}^{E-1}\sum_{b=1}^B\sum_{t=2}^{T+1}\mathbb{E}\left[\left\Vert\tilde{\by}_{e,b}^{(t)}-\hat{y}_{e,b}^{(t)}\right\Vert^2\right]\\
            &+4T\dfrac{\sigma^2}{B}\dfrac{(2-p)m-(1-p)m^2}{1-(1-p)(1-m)^2}+\dfrac{3}{m}\sum_{e=1}^E\mathbb{E}\left[\left\Vert\tilde{u}^{(1)}_e-\nabla_e\ell(\w^{(1)})\right\Vert^2\right].
        \end{aligned}
    \end{equation}
\begin{proof}
    For $t=2,3,\cdots,T$, we denote $\psi(t)$ as the last moment in which the sample is randomly obtained with $\mathcal{D}$ as of the $t$-th iteration. Specially,
    \begin{align*}
        \psi(t):=\max_{\tau\in\mathbb{S}_t}\tau,\quad\text{where }\mathbb{S}_t:=\{\tau=2,3,\cdots,t|\text{sampling randomly at iteration }\tau\}.
    \end{align*}
    Then, with the fact that the $p_2=1$, it holds for $\tau=2,\cdots,t$ that $\text{Pr}(\psi(t)=\tau)=\begin{cases}(1-p)^{t-2}, \text{ if }\tau=2 \\ p(1-p)^{t-\tau}, \text{ else}. \end{cases}$

    For $e=1,2,\cdots,E-1$ and $t=2,3,\cdots,T+1$, the error between the evaluated gradient and the true gradient satisfies:
    {\small
    \begin{equation}
    \label{u-nabla,new_large}
    \begin{aligned}
        &\tilde{u}^{(t)}_e-\nabla_e\ell(\w^{(t)})\\
        =&\underbrace{\sum_{\tau=\psi(t)}^t\dfrac{1}{B}\sum_{b=1}^Bm(1-m)^{t-\tau}\left(\nabla_2a_e(\tilde{y}_{e-1,b}^{(\tau)},w_e^{(\tau)})^{\tran}\tilde{v}_{e,b}^{(\tau)}-\nabla_2a_e(\hat{y}_{e-1,b}^{(\tau)},w_e^{(\tau)})^{\tran}\hat{v}_{e,b}^{(\tau)}\right)}_{:=\Xi_{e,1}}\\
        &+\sum_{\tau=\psi(t)}^t\dfrac{1}{B}\sum_{b=1}^Bm(1-m)^{t-\tau}\left(\nabla_2a_e(\hat{y}_{e-1,b}^{(\tau)},w_e^{(\tau)})^{\tran}\hat{v}_{e,b}^{(\tau)}-\nabla_e\ell(\w^{\tau})\right)\\
        &+\underbrace{\sum_{\tau=\psi(t)}^tm(1-m)^{t-\tau}\left(\nabla_e\ell(\w^{(\tau)})-\nabla_e\ell(\w^{(t)})\right)}_{:=\Xi_{e,2}}+\underbrace{(1-m)^{t+1-\psi(t)}\left(\nabla_e\ell(\w^{(\psi(t)-1)})-\nabla_e\ell(\w^{(t)})\right)}_{:=\Xi_{e,3}}\\
        &+(1-m)^{t+1-\psi(t)}\left(\tilde{u}_e^{(\psi(t)-1)}-\nabla_e\ell(\w^{(\psi(t)-1)})\right),
    \end{aligned}
    \end{equation}}
    where the first equation is from the momentum update rule. Moreover, we use $\Xi_{e,1},\Xi_{e,2},\Xi_{e,3}$ to denote some complex terms, which have been shown in Eq. \eqref{u-nabla,new_large}.

    Taking the $\ell_2$-norm and conditional expectation with respect to $\mathcal{F}^{(\psi(t))}$ on both sides of Eq. \eqref{u-nabla,new_large}, we can obtain:
    \begin{equation}
    \label{u-nabla,new21_large}
    \begin{aligned}
        &\mathbb{E}\left[\left\|\tilde{u}_e^{(t)}-\nabla_e\ell(\w^{(t)})\right\|^2\middle|\mathcal{F}^{(\psi(t))}\right]\\
        =&\mathbb{E}\left[\left\|\sum_{\tau=\psi(t)}^t\dfrac{1}{B}\sum_{b=1}^Bm(1-m)^{t-\tau}\left(\nabla_2a_e(\hat{y}_{e-1,b}^{(\tau)},w_{e}^{(\tau)})^{\tran}\hat{v}_{e,b}^{(\tau)}-\nabla_e\ell(\w^{\tau})\right)\right\|^2\middle|\mathcal{F}^{(\psi(t))}\right]\\
        &+\mathbb{E}\left[\left\|(1-m)^{t+1-\psi(t)}\left(\tilde{u}_e^{(\psi(t)-1)}-\nabla_e\ell(\w^{(\psi(t)-1)})\right)+\Xi_{e,1}+\Xi_{e,2}+\Xi_{e,3}\right\|^2\middle|\mathcal{F}^{(\psi(t))}\right]\\
        &+2\mathbb{E}\Bigg[\Bigg\langle\sum_{\tau=\psi(t)}^t\dfrac{1}{B}\sum_{b=1}^Bm(1-m)^{t-\tau}\left(\nabla_2a_e(\hat{y}_{e-1,b}^{(\tau)},w_e^{(\tau)})^{\tran}\hat{v}_{e,b}^{(\tau)}-\nabla_e\ell(\w^{\tau})\right),\\
        &\quad\quad\quad\quad(1-m)^{t+1-\psi(t)}\left(\tilde{u}_e^{(\psi(t)-1)}-\nabla_e\ell(\w^{(\psi(t)-1)})\right)+\Xi_{e,1}+\Xi_{e,2}+\Xi_{e,3}\Bigg\rangle\Bigg|\mathcal{F}^{(\psi(t))}\Bigg].
    \end{aligned}
    \end{equation}
    Thus, we can get:
    \begin{equation}
    \label{u-nabla,new2_large}
    \begin{aligned}
        &\mathbb{E}\left[\left\|\tilde{u}_e^{(t)}-\nabla_e\ell(\w^{(t)})\right\|^2\middle|\mathcal{F}^{(\psi(t))}\right]\\
        \leq&\dfrac{2}{B^2}\sum_{b=1}^B\mathbb{E}\left[\left\|\sum_{\tau=\psi(t)}^tm(1-m)^{t-\tau}\left(\nabla_2a_e(\hat{y}_{e-1,b}^{(\tau)},w_e^{(\tau)})^{\tran}\hat{v}_{e,b}^{(\tau)}-\nabla_e\ell(\w^{\tau})\right)\right\|^2\middle|\mathcal{F}^{(\psi(t))}\right]\\
        &+\mathbb{E}\left[\left\|\Xi_{e,1}+\Xi_{e,2}+\Xi_{e,3}\right\|^2\middle|\mathcal{F}^{(\psi(t))}\right]\\
        &+\mathbb{E}\left[\left\|(1-m)^{t+1-\psi(t)}\left(\tilde{u}_e^{(\psi(t)-1)}-\nabla_e\ell(\w^{(\psi(t)-1)})\right)+\Xi_{e,1}+\Xi_{e,2}+\Xi_{e,3}\right\|^2\middle|\mathcal{F}^{(\psi(t))}\right],
    \end{aligned}
    \end{equation}
    where the inequality is due to Cauchy-Schwarz inequality, Assumption \ref{assumption:unbiased} and the fact that samples in the batch are obtained independently.

    For the second term of the right-hand-side of Eq. \eqref{u-nabla,new2_large}, it holds that:
    \begin{equation}
    \label{eq123_large}
    \begin{aligned}
        &\mathbb{E}\left[\left\|\Xi_{e,1}+\Xi_{e,2}+\Xi_{e,3}\right\|^2\middle|\mathcal{F}^{(\psi(t))}\right]\\
        \leq&\dfrac{2}{B}\sum_{\tau=\psi(t)}^t\sum_{b=1}^Bm(1-m)^{t-\tau}\mathbb{E}\left[\left\|\nabla_2a_e(\tilde{y}_{e-1,b}^{(\tau)},w_e^{(\tau)})^{\tran}\tilde{v}_{e,b}^{(\tau)}-\nabla_2a_e(\hat{y}_{e-1,b}^{(\tau)},w_e^{(\tau)})^{\tran}\hat{v}_{e,b}^{(\tau)}\right\|^2\middle|\mathcal{F}^{(\psi(t))}\right]\\
        &+2\sum_{\tau=\psi(t)}^tm(1-m)^{t-\tau}\mathbb{E}\left[\left\|\nabla_e\ell(\w^{(\tau)})-\nabla_e\ell(\w^{(t)})\right\|^2\middle|\mathcal{F}^{(\psi(t))}\right]\\
        &+(1-m)^{t+1-\psi(t)}\mathbb{E}\left[\left\|\nabla_e\ell(\w^{(\psi(t)-1)})-\nabla_e\ell(\w^{(t)})\right\|^2\middle|\mathcal{F}^{(\psi(t))}\right],
    \end{aligned}
    \end{equation}
    where the inequality holds is due to the convexity of the $\ell_2$-norm.

    Moreover, for the last term, it also holds that:
    {\small
    \begin{equation}
    \label{psi_i+123_large}
    \begin{aligned}
        &\mathbb{E}\left[\left\|(1-m)^{t+1-\psi(t)}\left(\tilde{u}_e^{(\psi(t)-1)}-\nabla_e\ell(\w^{(\psi(t)-1)})\right)+\Xi_{e,1}+\Xi_{e,2}+\Xi_{e,3}\right\|^2\middle|\mathcal{F}^{(\psi(t))}\right]\\
        \leq&\dfrac{2}{B}\sum_{\tau=\psi(t)}^t\sum_{b=1}^Bm(1-m)^{t-\tau}\mathbb{E}\left[\left\|\nabla_2a_e(\tilde{y}_{e-1,b}^{(\tau)},w_e^{(\tau)})^{\tran}\tilde{v}_{e,b}^{(\tau)}-\nabla_2a_e(\hat{y}_{e-1,b}^{(\tau)},w_e^{(\tau)})^{\tran}\hat{v}_{e,b}^{(\tau)}\right\|^2\middle|\mathcal{F}^{(\psi(t))}\right]\\
        &+2\sum_{\tau=\psi(t)}^tm(1-m)^{t-\tau}\mathbb{E}\left[\left\|\nabla_e\ell(\w^{(\tau)})-\nabla_e\ell(\w^{(t)})\right\|^2\middle|\mathcal{F}^{(\psi(t))}\right]\\
        &+(1-m)^{t+1-\psi(t)}\mathbb{E}\left[\left\|\left(\tilde{u}_e^{(\psi(t)-1)}-\nabla_e\ell(\w^{(\psi(t)-1)})\right)+\left(\nabla_e\ell(\w^{(\psi(t)-1)})-\nabla_e\ell(\w^{(t)})\right)\right\|^2\middle|\mathcal{F}^{(\psi(t))}\right].
    \end{aligned}
    \end{equation}}

    With Assumption \ref{assumption:unbiased}, we can get:
    \begin{equation*}
    \begin{aligned}
        &\dfrac{1}{B^2}\mathbb{E}\left[\sum_{b=1}^B\sum_{e=1}^E\left\|\sum_{\tau=\psi(t)}^tm(1-m)^{t-\tau}\left(\nabla_2a_e(\hat{y}_{e-1,b}^{(\tau)},w_e^{(\tau)})^{\tran}\hat{v}_{e,b}^{(\tau)}-\nabla_e\ell(\w^{\tau})\right)\right\|^2\middle|\mathcal{F}^{(\psi(t))}\right]\\
        \leq&\left(\sum_{\tau=\psi(t)}^{t}m(1-m)^{t-\tau}\right)^2\dfrac{\sigma^2}{B}.
    \end{aligned}
    \end{equation*}

    Finally, with the same process of Lemma \ref{lemma3}, we can present the result that:
    \begin{equation*}
        \begin{aligned}
            &\sum_{e=1}^E\sum_{t=1}^{T+1}\mathbb{E}\left[\left\Vert\tilde{u}^{(t)}_e-\nabla_e\ell(\w^{(t)})\right\Vert^2\right]\\
            \leq&32L_{\nabla \ell}^2\left(\dfrac{p+m}{m^2(1-(1-p)(1-\frac{m}{2}))}+\dfrac{1}{m^2}\right)\sum_{e=1}^E\sum_{t=1}^{T}\mathbb{E}\left[\left\Vert w^{(t+1)}_e-w^{(t)}_e\right\Vert^2\right]\\
            &+\dfrac{8(L_{\nabla a}^\circ)^2}{B}\sum_{e=1}^{E-1}\sum_{b=1}^B\sum_{t=2}^{T+1}\mathbb{E}\left[\left\Vert\tilde{ v}_{e,b}^{(t)}- \hat{v}_{e,b}^{(t)}\right\Vert^2\right]+
            \dfrac{8(L_{\nabla a}')^2}{B}\sum_{e=1}^{E-1}\sum_{b=1}^B\sum_{t=2}^{T+1}\mathbb{E}\left[\left\Vert\tilde{\by}_{e,b}^{(t)}-\hat{y}_{e,b}^{(t)}\right\Vert^2\right]\\
            &+4T\dfrac{\sigma^2}{B}\dfrac{(2-p)m-(1-p)m^2}{1-(1-p)(1-m)^2}+\dfrac{3}{m}\sum_{e=1}^E\mathbb{E}\left[\left\Vert\tilde{u}^{(1)}_e-\nabla_e\ell(\w^{(1)})\right\Vert^2\right].
        \end{aligned}
    \end{equation*}
\end{proof}
\end{lemma}

Compared with Eq. \eqref{u-nablal}, the most difference of Eq. \eqref{u-nablal_large} is that the noise term $\sigma^2$ has been replaced to $\frac{\sigma^2}{B}$. And thus we can simply obtain the convergence of \ours in the large batch scenario.
\begin{lemma}[Convergence rate of \oursfu with large batch]
Suppose Assumption \ref{assumption:smoothness}, \ref{assumption:unbiased}, and \ref{assumption:compressor} hold. Then for Algorithm \ref{alg:clapping} there exist $\gamma,m,p>0$ such that: 
\begin{equation}
\label{descent lemma__GENERAL_e3}
\begin{aligned}
    \dfrac{1}{T}\sum_{t=1}^T\mathbb{E}\left[\left\Vert\nabla \ell(\w^{(t)})\right\Vert^2\right]\lesssim&\dfrac{\sigma}{\sqrt{BT}}+\dfrac{1}{T(1-\omega_B)^{E-1}(1-\omega_F)^{E-1}}.
\end{aligned}
\end{equation}
\end{lemma}
\begin{lemma}[Convergence rate of \oursfc with large batch]
Suppose Assumption \ref{assumption:smoothness}, \ref{assumption:unbiased}, and \ref{assumption:compressor} hold. Then for Algorithm \ref{alg:clapping} there exist $\gamma,m>0$ such that: 
\begin{equation}
\label{descent lemma__GENERAL_e3}
\begin{aligned}
    &\dfrac{1}{T}\sum_{t=1}^T\mathbb{E}\left[\left\Vert\nabla \ell(\w^{(t)})\right\Vert^2\right]\lesssim\dfrac{\sigma^{\frac{4}{3}}}{(BT)^{\frac{1}{3}}(1-\omega_B)^{\frac{4(E-1)}{3}}(1-\omega_F)^{\frac{4(E-1)}{3}}}\\
    &\quad\quad+\dfrac{1}{T}\left(\dfrac{1}{(1-\omega_B)^{E-1}(1-\omega_F)^{E-1}}+\dfrac{\omega_F^2+\omega_B}{(1-\omega_B)^{2(E-1)}(1-\omega_F)^2}+\dfrac{1}{(1-\omega_F)^{2(E-2)-1}}\right).
\end{aligned}
\end{equation}
\end{lemma}

\section{Additional discussions}
In this section, we present a series of additional discussions to the proposed algorithm, theoretical analysis and emprical efficiency.
\subsection{When should communication compression be enabled?}
\label{app:when_enable_compression}
\ours is designed for \emph{communication-bound} pipeline parallelism, where inter-stage transfers of activations/adjoints dominate iteration time under limited bandwidth.
In high-bandwidth intra-node deployments (e.g., NVLink/NVSwitch) where communication is not the bottleneck, compression may offer little benefit and can introduce overhead; in such settings it is reasonable to disable compression.
Our focus is therefore on regimes with constrained cross-node/over-the-internet links, where compression is often necessary for efficiency and \ours provides an optimization-stability mechanism under aggressive inter-stage compression (see Appendix~\ref{app:llama3b_bandwidth} for a bandwidth-sensitivity study).

\subsection{Additional details on multi-worker scenarios}
\label{appendix:aq-sgd}
In Section \ref{sec:aq-sgd}, we present the forward and backward process of AQ-SGD \citep{wang2022fine}. However, it is noteworthy that the convergence analysis of \citep[Appendix A.1.]{wang2022fine} assume that the machine $e$ computes the gradient of activation by
\begin{align}\label{aqsgd-update-wrong}
    v_{e-1}^{(t)}=\nabla_1 a_e(\tilde{y}^{(t)}_{x,e-1},w^{(t)}_e)\tran \mathcal{C}\left[\dfrac{\partial a_{e+1}\circ\cdots\circ a_{E}}{\partial \tilde{y}^{(t)}_{x,e}}(\tilde{y}^{(t)}_{x,e},w^{(t)}_{e+1},\cdots,w^{(t)}_E)\right],
\end{align}
where $a_{e+1}\circ\cdots\circ a_{E}$ denotes the composition of $a_{e+1},\cdots, a_{E}$. However, the parameters $w^{(t)}_{e+1},\cdots,w^{(t)}_E$ are held in different machines. Consequently, the computation of the partial gradient in the last term of \eqref{aqsgd-update-wrong} cannot avoids the communication between workers. Thus, it cannot be obtained losslessly through the communication compression during the propagation process. Thus, there causes a mismatch between the analysis in \citep{wang2022fine} and the reality. In fact, the gradient of activation $v_{e-1}^{(t)}$ can be obtained by the back-propagation as:
\begin{align}\label{aqsgd-update-right}
    v_{e-1}^{(t)}=\nabla_1 a_e(\tilde{y}^{(t)}_{x,e-1},w^{(t)}_e)\tran \mathcal{C}\left[\nabla_1 a_{e+1}(\tilde{y}^{(t)}_{x,e},w^{(t)}_{e+1})\tran\mathcal{C}\left[\cdots\mathcal{C}\left(\nabla_1 a_{E}(\tilde{y}^{(t)}_{x,E-1},w^{(t)}_{E})\right)^{\tran}\cdots\right]\right],
\end{align}
which calls for the analysis of error accumulation because of the multiple compression.

\subsection{On the relation to vanilla SGD and the $m=0$ case}
The convergence bounds of \oursfu in Lemma~\ref{lemma: FU_e_2} and Lemma~\ref{lemma: convergence_fu_e} are derived under the standing
assumption that the lazy--sampling probability is a fixed positive constant,
i.e., $p_t \equiv p \in (0,1]$, and that the stepsize $\gamma$ and the
momentum parameter $m$ are chosen accordingly. In this regime we obtain the
nonconvex rate
\[
  \frac1T\sum_{t=1}^T \mathbb{E}\bigl\|\nabla \ell(w^{(t)})\bigr\|^2
  \;\lesssim\;
  \frac{\sigma}{\sqrt{T}}
  \;+\;
  \frac{1}{T(1-\omega_F)(1-\omega_B)}
\]
up to logarithmic factors, where the variance term $\sigma^2$ comes from the
stochasticity of the activation gradients and the compression errors.

\smallskip
\noindent\textbf{(i) Behaviour as $p \to 0$.}
The proofs of Lemma~\ref{lemma: FU_e_2} and Lemma~\ref{lemma: convergence_fu_e} crucially require that $p$ is bounded away
from zero.  Indeed, the coefficient in front of the variance term,
\[
  \sigma^2 \cdot
  \frac{(2-p)m-(1-p)m^2}{1-(1-p)(1-m)^2},
\]
remains $\Theta(\sigma^2)$ as $p\to 0$.  In this limiting regime the algorithm keeps
reusing the same sample over a long horizon, and the bound only guarantees
convergence to a variance--dominated neighbourhood, rather than
$\|\nabla\ell(w^{(t)})\|\to 0$.  This behaviour is consistent with classical
SGD intuition when the data is not refreshed: the iterates converge to a
noise--dominated region determined by the fixed sample, instead of an exact
stationary point.  A similar phenomenon appears in the large--batch analysis
(see Eq.~\eqref{u-nablal_large}), where setting $p=0$ reduces the relevant ratio to~$1$ and
leaves a residual term of order $4\sigma^2/B$ that does not vanish with~$T$.
Therefore our rate statements should be interpreted for fixed $p>0$ and are
not intended to cover the degenerate case $p=0$.

\smallskip
\noindent\textbf{(ii) Why $m=0$ does not reduce to vanilla SGD.}
Even when $m=0$, Algorithm~\ref{alg:clapping} still performs forward and
backward passes with \emph{compressed} activations and adjoints, together with
\emph{error feedback} across pipeline stages.  The update direction used at
stage~$e$ is
\[
  \tilde u_e^{(t)}
  =
  \nabla_2 a_e\bigl(\tilde y_{e-1}^{(t)},w_e^{(t)}\bigr)^\top
  \tilde v_e^{(t)},
\]
where $(\tilde y,\tilde v)$ are generated by the compression and
error–feedback recursions.  Hence, the case $m=0$ in our framework corresponds
to SGD with a \emph{biased, compression–affected} gradient estimator, rather
than standard vanilla SGD on fresh, uncompressed samples.

The residual $\sigma^2$ term in Lemma~\ref{lemma: FU_e_2} and Lemma~\ref{lemma: convergence_fu_e} exactly captures the fact
that, in this setting, compression induces a persistent gradient bias that is
not completely averaged out when $m=0$.  This residual bias arises from the
pipeline compression structure and is present even without lazy sampling; it
is not a contradiction with the classical convergence guarantees for vanilla
SGD.

\smallskip
\noindent\textbf{(iii) Recovery of the standard nonconvex SGD rate.}
If we additionally disable compression ($\omega_F=\omega_B=0$) and turn off
lazy sampling ($p_t\equiv 1$), Algorithm~\ref{alg:clapping} reduces to
pipeline–parallel vanilla SGD.  In this special case the error--feedback
recursions collapse and our analysis recovers the standard nonconvex rate
\[
  \frac1T\sum_{t=1}^T \mathbb{E}\bigl\|\nabla \ell(w^{(t)})\bigr\|^2
  = \mathcal{O}\!\left(\frac{\sigma}{\sqrt{T}}\right)
\]
without additive residual terms.  Thus there is no inconsistency between our
bounds and the known behaviour of vanilla SGD; the additive $\sigma^2$ term
only appears in the genuinely compressed-and-lazy–sampled regime we target.

\subsection{Lazy sampling, convergence and generalization}
\label{app:discuss_pt}
Lazy sampling keeps the marginal data distribution unchanged: at each iteration we either reuse the previous sample $x^{(t-1)}$ with probability $1-p_t$ or draw a fresh sample $x^{(t)}\sim\mathcal{D}$ with probability $p_t$. Hence the sequence $\{x^{(t)}\}_{t\ge1}$ is no longer i.i.d., but every marginal law is still $\mathcal{D}$. Lazy sampling therefore introduces temporal correlations without systematically excluding any part of the data distribution.

The convergence bounds of \ours make the role of $\{p_t\}$ explicit. For the FC variant, the variance term
\[
  \sigma^2 \cdot \frac{(2-p)m-(1-p)m^2}{1-(1-p)(1-m)^2}
\]
remains of order $\Theta(\sigma^2)$ when $p_t$ is chosen too small, which prevents the averaged gradient norm from vanishing. Intuitively, as $p_t\to 0$ the algorithm repeatedly optimizes on almost the same sample, so the optimization dynamics are dominated by the noise induced by that sample and the method can only converge to a noise-dominated neighborhood of a stationary point.

On the other hand, if $p_t$ is taken very close to $1$, a new sample is seen at almost every step and the temporal correlation disappears, but the error–feedback iterates have fewer opportunities to cancel the compression bias. In this regime \oursfc may converge more slowly, whereas \oursfu effectively approaches standard (uncompressed) SGD and thus loses its communication advantages.

These observations suggest that the lazy–sampling schedule $\{p_t\}$ should be chosen in an \emph{intermediate} regime. It must be bounded away from zero to ensure that the effective number of distinct samples grows linearly with $T$, while it need not (and in general should not) be exactly one. In our analysis $p_t$ thus plays the role of a bias–variance trade-off parameter: moderate values of $p_t$ guarantee both sufficient exploration of the data distribution and enough sample reuse for error–feedback to attenuate compression bias.


\subsubsection{Adaptive refresh rates and scope of interpretation}
Only fixed refresh probabilities $p$ were evaluated in the additional experiments. Although a time-varying or adaptive $p$ is compatible with the algorithmic update, neither an adaptive-policy experiment nor an adaptive-schedule convergence guarantee is available. Compression-induced implicit regularization remains a possible explanation for occasional differences in held-out scores, rather than a causally established mechanism.

\section{Memory overhead analysis for LLMs pre-training with \ours}
\label{appendix: memory overhead}
In this section, we present the analysis for the memory overhead introduced by \ours in the pre-training tasks for LLMs.

\textbf{Basic setup. }We consider a pre-training task on LLaMA-2-based models with SwiGLU activation \citep{touvron2023LLaMA} running on Nvidia A100 80G GPUs. The dataset is C4-en \citep{raffel2020exploring}, which is primarily intended for pre-training language models and word representations on a large scale. We use the T5-base tokenizer with a sequence length of \(s=4096\), resulting in a total of approximately \(|\mathcal{D}|\approx45.6M\) training samples. The micro-batch size is \(B=16\), and the total batch size is set to 256. The optimizer used is Adam \citep{kingma2014adam}, and we employ BF16 precision, where each parameter requires 2 bytes for storage.

\begin{table}[t]
\centering
\caption{Hyperparameter configurations for LLaMA-2 models of different scales. `Num\_workers' denotes the total workers for a training pipeline.}
\label{table: pretraining setting}
\begin{threeparttable}
\begin{tabular}{cccccc}
\toprule
Parameters & Hidden size ($h$) & Heads ($a$) & Layers ($L$)   & Vocabulary size ($V$)   & Num\_workers   \\ \midrule
7B         & 4096       & 32    & 32     & 32000    & 2                \\  
13B         & 5120   & 40    & 40    & 32000   & 4                \\  
70B         & 8192          & 64    & 80   & 32000      & 16                \\ \bottomrule
\end{tabular}    
\end{threeparttable}
\end{table}

If we assume the intermediate size \(h_{\text{ff}} = \frac{8h}{3}\), we use vanilla gradient checkpointing (GCP) \cite{chen2016training} to save activations and do not use GQA for 70B models. The basic memory overhead for parameters, gradients, optimizer states, and activations is:
$$\underbrace{4Vh + 48Lh^2}_{\text{parameters, gradients, optimizer states}} + \underbrace{LBsh}_{\text{Activation memory with GCP}}$$
For \ours, it should cache the term \(\tilde{y}_e\) in workers $e$ and \(e+1\), as well as the term \(\tilde{v}_{e-1}\) in workers \(e-1\) and $e$. Since we often split the model at the end of some transformer block, both \(\tilde{y}_e\) and \(\tilde{v}_{e-1}\) have a size of Bsh. Thus, the memory overhead of \ours is:
$$4 \times (\text{num\_workers} - 1)Bsh.$$
Similarly, the memory overhead of AQ-SGD can be expressed as:
$$2 \times (\text{num\_workers} - 1)|\mathcal{D}|sh.$$

We present the memory overhead analysis for pre-training LLaMA-2 models with different communication compression algorithms under pipeline parallelism in Table \ref{table: pretraining memory llama2}. It can be observed that the memory overhead introduced by \ours is less than 11\% of the basic memory even when the model size scales to 70B, which is acceptable for practical pre-training tasks. Meanwhile, AQ-SGD requires thousands of GBs to store the sample-wise cache, making it unsuitable for pre-training tasks.

\begin{table}[t]
\centering
\caption{The memory overhead for pre-training LLaMA-2 models with different model size. `Basic M.' means the memory overhead for parameters, gradients, optimizer states, and activation memorys. `\ours M.' means the additional memory introduced by \ours. `AQ-SGD M.' means the additional memory introduced by AQ-SGD. `Ratio' denotes the ration between the additional memory of \ours and the basic memory.}
\label{table: pretraining memory llama2}
\begin{threeparttable}
\begin{tabular}{ccccc}
\toprule
Parameters & Basic M. (GB) & \ours M. (GB) & Ratio & AQ-SGD M. (GB)   \\ \midrule
7B         & 65.0       & 2.0    & 3.07\%     & 2850.0               \\  
13B         & 118.1   & 7.5    & 6.35\%    & 10687.5             \\  
70B         & 562.0          & 60.0   & 10.68\%   & 85500.0                  \\ \bottomrule
\end{tabular}    
\end{threeparttable}
\end{table}

\section{Experimental details}
\label{app:experiment}
In this section, we present the details of our numerical experiments, which were discussed in Section \ref{section:experiment}. Additionally, we provide additional experimental results that were not included in the main text due to space limitations.
\subsection{Synthetic Logistic Regression}
\label{app:toyexperiment}
\begin{figure*}[t]
\centering
\includegraphics[width=0.5\textwidth]{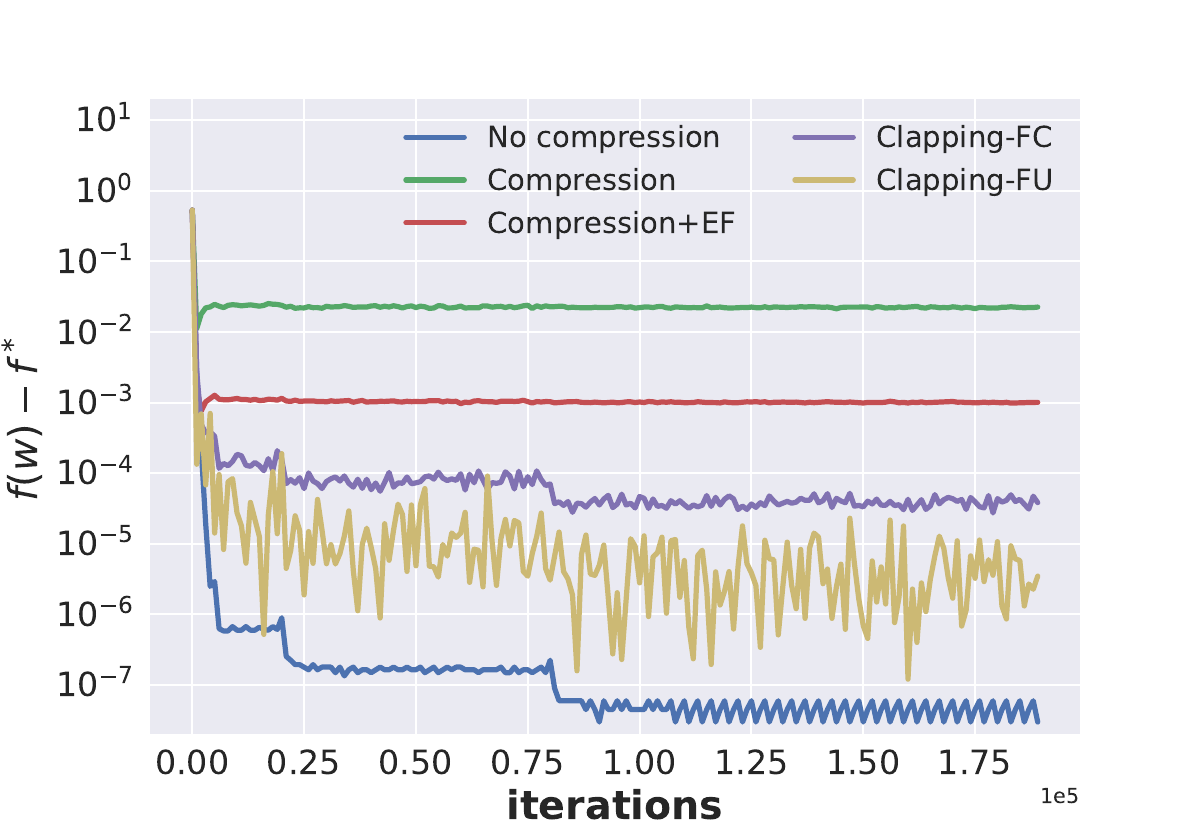}
\caption{\small The gap between the current loss and the optimal loss for different algorithms for the logistic regression problem.}
\label{fig:appendix_exp_toy_diff}
\end{figure*}

Herein, we consider the following logistic regression task with the objective function given by:
\begin{align}
\label{exp:toymodel_detail}
f(w,b)=\mathbb{E}_{(\xi,\zeta)\sim\mathcal{D}}\left[\ln(1 + \exp(-\zeta\cdot (\xi^\top w + b))) + C_r \Vert w \Vert^2 + C_r b^2\right],
\end{align}
where $w$ and $b$ denote the regression parameters, and $C_r$ is a regularization parameter. $\mathcal{D}$ represents the finite sample set. Here, we the sample size $|\mathcal{D}|$. For each stochastic sample $(\xi,\zeta)$ in $\mathcal{D}$, $\xi$ is a 200-dimensional vector generated as the standard $\xi^*+\varepsilon$, where each element of $\xi^*$ is independently normal distribution drawn from n $\mathcal{N}(0,0.5)$ and each element of $\varepsilon$ is independently drawn from $\mathcal{N}(0,0.3)$. 
Let $C_r=0.005$. Then, Equation \eqref{exp:toymodel_detail} is strongly convex. We employ gradient descent to find the minimum of Equation \eqref{exp:toymodel_detail}, denote as $f^*$. Subsequently, we split the model into two parts according to the formulation in \eqref{pipeline parallel}, where
\begin{align*}
y_1=-\zeta\cdot (\xi^\top w + b),\quad y_2=\ln(1+\exp(y_1)) + C_r \Vert w \Vert^2 + C_r b^2.
\end{align*}
After obtaining $y_1$ during the forward-propagation process, we introduce an error to $y_1$ as $\delta_1$ to simulate activation compression and use $\tilde{y}_1=y_1+\delta_1$ to compute $y_2$, where $\delta_1$ follows the uniform distribution $U(-0.2,0.2)$.
We set the batch size to 128. Then, we run the moving-average Stochastic Gradient Descent (SGD) algorithm for 200,000 iterations with different error compensation strategies, including direct compression without error feedback \textbf{(Compression)}, compression with error feedback but without lazy sampling \textbf{(Compression + EF)}, \oursfclast, and \oursfulast. We compare these strategies with the case of no compression \textbf{(No Compression)}. The moving-average term  is set to 0.9, and the step-size  is initialized to 0.1. When the iteration number , the step-size is multiplied by 0.5, and the historical gradient is cleared.
Figure \ref{fig:appendix_exp_toy_diff} depicts the gap between the current loss and the optimal loss for different algorithms. Both direct compression and compression with error feedback but without lazy sampling fail to converge as the step-size decreases. However, both \oursfc and \oursfu can converge as the step-size decreases, and \oursfu convergence faster than \oursfc to a smaller minimum. The error introduced by activation compression only slows down the convergence process, which validates our theoretical findings in Section \ref{sec:convergence}. 

\subsection{Training ResNet-18 on CIFAR-10}
We trained a ResNet-18 \citep{he2016deep} model on CIFAR-10 \citep{krizhevsky2009learning} dataset on a NVIDIA A100 80G GPU with communication compression algorithms for pipeline-parallel learning including directly compression and \texttt{Clapping}-\textbf{FC}. The basic setting is similar to that of \citep{rudakov2023activations}. We split the model into 4 parts and used 3 direct quantization \citep{alistarh2017qsgd} with different bits to simulate the communication compression. We set the batch size to 128. We trained 100 epochs with directly compression and trained with \oursfc for the same number of iterations. We used the SGD optimizer with momentum 0.9 and weight decay 5e-4. The learning rate is initialized to 0.01 and scheduled by a cosine annealing scheduler with $T_{\max}=200$.  During inference, we obtain the test accuracy with the sample compression as training and without any compression, respectively.

\begin{table}[t]
\caption{The best test accuracy in training ResNet-18 on CIFAR-10.}
\label{table:app_resnet_cifar10}
\setlength{\tabcolsep}{4pt}
\renewcommand{\arraystretch}{1.1}
\centering
\begin{threeparttable}
\begin{tabular}{cccccccc}
\toprule
& \multirow{2}{*}{Strategy} & 
\multirow{2}{*}{\begin{tabular}[c]{@{}c@{}} Direct \\ comp. \end{tabular}} & \multicolumn{5}{c}{\ours with $p_t$}                   \\ \cline{4-8} 
& &                          & 0.2 & 0.4 & 0.6 & 0.8 & 1.0 \\ \midrule
\multirow{6}{*}{C. on} & f5-b5                        &  \multicolumn{1}{c}{\begin{tabular}[c]{@{}c@{}} 0.8604 \\ (0.00170) \end{tabular}}   &   \multicolumn{1}{c}{\begin{tabular}[c]{@{}c@{}} 0.8719 \\ (0.01041) \end{tabular}}   &  \multicolumn{1}{c}{\begin{tabular}[c]{@{}c@{}} 0.8889 \\ (0.00257) \end{tabular}}    &  \multicolumn{1}{c}{\begin{tabular}[c]{@{}c@{}} \textbf{0.9034} \\ \textbf{(0.00117)} \end{tabular}}    &  \multicolumn{1}{c}{\begin{tabular}[c]{@{}c@{}} 0.9016 \\ (0.00066) \end{tabular}}    &  \multicolumn{1}{c}{\begin{tabular}[c]{@{}c@{}} 0.8963 \\ (0.00209) \end{tabular}}    \\
& f6-b6                        &   \multicolumn{1}{c}{\begin{tabular}[c]{@{}c@{}} 0.9024 \\ (0.00709) \end{tabular}}   &  \multicolumn{1}{c}{\begin{tabular}[c]{@{}c@{}} 0.8985 \\ (0.00997) \end{tabular}}    &  \multicolumn{1}{c}{\begin{tabular}[c]{@{}c@{}} 0.9119 \\ (0.00271) \end{tabular}}    &  \multicolumn{1}{c}{\begin{tabular}[c]{@{}c@{}} \textbf{0.9185} \\ \textbf{(0.00466)} \end{tabular}}    &  \multicolumn{1}{c}{\begin{tabular}[c]{@{}c@{}} 0.9150 \\ (0.00373) \end{tabular}}    &  \multicolumn{1}{c}{\begin{tabular}[c]{@{}c@{}} 0.9133 \\ (0.00397) \end{tabular}}  \\
& f6-b8                        &   \multicolumn{1}{c}{\begin{tabular}[c]{@{}c@{}} 0.9056 \\ (0.00558) \end{tabular}}   &  \multicolumn{1}{c}{\begin{tabular}[c]{@{}c@{}} {0.9014} \\ {(0.00199)} \end{tabular}}    &  \multicolumn{1}{c}{\begin{tabular}[c]{@{}c@{}} 0.9109 \\ (0.00496) \end{tabular}}    &  \multicolumn{1}{c}{\begin{tabular}[c]{@{}c@{}} 0.9140 \\ (0.00805) \end{tabular}}    &  \multicolumn{1}{c}{\begin{tabular}[c]{@{}c@{}} \textbf{0.9184} \\ \textbf{(0.00129)} \end{tabular}}    &  \multicolumn{1}{c}{\begin{tabular}[c]{@{}c@{}} 0.9160 \\ (0.00535) \end{tabular}}\\ \midrule
\multirow{6}{*}{C. off}          & f5-b5                        &  \multicolumn{1}{c}{\begin{tabular}[c]{@{}c@{}} 0.8623 \\ (0.00271) \end{tabular}}   &   \multicolumn{1}{c}{\begin{tabular}[c]{@{}c@{}} 0.8733 \\ (0.01010) \end{tabular}}   &  \multicolumn{1}{c}{\begin{tabular}[c]{@{}c@{}} 0.8907 \\ (0.00187) \end{tabular}}    &  \multicolumn{1}{c}{\begin{tabular}[c]{@{}c@{}} \textbf{0.9045} \\ \textbf{(0.00104)} \end{tabular}}    &  \multicolumn{1}{c}{\begin{tabular}[c]{@{}c@{}} 0.9022 \\ (0.00039) \end{tabular}}    &  \multicolumn{1}{c}{\begin{tabular}[c]{@{}c@{}} 0.8982 \\ (0.00197) \end{tabular}}    \\
& f6-b6                        &   \multicolumn{1}{c}{\begin{tabular}[c]{@{}c@{}} 0.9025 \\ (0.00643) \end{tabular}}   &  \multicolumn{1}{c}{\begin{tabular}[c]{@{}c@{}} 0.8986 \\ (0.00978) \end{tabular}}    &  \multicolumn{1}{c}{\begin{tabular}[c]{@{}c@{}} 0.9123 \\ (0.00282) \end{tabular}}    &  \multicolumn{1}{c}{\begin{tabular}[c]{@{}c@{}} \textbf{0.9186} \\ \textbf{(0.00438)} \end{tabular}}    &  \multicolumn{1}{c}{\begin{tabular}[c]{@{}c@{}} 0.9150 \\ (0.00412) \end{tabular}}    &  \multicolumn{1}{c}{\begin{tabular}[c]{@{}c@{}} 0.9137 \\ (0.00373) \end{tabular}}  \\
& f6-b8                        &   \multicolumn{1}{c}{\begin{tabular}[c]{@{}c@{}} 0.9058 \\ (0.00551) \end{tabular}}   &  \multicolumn{1}{c}{\begin{tabular}[c]{@{}c@{}} {0.9015} \\ {(0.00179)} \end{tabular}}    &  \multicolumn{1}{c}{\begin{tabular}[c]{@{}c@{}} 0.9109 \\ (0.00491) \end{tabular}}    &  \multicolumn{1}{c}{\begin{tabular}[c]{@{}c@{}} 0.9139 \\ (0.00798) \end{tabular}}    &  \multicolumn{1}{c}{\begin{tabular}[c]{@{}c@{}} \textbf{0.9184} \\ \textbf{(0.00121)} \end{tabular}}    &  \multicolumn{1}{c}{\begin{tabular}[c]{@{}c@{}} 0.9161 \\ (0.00533) \end{tabular}}  \\ \bottomrule
\end{tabular}
\end{threeparttable}
\begin{tablenotes}
     \item \textit{'f[A]-b[B]'} means compressing activation to A bits in forward propagation and compressing gradient to B bits in backward propagation. \textit{'C. on'} means taking the same compression during inference and \textit{'C. off'} means taking no compression.
\end{tablenotes}
\end{table}

Table \ref{table:app_resnet_cifar10} shows the average best test accuracy of directly compression and \ours with different lazy sampling coefficient $p_t$ with different compression strategies over 3 independent runs. We can observe that the error feedback technique and lazy sampling strategy can improve the prediction accuracy in the communication compression of pipeline parallelism.

\subsection{Fine-tuning LLMs}
\subsubsection{Fine-tuning on GLUE benchmark}
\label{appendix: finetune glue}
We fine-tune pre-trained RoBERTa-large \citep{liu2019roberta} dataset on GLUE benchmark \citep{wang2018glue} with communication compression algorithms including \ours and direct compression with EF21 \citep{richtarik2021ef21} and compare then with fine-tuning without compression on two Nvidia A800 GPUs. We use the TopK compression with 30\% elements at the middle of the networks. The batch size are all set 64 and the learning rate is 1e-5. For each dataset, we fine-tune the model for 10 epochs using both fine-tuning without compression and the EF21 algorithm. We also tune fine-tune the model for the same number of iterations as in 10 epochs by \ours with lazy sampling coefficient $p=0.5$. As Table \ref{table: fine-tuning result_GLUE} illustrates, \ours outperforms EF21 in majority of tasks and even outperforms fine-tuning without communication compression in tasks including MRPC and RTE. And \ours achieves the highest average score among all the algorithms.


\subsubsection{Additional Clapping-FU results on GLUE}
\label{app:rebuttal-glue-fu}
The RoBERTa-large GLUE results in the main text report FC (Table~\ref{table: fine-tuning result_GLUE}). We additionally record the corresponding FU scores for all eight tasks in Table~\ref{tab:rebuttal-glue-fu}. The FC column is repeated to facilitate direct comparison. Scores and their average follow the reporting scale used in the original GLUE experiment.

\begin{table}[t!]
\centering
\caption{RoBERTa-large GLUE task scores for Clapping-FC and Clapping-FU (higher is better).}
\label{tab:rebuttal-glue-fu}
\setlength{\tabcolsep}{4pt}
\renewcommand{\arraystretch}{1.12}
\begin{tabular}{@{}l*{9}{c}@{}}
\toprule
Method & MNLI & SST-2 & MRPC & CoLA & QNLI & QQP & RTE & STS-B & Avg. \\
\midrule
Clapping-FC & $89.94$ & $96.06$ & $90.50$ & $63.98$ & $94.04$ & $91.89$ & $84.19$ & $91.71$ & $87.79$ \\
Clapping-FU & $89.92$ & $96.18$ & $90.25$ & $64.20$ & $94.17$ & $91.97$ & $83.82$ & $92.27$ & $87.85$ \\
\bottomrule
\end{tabular}
\end{table}

\subsubsection{Fine-tuning LLaMA models with TopK compressor}
\label{appendix: finetune wiki}
We fine-tuned a pretrained LLaMA-2 7B \citep{touvron2023LLaMA} model and a LLaMA-3 8B \citep{grattafiori2024llama} model on Wikitext dataset (wikitext-2-raw-v1 version) \citep{merity2016pointer} on two NVIDIA A100 80G GPUs with communication compression algorithms for pipeline-parallel learning including directly compression, AQ-SGD, EF21 and \oursfclast. The block size was set to 1024. We used batch size 8, and we fine-tuned the model for 4 epochs for the competitive algorithms and we use the same iterations for \ours, respectively. We use the SGD optimizer and FP16 for fine-tuning. The learning rate is initialized from $2\times 10^{-5}$ and scheduled by a linear scheduler. We compare different compression algorithms with the compressor of Top-5\% \citep{wangni2018gradient}. For each model, we independently repeat \ours with different lazy sampling coefficient $p_t$ including $\{0.3,0.4,0.5\}$.

Table \ref{table: fine-tuning result_wiki} has present the evaluation accuracy of different approaches under Top-5\% compressor. It can be observed that \ours outperforms other algorithms, including direct compression, EF21, and AQ-SGD. By tuning the lazy-sampling coefficient $p$, \textbf{\oursfc achieves 95\% communication saving with less than 0.5\% error in practical fine-tuning tasks}.


\subsubsection{Independent runs and WikiText evaluation}
\label{app:rebuttal-seeds}
For the LLaMA-2 7B and LLaMA-3 8B WikiText-2 fine-tuning experiments, the original run used seed 555. Two independent repetitions with seeds 42 and 123 were added. The model, data, training budget, Top-5\% compressor, and, for FC and FU, the refresh probability $p=0.5$ were held fixed across these runs. Table~\ref{tab:rebuttal-three-seed} reports the mean and standard deviation of validation accuracy over the three seeds.

\begin{table}[t!]
\centering
\caption{WikiText-2 validation accuracy over three seeds (higher is better). FC and FU use $p=0.5$ and Top-5\% compression.}
\label{tab:rebuttal-three-seed}
\setlength{\tabcolsep}{3.8pt}
\renewcommand{\arraystretch}{1.12}
\begin{tabular}{@{}l*{6}{c}@{}}
\toprule
Model & No comp. & Direct & EF21 & AQ-SGD & Clapping-FC & Clapping-FU \\
\midrule
LLaMA-2 7B
& \multicolumn{1}{c}{\begin{tabular}[c]{@{}c@{}}0.5954\\(0.0015)\end{tabular}}
& \multicolumn{1}{c}{\begin{tabular}[c]{@{}c@{}}0.5526\\(0.0055)\end{tabular}}
& \multicolumn{1}{c}{\begin{tabular}[c]{@{}c@{}}0.5567\\(0.0102)\end{tabular}}
& \multicolumn{1}{c}{\begin{tabular}[c]{@{}c@{}}0.5684\\(0.0026)\end{tabular}}
& \multicolumn{1}{c}{\begin{tabular}[c]{@{}c@{}}0.5855\\(0.0038)\end{tabular}}
& \multicolumn{1}{c}{\begin{tabular}[c]{@{}c@{}}0.5960\\(0.0015)\end{tabular}} \\
LLaMA-3 8B
& \multicolumn{1}{c}{\begin{tabular}[c]{@{}c@{}}0.5953\\(0.0034)\end{tabular}}
& \multicolumn{1}{c}{\begin{tabular}[c]{@{}c@{}}0.5648\\(0.0025)\end{tabular}}
& \multicolumn{1}{c}{\begin{tabular}[c]{@{}c@{}}0.5649\\(0.0033)\end{tabular}}
& \multicolumn{1}{c}{\begin{tabular}[c]{@{}c@{}}0.5703\\(0.0028)\end{tabular}}
& \multicolumn{1}{c}{\begin{tabular}[c]{@{}c@{}}0.5893\\(0.0045)\end{tabular}}
& \multicolumn{1}{c}{\begin{tabular}[c]{@{}c@{}}0.5919\\(0.0008)\end{tabular}} \\
\bottomrule
\end{tabular}
\end{table}

The results were compared with direct compression using paired, two-sided $t$-tests on the three matched seeds. Table~\ref{tab:rebuttal-seed-significance} reports the stated unadjusted $p$-values. All four comparisons pass a Holm correction at family-wise significance level $0.05$; no significance claim against NoComp is made.

\begin{table}[t!]
\centering
\caption{Paired two-sided $t$-test $p$-values for WikiText-2 accuracy compared with direct compression; matched seeds 42, 123, and 555.}
\label{tab:rebuttal-seed-significance}
\begin{tabular}{lcc}
\toprule
Model & FC vs. Direct & FU vs. Direct \\
\midrule
LLaMA-2 7B & $0.0242$ & $0.0050$ \\
LLaMA-3 8B & $0.0230$ & $0.0037$ \\
\bottomrule
\end{tabular}
\end{table}


\subsubsection{FC/FU comparison across refresh rates and compressor strengths}
\label{app:rebuttal-fu-sweeps}
Table~\ref{tab:rebuttal-fu-p-sweep} gives the single-seed (seed 555) WikiText-2 validation accuracies when $p$ varies over $\{0.3,0.4,0.5\}$, using the Top-5\% compressor. The FC values also appear in the main experimental table (Table~\ref{table: fine-tuning result_wiki}); the FU values are additional results. These single-seed observations are distinct from the three-seed summary in Table~\ref{tab:rebuttal-three-seed}.

\begin{table}[t!]
\centering
\caption{WikiText-2 validation accuracy with Top-5\% compression for different refresh probabilities (higher is better).}
\label{tab:rebuttal-fu-p-sweep}
\setlength{\tabcolsep}{4.5pt}
\renewcommand{\arraystretch}{1.15}
\begin{tabular}{@{}l*{6}{c}@{}}
\toprule
\multirow{2}{*}{Model} & \multicolumn{3}{c}{Clapping-FC} & \multicolumn{3}{c}{Clapping-FU} \\
\cmidrule(lr){2-4}\cmidrule(l){5-7}
 & $p=0.3$ & $p=0.4$ & $p=0.5$ & $p=0.3$ & $p=0.4$ & $p=0.5$ \\
\midrule
LLaMA-2 7B & $0.5960$ & $0.5920$ & $0.5877$ & $0.5923$ & $0.5944$ & $0.5947$ \\
LLaMA-3 8B & $0.5865$ & $0.5969$ & $0.5887$ & $0.5862$ & $0.5910$ & $0.5912$ \\
\bottomrule
\end{tabular}
\end{table}

A separate, stronger-compression setting retains only the largest 1\% of elements, with $p=0.4$ for FC and FU. The corresponding WikiText-2 validation accuracies are given in Table~\ref{tab:rebuttal-top-one}.

\begin{table}[H]
\centering
\caption{WikiText-2 validation accuracy with Top-1\% compression ($p=0.4$ for FC/FU; higher is better).}
\label{tab:rebuttal-top-one}
\small
\begin{tabular}{lccc}
\toprule
Model & Direct & Clapping-FC & Clapping-FU \\
\midrule
LLaMA-2 7B & $0.5339$ & $0.5705$ & $0.5741$ \\
LLaMA-3 8B & $0.4627$ & $0.5745$ & $0.5761$ \\
\bottomrule
\end{tabular}
\end{table}

FC performs compressed error-feedback communication at every step, including refresh steps. FU uses uncompressed boundary transfers on refresh steps and compressed error-feedback transfers on reuse steps. For a constant refresh probability $p$, FC has $1/(1-p)$ times as many compressed-step opportunities as FU in expectation; this factor equals approximately $1.67$ at $p=0.4$. This counts compressed steps rather than specifying an end-to-end byte ratio, which additionally depends on message size and the compressor.

\subsubsection{Pipeline-depth sensitivity on LLaMA-2 13B}
\label{app:rebuttal-13b-depth}
\begin{table}[t]
\centering
\caption{The evaluation perplexity ($\downarrow$) and accuracy ($\uparrow$) for the fine-tuning task for LLaMA-2 13B under Wikitext with \oursfc with different number of workers ($E$).}
\label{table:leiji_ppl_acc}
\begin{threeparttable}
\begin{tabular}{rccc}
\toprule
Number of workers ($E$) & Validation loss $\downarrow$ & evaluation perplexity & Validation accuracy   \\ \midrule
1         & 1.813602 & 6.1325       & 0.5823              \\  
2         & 1.915392 & 6.7896 & 0.5654          \\  
4         & 2.077440 & 7.9840 & 0.5422                   \\ 
9         & 2.603682 & 13.5134 & 0.4606             \\16 & 2.842107 & 17.1519 & --- \\
32 & 3.227581 & 25.2186 & --- \\
\bottomrule
\end{tabular}    
\end{threeparttable}
\end{table}

\textbf{Ablations on the impact of the number of workers on performance.}
In Section~\ref{sec:convergence}, we showed that communication compression errors can accumulate exponentially with the number of workers $E$. To empirically examine this effect, we fine-tuned a pre-trained LLaMA-2 13B model \citep{touvron2023LLaMA} on WikiText for 4 epochs using \oursfc with a Top-20\% compressor and a fixed refresh probability of $p=0.5$. We varied the number of workers over $E\in\{1,2,4,9,16,32\}$ while keeping the other configurations unchanged, except for increasing the initial learning rate to $5\times10^{-5}$ to amplify the effect of error accumulation. Table~\ref{table:leiji_ppl_acc} reports validation loss and perplexity across all six configurations, together with validation accuracy for $E\leq9$. The results show that validation loss and perplexity increase as the number of workers grows, consistent with the accumulation of compression errors across pipeline stages.


\subsubsection{Sequential non-IID data-order stress test}
\label{app:rebuttal-sequential-noniid}
We consider a single sequential-order WikiText experiment on LLaMA-2 13B. Examples are traversed in their stored order with random shuffling disabled. The sequence length is 1,024, micro-batch size is one per GPU, global batch size is eight, and training lasts four epochs in FP16 using SGD with learning rate $5\times10^{-5}$. Lazy-only and Clapping-FU use $p=0.5$; the recorded methods are NoComp, Lazy-only (sample reuse without compression), and Clapping-FU (sample reuse with compression). All other settings are held fixed. Table~\ref{tab:rebuttal-noniid} reports the held-out results.

\begin{table}[t!]
\centering
\caption{LLaMA-2 13B WikiText sequential non-IID stress test.}
\label{tab:rebuttal-noniid}
\begin{tabular}{lcc}
\toprule
Method & Evaluation loss $\downarrow$ & Accuracy $\uparrow$ \\
\midrule
NoComp & $1.76885$ & $0.59511$ \\
Lazy-only & $1.76169$ & $0.59598$ \\
Clapping-FU & $1.76027$ & $0.59614$ \\
\bottomrule
\end{tabular}
\end{table}

This experiment tests one sequential data stream and one random seed; it does not establish generalization guarantees for arbitrary non-IID orderings.


\subsubsection{Generated-output quality: instruction tuning and LLM-as-a-Judge}
\label{app:rebuttal-judge}
For an output-quality stress test, LLaMA-2 7B is instruction-tuned for 500 steps on UltraChat 200K~\cite{ding2023enhancing} using NoComp, Direct, FC, and FU. The methods share the initialization, random seed, token budget, batch size, and decoding configuration. The learning rate is fixed to the intentionally aggressive value $5\times10^{-5}$ for all four methods; this is a matched stress test rather than a separate learning-rate search for each method.

The evaluation uses 200 shared held-out prompts and a separately fine-tuned LLaMA-2 7B judge. Each response pair is assessed in both A/B orientations. In Table~\ref{tab:rebuttal-judge}, ``X wins'' means that the judge favors the first method named in the comparison; ``Y wins'' favors the second method; ``Tie'' denotes equal preference; ``Abstain'' means no unambiguous judgment was extracted. Values are reported to one decimal place, with minor deviations from 100\% possible due to rounding.

\begin{table}[t!]
\centering
\caption{Position-balanced pairwise LLM-as-a-Judge judgments on 200 UltraChat prompts after 500-step LLaMA-2 7B instruction tuning at learning rate $5\times10^{-5}$. Each method pair is judged in both orderings.}
\label{tab:rebuttal-judge}
\setlength{\tabcolsep}{8pt}
\begin{tabular}{lrrrr}
\toprule
Comparison (X vs. Y) & X wins & Tie & Y wins & Abstain \\
\midrule
FC vs. Direct & $51.0\%$ & $7.0\%$ & $36.5\%$ & $5.5\%$ \\
FU vs. Direct & $59.8\%$ & $5.0\%$ & $29.0\%$ & $6.3\%$ \\
Direct vs. NoComp & $22.5\%$ & $5.0\%$ & $65.0\%$ & $7.5\%$ \\
FC vs. NoComp & $25.0\%$ & $3.3\%$ & $63.5\%$ & $8.3\%$ \\
FU vs. NoComp & $31.5\%$ & $2.8\%$ & $54.5\%$ & $11.3\%$ \\
FU vs. FC & $56.3\%$ & $5.3\%$ & $31.5\%$ & $7.0\%$ \\
\bottomrule
\end{tabular}
\end{table}

The results demonstrate that Clapping consistently outperforms direct compression in generated-output quality, with Clapping-FU achieving the strongest performance among the communication-compressed methods. In particular, Clapping-FU achieves a win rate of 59.8\% against direct compression, compared with a loss rate of 29.0\%, highlighting its effectiveness in preserving model quality under aggressive communication compression.

\subsubsection{Fine-tuning with multiple compression.}
Here we present the experimental result in fine-tuning tasks with multiple compression. Different from the fine-tuning tasks introduced in Appendix \ref{appendix: finetune glue} and \ref{appendix: finetune wiki}, the experimental setup we take here is more \textbf{STRICT} than the practical fine-tuning tasks to present a comparsion for different algorithms. 

\textbf{Fine-tuning GPT-2 on Wikitext. }
We fine-tuned a pretrained GPT-2 \citep{radford2019language} model on Wikitext dataset (wikitext-2-raw-v1 version) \citep{merity2016pointer} on a NVIDIA A100 80G GPU with communication compression algorithms for pipeline-parallel learning including directly compression, AQ-SGD, \ours without lazy sampling and \oursfclast. As the experiment is taken in a single GPU, we simulate the communication compression by adding the corresponding error to the activation and gradient. The block size was set to 1024. We used batch size 8, and we fine-tuned the model for 8 epochs for the competitive algorithms and we use the same iterations for \ours, respectively. We use the AdamW optimizer \citep{loshchilov2017decoupled} and FP16 for fine-tuning. The learning rate is initialized from $2\times 10^{-5}$ and scheduled by a linear scheduler. We compare different compression algorithms under three compression strategies. These strategies integrated different compressor including TopK \citep{wangni2018gradient}, direct quantization \citep{alistarh2017qsgd}, and Natural compression \citep{horvoth2022natural} on GPT-2 with different model size, which can be summarized as follow.
\begin{itemize}[leftmargin=18.8mm]
    \item[\textbf{Strategy S1}] Basic model: GPT-2 small;

    Compression position: At the end of layer 2, 5, 8;

    Forward compressor: Top 40\%;

    Backward compressor: Top 40\%.
    \item[\textbf{Strategy S2}] Basic model: GPT-2 medium;

    Compression position: At the end of layer 4, 10, 16;

    Forward compressor: Direct taking Natural Compression;

    Backward compressor: Quantizing to 8 bits, then taking Natural Compression.
    \item[\textbf{Strategy S3}] Basic model: GPT-2 medium;

    Compression position: At the end of layer 4, 10, 16;

    Compressor in layer 4: Direct taking Natural Compression;

    Compressor in layer 10 and 16: Quantizing to 8 bits, then taking Natural Compression.
\end{itemize}

For each compression strategies, we independently repeat \ours with different lazy sampling coefficient $p_t$ including $\{0.1,0.2,0.3,0.4,0.6,0.8\}$ for 3 times and select the $p_t$ with the best average evaluation accuracy.

\begin{figure}[t!]
\centering
	\subfigure{\centering
        \includegraphics[width=0.3\textwidth]{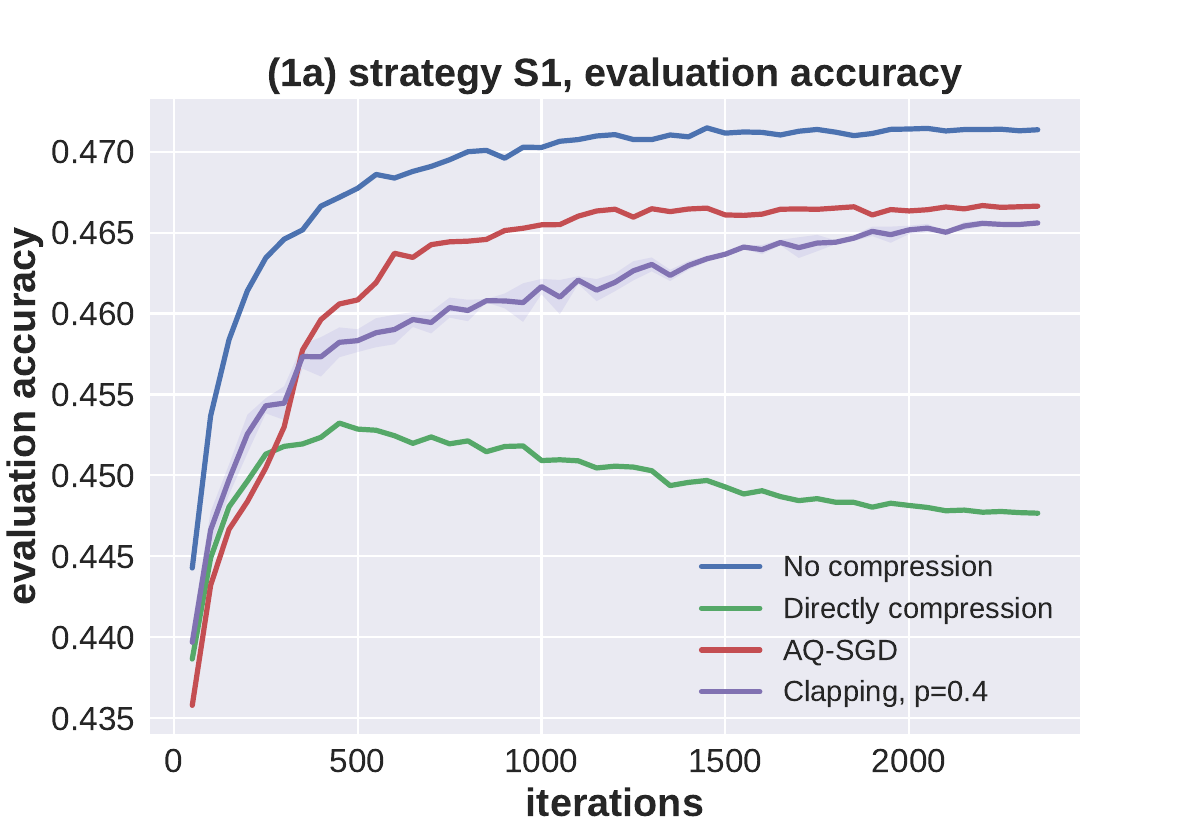}}
  \hspace{-15pt}
	\subfigure{\centering
		\includegraphics[width=0.3\textwidth]{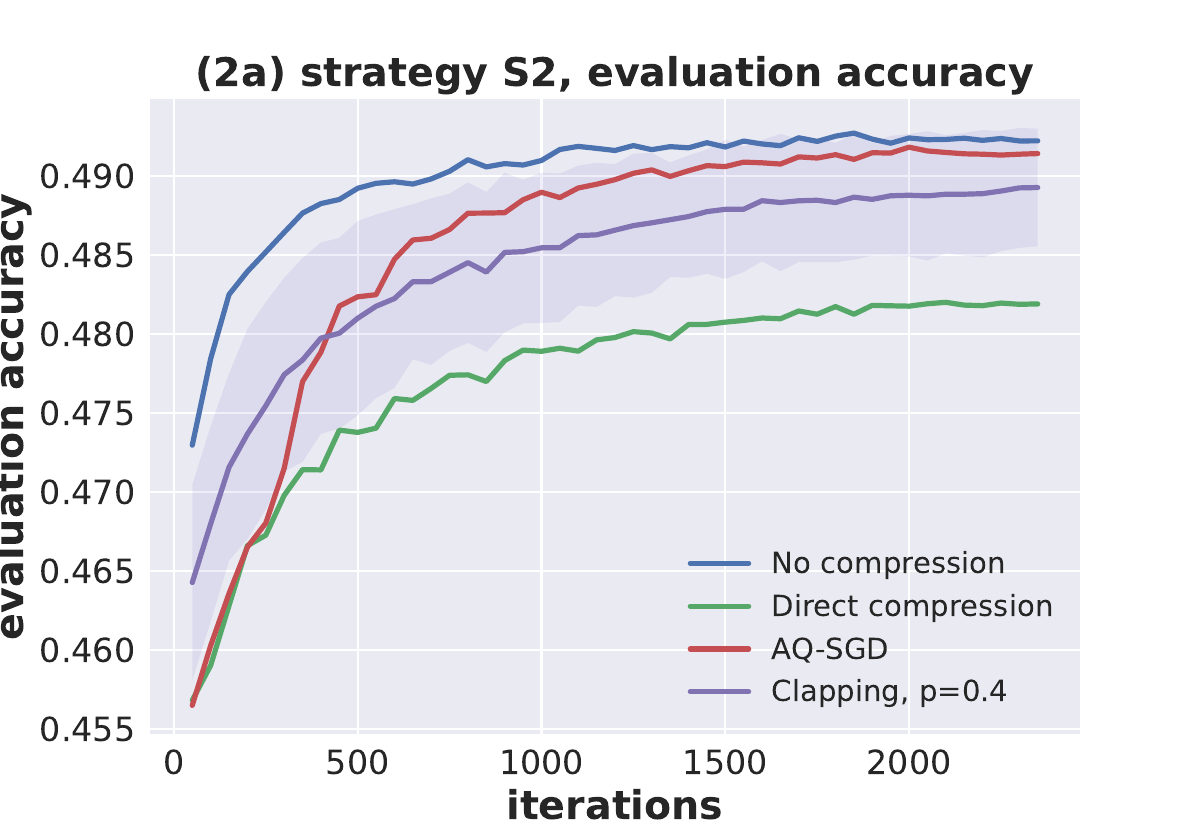}}
  \hspace{-15pt}
	\subfigure{\centering
		\includegraphics[width=0.3\textwidth]{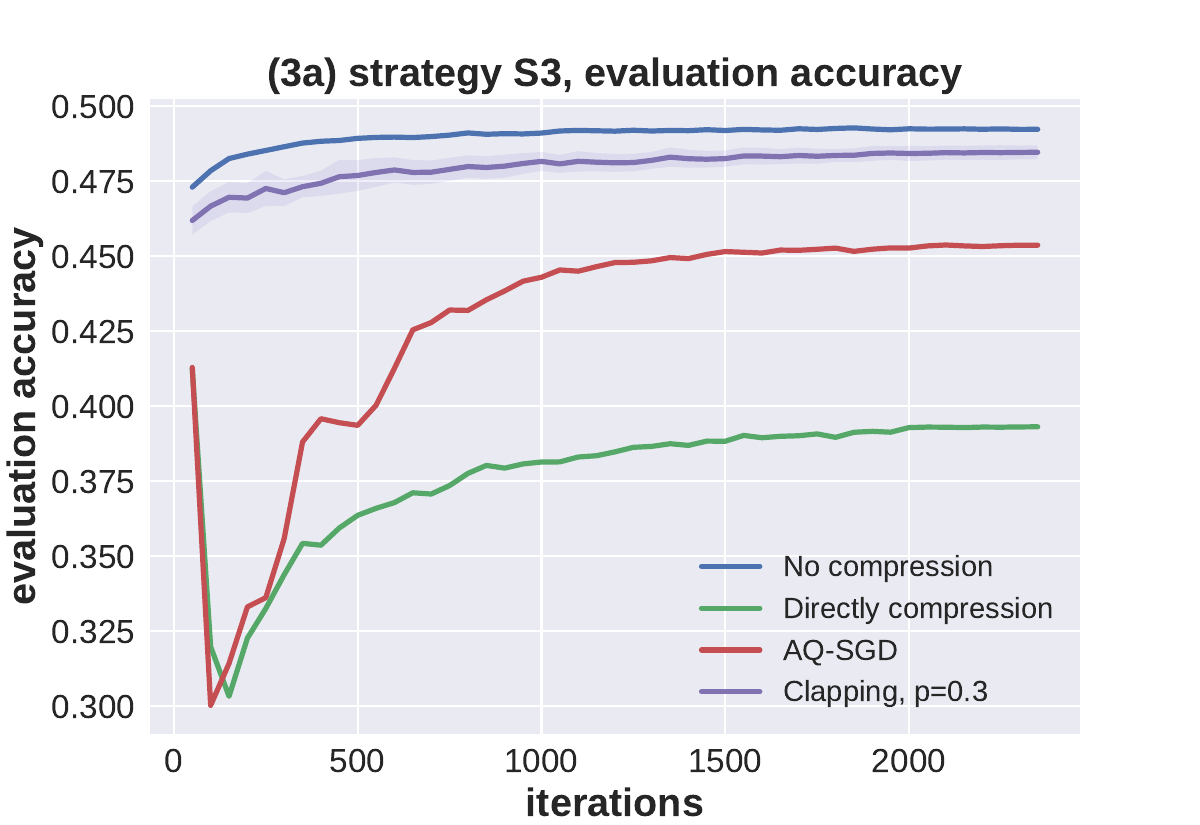}}
    \qquad
	\subfigure{\centering
        \includegraphics[width=0.3\textwidth]{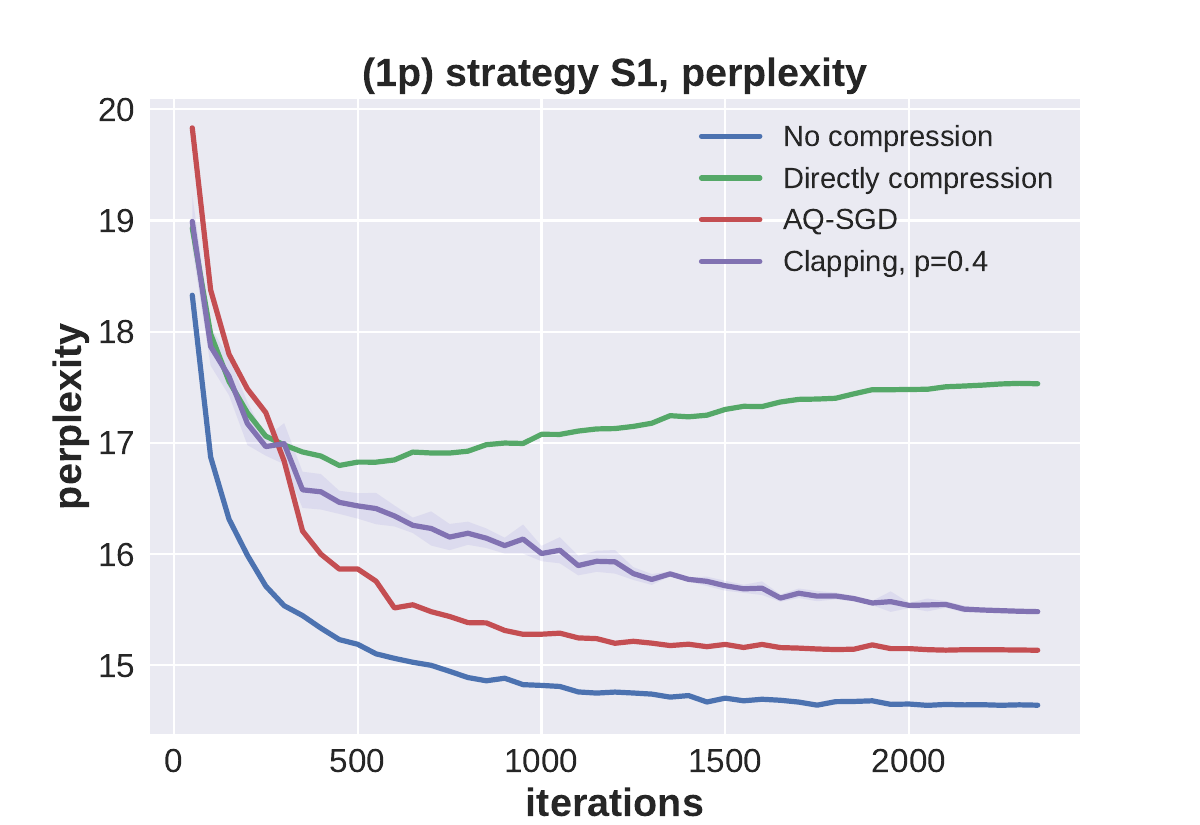}}
  \hspace{-15pt}
	\subfigure{\centering
		\includegraphics[width=0.3\textwidth]{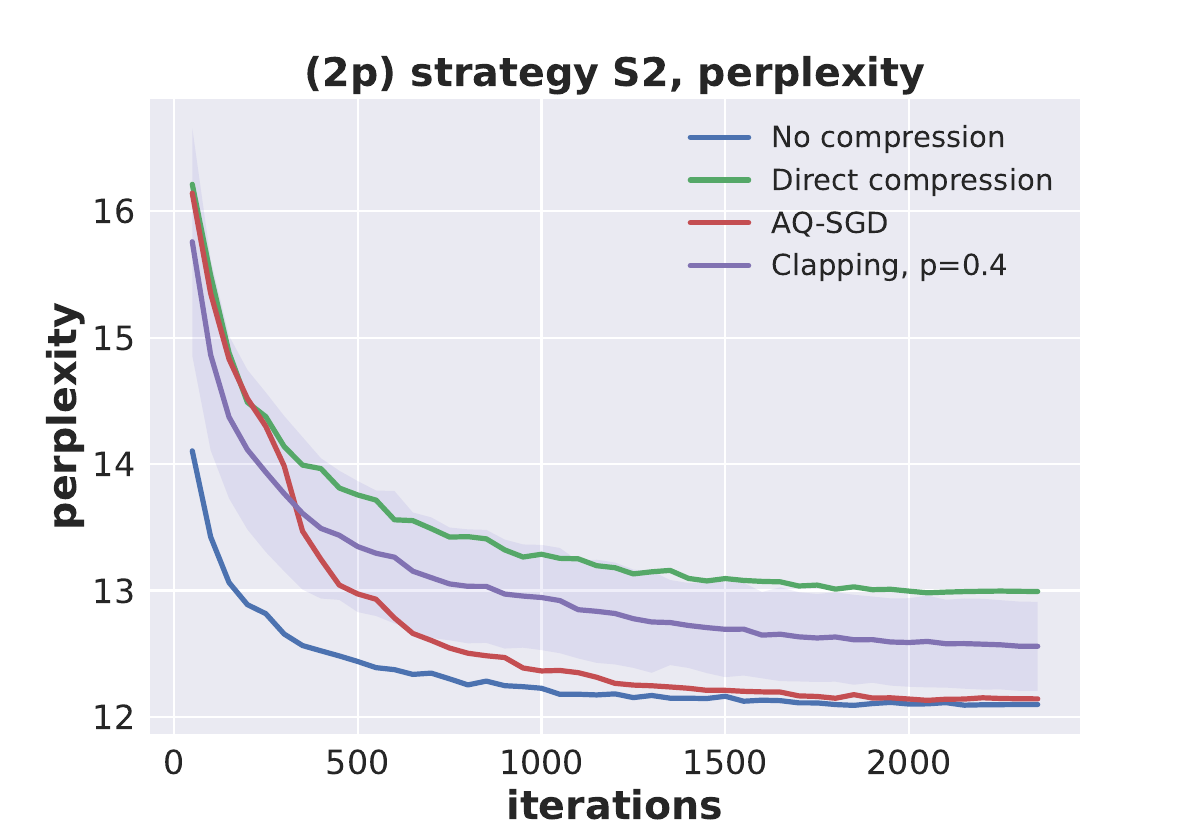}}
  \hspace{-15pt}
	\subfigure{\centering
		\includegraphics[width=0.3\textwidth]{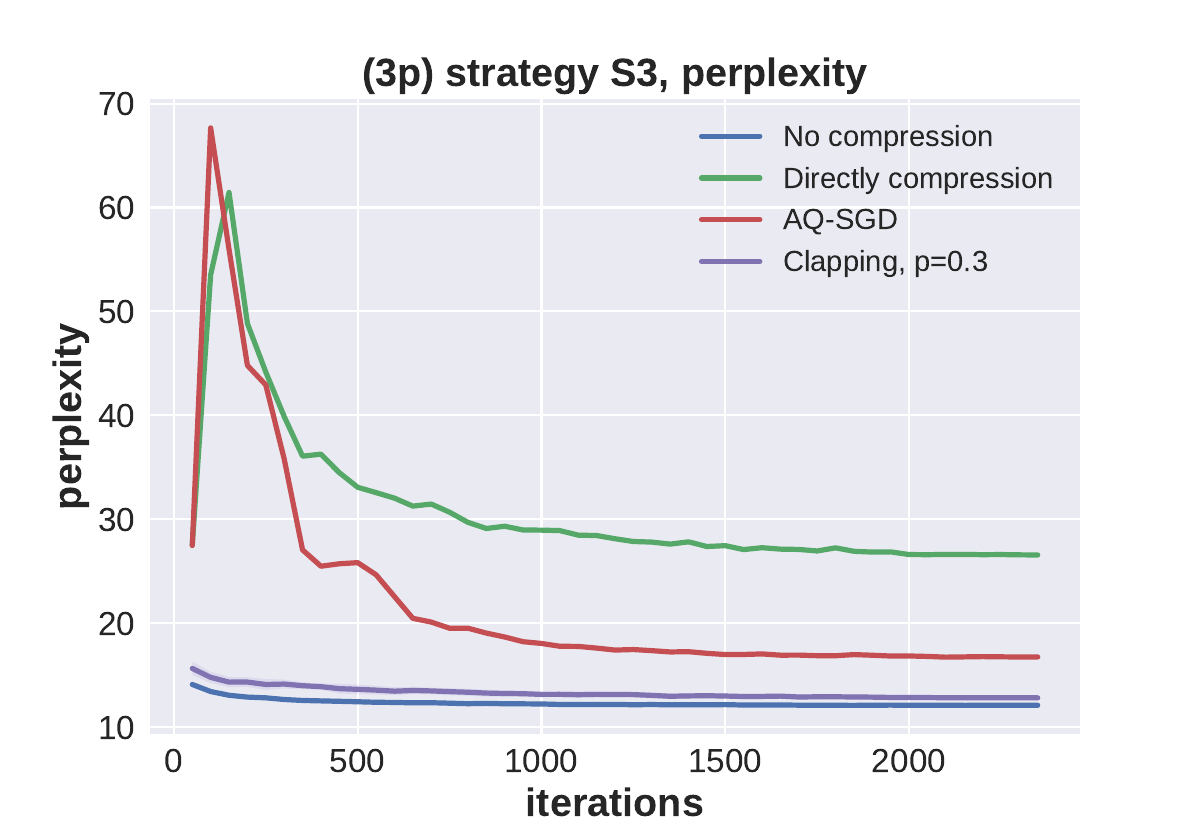}}
\caption{\small The evaluation accuracy and perplexity of GPT-2 fine-tuning with different compression strategies and difference compression algorithms. (Left: Strategy S1, Middle: Strategy S2, Right: Strategy S3.)}
\label{fig: gpt2 tuning-appendix}
\end{figure}

Figure \ref{fig: gpt2 tuning-appendix} illustrates the evaluation accuracy and perplexity of \ours with best $p_t$ as well as the other algorithms. It can be observed that \ours outperforms direct compression. Moreover, lazy sampling can improve the evaluation accuracy. In the case, \ours achieves a better accuracy than AQ-SGD because AQ-SGD suffers from a high compression error in the beginning.\

\begin{table}[t]
\caption{Evaluation perplexity of \ours in GPT-2 fine-tuning with different lazy coefficient $p_t$.}
\setlength{\tabcolsep}{9pt}
\renewcommand{\arraystretch}{1.1}
\label{table:app_different_p_t}
\centering
\begin{threeparttable}
\begin{tabular}{cccccccc}
\toprule
 \multirow{2}{*}{Strategy} & \multicolumn{7}{c}{$p_t$}                   \\ \cline{2-8} 
                           & \vspace{-1pt}0.1 & 0.2 & 0.3 & 0.4 & 0.6 & 0.8 & 1.0 \\ \midrule
 S1                        & \multicolumn{1}{c}{\begin{tabular}[c]{@{}c@{}} 16.118 \\ (0.095) \end{tabular}}   &   \multicolumn{1}{c}{\begin{tabular}[c]{@{}c@{}} 15.597 \\ (0.070) \end{tabular}}   &  \multicolumn{1}{c}{\begin{tabular}[c]{@{}c@{}} 15.523 \\ (0.013) \end{tabular}}    &  \multicolumn{1}{c}{\begin{tabular}[c]{@{}c@{}} \textbf{15.482} \\ \textbf{(0.014)} \end{tabular}}    &  \multicolumn{1}{c}{\begin{tabular}[c]{@{}c@{}} 15.675 \\ (0.015) \end{tabular}}    &  \multicolumn{1}{c}{\begin{tabular}[c]{@{}c@{}} 15.999 \\ (0.014) \end{tabular}}    &  \multicolumn{1}{c}{\begin{tabular}[c]{@{}c@{}} 16.582 \\ (0.417) \end{tabular}}  \\
 S2                        & \multicolumn{1}{c}{\begin{tabular}[c]{@{}c@{}} 13.432 \\ (0.118) \end{tabular}}   &   \multicolumn{1}{c}{\begin{tabular}[c]{@{}c@{}} 12.992 \\ (0.239) \end{tabular}}   &  \multicolumn{1}{c}{\begin{tabular}[c]{@{}c@{}} 12.713 \\ (0.343) \end{tabular}}    &  \multicolumn{1}{c}{\begin{tabular}[c]{@{}c@{}} \textbf{12.558} \\ \textbf{(0.352)} \end{tabular}}    &  \multicolumn{1}{c}{\begin{tabular}[c]{@{}c@{}} 12.630 \\ (0.306) \end{tabular}}    &  \multicolumn{1}{c}{\begin{tabular}[c]{@{}c@{}} 12.645 \\ (0.250) \end{tabular}}    &  \multicolumn{1}{c}{\begin{tabular}[c]{@{}c@{}} 12.756 \\ (0.182) \end{tabular}}  \\
S3                        & \multicolumn{1}{c}{\begin{tabular}[c]{@{}c@{}} 15.337 \\ (0.187) \end{tabular}}   &   \multicolumn{1}{c}{\begin{tabular}[c]{@{}c@{}} 14.330 \\ (0.372) \end{tabular}}   &  \multicolumn{1}{c}{\begin{tabular}[c]{@{}c@{}} \textbf{12.815} \\ \textbf{(0.222)} \end{tabular}}    &  \multicolumn{1}{c}{\begin{tabular}[c]{@{}c@{}} 14.106 \\ (0.208) \end{tabular}}    &  \multicolumn{1}{c}{\begin{tabular}[c]{@{}c@{}} 14.261 \\ (0.218) \end{tabular}}    &  \multicolumn{1}{c}{\begin{tabular}[c]{@{}c@{}} 14.466 \\ (0.231) \end{tabular}}    &  \multicolumn{1}{c}{\begin{tabular}[c]{@{}c@{}} 14.712 \\ (0.246) \end{tabular}}  \\ \bottomrule
\end{tabular}
\end{threeparttable}
\end{table}

And we also present the best evaluation accuracy and perplexity of \ours, which is shown in Table \ref{table:app_different_p_t}. From Table \ref{table:app_different_p_t}, we can observe that a suitable $p_t$ like 0.3 or 0.4 can benefit the convergence and generalization in the fine-tuning tasks. And a large $p_t$ may harm the generalization, which matches our discussion in Section \ref{sec:convergence}.

\textbf{Fine-tuning GPT-2 on arXiv abstracts. }
\label{sec:app_gpt_arxiv}
We also respectively fine-tuned pretrained GPT-2 small models and GPT-2 medium models \citep{radford2019language} on a dataset with 30K arXiv abstracts \citep{wang2022fine} on 8 NVIDIA RTX 4090 GPUs with communication compression algorithms for pipeline-parallel learning including directly compression, \oursfc without lazy sampling and \oursfc with lazy sampling. We also compared those approaches to the case with no compression. The block size was set to 1024. We set data parallel degree 8 with the macro-batch size 16 and the micro-batch size 2. As each micro-batch is computed in a single GPU, we simulate the communication compression by adding the corresponding error to the activation and gradient. And we fine-tuned the model for 8 epochs with other competitive algorithms and the same iterations with \ours, respectively. We simulated the communication compression for each GPU by adding the compression hook with Top-K \citep{wangni2018gradient} compressor that keep 50\% of the elements. For GPT-2 small, we add the hook at the end of the 2-nd, 5-th, and 8-th transformer layer. For GPT-2 medium, we add the hook at the end of the 4-th, 10-th, and 16-th transformer layer. 

\begin{figure}[t!]
\centering
	\subfigure{\centering
        \includegraphics[width=0.45\textwidth]{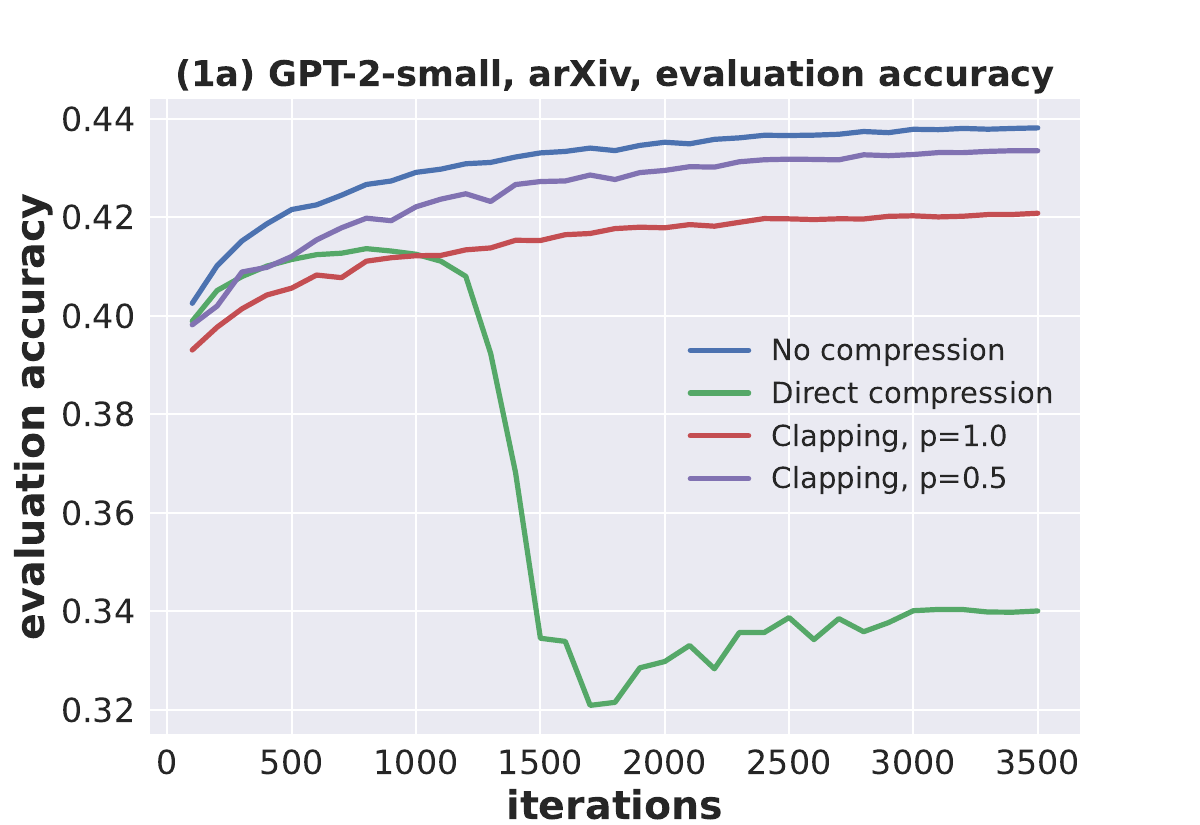}}
  \hspace{-20pt}
	\subfigure{\centering
		\includegraphics[width=0.45\textwidth]{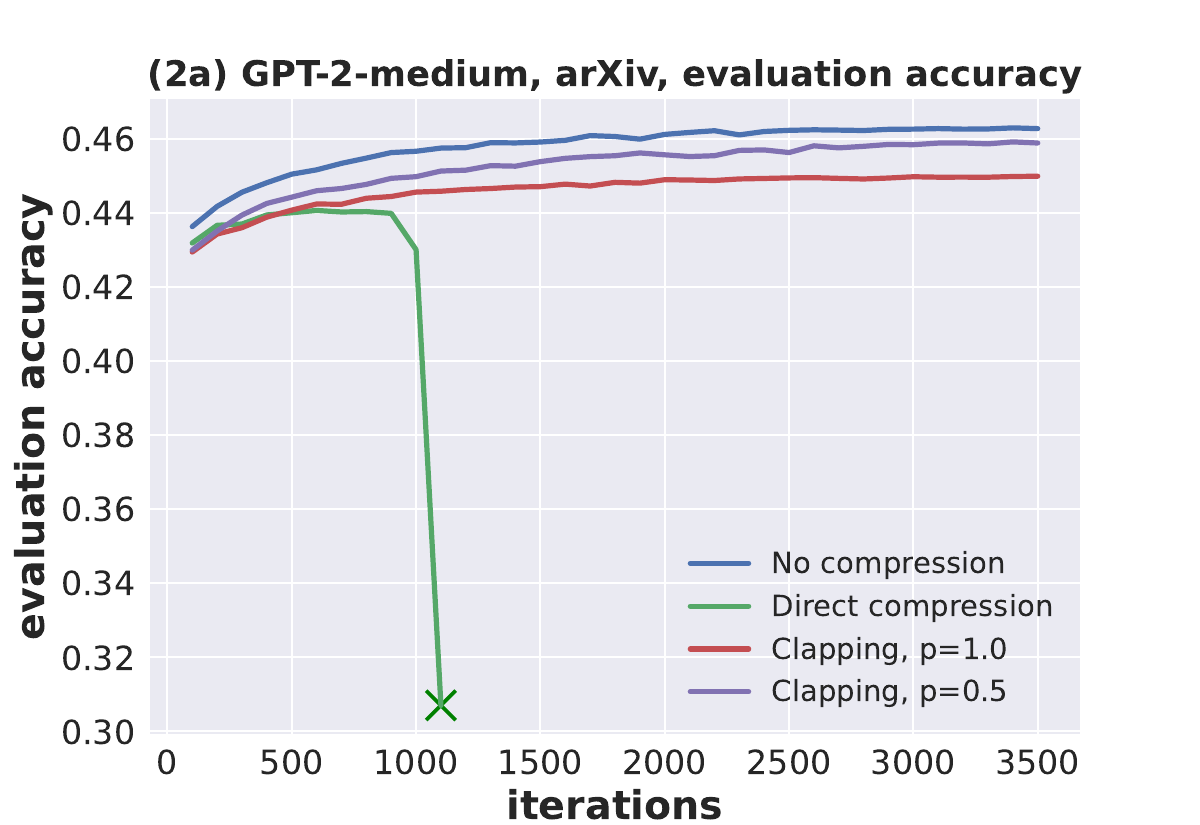}}
    \qquad
	\subfigure{\centering
        \includegraphics[width=0.45\textwidth]{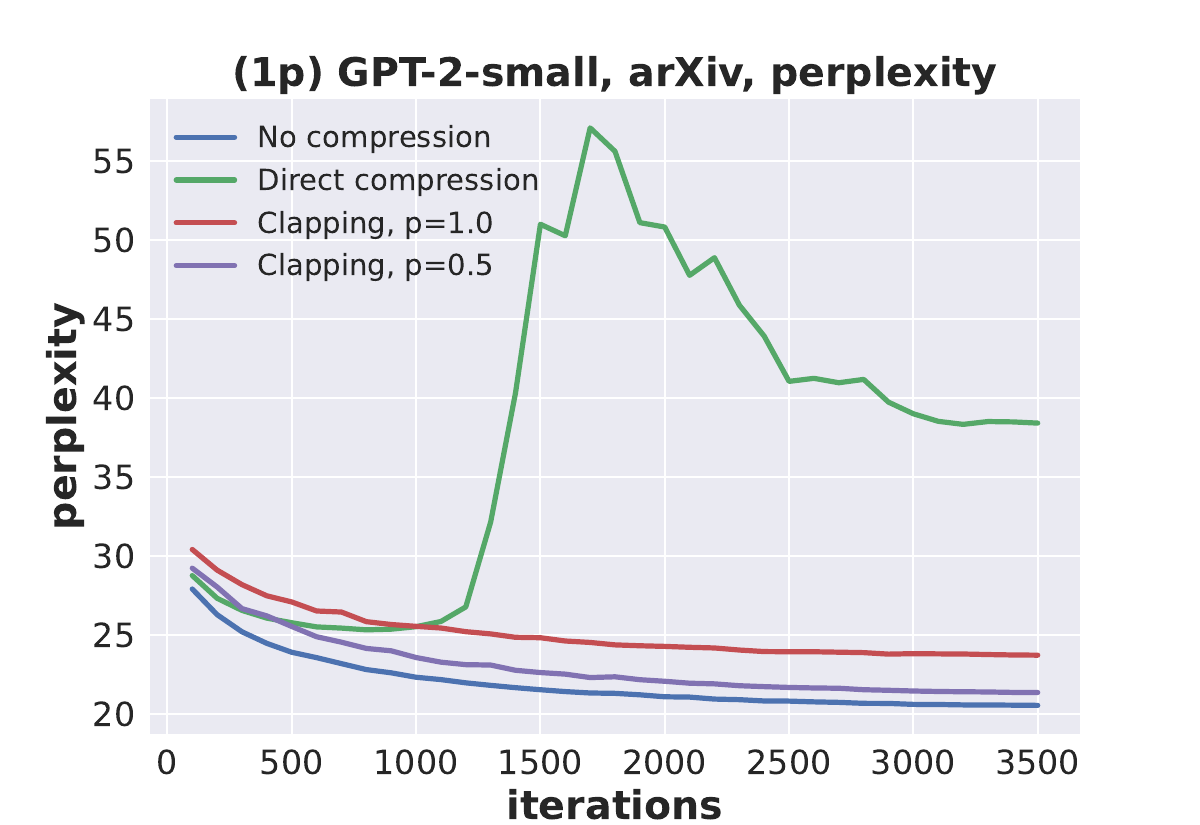}}
  \hspace{-20pt}
	\subfigure{\centering
		\includegraphics[width=0.45\textwidth]{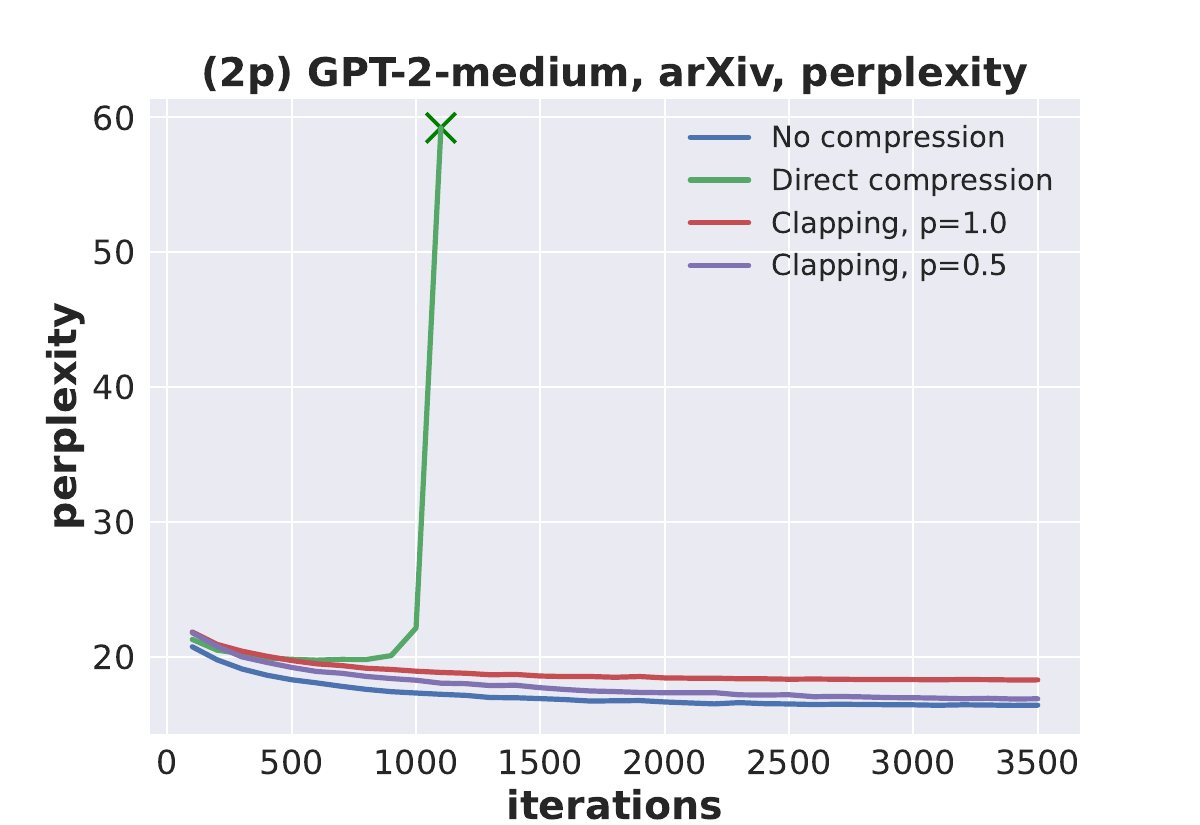}}
\caption{\small The evaluation accuracy and perplexity of GPT-2 small and GPT-2 medium fine-tuning by arXiv with Top50\%. (Left: GPT-2 small, Right: GPT-2 medium.)}
\label{fig: gpt2 tuning arxiv-appendix}
\end{figure}

We use the AdamW optimizer \citep{loshchilov2017decoupled} and FP16 for fine-tuning. The learning rate is initialized from $5\times 10^{-5}$ and scheduled by a cosine scheduler. The lazy sampling coefficient $p_t$ for \ours with lazy sampling was set to 0.5.  Figure \ref{fig: gpt2 tuning arxiv-appendix} illustrates the evaluation accuracy and perplexity. It can be observed that both \ours with lazy sampling and \ours without lazy sampling outperforms direct compression. Moreover, lazy sampling can make the evaluation accuracy and perplexity comparable to those of non-compressed case. Specifically, we can find that direct compression in GPT-2 medium suffers from the non-convergence but \ours not, which can illustrates the benefit of error feedback technique. 

\textbf{Fine-tuning LLaMA2-7B on Wikitext. }
We fine-tuned a pre-trained LLaMA2-7B model \citep{touvron2023LLaMA} on Wikitext dataset (wikitext-2-raw-v1 version) \citep{merity2016pointer} on 4 NVIDIA A100 40G GPUs with communication compression algorithms for pipeline-parallel learning including directly compression, \oursfc without lazy sampling and \oursfc with $p=0.3,0.5$. We also compared those approaches to the case with no compression. We set data parallel degree 4 with the macro-batch size 8 and the micro-batch size 2. As each micro-batch is computed in a single GPU, we simulate the communication compression by adding the corresponding error to the activation and gradient. We simulated the communication compression for each GPU by adding the compression hook at the end of 8-th, 16-th, and 24-th transformer layers. The block size was set to 1024. We fine-tuned the model for 6 epoch with the other competitive algorithms and the same iterations with \ours, respectively. We used the SGD optimizer with momentum 0.9 and FP16 for fine-tuning. The learning rate is $5\times 10^{-5}$. And we used TopK \citep{wangni2018gradient} that keeps $50\%$ of elements and natural compression \citep{horvoth2022natural} as the compressor, respectively. According to \citep{horvoth2022natural}, FP16 training with natural compression can compression the activations and gradients to 6-bit, thus the communication overhead is 37.5\% compared to the non-compressed scenario.
\begin{figure}[t!]
\centering
	\subfigure{\centering
        \includegraphics[width=0.45\textwidth]{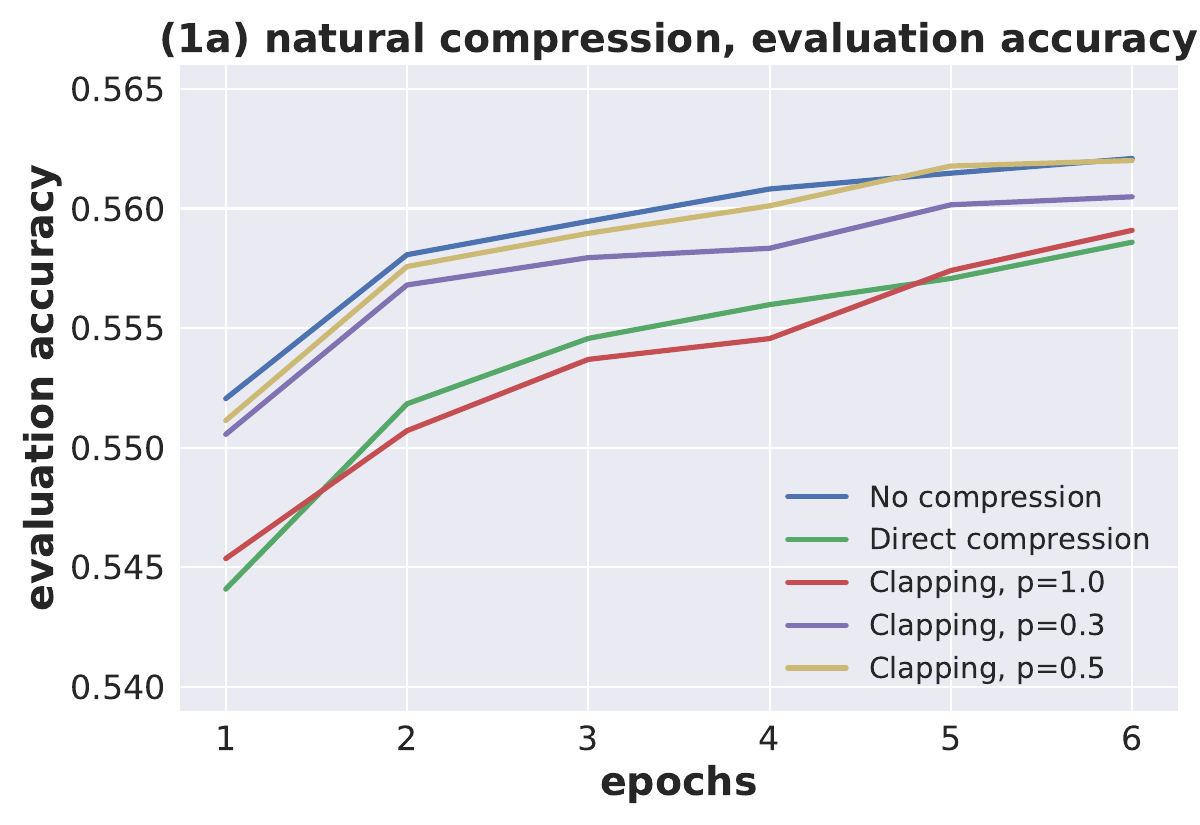}}
  \hspace{-5pt}
	\subfigure{\centering
		\includegraphics[width=0.45\textwidth]{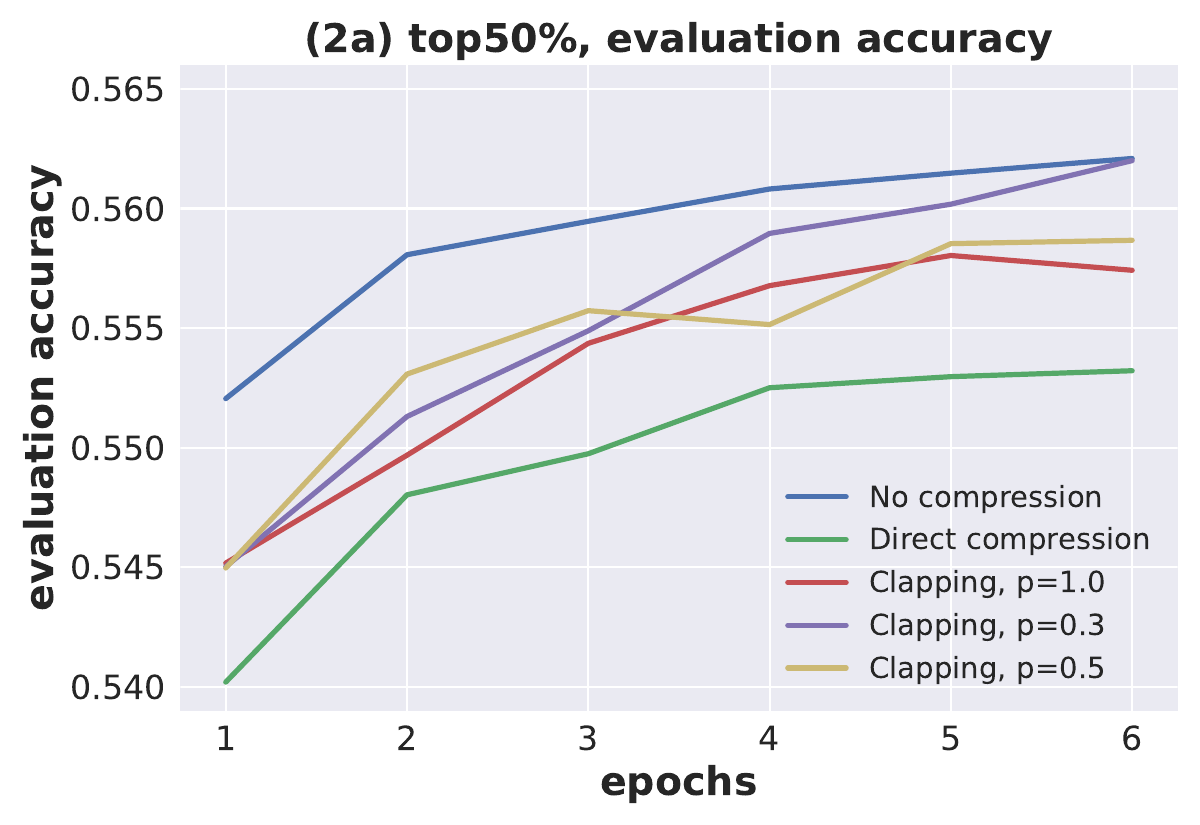}}
  \qquad
	\subfigure{\centering
        \includegraphics[width=0.45\textwidth]{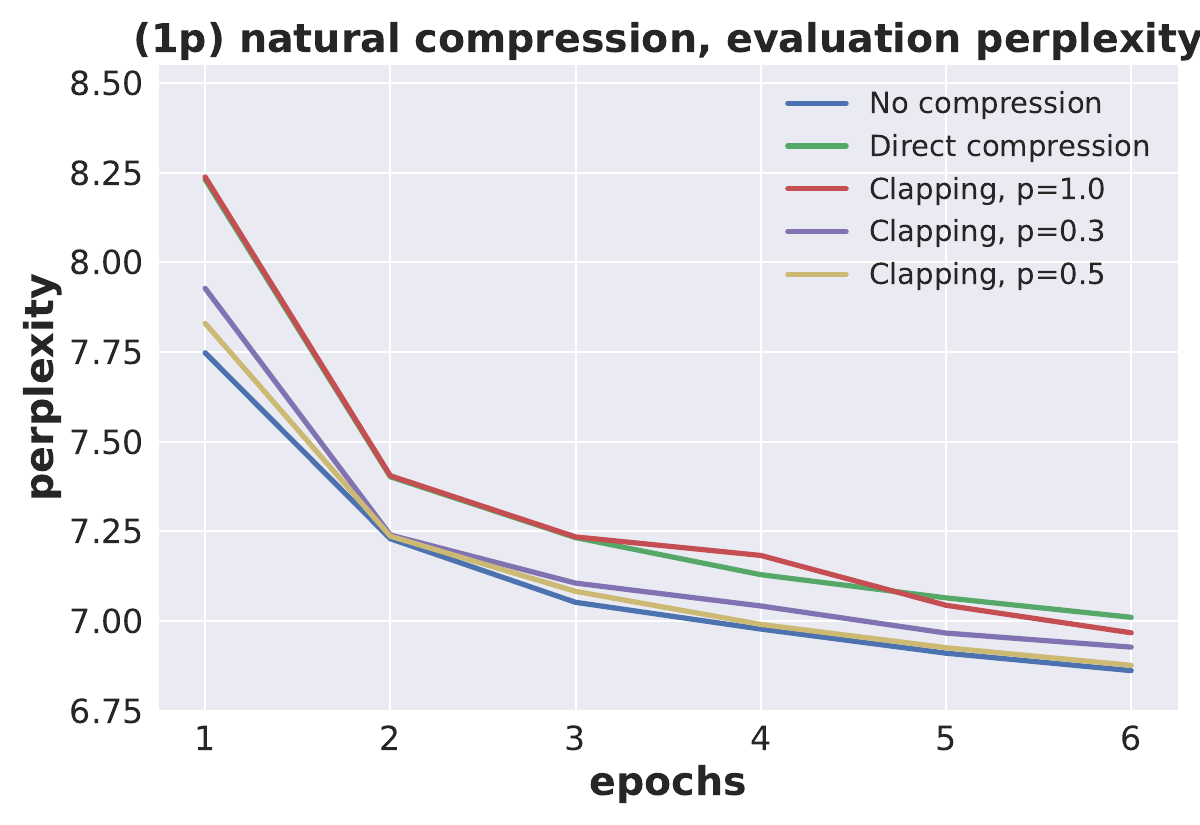}}
  \hspace{-5pt}
	\subfigure{\centering
		\includegraphics[width=0.45\textwidth]{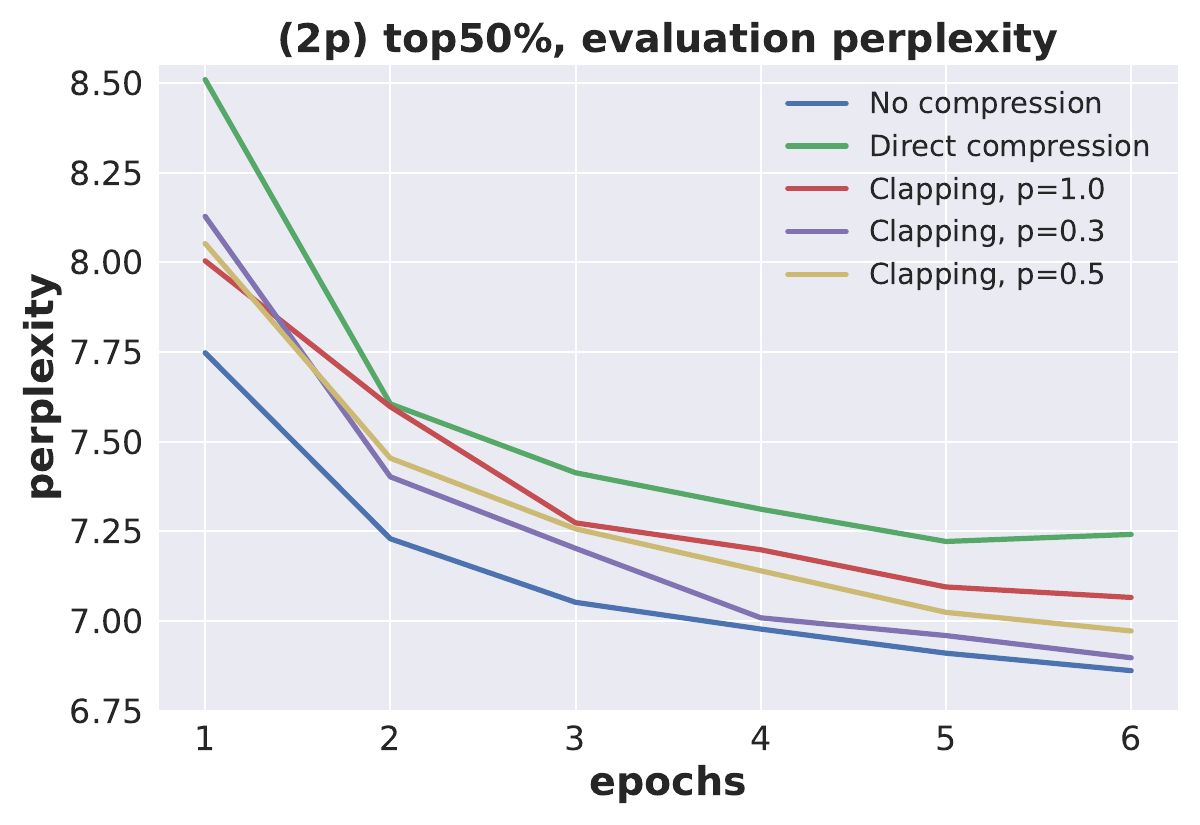}}
  \vspace{-10pt}

\caption{\small The evaluation accuracy and perplexity of LLaMA-2 fine-tuning with Wikitext-2 with different compressors and different compression algorithms. (Left: natural compression. Right: Top50\%)}
\label{fig: LLaMA2 finetuning-appendix}
\end{figure}

\begin{figure}[t!]
\centering
\vspace{-15pt}
	\subfigure{
    \includegraphics[width=0.45\textwidth]{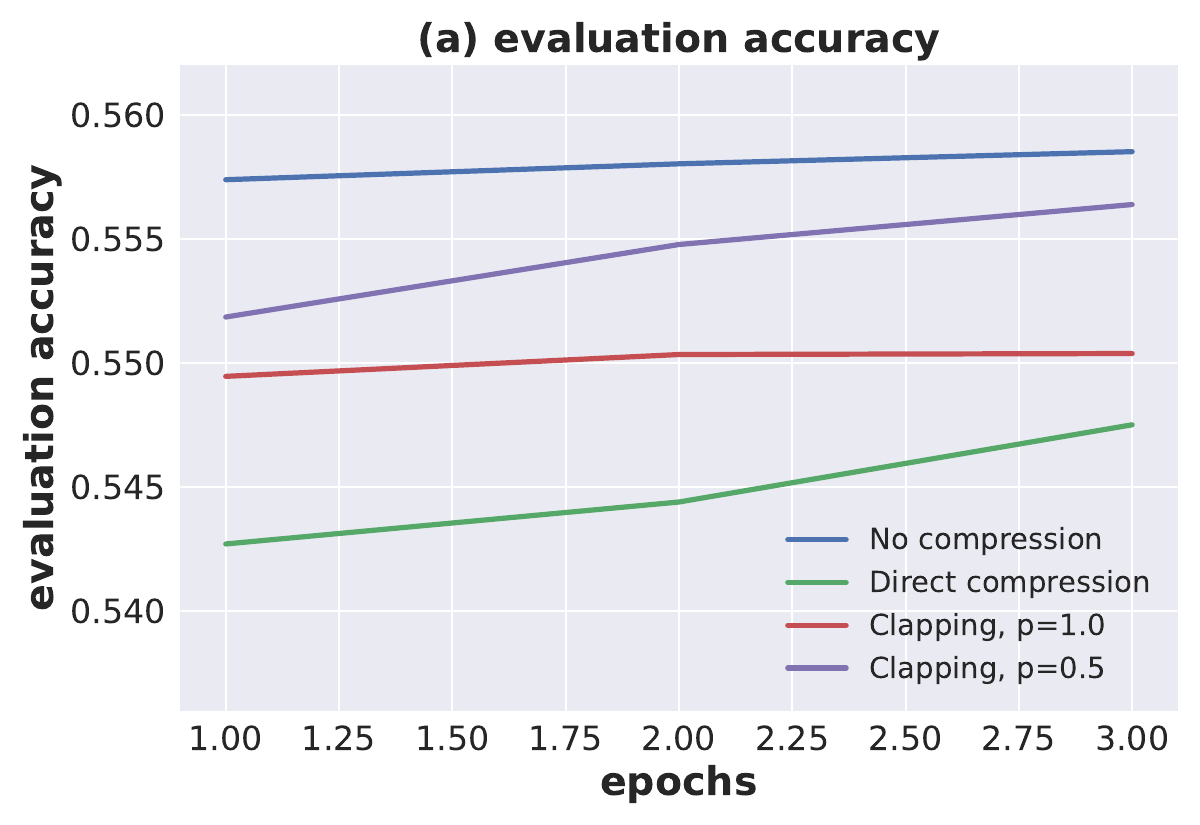}}\hspace{-5pt}
	\subfigure{
    \includegraphics[width=0.45\textwidth]{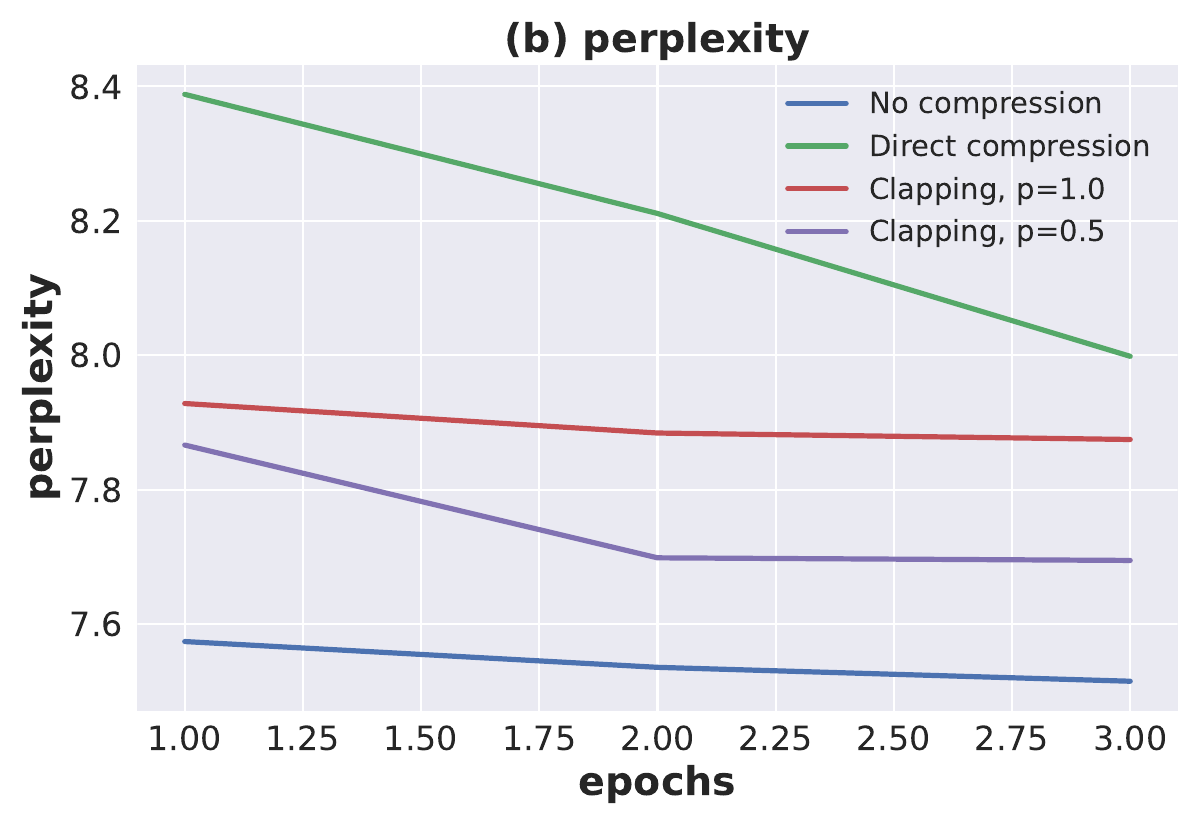}}
  \vspace{-10pt}
\hspace{-15pt}
\caption{The evaluation accuracy of various approaches in LLaMA-2 fine-tuning task by arXiv abstracts.}
\vspace{-8pt}
\label{fig: LLaMA_arxiv}
\end{figure}

Figure \ref{fig: LLaMA2 finetuning-appendix} shows the evaluation accuracy and perplexity of different algorithms in LLaMA-2 fine-tuning tasks. It can be observed that \ours outperforms directly compression in the evaluation scalability. It is also noteworthy that for both compressor, \ours can achieve nearly the \textbf{SAME} evaluation accuracy and perplexity as the non-compressed case. Thus one can tune a suitable $p_t$ to let \ours totally eliminate the negative impact of communication compression with more than \texttt{$2\times$} communication saving. Figure \ref{fig: LLaMA_arxiv} illustrate the evaluation accuracy and perplexity of different compression algorithms. It can be observed the benefit of error feedback technique and lazy sampling strategy.

\begin{figure}[t]
	\centering
	\begin{minipage}[c]{0.98\textwidth}
    \centering
    	\subfigure{
    		\includegraphics[width=0.45\textwidth]{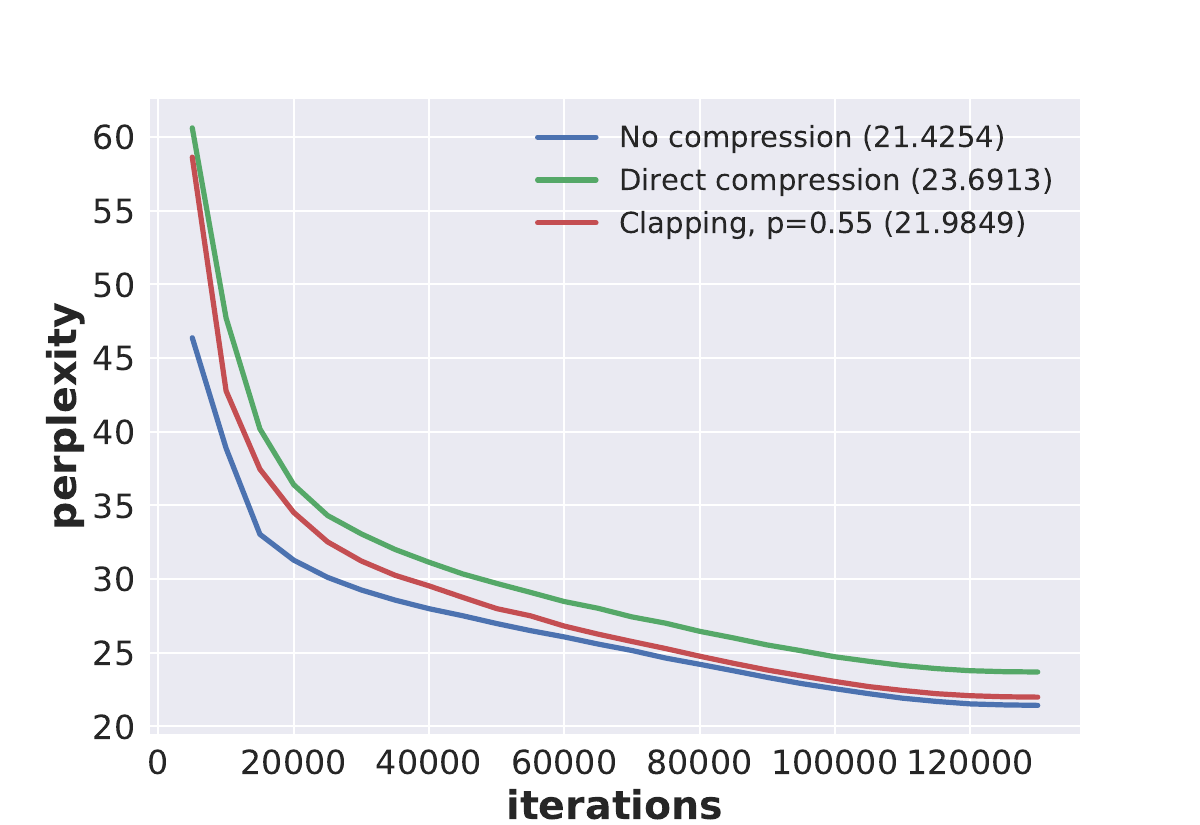}}
    	\subfigure{
    		\includegraphics[width=0.45\textwidth]{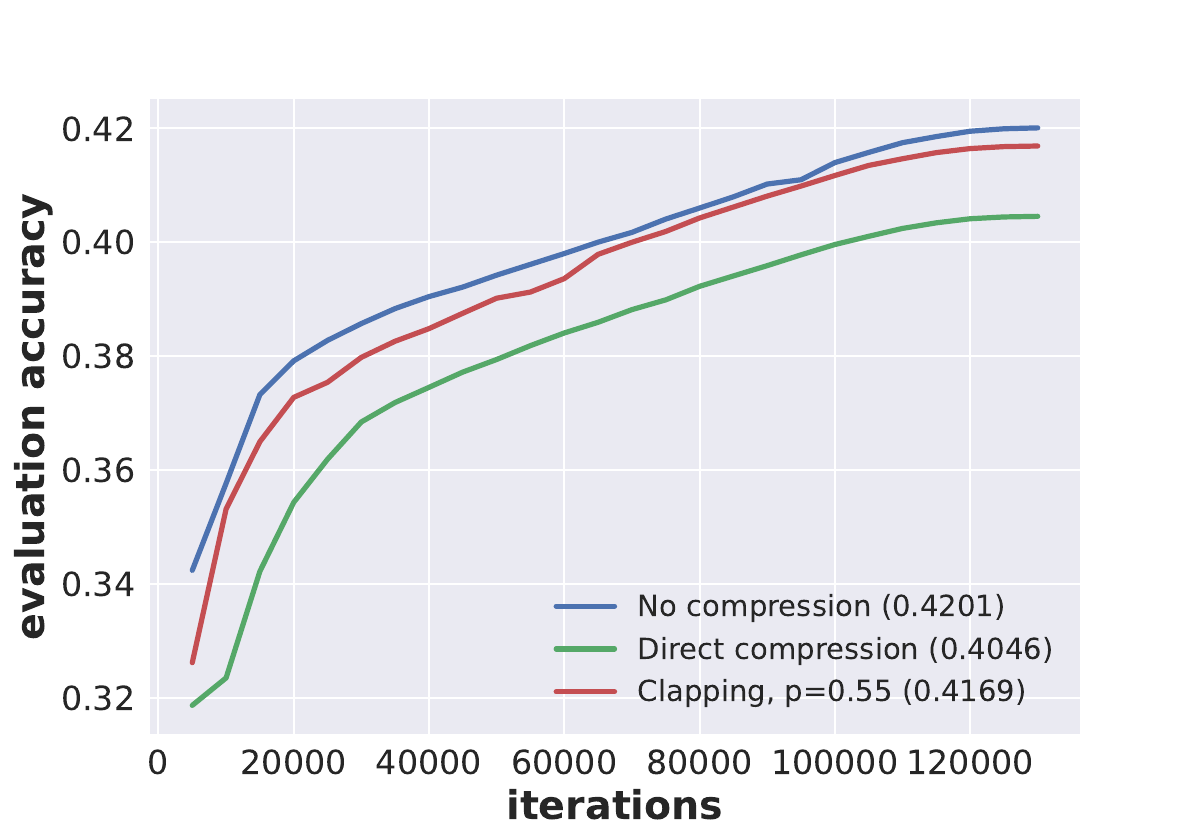}}
    \caption{\small The evaluation perplexity (left) and accuracy (right) with various communication compression algorithms for pre-training GPT-2.}
    \label{fig: gpt2_training-main}
	\end{minipage}  
\end{figure}

\begin{figure}[t]
\vspace{-10pt}
	\centering
	\begin{minipage}[c]{0.98\textwidth}
    \centering
      \hspace{-20pt}
    \vspace{-5pt}
    	\subfigure{
    		\includegraphics[width=0.72\textwidth]{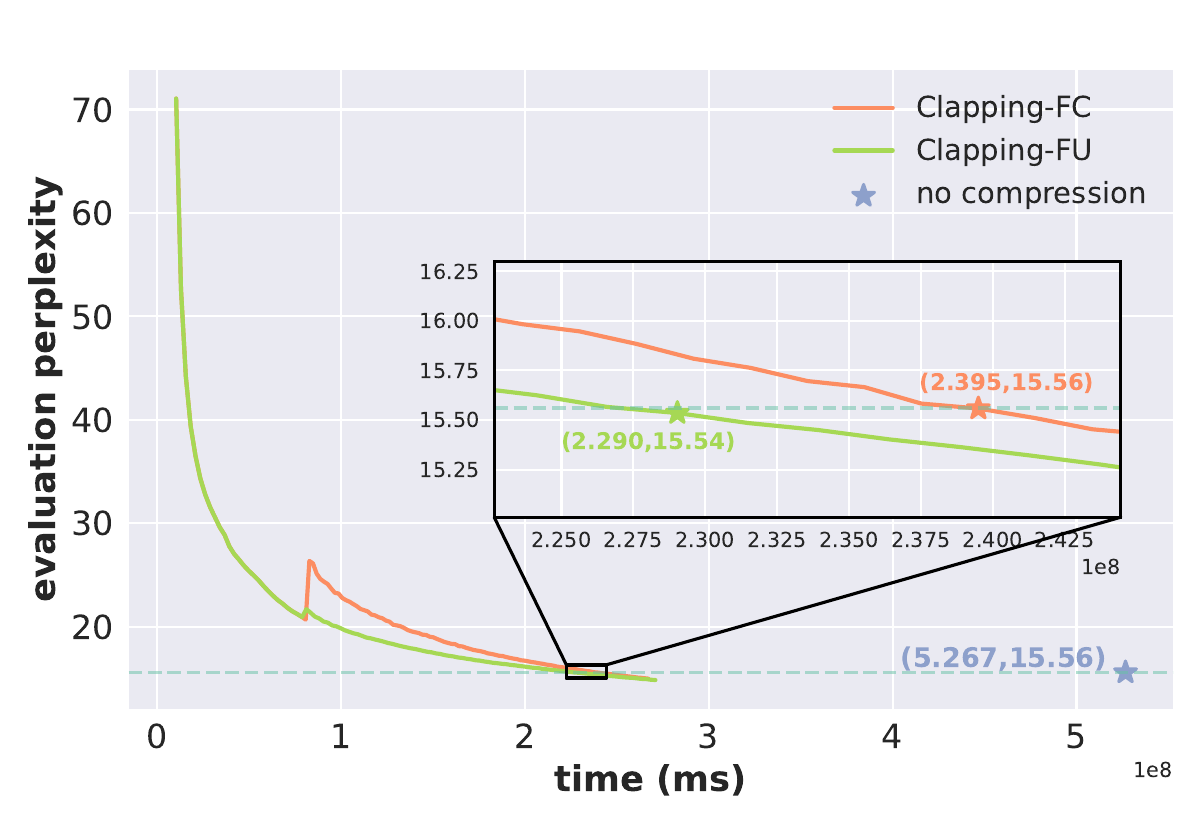}}
    \caption{\small The evaluation perplexity with various communication compression algorithms for pre-training LLaMA-2 1B.}
    \label{fig: LLaMA2 training-appendix2}
	\end{minipage} 
\end{figure}

\begin{table}[t]
\centering
\caption{Single-iteration time (ms) and speedup of Clapping-FC and non-compression across different network bandwidths}
\label{tab:speedup}
\begin{threeparttable}
\begin{tabular}{lccccc}
\toprule
\textbf{Method} & \textbf{100 MB/s} & \textbf{200 MB/s} & \textbf{300 MB/s} & \textbf{400 MB/s} & \textbf{500 MB/s} \\
\midrule
Clapping-FC     & 572.68           & 460.18           & 422.68           & 403.93           & 392.68          \\
No Compression  & 947.68           & 647.61           & 547.28           & 497.63           & 467.68          \\
\cmidrule(r){1-6}
Speedup & 1.65$\times$             & 1.41$\times$             & 1.30$\times$             & 1.23$\times$             & 1.19            \\
\bottomrule
\end{tabular}
\end{threeparttable}
\end{table}

\begin{table}[t]
\centering
\caption{Comparison of training time ($\times 10^5$ s) to reach target evaluation perplexity for LLaMA-2 3B pre-training under bandwidth constraint}
\begin{threeparttable}
\begin{tabular}{cccc}
\toprule
\textbf{Target evaluation perplexity} & \textbf{17.5} & \textbf{17.0} \\
\midrule
No compression & 2.34036 & 2.69330 \\
Clapping-FC & 2.08600 & 2.42791 \\
\bottomrule
\end{tabular}
\end{threeparttable}
\label{tab:training_time_comparison}
\end{table}

\begin{table}[t]
\centering
\caption{Comparison of evaluation loss across different compression methods}
\label{tab:compressor_comparison}
\begin{threeparttable}
\begin{tabular}{cccc}
\toprule
  \textbf{Top-50\%} & \textbf{8-bit Quantization} & \textbf{Natural Compression} & \textbf{Hybrid Compressor} \\
\midrule
 5.0260 & 3.1425 & 3.0990 & 2.8944 \\
\bottomrule
\end{tabular}
\end{threeparttable}
\end{table}

\textbf{Fine-tuning LLaMA2-7B on arXiv abstracts. }
We fine-tuned a pre-trained LLaMA2-7B model \citep{touvron2023LLaMA} on a dataset with 30K arXiv abstracts \citep{merity2016pointer} on 4 NVIDIA A100 40G GPUs with communication compression algorithms for pipeline-parallel learning including directly compression, \oursfc without lazy sampling and \oursfc with $p=0.5$. We also compared those approaches to the case with no compression. We set data parallel degree 4 with the macro-batch size 2 and the micro-batch size 8. We simulated the communication compression for each GPU by adding the corresponding error to the activation and gradient at the end of 8-th, 16-th, and 24-th transformer layers. The block size was set to 1024. We fine-tuned the model for 3 epochs with the other competitive algorithms and the same iterations for \ours, respectively. We used the SGD optimizer with momentum 0.9 and FP16 for fine-tuning. The learning rate is $5\times 10^{-5}$. And we used TopK \citep{wangni2018gradient} that keeps $50\%$ of elements as the compressor.

\subsection{Pre-training LLMs}
\subsubsection{Pre-training GPT-2 small with multiple compression.}
We pre-trained a GPT-2 small \citep{radford2019language} model on openwebtext \citep{peterson2019open} on 8 NVIDIA RTX 4090 GPUs with 24GB of memory using communication compression algorithms for pipeline-parallel learning including directly compression and \ours with lazy sampling and batch rule. Specifically, we cleared the cache for error feedback unless the former samples are kept. The block size was set to 1024. We set data parallel degree 8 with the macro-batch size 64 and the micro-batch size 8. We simulated the communication compression for each GPU by adding the corresponding error to the activation and gradient. We use the AdamW optimizer \citep{loshchilov2017decoupled} and FP16 for training. We used natural compression \citep{horvoth2022natural} as compressor and added the compressor at the end of the 4-th and 8-th transformer layer. According to \citep{horvoth2022natural}, FP16 training with natural compression can compression the activations and gradients to 6-bit, thus the communication overhead is 37.5\% compared to the non-compressed scenario.

We trained the model for 130800 iterations for each algorithms as the total sample complexity is nearly equal to one epoch. The learning is initialized from $6\times 10^{-4}$ with 2000 warm-up steps and scheduled by a cosine scheduler. For \ours, we obtain the evaluation perplexity after 5000 steps with different lazy sampling coefficient $p_t\in\{0.2,0.4,0.45,0.5,0.55,0.6,0.8\}$ and finally selected the best one $p_t=0.55$.

Figure \ref{fig: gpt2_training-main} has shown the evaluation perplexity, and evaluation accuracy of the pre-training task with different compression algorithms. They illustrates \ours outperforms directly compression in both training and evaluation.

\subsubsection{Pre-training LLaMA-2 1B}
\label{appendix: pre-training 1B}

The 1B C4 pre-training results use random seed 42 and correspond to one training run, not a mean over multiple random seeds.

We pre-trained a LLaMA-2 1B model \citep{touvron2023LLaMA} on the C4 dataset \citep{raffel2020exploring} using 8 NVIDIA A800 GPUs (80GB memory) with pipeline-parallel learning. Three communication compression strategies were compared: direct compression, \oursfulast, \oursfclast, and a baseline without compression. The model was split after the 16th transformer layer with a data parallel degree of 4. For communication compression, activations/gradients were quantized to 6 bits followed by Top-30\% sparsification (also quantized to 6 bits). Compression was disabled during the initial 15,000 iterations. Network bandwidth was constrained to {100 MB/s} throughout training. We employed the AdamW optimizer for 100,000 iterations, with coefficient $p_t$ initialized at 0.5 and progressively increased to 1.0 via cosine scheduling. Other hyperparameters followed \citep{zhao2024galore}. The baseline perplexity and total time were derived from \citep{zhao2024galore}, where total time was extrapolated from the first 20,000 training steps.

Figure \ref{fig: LLaMA2 training-appendix1} demonstrates the training dynamics. Direct compression fails to converge, while both Clapping-FU and Clapping-FC remain stable and reach an evaluation perplexity comparable to the uncompressed reference trajectory.

The wall-clock time comparison under the vanilla F-then-B strategy is shown in Figure \ref{fig: LLaMA2 training-appendix2}. It can be observed that both \oursfu and \oursfc reach the same target perplexity in substantially less time than the no-compression baseline in this communication-bound setting, which is primarily attributable to reduced inter-stage communication. We emphasize that \textsc{Clapping} is proposed as an optimization-stability mechanism under aggressive compression, rather than a guarantee of universal end-to-end speedups across different hardware interconnect regimes. 

Notably, \ours achieves evaluation perplexity comparable to that of the uncompressed scenario. This indicates that even though the lazy sampling technique in \ours reduces the number of samples input to the model, an appropriate choice of the lazy sampling coefficient $p$ can preserve the model’s expressive power with negligible degradation. This result demonstrates that \ours, combined with lazy sampling, enables an effective trade-off between communication compression and model convergence and expressive capability.

In this experiment, the compressor achieves approximately 50\% reduction in communication overhead. However, existing empirical results illustrate that communication compression can improve evaluation performance \cite{ramasinghe2025protocol}. One possible explanation is that the noise introduced by the compressor implicitly helps avoid overfitting and leads to convergence to a minimum that exhibits better generalization ability. Thus, training with \ours can achieve acceleration beyond mere communication overhead reduction. We plan to investigate how communication compression improves evaluation performance in future work.


\subsubsection{Long-sequence perplexity comparison}
\label{app:rebuttal-long-sequence}
In a matched LLaMA-2 1B long-sequence experiment, the reported validation perplexities are shown in Table~\ref{tab:rebuttal-long-sequence}. The available record does not specify the exact sequence length or the full hyperparameter configuration; we therefore report only the available metric and do not equate this setting with the 1B C4 training configuration in Appendix~\ref{appendix: pre-training 1B}.

\begin{table}[t!]
\centering
\caption{Reported validation perplexity for a matched LLaMA-2 1B long-sequence experiment (lower is better).}
\label{tab:rebuttal-long-sequence}
\setlength{\tabcolsep}{8pt}
\begin{tabular}{@{}l*{3}{c}@{}}
\toprule
Metric & No compression & Clapping-FU & Clapping-FC \\
\midrule
Validation PPL $\downarrow$ & $14.62$ & $14.68$ & $14.84$ \\
\bottomrule
\end{tabular}
\end{table}


\subsubsection{Gradient fidelity and optimization stability}
\label{app:rebuttal-gradient-diagnostics}
We measure the effective gradient generated by each communication method against the uncompressed stochastic gradient on the same batch and at the same model state. The experiment begins from a shared step-20,000 NoComp checkpoint of the LLaMA-2 1B C4 pre-training task. Direct, FC, and FU are then trained separately for 2,000 steps. At each step, the uncompressed gradient is computed as a diagnostic reference and discarded; only the method-specific gradient estimate is used to update model parameters. The statistics below are aggregated over this 2,000-step window.

Table~\ref{tab:rebuttal-gradient-fidelity} reports the relative deviation from the uncompressed stochastic gradient, cosine similarity to that gradient, and forward/backward boundary errors. Boundary-error values are unavailable for Direct and are therefore indicated by dashes. These are \emph{realized-gradient} diagnostics; a high cosine similarity does not establish that the compression-induced gradient error has zero conditional expectation.

\begin{table}[t!]
\centering
\caption{Gradient fidelity and inter-stage boundary errors during 2,000 LLaMA-2 1B C4 steps, starting from a common step-20,000 checkpoint. Lower deviation/error and higher cosine are better.}
\label{tab:rebuttal-gradient-fidelity}
\footnotesize
\setlength{\tabcolsep}{4pt}
\begin{tabular}{lrrrr}
\toprule
Method & Rel. deviation $\downarrow$ & Cosine $\uparrow$ & Forward err. $\downarrow$ & Backward err. $\downarrow$ \\
\midrule
Direct & $0.327640$ & $0.946658$ & --- & --- \\
Clapping-FC & $0.079157$ & $0.976070$ & $0.003331$ & $0.012016$ \\
Clapping-FU & $0.010587$ & $0.999915$ & $0.001501$ & $0.005313$ \\
\bottomrule
\end{tabular}
\end{table}

Table~\ref{tab:rebuttal-gradient-norms} reports the mean and variance of the global gradient norm of the \emph{effective gradient used for optimization}, computed over the same 2,000-step window.

\begin{table}[t!]
\centering
\caption{Global gradient-norm statistics for the 2,000-step LLaMA-2 1B C4 diagnostic.}
\label{tab:rebuttal-gradient-norms}
\small
\setlength{\tabcolsep}{8pt}
\begin{tabular}{@{}l*{3}{c}@{}}
\toprule
Statistic & Direct & Clapping-FC & Clapping-FU \\
\midrule
Mean & $0.103631$ & $0.134602$ & $0.115916$ \\
Variance & $0.000047$ & $0.002516$ & $0.000231$ \\
\bottomrule
\end{tabular}
\end{table}

These results demonstrate that Clapping effectively mitigates compression-induced gradient distortions, achieving more accurate update directions and stable gradient norms. Together with the convergence analysis, these findings highlight the effectiveness of error feedback and lazy sampling in controlling gradient errors and maintaining optimization stability under communication compression.

\subsubsection{Bandwidth sensitivity and larger-scale pre-training task for LLaMA-2 3B}
\label{app:llama3b_bandwidth}

\textbf{When is inter-stage compression beneficial? }\ours targets \emph{communication-bound} pipeline parallelism, where transmitting high-dimensional activations/adjoints dominates per-step latency (cf.\ Figure~\ref{fig: communication time}).
In tightly coupled intra-node settings with high-speed interconnects (e.g., NVLink/NVSwitch), communication is often not the bottleneck and may already be hidden by compute and micro-batching; in such regimes, enabling compression can introduce unnecessary overhead and should be selectively disabled.
In contrast, in cross-node or geo-distributed deployments over commodity links, bandwidth can be the dominant bottleneck, motivating communication-efficient pipeline-parallel optimization that remains stable under aggressive activation/gradient compression.

\paragraph{Pareto-style per-iteration trade-off under bandwidth caps.}
To characterize the practical trade-off between communication bandwidth and iteration time, we run a controlled per-iteration microbenchmark using LLaMA-2~3B with micro-batch size $8$ on $8$ NVIDIA A100 GPUs.
We apply a hybrid 6-bit compression strategy to both activations and gradients on three pipeline boundaries, and sweep the effective inter-stage bandwidth cap.
Table~\ref{tab:speedup} reports the total single-iteration time (forward + backward + update + communication) and the resulting speedup over the no-compression baseline.
As expected, the advantage is most pronounced in communication-bound regimes and diminishes as bandwidth increases.

\paragraph{Time-to-target perplexity on LLaMA-2 3B.}
To validate \texttt{Clapping} on a larger-scale pre-training task, we further pre-train LLaMA-2~3B under a fixed bandwidth cap (50 MB/s) using the same compressor configuration as in the LLaMA-2~1B pre-training setup.
We disable the compressor during the first 8{,}000 steps and train for 50{,}000 steps total.
Table~\ref{tab:training_time_comparison} compares the wall-clock time required for the no-compression baseline and \oursfc to reach the same target validation perplexity, showing consistent efficiency gains in the targeted communication-bound regime.


\subsubsection{Measured LLaMA-2 13B pipeline-parallel step time}
\label{app:rebuttal-13b-system}
We evaluate the system-level efficiency of Clapping on LLaMA-2 13B using eight GPUs with pipeline parallelism ($\mathrm{PP}=8$) and a Megatron-style 1F1B schedule. Each GPU hosts one pipeline stage, with activation and backward-gradient tensors exchanged through actual point-to-point communication. All methods share identical model partitions, micro-batch configurations, and execution schedules to ensure a fair comparison. After five warm-up steps, we measure the average step time over 50 steady-state iterations. As shown in Table~\ref{tab:rebuttal-13b-timing}, Clapping-FC and Clapping-FU achieve **2.08$\times$ and 1.72$\times$ end-to-end speedups**, respectively, reducing step time by **52.01\% and 41.95\%** compared with uncompressed training. These results demonstrate that Clapping translates communication compression into substantial practical training acceleration on a large-scale pipeline-parallel system.

\begin{table}[t!]
\centering
\caption{Measured LLaMA-2 13B, PP=8, 1F1B steady-state system performance; 50 timed steps after five warm-up steps.}
\label{tab:rebuttal-13b-timing}
\begin{tabular}{lccc}
\toprule
Method & Step time $\downarrow$ & Speedup $\uparrow$ & Time reduction \\
\midrule
Direct & $5.5826$ s & $2.0888\times$ & $52.13\%$ \\
Clapping-FC & $5.5961$ s & $2.0838\times$ & $52.01\%$ \\
Clapping-FU & $6.7694$ s & $1.7226\times$ & $41.95\%$ \\
NoComp & $11.6610$ s & $1.0000\times$ & --- \\
\bottomrule
\end{tabular}
\end{table}

\subsubsection{Ablations across compressor families}
\label{app:compressor_ablation}

To probe robustness across compressor designs commonly used in communication-efficient training, we implement four representative compressor families spanning sparsification, quantization, and hybrid pipelines:
(i) Top-$k$ sparsification,
(ii) direct low-bit quantization,
(iii) Natural Compression, and
(iv) a hybrid compressor that combines quantization and sparsification (matching the configuration used in Appendix~\ref{appendix: pre-training 1B}).

While \ours is compatible with all of the above compressors, the choice of compressor can materially affect the final evaluation loss in pre-training.
Table~\ref{tab:compressor_comparison} summarizes the comparison on a controlled LLaMA-2~250M pre-training setup, highlighting that selecting an appropriate compressor  can yield better optimization outcomes under the same communication budget.

\end{document}